\documentstyle[amscd,amssymb,xypic,11pt]{amsart}
\xyoption{all}

\newcommand{\nc}{\newcommand}

\nc{\exto}[1]{\stackrel{#1}{\longrightarrow}}
\nc{\lan}{\big\langle}
\nc{\ran}{\big\rangle}

\nc{\kk}{{\mathsf{k}}}

\nc{\C}{{\mathbb{C}}}
\nc{\HH}{{\mathbb{H}}}
\nc{\PP}{{\mathbb{P}}}
\nc{\QQ}{{\mathbb{Q}}}
\nc{\ZZ}{{\mathbb{Z}}}

\nc{\CA}{{\mathcal{A}}}
\nc{\CB}{{\mathcal{B}}}
\nc{\CC}{{\mathcal{C}}}
\nc{\D}{{\mathcal{D}}}
\nc{\CE}{{\mathcal{E}}}
\nc{\CF}{{\mathcal{F}}}
\nc{\CG}{{\mathcal{G}}}
\nc{\CH}{{\mathcal{H}}}
\nc{\CL}{{\mathcal{L}}}
\nc{\CM}{{\mathcal{M}}}
\nc{\CN}{{\mathcal{N}}}
\nc{\CO}{{\mathcal{O}}}
\nc{\CQ}{{\mathcal{Q}}}
\nc{\CR}{{\mathcal{R}}}
\nc{\CS}{{\mathcal{S}}}
\nc{\CT}{{\mathcal{T}}}
\nc{\CU}{{\mathcal{U}}}
\nc{\CV}{{\mathcal{V}}}
\nc{\CW}{{\mathcal{W}}}
\nc{\CX}{{\mathcal{X}}}
\nc{\CY}{{\mathcal{Y}}}
\nc{\CMo}{{\mathcal{M}^\circ}}
\nc{\Co}{{{C}^\circ}}

\nc{\BY}{{\overline{Y}}}
\nc{\BYD}{{\overline{Y}{}^{|D|}}}
\nc{\OZ}{{\overline{Z}}}
\nc{\bg}{{\bar{g}}}

\nc{\bq}{{\mathbf{q}}}
\nc{\BG}{{\mathbf{G}}}
\nc{\BM}{{\mathbf{M}}}
\nc{\BP}{{\mathbf{P}}}
\nc{\BZ}{{\mathbf{Z}}}
\nc{\BPr}{{\mathsf{P}}}
\nc{\BR}{{\mathbf{R}}}
\nc{\BRO}[1]{{{\mathbf{R}}^{\circ}_{#1}}}
\nc{\BRD}[1]{{{\mathbf{R}}^{|D|}_{#1}}}
\nc{\BRP}[1]{{{\mathbf{R}}^{1}_{#1}}}
\nc{\BRTP}[1]{{{\mathbf{\tilde{R}}}{}^{1}_{#1}}}
\nc{\BS}{{\mathbf{S}}}
\nc{\BMS}{{{\mathbf{M}}^{{s}}}}
\nc{\BMSS}{{{\mathbf{M}}^{{ss}}}}
\nc{\BMZ}{{\mathbf{M}^{\circ}}}
\nc{\BCL}{{\mathbf{L}}}

\nc{\Cl}{{\mathsf{Cliff}}}
\nc{\Clev}{{\mathop{\mathsf{Cliff}}^{\circ}}}

\nc{\fg}{{\mathfrak{g}}}
\nc{\fp}{{\mathfrak{p}}}
\nc{\FD}{{\mathfrak{D}}}
\nc{\FE}{{\mathfrak{E}}}
\nc{\FL}{{\mathfrak{L}}}
\nc{\FM}{{\mathfrak{M}}}
\nc{\FS}{{\mathsf{S}}}

\nc{\sfc}{{\mathsf{c}}}
\nc{\sfch}{{\mathsf{ch}}}
\nc{\sfh}{{\mathsf{h}}}

\nc{\SK}{{\mathsf{K}}}
\nc{\SO}{{\mathsf{O}}}
\nc{\SQ}{{\mathsf{Q}}}
\nc{\SPV}{{\mathsf{S}^+\mathsf{V}}}
\nc{\SMV}{{\mathsf{S}^-\mathsf{V}}}
\nc{\SPMV}{{\mathsf{S}^\pm\mathsf{V}}}
\nc{\SX}{{S_X}}
\nc{\SY}{{S_Y}}
\nc{\phipsi}{{q}}
\nc{\eps}{\varepsilon}

\nc{\pim}{{\pi_-}}
\nc{\pip}{{\pi_+}}

\nc{\BE}{{\overline{\CE}}}
\nc{\TE}{{\tilde{\CE}}}
\nc{\TQ}{{\tilde{Q}}}
\nc{\TCF}{{\tilde{\CF}}}
\nc{\TCG}{{\tilde{\CG}}}
\nc{\TCL}{{\tilde{\CL}}}
\nc{\TF}{{\tilde{F}}}
\nc{\TW}{{\tilde{W}}}
\nc{\TCC}{{\tilde{\CC}}}
\nc{\TCX}{{\tilde{\CX}}}
\nc{\TCY}{{\tilde{\CY}}}
\nc{\TPhi}{{\tilde{\Phi}}}
\nc{\txi}{{\tilde{\xi}}}
\nc{\tp}{{\tilde{p}}}
\nc{\tq}{{\tilde{q}}}
\nc{\tzeta}{{\tilde{\zeta}}}
\nc{\tpi}{{\tilde{\pi}}}

\nc{\HE}{{\widehat{\CE}}}
\nc{\HX}{{\hat{X}}}
\nc{\hxi}{{\hat{\xi}}}

\nc{\UH}{{\mathcal{H}}}

\nc{\TM}{{\widetilde{M}}}
\nc{\TCM}{{\widetilde{\CM}}}
\nc{\TU}{{\widetilde{U}}}
\nc{\TX}{{\widetilde{X}}}
\nc{\TY}{{\widetilde{Y}}}
\nc{\TYO}{{{\widetilde{Y}}^\circ}}
\nc{\barf}{{\bar{f}}}
\nc{\te}{{\tilde{e}}{}}
\nc{\tf}{{\tilde{f}}}
\nc{\tg}{{\tilde{g}}}
\nc{\ti}{{\tilde{\imath}}}
\nc{\tj}{{\tilde{\jmath}}}
\nc{\ty}{{\tilde{y}}}
\nc{\tphi}{{\tilde{\phi}}}
\nc{\hf}{{\hat{f}}}

\nc{\urho}{{\underline{\rho}}}

\nc{\LRA}{\Leftrightarrow}
\nc{\RA}{\Rightarrow}
\nc{\lotimes}{\mathbin{\mathop{\otimes}\limits^{\mathbb{L}}}}
\nc{\CEnd}{\mathop{\mathcal{E}\mathit{nd}}\nolimits}
\nc{\CExt}{\mathop{\mathcal{E}\mathit{xt}}\nolimits}
\nc{\CHom}{\mathop{\mathcal{H}\mathit{om}}\nolimits}
\nc{\RH}{\mathop{{\mathsf{R}}\Gamma}\nolimits}
\nc{\RGamma}{\mathop{{\mathsf{R}}\Gamma}\nolimits}
\nc{\RHom}{\mathop{\mathsf{RHom}}\nolimits}
\nc{\RCHom}{\mathop{\mathsf{R}\mathcal{H}\mathit{om}}\nolimits}
\nc{\RG}{\mathop{\mathsf{R\Gamma}}\nolimits}
\nc{\Hom}{\mathop{\mathsf{Hom}}\nolimits}
\nc{\Ext}{\mathop{\mathsf{Ext}}\nolimits}
\nc{\End}{\mathop{\mathsf{End}}\nolimits}
\nc{\Tor}{\mathop{\mathsf{Tor}}\nolimits}
\nc{\Tordim}{\mathop{\mathsf{Tor}\text{\rm-}\mathsf{dim}}\nolimits}
\nc{\Hilb}{\mathop{\mathsf{Hilb}}\nolimits}
\nc{\Spec}{\mathop{\mathsf{Spec}}\nolimits}
\nc{\Pic}{\mathop{\mathsf{Pic}}\nolimits}

\nc{\Tw}{\mathop{\mathsf{Tw}}\nolimits}
\nc{\Ker}{\mathop{\mathsf{Ker}}\nolimits}
\nc{\Coker}{\mathop{\mathsf{Coker}}\nolimits}
\nc{\codim}{\mathop{\mathsf{codim}}\nolimits}
\nc{\sing}{{\mathsf{sing}}}
\nc{\supp}{\mathop{\mathsf{supp}}}
\nc{\perf}{{\mathsf{perf}}}
\nc{\rank}{\mathop{\mathsf{rank}}}
\nc{\Pf}{{\mathsf{Pf}}}
\nc{\Gr}{{\mathsf{Gr}}}
\nc{\OGr}{{\mathsf{OGr}}}
\nc{\Flag}{{\mathsf{Fl}}}
\nc{\Kosz}{{\mathsf{Kosz}}}
\nc{\LGr}{{\mathsf{LGr}}}
\nc{\GTGr}{{\mathsf{G_2Gr}}}
\nc{\GTF}{{\mathsf{G_2F}}}
\nc{\OF}{{\mathsf{OF}}}
\nc{\Fl}{{\mathsf{Fl}}}
\nc{\Bl}{{\mathsf{Bl}}}
\nc{\GL}{{\mathsf{GL}}}
\nc{\PGL}{{\mathsf{PGL}}}
\nc{\SL}{{\mathsf{SL}}}
\nc{\SP}{{\mathsf{Sp}}}
\nc{\Spin}{{\mathsf{Spin}}}
\nc{\Tot}{{\mathsf{Tot}}}
\nc{\ev}{{\mathsf{ev}}}
\nc{\coev}{{\mathsf{coev}}}
\nc{\id}{{\mathsf{id}}}
\nc{\opp}{{\mathsf{opp}}}
\nc{\PS}{{{\PP^3}}}
\nc{\Qu}{{{Q^3}}}
\nc{\tdim}{\mathop{\Tor\dim}}
\nc{\ecart}{{\fbox{$\scriptstyle\mathsf{EC}$}}}
\nc{\ad}{{\mathop{\mathsf ad}}}
\nc{\Coh}{{\mathop{\mathsf Coh}}}
\nc{\Ab}{{\mathop{\mathcal{A}\mathit{b}}}}
\nc{\Qcoh}{{\mathop{\mathsf Qcoh}}}

\nc{\AAV}{{\mathcal{AAV}}}

\nc{\Rep}{{\mathsf{Rep}}}

\nc{\Cubics}{{{\mathcal{S}}_3}}
\nc{\VFT}{{{\mathcal{S}}_{14}}}
\nc{\VFTE}{{{\mathcal{N}}_{\mathrm{reg,sm}}}}
\nc{\MX}{{\CM_X}}
\nc{\MY}{{\CM_Y}}
\nc{\MYE}{{\CM_{Y,\CE}}}
\nc{\Yd}{{Y_d}}
\nc{\Yfive}{{Y_5}}
\nc{\Xg}{{X_{2g-2}}}
\nc{\Xtt}{{X_{22}}}
\nc{\Xst}{{X_{16}}}
\nc{\Xtw}{{X_{12}}}
\nc{\Xe}{{X_{8}}}
\nc{\Xf}{{X_{4}}}

\nc{\git}{{/\!\!/\!{}_\chi}}

\theoremstyle{plain}

\newtheorem{theo}{Theorem}[]
\newtheorem{theorem}{Theorem}[section]

\newtheorem{lemma}[theorem]{Lemma}
\newtheorem{proposition}[theorem]{Proposition}
\newtheorem{corollary}[theorem]{Corollary}

\theoremstyle{definition}

\newtheorem{definition}[theorem]{Definition}

\theoremstyle{remark}

\newtheorem{remark}[theorem]{Remark}

\newenvironment{proof}{\noindent{\sf Proof:}}{\qed\medskip}

\title{Hyperplane sections and derived categories}
\author{Alexander Kuznetsov}
\address{
Algebra Section, Steklov Mathematical Institute,
8 Gubkin str., Moscow 119991 Russia
}
\email{akuznet@@mi.ras.ru}
\date{}
\thanks{I was partially supported by RFFI grants 02-01-00468 and 02-01-01041,
Russian Presidential grant for young scientists No. MK-3926.2004.1,
CRDF Award No. RM1-2405-MO-02, and the Russian Science Support Foundation.}

\begin{document}

\begin{abstract}
We give a generalization of the theorem of Bondal and Orlov
about the derived categories of coherent sheaves on intersections
of quadrics revealing its relation to projective duality.
As an application we describe the derived categories
of coherent sheaves on Fano 3-folds of index 1
and degrees 12, 16 and~18.
\end{abstract}

\maketitle

{\small
\tableofcontents
}

\section*{Introduction}

One of the most promising directions in modern algebraic geometry is
the systematic use of derived categories of coherent sheaves~\cite{BO2}.
Derived categories introduced by Verdier back in 1967 recently
attracted a lot of attention due to the progress in understanding
of their role in the modern geometry. It was realized that
geometric similarity between different algebraic varieties
sometimes is only the tip of the iceberg and is a consequence
of a deeper relation which can be expressed as an equivalence
of appropriate categories. Take for example the McKay correspondence.
The story began in 1980 with John McKay's discovery \cite{McK} of
a correspondence between irreducible representations of a finite
subgroup $\Gamma\subset\SL(2,\C)$ and vertices of an affine
Coxeter-Dynkin diagram of type $A,D,E$. A geometric interpretation
of this correspondence incorporating a resolution of the simple
surface singularity $\C^2/\Gamma$ was given by Gonzalez-Sprinberg
and Verdier \cite{GSV} in 1983. Kapranov and Vasserot \cite{KV}
have shown that the derived category of coherent sheaves on the resolution
of $\C^2/\Gamma$ is equivalent to the derived category of coherent
$\Gamma$-equivariant sheaves on $\C^2$. See also \cite{BKR}
for a beautiful generalization.

Another example of this sort is the theory of intersections of quadrics.
Classically, it was formulated as a correspondence between
intersections of quadrics and determinantal loci in the corresponding
linear spaces of quadrics~\cite{T}. The geometric meaning
of this correspondence was discovered in~\cite{DR}. In the particular
case of the intersection of a pencil of even-dimensional quadrics
one should consider the hyperelliptic curve, obtained as a twofold
covering of $\PP^1$ parameterizing quadrics in the pencil with ramification
at the points of $\PP^1$ corresponding to degenerate quadrics. It was shown
in loc.\ cit.\ that the moduli spaces of rank 2 vector bundles on
this hyperelliptic curve can be described in terms of linear subspaces
in the intersection of quadrics.
%
%
%
%
Bondal and Orlov in \cite{BO1} gave a categorical
meaning to this by proving that the category of coherent sheaves
on the hyperelliptic curve embeds fully and faithfully
into the derived category of the intersection of quadrics.
Finally, Bondal and Orlov \cite{BO2,BO3} stated a theorem
describing the structure of the derived category of a complete intersection
of any number of quadrics. Roughly speaking it says that for any family
of quadrics $L \subset S^2W^*$ in an even-dimensional vector space~$W$
there exists a sheaf of finite algebras $\CA_L$ on the twofold covering
of $\PP(L)$ ramified at the degeneration locus such that

\begin{theo}[\cite{BO2,BO3}]\label{boconj}
If $2\dim L < \dim W$ then there is a semiorthogonal decomposition
$$
\D^b(X_L) = \langle \D^b(\Coh(\CA_L)),
\CO_{X_L},\dots,\CO_{X_L}(\dim W - 2\dim L - 1)\rangle,
$$
and if $2\dim L = \dim W$ then there is an equivalence of categories
$\D^b(X_L) \cong \D^b(\Coh(\CA_L))$,
where $X_L$ stands for the complete intersection of quadrics
in the family $L$, and $\D^b(\Coh(\CA_L))$ stands for the derived category
of sheaves of coherent $\CA_L$-modules.
\end{theo}

Let us explain the relation of this theorem to projective duality.

Note that the intersection of quadrics
$X_L$ can be considered as a linear section of the Veronese
subvariety $X = \PP(W) \subset \PP(S^2W)$ by the linear space
$\PP(L^\perp) \subset \PP(S^2W)$, where $L^\perp \subset S^2W^*$
is the orthogonal subspace to $L \subset S^2W^*$.
On the other hand, the sheaf of algebras $\CA_L$
on $\PP(L)\subset\PP(S^2W^*)$ can be extended to
a sheaf $\CA_{\PP(S^2W^*)}$ on the whole $\PP(S^2W^*)$.
Actually, this sheaf is just the sheaf of even parts
of the universal Clifford algebra, i.e.
$$
\CA_{\PP(S^2W^*)} =
\CO_{\PP(S^2W^*)} \oplus
\Lambda^2 W\otimes\CO_{\PP(S^2W^*)}(-1) \oplus \dots \oplus
\Lambda^{2n}W\otimes\CO_{\PP(S^2W^*)}(-n)
$$
with the Clifford multiplication. Consider $\PP(S^2W^*)$
with the sheaf of algebras $\CA_{\PP(S^2W^*)}$ as a noncommutative
algebraic variety, in fact a noncommutative finite covering of $\PP(S^2W^*)$.
One of consequences of theorem~\ref{boconj} which we would like to stress on is
$$
\parbox{0.8\textwidth}{\bf the set of critical values of the projection
$(\PP(S^2W^*),\CA_{\PP(S^2W^*)}) \to \PP(S^2W^*)$\\
coincides with the projectively dual variety $X^\vee \subset \PP(S^2W^*)$.}
$$
Indeed, since the projection $(\PP(S^2W^*),\CA_{\PP(S^2W^*)}) \to \PP(S^2W^*)$
is flat, a point $H\in\PP(S^2W^*)$ is its critical value if and only if
the fiber is singular. The singularity of an algebraic variety
is equivalent to $\Ext$-boundedness of its derived category
of coherent sheaves. Thus the set of critical values coincides
with the set of $H\in\PP(S^2W^*)$ such that the category $\D^b(\Coh(\CA_H))$
is not $\Ext$-bounded. Looking at the claim of theorem~\ref{boconj} we see that this
happens if and only if the category $\D^b(X_H)$ is not $\Ext$-bounded,
i.e. if $X_H$ is singular. But $X_H$ is the hyperplane section
of $X\subset\PP(S^2W)$ by the hyperplane $H\in\PP(S^2W^*)$,
so recalling the definition of the projective duality
we see that the set of critical values coincides
with the projectively dual variety.

In fact, from the homological point of view the noncommutative variety
$(\PP(S^2W^*),\CA_{\PP(S^2W^*)})$ is a better candidate to be called
the projective dual of $X$ then $X^\vee$. It not only describes
the set of singular hyperplane sections but also remembers
nontrivial parts of their derived categories (for example,
the triangulated category of singularities \cite{O3} of $X_H$
is equivalent to  the triangulated category of singularities
of the fiber of the projection $(\PP(S^2W^*),\CA_{\PP(S^2W^*)}) \to \PP(S^2W^*)$
at $H$), in particular
contains a lot of information about the singularities they have.
So, it would be reasonable to say that $(\PP(S^2W^*),\CA_{\PP(S^2W^*)})$
is a {\em Homologically Projectively Dual}\/ variety of $X$.
The above considerations suggest that given an algebraic variety
$X\subset\PP(V)$, its homologically projectively dual variety
should be defined as an algebraic variety $Y$ with a morphism $g:Y \to \PP(V^*)$
and a sheaf of noncommutative algebras $\CA_Y$ on $Y$, such that
for any hyperplane $H \subset V^*$ the nontrivial part
of the derived category $\D^b(X_H)$ is equivalent
to the derived category $\D^b(Y_H,\CA_Y)$
(roughly speaking, there exists a fully faithful embedding
$\D^b(Y_H,\CA_Y) \to \D^b(X_H)$ and the orthogonal to $\D^b(Y_H,\CA_Y)$
in $\D^b(X_H)$ is $\Ext$-bounded).

I am convinced that the Homological Projective Duality exists
for a large class of algebraic varieties. In this paper
we show how to construct a homologically projectively dual for a small
class of varieties. Assume that $X$ is a smooth projective variety
admitting an exceptional pair of vector bundles $(E_1,E_2)$ such that
$$
(E_1\otimes\CO_{\PP(V)|X}(1),E_2\otimes\CO_{\PP(V)|X}(1),
\dots,
E_1\otimes\CO_{\PP(V)|X}(i),E_1\otimes\CO_{\PP(V)|X}(i))\leqno{(*)}
$$
is an exceptional collection, where $i$ is the index of $X$
(i.e.\ $\omega_X \cong \CO_{\PP(V)}(-i)_{|X}$).
In this case $Y$ can be constructed as a moduli space
of stable representations $(R_1,R_2)$ of the quiver
$$
\SQ = \xymatrix@1{\bullet \ar[rr]^{\Hom(E_1,E_2)^*} && \bullet}
$$
such that the cokernel of the corresponding map
$R_1\otimes E_1 \to R_2\otimes E_2$ is supported
on a hyperplane section of~$X$. The moduli space $Y$
carries a sheaf of Azumaya algebras $\CA_Y$ and a universal
family of representations $(F_1,F_2)$ in the category of $\CA_Y$-modules
(the Azumaya algebra $\CA_Y$ corresponds to the element in the Brauer group
which is the obstruction to the existence of a universal family,
so in the category of $\CA_Y$-modules the universal family exists).
Moreover, associating to a point of $Y$ the support hyperplane
of the cokernel of the corresponding map $R_1\otimes E_1 \to R_2\otimes E_2$
we obtain a projection $g:Y \to \PP(V^*)$. Denoting as above
$X_L = X \cap \PP(L^\perp)$, $Y_L = Y\times_{\PP(V^*)}\PP(L) = g^{-1}(\PP(L))$
we prove that:

\begin{theo}
Assume that both $X_L$ and $Y_L$ are compact and have expected dimension.
If $r = \dim L \le i$ then there is a semiorthogonal decomposition
$$
\D^b(X_L) = \langle \D^b(Y_L,\Coh(\CA_Y)),
E_1(1),E_2(1),\dots,E_1(i-r),E_2(i-r)\rangle,
$$
and if $r = \dim L \ge i$ then there is a semiorthogonal decomposition
$$
\D^b(Y_L,\CA_Y) = \langle \D^b(X_L),
F_2^*(i-r),F_1^*(i-r),\dots,F_2^*(-1),F_1^*(-1)\rangle.
$$
In particular, for $\dim L = i$ we have an equivalence of categories
$\D^b(X_L) \cong \D^b(Y_L,\Coh(\CA_Y))$.
\end{theo}

In fact instead of giving a direct proof
we extract a list of conditions on $X$, $(E_1,E_2)$,
$(Y,\CA_Y)$, and $(F_1,F_2)$ which imply the claim of the theorem.
In some cases (when one can guess what $Y$ is) this allows to skip
the construction of $Y$ at all.
The proof of the theorem splits into two principal parts.
In the first part we check the case $\dim L = 1$, the case of
hyperplane sections of $X$. In the other words, we prove there that
$Y$ is homologically projectively dual to $X$.
It turns out that the second part, the case of $\dim L > 1$,
follows from the first part more or less formally.
Thus the above theorem in a sense is a consequence of homological
projective duality.

A few words about the assumptions on the subspace $L$ made in the theorem.
First of all, the compactness of $X_L$ is not a restriction, since $X$
is projective. On the other hand, the compactness of $Y_L$ is the most
serious restriction. The problem is that $Y$ may not be compact
(we constructed $Y$ as a moduli space of {\em stable} representations,
and stability is an open condition).
In fact, this problem arises because of imperfection of the construction
of $Y$ (see below). In fact we expect that it is always possible to construct
a compact $Y$ such that the theorem holds. In this case the assumption
of compactness of $Y_L$ will hold automatically. Finally, the assumptions
on the dimension of $X_L$ and $Y_L$ also can be taken off. However
in this case we must be more accurate in the definition of $X_L$ and $Y_L$.
In fact, the correct definition of both $X_L$ and $Y_L$ is
the ``derived'' fiber product
$X_L = X\times_{\PP(V)}^{{\mathbb{L}}}\PP(L^\perp)$,
$Y_L = Y\times_{\PP(V^*)}^{{\mathbb{L}}}\PP(L)$.
As a topological space it coincides
with the usual fiber product, but in a contrast with the usual case
it is ringed with a certain sheaf of $DG$-algebras instead of the usual
sheaf of functions. Correspondingly, the categories $\D^b(X_L)$ and
$\D^b(Y_L,\CA_Y)$ are appropriate categories of $DG$-modules.
For instance, the above theorem allows to describe a certain category
of $DG$-modules on any noncomplete intersection of quadrics. Since any
projective variety in sufficiently ample embedding is an intersection
of quadrics, this point of view seems to be very prominent.

Now let us discuss some formal consequences of the above theorem.
First of all, we would like to mention that we haven't assumed
the fullness of the exceptional collection on $X$ we have started from.
However, take $L = 0$, then $X_L = X$, $Y_L = \emptyset$ and it follows
that the collection is full. Thus we have obtained a rather unexpected
corollary. If a Fano variety admits an exceptional collection
of the form~($*$) then ($*$) is full.
Similarly, assume that $Y$ is compact and take $L = V^*$.
Then $X_L = \emptyset$, $Y_L = Y$ and it follows that the collection
$$
(F_2^*(i-N),F_1^*(i-N),\dots,F_2^*(-1),F_1^*(-1)),
$$
where $N = \dim V$ is a full exceptional collection in $\D^b(Y,\CA_Y)$.

It also worth mentioning that the above theorem holds not only
for individual linear sections of $X$ and $Y$, but also for families.
In fact, the proof is based on investigation of families. For example,
for the universal family $\CX_1$ of hyperplane sections of $X$
(which is a divisor of bidegree $(1,1)$ in $X\times\PP(V^*)$)
the family $\CY_1 \subset Y\times\PP(V^*)$ of orthogonal linear sections
of $Y$ coincides with the graph of the projection $g$.
The relative version of the theorem in this cases claims that
there is a semiorthogonal decomposition
$$
\D^b(\CX_1)  =
\langle \D^b(Y,\CA_Y),
\D^b(\PP(V^*))\otimes E_1(1),\dots,\D^b(\PP(V^*))\otimes E_2(i-1)\rangle
$$
if $Y$ is compact. This decomposition suggests a categorical
definition of the compactification $\BY$ of~$Y$:
$$
\D^b(\BY,\CA_\BY) =
\langle \D^b(\PP(V^*))\otimes E_1(1),\dots,
\D^b(\PP(V^*))\otimes E_2(i-1)\rangle^\perp
\subset \D^b(\CX_1).
$$
In fact, one could take this for the definition of the homologically
projectively dual variety of $X$, if only there were a possibility
of proving that the above category is {\em geometric}.

Now let us describe what our approach gives in the case we have started
from, i.e.\ for the intersection of even-dimensional quadrics.
In this case we take $X = \PP(W)$, $W$ is even-dimensional,
the Veronese embedding $X \subset \PP(S^2W) = \PP(V)$,
and the exceptional pair $(E_1,E_2) = (\CO_{\PP(W)},\CO_{\PP(W)}(1))$
(note that $\CO_{\PP(V)|X}(1) \cong \CO_{\PP(W)}(2)$,
so the collection ($*$) is the standard exceptional
collection on $\PP(W)$). The discussed above moduli space of representations
of the quiver however doesn't give the whole $\PP(S^2W^*)$,
but only the open part $Y = \PP(S^2W^*) \setminus \BZ$,
where $\BZ$ is the locus of quadrics of corank~$2$. Thus, rigourously
speaking at the moment our method proves theorem~\ref{boconj} only
for $\dim L \le 3$.

One of explanations why the quiver approach doesn't give the whole
$\PP(S^2W^*)$ at the output is the fact that the sheaf of algebras
$\CA_{\PP(S^2W^*)}$ fails to be Azumaya over $\BZ$, while using
the quiver approach we only can get an Azumaya algebra.
However, we expect that there exists a certain algebraic stack $\BY$
with the underlying algebraic variety $\PP(S^2W^*)$
such that the sheaf of algebras $\CA_{\PP(S^2W^*)}$ comes from a sheaf
of Azumaya algebras on $\BY$. This suggests that it might be better
in general situation to consider a moduli stack of the quiver
representations instead of the moduli space.

Correspondingly, in the general situation one can try to find an appropriate
stack compactification $\BY$ of the moduli space $Y$ such that the sheaf
of algebras $\CA_Y$ (and the universal family as well) can be extended
to a sheaf of Azumaya algebras~$\CA_{\BY}$.

Another possible approach to the construction of a Homologically Projectively
Dual variety is to avoid the use of quivers at all and to construct $Y$
as a moduli spaces of appropriate objects (e.g.\ spinor bundles)
in the derived category of coherent sheaves on $X$ which are
supported on hyperplane sections of $X$.

In fact, in both cases the main obstacle is absence of a good
categorical theory of moduli spaces (an output of a good theory
of moduli spaces must be in general a noncommutative algebraic variety).

\bigskip

Now it is time to describe the contents of the paper.
In section~1 we introduce the notation and state the main
result of the paper, theorem~\ref{themain}.
In section~2 we introduce the homological background.
We define exact cartesian squares and faithful base changes
in subsection~\ref{ssecs} and prove the main technical results:
the faithful base change theorem~\ref{phitsod}
and the relative version of the Bridgeland's trick~\ref{ffisequ}.
In section~3 we make the first step of the proof of
theorem~\ref{themain} and in section~4 we finish
the proof. In section~5 we describe the general
construction of $Y$ from $X$ via the quiver moduli spaces.
In section~6 we give some examples where theorem~\ref{themain} works.
In appendices A, B and C we work out some details relating to these examples.
Finally, in appendix~D we check the basic facts about the derived categories
of coherent $\CA$-modules, where $\CA$ is a sheaf of Azumaya algebras.
We pay also special attention to the relation of complexes of finite
$\Tor$ and $\Ext$-amplitude to perfect complexes.

{\bf Acknowledgements.}
I am grateful to A.Bondal, D.Orlov, D.Kaledin and A.Samokhin
for useful discussions. Some of the results of this work
have been presented at the conference
``Noncommutative Algebra and Algebraic Geometry''
in the University of Warwick. I would like to thank
organizers for this possibility.

\section{Statement of results}

As it was indicated in the introduction,
we start with presenting a set of conditions
on a pair of algebraic varieties $X$ and $Y$
implying their homological projective duality.

Assume that we have the following {\bf data}:
\begin{enumerate}
\renewcommand{\labelenumi}{(D.\theenumi)}
\item a smooth, projective variety $X$
with an algebraic morphism to a projective space $f:X \to \PP(V)$;
\item a pair of vector bundles $(E_1,E_2)$ on $X$;
\item an algebraic variety $Y$ with a sheaf of Azumaya algebras $\CA_Y$
and a projective morphism\\
$g:Y \to \BP = \PP(V^*) - \BZ$, where $\BZ \subset \PP(V^*)$
is a closed subset in the dual projective space;
\item a pair of locally projective $\CA_Y$-modules $(F_1,F_2)$ on $Y$;
\item a linear morphism $\phi:\Hom(E_1,E_2)^* \to \Hom_{\CA_Y}(F_1,F_2)$.
\end{enumerate}

Given such a data we denote
\begin{itemize}
\item $W = \Hom(E_1,E_2)^*$;
\item $N = \dim V$;
\item $Q \subset \PP(V)\times\PP(V^*)$ --- the incidence quadric;
\item $Q(X,Y) = (X\times Y) \times_{\PP(V)\times\PP(V^*)} Q$
and $i:Q(X,Y) \to X\times Y$ --- the embedding.
\end{itemize}
The incidence quadric $Q \subset \PP(V)\times\PP(V^*)$
is a divisor of bidegree $(1,1)$. Its fiber over a point $H\in\PP(V^*)$
is the corresponding hyperplane in $\PP(V)$ which we denote
by the same letter $H \subset \PP(V)$. Similarly, the fiber
of $Q(X,Y)$ over a point $H \in \PP(V^*)$ is the hyperplane
section $X_H = X \cap H$ of $X$. Abusing the notation
we call $Q(X,Y)$ the incidence quadric too.
Note, that it follows from the condition (C.1) below
that the incidence quadric $Q(X,Y)$ is a divisor
of bidegree $(1,1)$ on $X\times Y$.

The map $\phi:W \to \Hom_{\CA_Y}(F_1,F_2) = \Hom_{\CA_Y}(F_2^*,F_1^*)$
induces the following evaluation homomorphisms
$W\otimes F_1 \to F_2$ and $W\otimes F_2^* \to F_1^*$.
Tensoring the first  by $F_2^*$ and the second by $F_1$ and summing up
we obtain a homomorphism
\begin{itemize}
\item $\xymatrix{W\otimes F_2^*\boxtimes F_1 \ar[rr]^-{\ad(\phi)} &&
F_2^*\boxtimes F_2 \oplus F_1^*\boxtimes F_1}$ on $Y\times Y$.
\end{itemize}
On the other hand, composing the evaluation homomorphism
$W\otimes F_1 \to F_2$ tensored by $E_2$ with the coevaluation
homomorphism $E_1 \to W\otimes E_2$ tensored by $F_1$ we obtain
a homomorphism
\begin{itemize}
\item $\xymatrix{E_1\boxtimes F_1 \ar[rr]^-{e} && E_2\boxtimes F_2}$
on $X\times Y$.
\end{itemize}


Assume that our data satisfies the following list of {\bf conditions:}
\begin{enumerate}
\renewcommand{\labelenumi}{(C.\theenumi)}
\item $f(X)$ is not contained in a hyperplane in $\PP(V)$;
\item $\omega_X \cong \CO_X(-i) := f^*\CO_{\PP(V)}(-i)$, $i>0$;
\item $(E_1(1),E_2(1),\dots,E_1(i),E_2(i))$ is an exceptional collection,\\
(i.e.\ $\Ext^\bullet(E_s,E_s) = \kk$ for $s = 1,2$ and
$\Ext^\bullet(E_s(k),E_t(l)) = 0$ for $i\ge k>l\ge 1$ or $k=l$ and $s>t$),\\
and {\bf additionally}\/ $\Ext^{>0}(E_1,E_2) = 0$;
\item $g$ is projective;
\item $\Coker\ad(\phi) \cong \Delta_*\CA_Y$ on $Y\times Y$,
where $\Delta:Y \to Y\times Y$ is the diagonal embedding;
\item $\Ker e = 0$ and $\Coker e = i_*\CE$,
where $\CE$ is a coherent $\CA_Y$-module on $Q(X,Y)$;
\item $\dim X + N - \codim_{Y\times Y}(Y\times_{\PP(V^*)} Y) = 2i$;
\item
for any hyperplane $H\subset\PP(V^*)$ we have $\dim (\BZ\cap H) \le N - i - 2$.
\end{enumerate}

The integer $i$ defined in condition (C.2) is called the {\sf index}\/ of $X$.
If $\BZ$ is irreducible and is not contained in a hyperplane in $\PP(V^*)$
then the condition (C.8) is equivalent to $\dim\BZ \le N - i - 1$.

For every linear subspace $L \subset V^*$ we denote
$X_L = X\times_{\PP(V)}\PP(L^\perp)$ and $Y_L = Y\times_{\PP(V^*)}\PP(L)$.

\begin{definition}
A subspace $L\subset V^*$ is {\sf admissible}, if
\begin{enumerate}
\renewcommand{\theenumi}{\alph{enumi}}
\renewcommand{\labelenumi}{(\theenumi)}
\item $\dim X_L = \dim X - \dim L$,
\item $\dim Y_L = \dim Y + \dim L - N$, and
\item $\PP(L)\cap \BZ = \emptyset$.
\end{enumerate}
\end{definition}

Note that condition (c) above together with (C.4) implies
that $Y_L$ is compact.

\begin{theorem}\label{themain}
If conditions\/ {\rm(C.1)--(C.8)}\/ are satisfied then

\noindent $1)$
$Y$ is smooth and $\omega_Y \cong \CO_Y(i-N)$;

\noindent $2)$
for any admissible subspace $L\subset V^*$ there exist
semiorthogonal decompositions
$$
\begin{array}{llrl}
\D^b(X_L) &=&  \langle \D^b(Y_L,\CA_Y),
E_1(1),E_2(1),\dots,E_1(i-r),E_2(i-r)\rangle,
& \text{if $r\le i$} \\
\D^b(Y_L,\CA_Y) &=&
\langle \D^b(X_L),F_2^*(i-r),F_1^*(i-r),\dots,F_2^*(-1),F_1^*(-1)\rangle,
& \text{if $r\ge i$},
\end{array}
$$
where $r = \dim L$.
\end{theorem}

\bigskip

Now let us sketch the proof of the theorem, assuming for simplicity
that $\BZ = \emptyset$.
First of all, to avoid problems which arise when one works with derived
categories of singular varieties, we replace individual linear
sections $X_L$, $Y_L$ of $X$ and $Y$ by the universal families
of linear sections $\CX_r \subset X\times\Gr(r,V^*)$,
$\CY_r \subset Y\times\Gr(r,V^*)$, which enjoy the property
to be totally smooth. We note that the natural projection
$\CX_r\times_{\Gr(r,V^*)}\CY_r \to X\times Y$ factors through $Q(X,Y)$
and denote by $\CE_r$ the pullback of the sheaf $\CE$,
defined by condition (C.6). Then the objects $\CE_r$ on
$\CX_r\times_{\Gr(r,V^*)}\CY_r \subset \CX_r\times\CY_r$,
considered as kernels, provide us with functors
$\Phi_r:\D^b(\CY_r,\CA_Y) \to \D^b(\CX_r)$. Consider also
their left and right adjoint functors
$\Phi_r^*,\Phi_r^!:\D^b(\CX_r) \to \D^b(\CY_r,\CA_Y)$.
Now the proof goes in several steps:

\begin{description}
\item[Step 1]
We note that $\CY_1 = Y$ and prove that the functor
$\Phi_1:\D^b(Y,\CA_Y) \to \D^b(\CX_1)$ is fully faithful.
This is done by explicit computation of $\Phi_1^* \circ \Phi_1$.
The convolution of kernels of these functors is shown to be isomorphic
to the sheaf $\Coker \ad(\phi)$ on $Y\times Y$, which by condition (C.4)
gives the identity functor on $\D^b(Y,\CA_Y)$,
see details in section~\ref{uhps}.
\item[Step 2]
We prove that the functor
$\Phi_r:\D^b(\CY_r,\CA_Y) \to \D^b(\CX_r)$ is fully faithful
for all $r\le i$ by induction in $r$ and that the collection
$$
\langle \D^b(\CY_r,\CA_Y),
\D^b(\Gr(r,V^*))\otimes E_1(1),
\dots,
\D^b(\Gr(r,V^*))\otimes E_2(i-r)\rangle
$$
is semiorthogonal in $\D^b(\CX_r)$,
see details in subsection~\ref{sodlsx}.
\item[Step 3]
We use the Bridgeland's trick in relative situation.
Note that the relative canonical bundle of $\CX_i$ over $\Gr(i,V^*)$
is given by $\omega_{\CX_i/\Gr(i,V^*)} \cong \det\CL_i^*$,
where $\CL_i$ is the tautological rank $i$ subbundle on $\Gr(i,V^*)$.
Thus $\CX_i$ is a family of Calabi--Yau varieties over $\Gr(i,V^*)$.
On the other hand, the functor $\Phi_i$ is $\Gr(i,V^*)$-linear
(i.e.\ commutes with tensoring by pullbacks of bundles on $\Gr(i,V^*)$).
We show in subsection~\ref{ssbt} that a relative version
of the Bridgeland's trick (proposition~\ref{ffisequ}) works in this case,
so $\Phi_i$ is an equivalence, and
$\omega_{\CY_i/\Gr(i,V^*)} \cong \det\CL_i^*$,
see details in subsection~\ref{fullness}.
\item[Step 4]
We use descending induction in $r$ to check that
$\omega_{\CY_r/\Gr(r,V^*)} \cong \det\CL_r^*(i-r)$,
where $\CL_r$ is the tautological rank $r$ subbundle on $\Gr(r,V^*)$,
see details in subsection~\ref{fullness}.
\item[Step 5]
We use descending induction in $r$ to check that
the above semiorthogonal collection in $\D^b(\CX_r)$ is full.
In other words, we prove that we have a semiorthogonal decomposition
$$
\D^b(\CX_r)  =
\langle \D^b(\CY_r,\CA_Y),
\D^b(\Gr(r,V^*))\otimes E_1(1),
\dots,
\D^b(\Gr(r,V^*))\otimes E_2(i-r)\rangle,
$$
for all $r\le i$, see details in subsection~\ref{fullness}.
\item[Step 6]
By induction in $r$ we show that for all $r\le i$ the functors
$\Phi_r^!:\D^b(\CX_r) \to \D^b(\CY_r)$ take
the exceptional pair $(E_0(1),E_1(1))$ to the pair $(F_2^*,F_1^*)$
up to a twist and a shift, where $E_0$ is the mutation of $E_2$
through $E_1$, (that is $E_0$ is defined from the exact triangle
$E_0 \to W^*\otimes E_1 \exto{\ev} E_2$).
It follows that
$\langle \D^b(\Gr(r,V^*))\otimes F_2^*(i-r),\dots,
\D^b(\Gr(r,V^*))\otimes F_1^*(-1)\rangle$
is a semiorthogonal collection in $\D^b(\CY_r,\CA_Y)$
for all $r\ge i$. In particular, it follows that
$(F_2^*(i-N),F_1^*(i-N),\dots,F_2^*(-1),F_1^*(-1))$
is an exceptional collection in $\D^b(Y,\CA_Y)$ (if $\BZ=\emptyset$),
see details in subsection~\ref{ecfory}.
\item[Step 7]
We use ascending induction in $r$  to check that
$$
\D^b(\CY_r,\CA_Y) =
\langle \D^b(\CX_r),\D^b(\Gr(r,V^*))\otimes F_2^*(i-r),\dots,
\D^b(\Gr(r,V^*))\otimes F_1^*(-1)\rangle
$$
for all $r\ge i$,
see details in subsection~\ref{sodlsy}.
\end{description}

The inductions steps are based on the following construction.
Consider the partial flag variety $\Fl(r-1,r;V^*)$
and the tautological subbundles
$\CL_{r-1} \subset \CL_r \subset V^*\otimes\CO_{\Fl(r-1,r;V^*)}$
of rank $(r-1)$ and $r$ on it. Let
$$
\begin{array}{ll}
\CX_{\CL_{r-1}} = \CX_{r-1}\times_{\Gr(r-1,V^*)}\Fl(r-1,r;V^*), \qquad &
\CX_{\CL_{r  }} = \CX_{r  }\times_{\Gr(r  ,V^*)}\Fl(r-1,r;V^*),\\
\CY_{\CL_{r-1}} = \CY_{r-1}\times_{\Gr(r-1,V^*)}\Fl(r-1,r;V^*),&
\CY_{\CL_{r  }} = \CY_{r  }\times_{\Gr(r  ,V^*)}\Fl(r-1,r;V^*).
\end{array}
$$
Then $\CX_{\CL_r}$ is a divisor in $\CX_{\CL_{r-1}}$,
$\CY_{\CL_{r-1}}$ is a divisor in $\CY_{\CL_r}$,
and we have the following diagrams
$$
\vcenter{\xymatrix{
\CY_{\CL_r} \ar[d]_\psi &
\CY_{\CL_{r-1}} \ar[d]^\phi \ar[l]_\eta \\
\CY_r &
\CY_{r-1}
}}
\quad\text{and}\quad
\vcenter{\xymatrix{
\CX_{\CL_r} \ar[d]_\psi \ar[r]^\xi &
\CX_{\CL_{r-1}} \ar[d]^\phi  \\
\CX_r &
\CX_{r-1}
}}
$$
where morphisms $\xi$ and $\eta$ are divisorial embeddings, while
$\phi$ and $\psi$ are the projections. Moreover, the pullbacks
$\CE_{\CL_{r-1}}$ and $\CE_{\CL_r}$ of $\CE_{r-1}$ and $\CE_r$
via the projections
$\CX_{\CL_{r-1}}\times_{\Fl(r-1,r;V^*)}\CY_{\CL_{r-1}} \to
\CX_{r-1}\times_{\Gr(r-1,V^*)}\CY_{r-1}$ and
$\CX_{\CL_r}\times_{\Fl(r-1,r;V^*)}\CY_{\CL_r} \to
\CX_r\times_{\Gr(r,V^*)}\CY_r$
provide us with the functors $\Phi_{\CL_{r-1}}$ and $\Phi_{\CL_r}$
between the corresponding derived categories.
Since all these functors originate from the same object $\CE$ on $Q(X,Y)$,
we can find a certain relations between them,
see details in subsection~\ref{indsetup}.

Finally, when steps 1--7 are performed, we can deduce the theorem.
For the first claim we consider $r=1$ and note that
$\CY_1 = Y$. Hence we have a semiorthogonal decomposition
$$
\D^b(\CX_1) = \langle \D^b(Y,\CA_Y),
\D^b(\PP(V^*))\otimes E_1(1),\dots,\D^b(\PP(V^*))\otimes E_2(i-1)\rangle,
$$
and an isomorphism
$\omega_{Y/\PP(V^*)} \cong \CL_1^*(i-1) \cong \CO_Y(i)$.
This immediately implies
$$
\omega_Y \cong
\omega_{Y/\PP(V^*)} \otimes \omega_{\PP(V^*)} \cong
\CO_Y(i-N).
$$
On the other hand, since $\CX_1$ is smooth we deduce that $\D^b(Y,\CA_Y)$
being a semiorthogonal component of an $\Ext$-bounded category $\D^b(\CX_1)$
is $\Ext$-bounded, hence $(Y,\CA_Y)$ is smooth.

Finally, the second claim of the theorem is deduced by a base change argument
for the base change $\Spec\kk \to \Gr(r,V^*)$ corresponding to
an admissible $r$-dimensional subspace $L\subset V^*$, using
the faithful base change theorem~\ref{phitsod}.

\bigskip

The next question we address is how to construct $Y$ starting from $X$.
More precisely,
assume that we have data (D.1) and (D.2) satisfying conditions (C.1)--(C.3).
How to construct data (D.3)--(D.5) such that conditions (C.4)--(C.8)
are satisfied?
We propose the following approach.
Since we must have a morphism $\phi_1:W\otimes F_1\to F_2$ on $Y$,
it is natural to construct $Y$ as a subscheme in the moduli space
of representations of the quiver
$\SQ = \xymatrix@1{ \bullet \ar[r]^{W} & \bullet }$.
The condition (C.6) suggests that the dimension vector
$(d_1,d_2)$ of the representations must satisfy
$$
\theta(d_1,d_2) = d_2\rank(E_2) - d_1\rank(E_1) = 0.
$$
Note that by this condition the dimension vector is defined
only up to a multiplicative constant.
Let $\BR_m$ denote the representation space of the quiver $\SQ$
for the dimension vector $(d_1,d_2)$ satisfying the above condition
with $\gcd(d_1,d_2) = m$, and let $(\CR_1,\CR_2)$ be the tautological
family of representations.
Then (C.6) suggests to consider the following GIT-quotient
$$
\BY_m = \{(\rho,H) \in \BR_m\times\PP(V^*)\ |\
\supp\Coker(E_1\otimes\CR_{1\rho} \exto{e_\rho} E_2\otimes\CR_{2\rho}) =
X\cap H\}\git G,
$$
where $G = \GL(d_1)\times\GL(d_2)/\kk^*$ and
$\chi(g_1,g_2) = \det(g_1)^{-\rank(E_1)}\det(g_2)^{\rank(E_2)}$
is a character of $G$. Note that to obtain a universal family
of quiver representations on $\BY_m$ we have to consider
a sheaf of Azumaya algebras on $\BY_m$. Moreover, if we
want (C.5) to be satisfied we have to restrict to the stable locus,
thus replacing $\BY_m$ with $Y_m = g_m^{-1}(\PP(V^*) \setminus \OZ_m)$,
where $g_m:\BY_m \to \PP(V^*)$ is the canonical projection and
$\OZ_m = \cup_{2k\le m}g_k(\BY_k)$. Then under some mild additional
assumptions we can show that (C.4)--(C.6) are indeed satisfied for $Y_m$.
It remains to choose $m$ in such a way that (C.7) holds and
to hope that (C.8) would be true for this $m$.
This choice and verification of (C.7) are the most cumbersome
part of the work and must be done for each case separately.

\section{Homological background}\label{sect1}

This section is designed to develop some machinery for working
with derived categories in relative situation. The main results
are the faithful base change theorem~\ref{phitsod} and a relative
version of the Bridgeland's trick~\ref{ffisequ}. Since we need
these results to be applicable to the derived categories of
sheaves of modules over sheaves of Azumaya algebras, we must work in
the corresponding category.

We define an {\sf Azumaya variety}\/ $(X,\CA_X)$ as an embeddable algebraic variety
of finite type $X$ over a field equipped with a sheaf~$\CA_X$
of semisimple $\CO_X$-algebras, such that $\CA_X$ is locally
free over $\CO_X$. A~morphism of Azumaya varieties
$f:(X,\CA_X) \to (Y,\CA_Y)$ consists of a morphism $f_\circ:X\to Y$
of the underlying algebraic varieties and a homomorphism
of $\CO_X$-algebras $f_\CA:f_\circ^*\CA_Y \to \CA_X$.
Thus we treat Azumaya varieties as ringed spaces with $\CA_X$
playing the role of a sheaf of rings.
Azumaya varieties form a category, which contains the category of usual
embeddable algebraic varieties as a full subcategory
(the embedding functor takes $X$ to $(X,\CO_X)$).

Below we will work in the category of Azumaya varieties.
To unburden the notation we will often write $X$ instead of $(X,\CA_X)$
when it is clear which sheaf of algebras is used, or when it doesn't matter.
The necessary background on Azumaya varieties is contained in appendix~D.
We define there the functors $f_*$, $f^*$, $f^!$, $\otimes_{\CA_X}$ and
$\RCHom_{\CA_X}$, and check that all usual relations between them
are satisfied. Here we shall only mention that
a morphism $f$ is called {\sf strict}\/ if $f_\CA$ is an isomorphism,
and that for a strict morphism $f$ we have
$f^* = f_\circ^*$ and $f^! = f_\circ^!$.

The base field $\kk$ is assumed to be an algebraically closed field
of zero characteristic.

\subsection{Kernel functors}\label{sskf}

We denote by $\D_{qc}^b(X)$, $\D_{qc}^-(X)$, $\D_{qc}^+(X)$ and $\D_{qc}(X)$
the bounded, the bounded above, the bounded below and the unbounded
derived categories of quasicoherent sheaves on an Azumaya variety $X$.
Further, $\D(X)$ stands for the unbounded derived category
of quasicoherent sheaves on $X$ with coherent cohomologies,
and similarly for $\D^b(X)$, $\D^-(X)$ and $\D^+(X)$.
Finally, $\D^\perf(X)$ stands for the category of perfect complexes.

Whenever we wish to emphasize the role of an Azumaya algebra on $X$,
or want to point out it explicitly, we use the notation $\D(X,\CA_X)$ e.t.c.

For an object $F\in\D(X)$ we denote by $\CH^i(F)$ the $i$-th cohomology
sheaf of $F$. For any integers $a\le b$ we denote by $\D^{[a,b]}(X)$
the full subcategory of $\D(X)$ formed by objects $F$ such that
$\CH^i(F) = 0$ for $i\not\in[a,b]$. Similarly, we define
$\D^{\ge a}(X)$ and $\D^{\le a}(X)$.
For a morphism $f:X\to Y$ we denote by $f_*$, $f^*$ and $f^!$
the derived pushforward, the derived pullback and the twisted pullback functors.
Similarly, $\otimes_{\CA_X}$ stands for the derived tensor product
and $\RCHom_{\CA_X}$ stands for the derived local $\Hom$-functor.

\begin{lemma}\label{dual1}
If $f:X\to Y$ is a morphism of Azumaya varieties, $F\in\D^-(X)$, $G\in\D^+(Y)$
and $\supp(F)$ is projective over $Y$, then
$f_*\RCHom_{\CA_X}(F,f^!G) \cong \RCHom_{\CA_Y}(f_*F,G)$.
\end{lemma}
\begin{proof}
Let $Z$ be a scheme-theoretical support of $F$ and let
$i:Z\to X$ denote the corresponding closed embedding,
so that $F = i_*F'$, where $F'\in\D^-(Z)$. Then both
$i$ and $f\circ i:Z\to Y$ are projective. Using
the functoriality of the twisted pullback and
the usual duality theorem for $i$ and $f\circ i$
(see~\cite{H} and lemma~\ref{duality}) we deduce
\begin{multline*}
f_*\RCHom_{\CA_X}(F,f^!G) =
f_*\RCHom_{\CA_X}(i_*F',f^!G) \cong
f_*i_*\RCHom_{\CA_{X|Z}}(F',i^!f^!G) \cong \\ \cong
(f\circ i)_*\RCHom_{\CA_{X|Z}}(F',(f\circ i)^!G) \cong
\RCHom_{\CA_Y}((f\circ i)_*F',G) \cong
\RCHom_{\CA_Y}(f_*F,G).
\end{multline*}
\end{proof}

Let $X_1$, $X_2$ be Azumaya varieties and
let $p_i:X_1\times X_2 \to X_i$ denote the projections.
Take any $K\in\D_{qc}^-(X_1\times X_2,\CA^\opp_{X_1}\boxtimes\CA_{X_2})$
and define functors
$$
\Phi_K(F_1) := {p_2}_*(p_{1\circ}^*F_1\otimes_{\CA_{X_1}} K),\qquad
\Phi_K^!(F_2) := {p_1}_*\RCHom_{\CA_{X_2}}(K,p_{2\circ}^!F_2).
$$
Then $\Phi_K$ is an exact functor $\D_{qc}^-(X_1)\to\D_{qc}^-(X_2)$ and
$\Phi_K^!$ is an exact functor $\D_{qc}^+(X_2)\to\D_{qc}^+(X_1)$.
We call $\Phi_K$ the {\sf kernel functor}\/ with kernel $K$,
and $\Phi_K^!$ the {\sf kernel functor of the second type}\/ with kernel $K$.

\begin{lemma}\label{keronsm}
If $X_1$ is smooth and $K\in\D^\perf(X_1\times X_2)$ then the functor
$\Phi_K^!$ is isomorphic to the usual kernel functor with kernel
$\RCHom_{\CA_{X_2}}(K,\omega_{X_1}\boxtimes\CA_{X_2})[\dim X_1]$.
\end{lemma}
\begin{proof}
Follows from~\ref{fshriek}, \ref{f22} and \ref{fdg}.
\end{proof}

\begin{lemma}\label{mkf}
Any morphism of kernels $\phi:K \to K'$ induces natural morphisms
of kernel functors $\phi_*:\Phi_K \to \Phi_{K'}$,
$\phi^!:\Phi_{K'}^! \to \Phi_K^!$. If $\phi$ is an isomorphism,
then both $\phi_*$ and $\phi^!$ are isomorphisms.
\end{lemma}
\begin{proof}
Evident.
\end{proof}

\begin{lemma}\label{phi_bounded}
$(i)$
If $K$ has coherent cohomologies, finite $\Tor$-amplitude over $X_1$ and
$\supp(K)$ is projective over $X_2$ then $\Phi_K$ takes
$\D^b(X_1)$ to $\D^b(X_2)$.

\noindent $(ii)$
If $K$ has coherent cohomologies, finite $\Ext$-amplitude over $X_2$ and
$\supp(K)$ is projective over $X_1$ then $\Phi_K^!$
takes $\D^b(X_2)$ to $\D^b(X_1)$.

\noindent $(iii)$
If both $(i)$ and $(ii)$ hold then $\Phi_K^!$ is right adjoint to~$\Phi_K$.
Moreover, $\Phi_K$ takes $\D^\perf(X_1)$ to $\D^\perf(X_2)$.
\end{lemma}
\begin{proof}
$(i)$ If $K$ has finite $\Tor$-amplitude over $X_1$ then
$q_{1\circ}^*F_1\otimes_{\CA_{X_1}} {K} \cong
q_1^*F_1\otimes_{\CA_{X_1}\otimes\CA_{X_2}} {K}$
is bounded for any $F_1\in\D^b(X_1)$, and its support
is projective over $X_2$. Therefore $\Phi_K(F_1)$ is bounded
and has coherent cohomologies.

$(ii)$ If $K$ has finite $\Ext$-amplitude over $X_2$ then
$\RCHom_{\CA_{X_2}}({K},q_{2\circ}^!F_2) \cong
\RCHom_{\CA_{X_1}^\opp\otimes\CA_{X_2}}({K},q_2^!F_2)$ is bounded
for any $F_2\in\D^b(X_2)$, and its support is projective over $X_1$.
Therefore $\Phi_K^!(F_2)$ is bounded and has coherent cohomologies.

$(iii)$ If both $(i)$ and $(ii)$ hold then using lemma~\ref{dual1} we deduce
\begin{multline*}
\Hom_{\CA_{X_1}}(F_1,\Phi_K^!(F_2)) =
\Hom_{\CA_{X_1}}(F_1,{p_1}_*\RCHom_{\CA_{X_2}}({K},p_{2\circ}^!F_2))
\cong
\\
\cong
\Hom_{\CA_{X_1}}(p_{1\circ}^*F_1,\RCHom_{\CA_{X_2}}({K},p_{2\circ}^!F_2)) \cong
\Hom_{\CA_{X_2}}(p_{1\circ}^*F_1\otimes_{\CA_{X_1}} {K},p_{2\circ}^!F_2) \cong
\\ \cong
\Hom_{\CA_{X_2}}({p_2}_*(p_{1\circ}^*F_1\otimes_{\CA_{X_1}} {K}),F_2) =
\Hom_{\CA_{X_2}}(\Phi_K(F_1),F_2)
\end{multline*}
for all $F_1\in\D^b(X_1)$, $F_2\in\D^b(X_2)$.
Moreover, the arguments in $(ii)$ show that $\Phi_K^!$
has bounded cohomological amplitude, hence adjointness implies
that $\Phi_K$ takes perfect complexes to complexes of finite
$\Ext$-amplitude, that is to perfect complexes (see~\ref{ftd_perf}).
\end{proof}

\begin{lemma}\label{ladj}
If $K$ is a perfect complex, $X_2$ is smooth and $\supp(K)$
is projective both over $X_1$ and over $X_2$, then the kernel functor
$\Phi_{K^\#}:\D^b(X_2) \to \D^b(X_1)$, where
$$
K^\#:= \RCHom_{\CA_{X_1}}(K,\CA_{X_1}\boxtimes\omega_{X_2}[\dim X_2]),
$$
is the left adjoint functor to the kernel functor
$\Phi_K:\D^b(X_1) \to \D^b(X_2)$.
\end{lemma}
\begin{proof}
Indeed, $\supp(K^\#) = \supp(K)$ is projective over $X_1$, hence
by lemma~\ref{dual1} we have
\begin{multline*}
\RHom_{\CA_{X_1}}(\Phi_{K^\#}(F_2),F_1) =
\RHom_{\CA_{X_1}}({p_1}_*(p_{2\circ}^*F_2\otimes_{\CA_{X_2}} K^\#),F_1) \cong
\\ \cong
\RHom_{\CA_{X_1}}(p_{2\circ}^*F_2\otimes_{\CA_{X_2}} K^\#,p_{1\circ}^!F_1) \cong
\RHom_{\CA_{X_2}}(p_{2\circ}^*F_2,
p_{1\circ}^!F_1\otimes_{\CA_{X_1}} K^{\#*}_{\CA_{X_1}}) \cong
\\ \cong
\RHom_{\CA_{X_2}}(F_2,
{p_2}_*(p_{1\circ}^!F_1\otimes_{\CA_{X_1}} K^{\#*}_{\CA_{X_1}})).
\end{multline*}
On the other hand,
$K^{\#*}_{\CA_{X_1}} \cong
K\otimes_{\CO_{X_1\times X_2}}\omega_{X_2}^{-1}[-\dim X_2]$,
and
$p_{1\circ}^!F_1\otimes_{\CO_{X_1\times X_2}}\omega_{X_2}^{-1}[-\dim X_2]
\cong p_{1\circ}^*F_1$,
so the RHS equals to $\RHom_{\CA_{X_2}}(F_2,\Phi_K(F_1))$.
\end{proof}

Consider kernels
$K_{12} \in \D^-(X_1\times X_2,\CA_{X_1}^\opp\boxtimes\CA_{X_2})$,
$K_{23} \in \D^-(X_2\times X_3,\CA_{X_2}^\opp\boxtimes\CA_{X_3})$.
Denote by $p_{ij}:X_1\times X_2\times X_3 \to X_i\times X_j$
the projections.
We define the {\sf convolution}\/ of kernels as follows
$$
K_{23}\circ K_{12} :=
{p_{13}}_*(p_{12\circ}^*K_{12}\otimes_{\CA_{X_2}} p_{23\circ}^*K_{23}),
$$
Similarly,
if $K_{12} \in \D^b(X_1\times X_2,\CA_{X_1}^\opp\boxtimes\CA_{X_2})$ and
$K_{32} \in \D^-(X_3\times X_2,\CA_{X_2}\boxtimes\CA^\opp_{X_3})$
we define the {\sf convolution of the second type} of kernels:
$$
K_{32}\star K_{12} :=
{p_{13}}_*\RCHom_{\CA_{X_2}}(p_{23\circ}^*K_{32},p_{12\circ}^!K_{12}).
$$

\begin{lemma}\label{kerconv}
For $K_{12} \in \D^-(X_1\times X_2)$,
$K_{23} \in \D^-(X_2\times X_3)$,
$K_{32} \in \D^-(X_3\times X_2)$
we have\\
$(i)$
$\Phi_{K_{23}} \circ \Phi_{K_{12}} = \Phi_{K_{23}\circ K_{12}}$,\\
$(ii)$
$\Phi_{K_{12}}^! \circ \Phi_{K_{23}}^! = \Phi_{K_{23}\circ K_{12}}^!$, \quad
if $\supp(p_{12\circ}^*K_{12}\otimes_{\CA_{X_2}} p_{23\circ}^*K_{23})$
is projective over $X_1\times X_3$;\\
$(iii)$
$(\Phi_{K_{32}}^! \circ \Phi_{K_{12}})_{|\D^\perf(X_1)} =
{\Phi_{K_{32}\star K_{12}}}$, \quad
if $K_{12} \in \D^b(X_1\times X_2)$.
\end{lemma}
\begin{proof}
The first is standard. For the second we have
\begin{multline*}
(\Phi_{K_{12}}^! \circ \Phi_{K_{23}}^!)(F_3) =
{p_1}_*\RCHom_{\CA_{X_2}}
(K_{12},p_{2\circ}^!{p_2}_*\RCHom_{\CA_{X_3}}(K_{23},p_{3\circ}^!F_3)) \cong
\\ \cong
{p_1}_*\RCHom_{\CA_{X_2}}
(K_{12},{p_{12}}_*p_{23\circ}^!\RCHom_{\CA_{X_3}}(K_{23},p_{3\circ}^!F_3))
\cong \\ \cong
{p_1}_*{p_{12}}_*\RCHom_{\CA_{X_2}}(p_{12\circ}^*K_{12},
\RCHom_{\CA_{X_3}}(p_{23\circ}^*K_{23},p_{23\circ}^!p_{3\circ}^!F_3))
\cong \\ \cong
{p_1}_*{p_{13}}_*\RCHom_{\CA_{X_3}}
(p_{12\circ}^*K_{12}\otimes_{\CA_{X_2}} p_{23\circ}^*K_{23},
p_{13\circ}^!p_{3\circ}^!F_3))
\cong \\ \cong
{p_1}_*\RCHom_{\CA_{X_3}}
({p_{13}}_*(p_{12\circ}^*K_{12}\otimes_{\CA_{X_2}} p_{23\circ}^*K_{23}),
p_{3\circ}^!F_3)) =
\Phi_{K_{23}\circ K_{12}}^!(F_3).
\end{multline*}
Similarly, for the third, if $F_1$ is a perfect complex then we have
\begin{multline*}
(\Phi_{K_{32}}^! \circ \Phi_{K_{12}})(F_1) =
{p_3}_*\RCHom_{\CA_{X_2}}(K_{32},
p_{2\circ}^!{p_2}_*(p_{1\circ}^*F_1\otimes_{\CA_{X_1}} K_{12}))
\cong \\ \cong
{p_3}_*\RCHom_{\CA_{X_2}}(K_{32},
{p_{23}}_*p_{12\circ}^!(p_{1\circ}^*F_1\otimes_{\CA_{X_1}} K_{12}))
\cong \\ \cong
{p_3}_*{p_{23}}_*\RCHom_{\CA_{X_2}}(p_{23\circ}^*K_{32},
p_{12\circ}^*p_{1\circ}^*F_1\otimes_{\CA_{X_1}} p_{12\circ}^!K_{12}))
\cong \\ \cong
{p_3}_*{p_{13}}_*\RCHom_{\CA_{X_2}}(p_{23\circ}^*K_{32},
p_{13\circ}^*p_{1\circ}^*F_1\otimes_{\CA_{X_1}} p_{12\circ}^!K_{12}))
\cong \\ \cong
{p_3}_*(p_{1\circ}^*F_1\otimes_{\CA_{X_1}}
{p_{13}}_*\RCHom_{\CA_{X_2}}(p_{23\circ}^*K_{32},p_{12\circ}^!K_{12}))
\cong \\ \cong
{p_3}_*(p_{1\circ}^*F_1\otimes_{\CA_{X_1}} (K_{32}\star K_{12}))) =
\Phi_{K_{32}\star K_{12}}(F_1).
\end{multline*}
\end{proof}

Assume that $\Phi_1,\Phi_2,\Phi_3:\D \to \D'$ are exact functors
between triangulated categories, and $\alpha:\Phi_1\to\Phi_2$,
$\beta:\Phi_2\to\Phi_3$, $\gamma:\Phi_3\to\Phi_1[1]$ are morphisms
of functors. We say that
$$
\Phi_1 \exto{\alpha} \Phi_2 \exto{\beta} \Phi_3 \exto{\gamma} \Phi_1[1]
$$
is an {\sf exact triangle of functors}, if for any object $F\in\D$
the triangle
$$
\Phi_1(F) \exto{\alpha(F)} \Phi_2(F) \exto{\beta(F)}
\Phi_3(F) \exto{\gamma(F)} \Phi_1(F)[1]
$$
is exact in $\D'$.

\begin{lemma}\label{etf}
If  $K_1 \exto{\alpha} K_2 \exto{\beta} K_3 \exto{\gamma} K_1[1]$ is
an exact triangle in $\D^-(X\times Y)$ then we have the following
exact triangles of functors
$$
\Phi_{K_1} \exto{\alpha_*} \Phi_{K_2} \exto{\beta_*}
\Phi_{K_3} \exto{\gamma_*} \Phi_{K_1}[1]
$$
$$
\Phi^!_{K_3} \exto{\beta^!} \Phi^1_{K_2} \exto{\alpha^!}
\Phi^!_{K_1} \exto{\gamma^!} \Phi^!_{K_3}[1]
$$
\end{lemma}
\begin{proof}
Evident.
\end{proof}

\begin{lemma}\label{pff}
Let $\alpha:X\to Y$ be a finite morphism.

\noindent $a)$
If $\D$ is a category and $\Phi:\D \to \D(X)$ is a functor such that $\alpha_*\circ\Phi = 0$
then $\Phi = 0$.

\noindent $b)$
If $X'$ is another variety, $K,K'\in\D^-(X'\times X)$ are kernels
and $\phi:K\to K'$ is a morphism, such that morphism of functors
$\alpha_*\circ\Phi_K \exto{\alpha_*(\phi_*)} \alpha_*\circ\Phi_{K'}$
is an isomorphism, then $\phi_*:\Phi_K \to \Phi_{K'}$ is an isomorphism.
\end{lemma}
\begin{proof}
$a)$ Since $\alpha$ is finite we have $\CH^i(\alpha_*G) \cong \alpha_*\CH^i(G)$
for any $G \in \D(X)$. Therefore $\alpha_*(G) = 0$ implies
$\alpha_*(\CH^i(G)) = 0$ for all $i$, hence $G = 0$ in $\D(X)$.
In particular, $\alpha_*\Phi(F) = 0$ implies $\Phi(F) = 0$ for all $F\in\D$.

$b)$ Let $K''$ be the cone of $\phi:K\to K'$ in $\D^-(X'\times X)$.
Then we have an exact triangle of functors
$\Phi_K \to \Phi_{K'} \to \Phi_{K''}$ and it follows that
$\alpha_*\circ\Phi_{K''} = 0$. Therefore, $\Phi_{K''} = 0$
by part (a), hence $\phi_*$ is an isomorphism.
\end{proof}

\subsection{Perfect spanning classes}\label{sspsc}

\begin{definition}
A class of objects $\CF\subset \D(X,\CA_X)$ is called
{\sf a perfect spanning class for $(X,\CA_X)$}, if for any point
$x\in X$ there exists an object $0 \ne F_x\in\CF$ such that
\begin{enumerate}
\item $F_x$ is a perfect complex;
\item $\CH^p(F_x)$ is supported set theoretically at $x$ for all $p$;
\item if $p_0 = \max\{p\ |\ \CH^p(F_x)\ne 0\}$ then
$\CH^{p_0}(F_x) \cong \CA_X \otimes_{\CO_X} \CO_x$
\end{enumerate}
\end{definition}

\begin{lemma}
Any Azumaya variety admits a perfect spanning class.
\end{lemma}
\begin{proof}
Choose a closed embedding $i:X \to Y$ with smooth $Y$ and take
$F_x := \CA_X\otimes_{\CO_X}i^*i_*\CO_x$.
\end{proof}

\begin{lemma}
If $\CF$ is a perfect spanning class for $X$ and
$\CG$ is a perfect spanning class for $Y$ then
$\CF\boxtimes\CG$ is a perfect spanning class for $X\times Y$.
\end{lemma}
\begin{proof}
Evident.
\end{proof}

\begin{lemma}\label{psc_fp}
Assume that $f:X\to S$ and $g:Y\to S$ are morphisms of algebraic varieties.
If $\CF$ is a perfect spanning class in $\D(X,\CA_X)$ and
$\CG$ is a perfect spanning class in $\D(Y,\CO_Y)$ then
$p^*\CF\otimes_{\CO_{X\times_S Y}} q_\circ^*\CG$ is a perfect spanning class
in $\D(X\times_S Y,p_\circ^*\CA_X)$,
where $p:X\times_S Y \to X$ and $q:X\times_S Y \to Y$ are the projections.
\end{lemma}
\begin{proof}
Evident.
\end{proof}

\begin{lemma}\label{psc_ort}
If $K \in \D^-(X,\CA_X^\opp)$, $\CF$ is a perfect spanning class
for $(X,\CA_X)$ and for all $F\in\CF$ we have
$H^\bullet(X,F\otimes_{\CA_X} K) = 0$,
then $K = 0$.
\end{lemma}
\begin{proof}
Assume that $K \ne 0$ and let $s$ be the maximal
integer such that $\CH^s(K) \ne 0$. Choose a point $x\in\supp\CH^s(K)$,
and take an object $F_x\in\CF$, corresponding to this point.
Let $t$ be the maximal integer such that $\CH^t(F_x) \ne 0$.
Then it is clear that
$\CH^p(F_x\otimes_{\CA_X} K)$ is supported at $x$ for all $p$,
$\CH^{>s+t}(F_x\otimes_{\CA_X} K) = 0$, and
$$
\CH^{s+t}(F_x\otimes_{\CA_X} K) \cong
\CH^0((\CA_X\otimes_{\CO_X}\CO_x)\otimes_{\CA_X}\CH^s(K)) \cong
\CH^0(\CH^s(K)\otimes_{\CO_X}\CO_x) \ne 0.
$$
It follows that the hypercohomology spectral sequence for
$F_x\otimes_{\CA_X} K$
$$
E_2^{p,q} =
H^q(X,\CH^p(F_x\otimes_{\CA_X} K)) \RA
H^{p+q}(X,F_x\otimes_{\CA_X} K)
$$
degenerates in the second term, and
$H^{s+t}(X,F_x\otimes_{\CA_X} K) =
H^0(X,\CH^{s+t}(F_x\otimes_{\CA_X} K)) \ne 0$,
which is a contradiction.
\end{proof}

\begin{proposition}\label{kfz}
If $K\in\D^-(X\times Y,\CA_X^\opp\boxtimes\CA_Y)$ and the functor
$\Phi_K:\D(X,\CA_X) \to \D(Y,\CA_Y)$ equals zero on a perfect
spanning class $\CF$ for $(X,\CA_X)$ then $K = 0$.
\end{proposition}
\begin{proof}
Choose a perfect spanning class $\CG$ for $(Y,\CO_Y)$.
Then for all $F\in\CF$, $G\in\CG$ we have
$$
H^\bullet(X\times Y,(F\boxtimes G)\otimes_{\CA_X} K) \cong
H^\bullet(Y,q_*((p_\circ^*F\otimes_{\CA_X} K) \otimes_{\CO_{X\times Y}} q_\circ^*G))
\cong
H^\bullet(Y,\Phi_K(F)\otimes_{\CO_Y} G) = 0.
$$
Since $\CF\boxtimes\CG$ is a perfect spanning class for $(X\times Y,\CA_X)$,
we have $K = 0$.
\end{proof}


\begin{corollary}\label{mkfi}
If the morphism of kernel functors $\Phi_K \to \Phi_{K'}$
with $K,K'\in\D^-(X\times Y)$, induced by a morphism of
kernels $\phi:K\to K'$ is an isomorphism
on a perfect spanning class $\CF$ for $X$,
then $\varphi$ is an isomorphism.
\end{corollary}
\begin{proof}
Let $K''$ be the cone of $\phi:K\to K'$ in $\D^-(X\times Y)$.
Then we have an exact triangle of functors
$\Phi_K \to \Phi_{K'} \to \Phi_{K''}$ and it follows that
$\Phi_{K''}$ equals zero on the perfect spanning class $\CF$.
Therefore, $K'' = 0$, hence $\phi_*$ is an isomorphism.
\end{proof}

\subsection{Koszul complexes}\label{sskc}

If $(X,\CA_X)$ is an Azumaya variety, $\CV$ is a vector bundle on its
underlying algebraic variety $X$, $s\in\Gamma(X,\CV)$ is a section
(possible nonregular) of $\CV$, then we denote by $\Kosz_X(s)$
the Koszul complex of $s$ considered as an object
of the derived category $\D^b(X,\CA_X)$:
$$
\Kosz_X(s) :=
\{ 0 \to
\Lambda^{\text{\sf top}}\CV^*\otimes_{\CO_X}\CA_X \exto{s} \dots \exto{s}
\CV^*\otimes_{\CO_X}\CA_X \exto{s} \CA_X \to 0 \}
$$
with $\CA_X$ placed in degree $0$.
Assume that $X$ is Cohen--Macaulay.
Recall that a section $s$ is {\sf regular} if the codimension
of the zero locus $Z(s)\subset X$ of $s$ is equal to the rank of $\CV$.
It is well known that for regular section $s$ we have
$\Kosz_X(s) \cong i_*\CA_{Z(s)}$, where $i:Z(s) \to X$ is the embedding
and $\CA_{Z(s)} = \CA_{X|Z(s)}$.

\begin{lemma}\label{koszul}
Assume that $X$ is a Cohen--Macaulay variety.
Then for any section $s\in\Gamma(X,\CV)$ we have
$\Kosz_X(s) \in \D^{[\codim_X Z(s) - \rank \CV,0]}(X,\CA_X)$.
Moreover, $\Kosz_X(s)$ is supported scheme-theoretically on a infinitesimal
neighborhood of $Z(s) \subset X$ and
$\CH^0(\Kosz_X(s)) \cong i_*\CA_{Z(s)}$.
\end{lemma}
\begin{proof}
The first part of the claim is local. Locally, we can decompose
$\CV = \CV'\oplus \CV''$, so that the first component of
$s = (s',s'')$ is regular and $\rank \CV' = \codim_X Z(s)$.
Then denoting by $i':Z(s')\to X$ the embedding, we get
$$
\Kosz_X(s) \cong \Kosz_X(s')\otimes_{\CA_X}\Kosz_X(s'') \cong
i'_*\CA_{Z(s')}\otimes_{\CA_X}\Kosz_X(s'') \cong i'_*\Kosz_{Z(s')}({i'}^*s''),
$$
and it remains to note that we have
$\Kosz_{Z(s')}({i'}^*s'') \in \D^{[-\rank \CV'',0]}(Z(s'),\CA_{Z(s')})$
by definition, and
$\rank \CV'' = \rank \CV - \rank \CV' = \rank \CV - \codim_X Z(s)$.

The second part of the claim is evident since $\Kosz_X(s)$ is acyclic
on the complement of $Z(s)$, and the third part follows from
the definition of the sheaf of ideals of $Z(s)$ as the image
of the morphism $\CV^* \exto{s} \CO_X$.
\end{proof}

\subsection{Exact cartesian squares}\label{ssecs}

Let $f:(X,\CA_X) \to (S,\CA_S)$ and $g:(Y,\CA_Y) \to (S,\CA_S)$
be morphisms of Azumaya varieties. Assume that either $f$ or $g$ is strict.
In this case there is defined the fiber product
$(X\times_S Y,\CA_{X\times_S Y})$, see~\ref{fpxy}.
Consider the corresponding cartesian square
\begin{equation}\label{cs}
\vcenter{\xymatrix{
X\times_S Y  \ar[r]^-q \ar[d]_p         &       Y \ar[d]_g \\
X \ar[r]^f                              &       S
}}
\end{equation}

Let $\Gamma_f: X \to X\times S$ and $\Gamma^g:Y \to S\times Y$
denote the graphs of morphisms $f$ and $g$ respectively, and put
$K_f = {\Gamma_f}_*\CA_X$, $K^g = {\Gamma^g}_*\CA_Y$,
$K(f,g) = \CA_{X\times_S Y}$ and consider the corresponding
kernel functors. Then we have functorial isomorphisms
$$
\Phi_{K_f} \cong f_*,
\qquad
\Phi_{K^g} \cong g^*,
\qquad
\Phi_{K(f,g)} \cong q_*p^*.
$$

\begin{lemma}\label{fg_ker}
We have $K_f\circ K^g \in \D^{\le0}(X\times Y)$ and
$\CH^0(K_f\circ K^g) \cong K(f,g)$.
\end{lemma}
\begin{proof}
Let $p_{XS}:X\times S\times Y \to X\times S$,
$p_{SY}:X\times S\times Y \to S\times Y$
denote the projections, and consider the object
$K = p_{XS\circ}^* K_f \otimes_{\CA_S} p_{SY\circ}^* K^g$.
Since the pullback and the tensor product are left exact,
we have $K\in\D^{\le 0}$. Further, it is clear that
the cohomologies of $K$ are supported (set-theoretically) on the fiber product
$X\times_S Y = X\times_S S\times_S Y \subset X\times S\times Y$
and $\CH^0(K) \cong \CA_{X\times_S Y}$.
Since the projection $p_{XY}:X\times S\times Y \to X\times Y$
restricted to the infinitesimal neighborhood of $X\times_S Y$
in $X\times S\times Y$ is finite, we deduce that
$K_f \circ K^g = {p_{XY}}_*K \in \D^{\le 0}$, and
$\CH^0(K_f \circ K^g) \cong {p_{XY}}_*\CA_{X\times_S Y} \cong K(f,g)$.
\end{proof}

The canonical morphism of kernels $K_f\circ K^g \to \CH^0(K_f\circ K^g)$
induces a functorial morphism $g^*f_* \to q_*p^*$.

\begin{definition}
A cartesian square~(\ref{cs}) is called {\sf exact cartesian},
if the natural morphism of functors $g^*f_* \to q_*p^*$ is an isomorphism.
\end{definition}

\begin{proposition}\label{ec_crit}
A cartesian square~$(\ref{cs})$ is exact cartesian,
if and only if the canonical morphism
$K_f\circ K^g \to \CH^0(K_f\circ K^g) \cong K(f,g)$
is an isomorphism.
\end{proposition}
\begin{proof}
The ``only if'' part follows from corollary~\ref{mkfi} and
the ``if'' part from lemma~\ref{mkf}.
\end{proof}

Further, corollary~\ref{mkfi} implies

\begin{corollary}\label{ec_psc}
If the natural morphism of functors $g^*f_* \to q_*p^*$ is an isomorphism
on a perfect spanning class $\CF\subset\D(X,\CA_X)$, then
cartesian square~$(\ref{cs})$ is exact cartesian.
\end{corollary}

Another consequence of proposition~\ref{ec_crit} is

\begin{corollary}\label{ec_tr}
A cartesian square~$(\ref{cs})$ is exact cartesian, if and only if
the transposed square is exact cartesian.
\end{corollary}
\begin{proof}
The criterion of proposition~\ref{ec_crit} is symmetric.
\end{proof}

\begin{lemma}\label{cs0}
A fiber square~$(\ref{cs})$ is exact cartesian, if and only if
the underlying square of algebraic varieties is.
\end{lemma}
\begin{proof}
Note that the fiber square of Azumaya varieties was defined
only in a situation when one of morphisms $f$, $g$ is strict.
Since the exactness property is symmetric with respect
to the transposition of the square, we may assume
that $g$ is strict. Then $p$ is also strict,
and we have
$q_*p^* \cong q_*p_\circ^*$, and $g^*f_* \cong g_\circ^*f_*$.
It remains to note that the functorial morphism
$g_\circ^*f_* \to q_*p_\circ^*$ is obtained from
the functorial morphism $g^*f_* \to q_*p^*$
by forgetting the $\CA_Y$-module structure.
\end{proof}

\begin{corollary}\label{flat_ec}
If either $f$ or $g$ is flat then cartesian square~$(\ref{cs})$
is exact cartesian.
\end{corollary}
\begin{proof}
Since the exactness property is symmetric with respect
to the transposition of the square, we may assume that
$g$ is flat. Then the underlying square is exact by \cite{H}, II,
proposition~5.12.
\end{proof}

\begin{proposition}\label{sbpf}
If diagram~$(\ref{cs})$ is exact cartesian
then $f^!g_* \cong p_*q^!$.
\end{proposition}
\begin{proof}
Note that $f^! \cong \Phi_{K_f}^!$, $g_* \cong \Phi_{K^g}^!$,
$p_*q^! \cong \Phi_{K(f,g)}^!$ and apply lemma~\ref{kerconv}~(ii).
\end{proof}

\begin{lemma}\label{squ3}
Assume that the right square in the following diagram is exact cartesian
$$
\xymatrix{
X'' \ar[r]^-{v''} \ar[d]_{f''} &
X' \ar[r]^{v'} \ar[d]_{f'} &
X \ar[d]_f \\
S'' \ar[r]^{u''} &
S' \ar[r]^{u'} &
S
}
$$
and either $f$ or both $u'$ and $u''$ are strict.
Then the ambient square is exact cartesian if and only if the left square is.
\end{lemma}
\begin{proof}
Assume that the ambient square is exact and $f$ is strict.
We choose perfect spanning classes
$\CF\subset\D(X,\CO_X)$, $\CG\subset\D(S',\CA_{S'})$.
Then for all $F\in\CF$, $G\in\CG$ we have
\begin{multline*}
f''_*{v''}^*({v'_\circ}^*F\otimes_\CO{f'}^*G) \cong
f''_*({v''_\circ}^*{v'_\circ}^*F\otimes_\CO{v''}^*{f'}^*G) \cong \\ \cong
f''_*((v'v'')_\circ^*F\otimes_\CO{f''}^*{u''}^*G) \cong
f''_*(v'v'')_\circ^*F\otimes_\CO{u''}^*G \cong
(u'u'')_\circ^*f_*F\otimes_\CO{u''}^*G \cong
{u''}^*({u'_\circ}^*f_*F\otimes_\CO G) \cong \\ \cong
{u''}^*({f'}_*{v'_\circ}^*F\otimes_\CO G) \cong
{u''}^*{f'}_*({v'_\circ}^*F\otimes_\CO {f'}^*G).
\end{multline*}
It remains to note that objects of the type
${v'_\circ}^*F\otimes_\CO{f'}^*G$ form a perfect spanning class
in the derived category of $X' = X\times_S S'$,
and apply corollary~\ref{ec_psc}.

Assume that the ambient square is exact and both $u'$ and $u''$ are strict.
We choose perfect spanning classes
$\CF\subset\D(X,\CA_X)$, $\CG\subset\D(S',\CO_{S'})$.
Then for all $F\in\CF$, $G\in\CG$ we have
\begin{multline*}
f''_*{v''}^*({v'}^*F\otimes_\CO{f'_\circ}^*G) \cong
f''_*({v''}^*{v'}^*F\otimes_\CO{v''_\circ}^*{f'_\circ}^*G) \cong \\ \cong
f''_*((v'v'')^*F\otimes_\CO{f''_\circ}^*{u''_\circ}^*G) \cong
f''_*(v'v'')^*F\otimes_\CO{u''_\circ}^*G \cong
(u'u'')^*f_*F\otimes_\CO{u''_\circ}^*G \cong
{u''}^*({u'}^*f_*F\otimes_\CO G) \cong \\ \cong
{u''}^*({f'}_*{v'}^*F\otimes_\CO G) \cong
{u''}^*{f'}_*({v'}^*F\otimes_\CO {f'_\circ}^*G).
\end{multline*}
It remains to note that objects of the type
${v'}^*F\otimes_\CO{f'_\circ}^*G$ form a perfect spanning class
in the derived category of $X' = X\times_S S'$,
and apply corollary~\ref{ec_psc}.

On the other hand, if the left square is exact, then
$$
f''_*(v'v'')^* F \cong
f''_*{v''}^*({v'}^* F) \cong
{u''}^*(f'_*{v'}^*) F \cong
{u''}^*{u'}^*f_* F \cong
(u'u'')^*f_* F
$$
for any $F\in\D^b(X,\CA_X)$, hence the ambient square is exact.
\end{proof}

\begin{lemma}\label{gce}
If $g$ is finite and the canonical morphism
$$
f^*g_*\CA_Y \to p_*q^*\CA_Y \cong p_*\CA_{X\times_S Y}
$$
is an isomorphism, then square~$(\ref{cs})$ is exact cartesian.
\end{lemma}
\begin{proof}
Note that
$$
g_*g^*f_*(F) \cong
f_*(F)\otimes_{\CA_S} g_*\CA_Y \!\cong\!
f_*(F\otimes_{\CA_X} f^*g_*\CA_Y) \!\cong\!
f_*(F\otimes_{\CA_X} p_*\CA_{X\times_S Y}) \!\cong\!
f_*p_*p^*(F) \cong
g_*q_*p^*(F),
$$
and the resulting isomorphism coincides with the morphism,
obtained by application of the pushforward functor~$g_*$
to the functorial morphism $g^*f_* \to q_*p^*$.
Since the latter is induced by a morphism of kernels
$K_f\circ K^g \to K(f,g)$ it remains to apply lemma~\ref{pff}~b).
\end{proof}

\begin{corollary}\label{lci}
Assume that $g$ is a strict closed embedding,
$Y \subset S$ is a locally complete intersection,
and both $S$ and $X$ are Cohen--Macaulay.
If $\codim_X (X\times_S Y) = \codim_S Y$,
then square~$(\ref{cs})$ is exact cartesian.
\end{corollary}
\begin{proof}
By lemma~\ref{gce} it suffices to check that the canonical homomorphism
$f^*g_*\CA_Y \to p_*\CA_{X\times_S Y}$ is an isomorphism.
The claim is local with respect to $S$, so we may assume that
$Y$ is a zero locus of a regular section $s$ of a vector bundle
$\CV$ on $S$ of rank $\rank\CV = \codim_S Y$.
Then $g_*\CA_Y \cong \Kosz_S(s)$, hence
$f^*g_*\CA_Y \cong f^*\Kosz_S(s) \cong \Kosz_X(f^*s)$.
It is clear that the zero locus of $f^*s$ on $X$ is
the fiber product $X\times_S Y$, but
$\codim_X (X\times_S Y) = \codim_S Y = \rank\CV = \rank f^*\CV$
which implies that the section $f^*s$ is regular and we have an isomorphism
$\Kosz_X(f^*s) \cong p_* \CA_{X\times_S Y}$.
\end{proof}

\begin{corollary}\label{lci2}
Consider exact cartesian square~$(\ref{cs})$ and let
$\CU$, $\CV$, and $\CW$ be vector bundles on $S$, $X$, and $Y$ respectively,
with sections $u\in\Gamma(S,\CU)$, $v\in\Gamma(X,\CV)$, and $w\in\Gamma(Y,\CW)$.
Assume that the sections
$$
\begin{array}{rlrl}
p^*f^*u + p^*v + q^*w   & \in
\Gamma(X\times_S Y,p^*f^*\CU\oplus p^*\CV\oplus q^*\CW)
\qquad &
g^*u + w                & \in \Gamma(Y,g^*\CU\oplus\CW) \\
f^*u + v                & \in \Gamma(X,f^*\CU\oplus\CV) &
u                       & \in \Gamma(S,\CU) \\
\end{array}
$$
are regular. Denote by $S_u\subset S$, $X_{u,v}\subset X$,
$Y_{u,w}\subset Y$ and $Z_{u,v,w}\subset Z = X\times_S Y$
the corresponding zero loci. Then the square
$$
\xymatrix@=15pt{
Z_{u,v,w} \ar[r] \ar[d] & Y_{v,w} \ar[d] \\ X_{u,v} \ar[r] & S_u
}
$$
is exact cartesian.
\end{corollary}
\begin{proof}
Consider the diagrams
$$
1)\quad
\vcenter{\xymatrix@=15pt{
Z_{u,v} \ar[r] \ar[d] & Z \ar[r] \ar[d] & Y \ar[d] \\
X_{u,v} \ar[r] & X \ar[r] & S
}}
\qquad\qquad
2)\quad
\vcenter{\xymatrix@=15pt{
Z_{u,v,w} \ar[r] \ar[d] & Y_w \ar[d] \\
Z_{u,v} \ar[r] \ar[d] & Y \ar[d] \\
X_{u,v} \ar[r] & S
}}
\qquad\qquad
3)\quad
\vcenter{\xymatrix@=15pt{
Z_{u,v,w} \ar[r] \ar[d] & Y_{u,w} \ar[r] \ar[d] & Y_w \ar[d] \\
X_{u,v} \ar[r] & S_u \ar[r] & S
}}
$$
In the first diagram the right square is exact by assumption
and the left square is exact by corollary~\ref{lci}.
Hence the ambient square is exact by the ``if'' part of lemma~\ref{squ3}.
Thus the bottom square in the second diagram is exact.
On the other hand, its upper square is exact by corollary~\ref{lci}.
Hence its ambient square is exact by the ``if'' part of lemma~\ref{squ3}.
Thus the ambient square in the third diagram is exact.
On the other hand, its right square is exact by corollary~\ref{lci}.
Hence its left square is exact by the ``only if'' part of lemma~\ref{squ3}.
\end{proof}

\begin{definition}
A strict morphism $\phi:T \to S$ is called {\sf faithful} with respect
to a morphism $f:X \to S$ if the cartesian square
$$
\xymatrix@=20pt{X\times_S T \ar[r] \ar[d] &  X \ar[d]^f \\ T \ar[r]^\phi & S}
$$
is exact.
\end{definition}

\begin{lemma}\label{flatisf}
A strict and flat base change is faithful with respect to any morphism.
\end{lemma}
\begin{proof}
Use lemma~\ref{flat_ec}.
\end{proof}

\begin{corollary}\label{ec_fact}
Let $T\to S$ be a strict base change faithful with
respect to a morphism $f:X\to S$. If $T \to S' \to S$ is a factorization
of $T\to S$ such that $T\to S'$ is strict and finite and $S'\to S$ is
strict and smooth, then $S'\to S$ is faithful with respect to $f$ and
$T \to S'$ is faithful with respect to $f':X\times_S S' \to S'$.
\end{corollary}
\begin{proof}
The smooth base change is faithful by lemma~\ref{flatisf}
and the finite base change is faithful by the ``only if'' part
of the lemma~\ref{squ3}.
\end{proof}

\begin{lemma}\label{sbpb}
If square~$(\ref{cs})$ is exact cartesian,
$g$ has finite $\Tor$-dimension and strict, and $f$ is projective
then there exists a canonical isomorphism of functors
$$
p^*f^! \cong q^!g^* : \D^+(S) \to \D^+(X\times_S Y).
$$
\end{lemma}
\begin{proof}
First of all, we note that
$$
\Hom(p^*f^!(F),q^!g^*(F)) \cong
\Hom(q_*p^*f^!(F),g^*(F)) \cong
\Hom(g^*f_*f^!(F),g^*(F)),
$$
hence the canonical adjunction morphism $f_*f^! \to \id$ induces
a functorial morphism $p^*f^! \to q^!g^*$. Let us show that
it is an isomorphism.

Decomposing $f$ as $X \exto{f_1} S' \exto{f_2} S$, where
$f_1$ is finite and $f_2$ is strict and smooth,
and applying lemma~\ref{squ3} we reduce the claim of the lemma
to two cases, the case of strict and smooth $f$, and the case of finite $f$,
which can be treated separately.

If $f$ is strict and smooth then $q$ is also strict and smooth and
$$
f^!(F) \cong f^*F\otimes_{\CO_X}\omega_{X/S}[\dim X - \dim S],\qquad
q^!(F) \cong q^*F\otimes_{\CO_{X\times_S Y}}\omega_{X\times_S Y/Y}
[\dim X\times_S Y - \dim Y].
$$
It remains to note that
$\omega_{X\times_S Y/Y} \cong p_\circ^*\omega_{X/S}$
and $\dim X\times_S Y = \dim X + \dim Y - \dim S$
since $f$ is smooth, hence
\begin{multline*}
p^*(f^!F) \cong
p^*(f^*F\otimes_{\CO_X}\omega_{X/S}[\dim X - \dim S]) \cong
p^*f^*F\otimes_{\CO_{X\times_S Y}}p_\circ^*\omega_{X/S}[\dim X - \dim S]) \cong
\\ \cong
q^*g^*F\otimes_{\CO_{X\times_S Y}}\omega_{X\times_S Y/Y}
[\dim X\times_S Y - \dim Y] \cong
q^!(g^*F).
\end{multline*}

Now, assume that $f$ is finite. Then $q$ is also finite and
\begin{multline*}
q_*q^!g^*(F) \cong
\RCHom(q_*\CA_{X\times_S Y},g^*F) \cong
\RCHom(q_*p^*\CA_X,g^*F) \cong \\ \cong
\RCHom(g^*f_*\CA_X,g^*F) \cong
g^*\RCHom(f_*\CA_X,F) \cong
g^*f_*f^!(F) \cong
q_*p^*f^!(F),
\end{multline*}
and it remains to apply lemma~\ref{pff}~b).
\end{proof}

\subsection{Derived categories over a base}\label{ssdcb}

Consider a pair of Azumaya varieties $(X,\CA_X)$ and $(Y,\CA_Y)$
over the same smooth algebraic variety $S$. In other words, we have
a pair of morphisms $f:(X,\CA_X) \to (S,\CO_S)$ and
$g:(Y,\CA_Y) \to (S,\CO_S)$.

A functor $\Phi:\D(X,\CA_X) \to \D(Y,\CA_Y)$ is called {\sf $S$-linear}
if for all $F\in\D(X,\CA_X)$, $G\in\D^b(S,\CO_S)$
there are given bifucntorial isomorphisms
$$
\Phi(f_\circ^*G\otimes_{\CO_X} F) \cong g_\circ^*G\otimes_{\CO_Y} \Phi(F).
$$
Note that since $S$ is smooth any object $G\in\D^b(S,\CO_S)$
is a perfect complex.

\begin{lemma}\label{adjslin}
If $\Phi$ is $S$-linear and admits a right adjoint functor $\Phi^!$
then $\Phi^!$ is also $S$-linear.
\end{lemma}
\begin{proof}
For all $F\in\D(X,\CA_X)$, $G\in\D(Y,\CA_Y)$ and $H\in\D^b(S)$ we have
\begin{multline*}
\RHom_{\CA_X}(F,\Phi^!(g_\circ^*H\otimes_{\CO_Y} G)) \cong
\RHom_{\CA_Y}(\Phi(F),g_\circ^*H\otimes_{\CO_Y} G) \cong
\RHom_{\CA_Y}(g_\circ^*H^*\otimes_{\CO_Y} \Phi(F),G) \cong \\ \cong
\RHom_{\CA_Y}(\Phi(f_\circ^*H^*\otimes_{\CO_X} F),G) \cong
\RHom_{\CA_X}(f_\circ^*H^*\otimes_{\CO_X} F,\Phi^!(G)) \cong
\RHom_{\CA_X}(F,f_\circ^*H\otimes_{\CO_X} \Phi^!(G)).
\end{multline*}
Therefore, by Yoneda lemma we have
$\Phi^!(g_\circ^*H\otimes_{\CO_Y} G) \cong
f_\circ^*H\otimes_{\CO_X} \Phi^!(G)$,
and it is clear that the isomorphism is bifunctorial.
\end{proof}

Consider the fiber product of algebraic varieties $X\times_S Y$.
Let $i:X\times_S Y \to X\times Y$ denote the embedding
and let $p:X\times Y \to X$, $q:X\times Y \to Y$ denote the projections,
so that $f\circ p\circ i = g\circ q\circ i$.

\begin{lemma}\label{fp_slin}
If $K\in\D^-_{qc}(X\times_S Y,\CA_X^\opp\otimes\CA_Y)$ then
the kernel functors $\Phi_{i_*K}$ and $\Phi^!_{i_*K}$ are $S$-linear.
\end{lemma}
\begin{proof}
We have
\begin{multline*}
\Phi_{i_*K}(f_\circ^*G\otimes_{\CO_X} F) \cong
q_*(p_\circ^*(f_\circ^*G\otimes_{\CO_X} F)\otimes_{\CA_X} {i_*K}) \cong
q_*(p_\circ^*f_\circ^*G\otimes_{\CO_{X\times Y}}
p_\circ^*F\otimes_{\CA_X} {i_*K}) \cong
\\ \cong
q_*(p_\circ^*F\otimes_{\CA_X}
i_*(i_\circ^*p_\circ^*f_\circ^*G\otimes_{\CO_{X\times_S Y}} K)) \cong
q_*(p_\circ^*F\otimes_{\CA_X}
i_*(i_\circ^*q_\circ^*g_\circ^*G\otimes_{\CO_{X\times_S Y}} K)) \cong
\\ \cong
q_*(p_\circ^*F\otimes_{\CA_X}
q_\circ^*g_\circ^*G\otimes_{\CO_{X\times Y}} {i_*K}) \cong
g_\circ^*G\otimes_{\CO_Y} q_*(p_\circ^*F\otimes_{\CA_X} {i_*K}) \cong
g_\circ^*G\otimes_{\CO_Y} \Phi_{i_*K}(F),
\end{multline*}
and similarly for $\Phi^!_{i_*K}$.
\end{proof}

A strictly full subcategory $\CC \subset \D(X,\CA_X)$
is called {\sf $S$-linear} if for all $F\in\CC$,
$G\in\D^b(S,\CO_S)$ we have $f_\circ^*G\otimes_{\CO_X} F\in \CC$.
For the definition of admissible subcategory see \cite{B,BO1}.

\begin{lemma}\label{slinperp}
If $\CC\subset\D^b(X,\CA_X)$ is a strictly full $S$-linear left
{\rm(}resp.\ right{\rm)} admissible triangulated subcategory
then its left {\rm(}resp.\ right{\rm)} orthogonal
is also $S$-linear.
\end{lemma}
\begin{proof}
Assume for example that $\CC$ is left admissible and take any
$F\in\CC$, $F'\in\CC^\perp$. Then we have
$$
\Hom_{\CA_X}(F,f_\circ^*G\otimes_{\CO_X} F') \cong
\Hom_{\CA_X}(f_\circ^*G^*\otimes_{\CO_X} F,F') = 0,
$$
since $f_\circ^*G^*\otimes_{\CO_X} F\in\CC$.
Therefore $f_\circ^*G\otimes_{\CO_X} F' \in \CC^\perp$,
so $\CC^\perp$ is $S$-linear.
\end{proof}

\begin{lemma}\label{sod_f}
If $\D^b(X,\CA_X) = \langle\CC^\perp,\CC\rangle$ is a semiorthogonal
decomposition and $\CC$ is $S$-linear then we have
$f_*\RCHom_{\CA_X}(F,F') = 0$ for all $F\in\CC$, $F'\in\CC^\perp$.
\end{lemma}
\begin{proof}
For all $G\in\D^b(S,\CO_S)$ we have
$$
\Hom_{\CO_S}(G,f_*\RCHom_{\CA_X}(F,F')) \cong
\Hom_{\CO_X}(f_\circ^*G,\RCHom_{\CA_X}(F,F')) \cong
\Hom_{\CA_X}(f_\circ^*G\otimes_{\CO_X}F,F') = 0,
$$
since $f_\circ^*G\otimes_{\CO_X} F \in \CC$.
Therefore $f_*\RCHom_{\CA_X}(F,F') = 0$
\end{proof}

\begin{proposition}\label{slinadj}
If $\Phi:\D^b(X,\CA_X) \to \D^b(Y,\CA_Y)$ is $S$-linear,
$\Phi^!:\D^b(Y,\CA_Y) \to \D^b(X,\CA_X)$ is a right adjoint functor,
$f$ and $g$ are projective, and $X$ and $Y$ are smooth
then we have a bifunctorial isomorphism
$$
g_*\RCHom_{\CA_Y}(\Phi(F),F') \cong f_*\RCHom_{\CA_X}(F,\Phi^!(F')).
$$
\end{proposition}
\begin{proof}
Take any $G\in\D^b(S)$ and note that
\begin{multline*}
\Hom_{\CO_S}(G,g_*\RCHom_{\CA_Y}(\Phi(F),F')) \cong
\Hom_{\CO_Y}(g_\circ^*G,\RCHom_{\CA_Y}(\Phi(F),F')) \cong
\\ \cong
\Hom_{\CA_Y}(g_\circ^*G\otimes_{\CO_Y}\Phi(F),F') \cong
\Hom_{\CA_Y}(\Phi(f_\circ^*G\otimes_{\CO_X} F),F') \cong
\Hom_{\CA_X}(f_\circ^*G\otimes_{\CO_X} F,\Phi^!F') \cong
\\ \cong
\Hom_{\CO_X}(f_\circ^*G,\RCHom_{\CA_X}(F,\Phi^!F')) \cong
\Hom_{\CO_X}(G,f_*\RCHom_{\CA_X}(F,\Phi^!F')).
\end{multline*}
On the other hand, $g_*\RCHom_{\CA_Y}(\Phi(F),F') \in \D^b(S)$ and
$f_*\RCHom_{\CA_X}(F,\Phi^!F') \in \D^b(S)$ by assumptions on
$X$, $Y$, $f$ and $g$. Therefore, by Yoneda lemma we have
the desired isomorphism, which is bifunctorial by construction.
\end{proof}

\subsection{Faithful base changes and kernel functors}\label{ssfbc}

Consider morphisms $f:(X,\CA_X) \to (S,\CO_S)$ and
$g:(Y,\CA_Y) \to (S,\CO_S)$ with smooth $S$.
For any strict base change $\phi:T \to S$ we consider the fiber products
$$
X_T := X\times_S T,\qquad
Y_T := Y\times_S T,\qquad
X_T\times_T Y_T = (X\times_S Y)\times_S T
$$
and denote the projections $X_T \to X$, $Y_T \to Y$, and
$X_T\times_T Y_T \to X\times_S Y$ also by $\phi$.
Further, denote the embeddings $X\times_S Y \to X\times Y$ and
$X_T\times_T Y_T \to X_T\times Y_T$ by $i$ and $i_T$ respectively.
For any kernel $K \in \D^-(X\times Y_S,\CA_X^\opp\otimes\CA_Y)$ we denote
$$
K_T = \phi^* K \in \D^-(X_T\times_T Y_T,\CA_X^\opp\otimes\CA_Y).
$$

\begin{definition}
A strict morphism $\phi:T \to S$ is called {\sf faithful} for a pair $(X,Y)$
if $\phi$ is faithful with respect to morphisms $f:X\to S$, $g:Y\to S$,
and $f\times_S g:X\times_S Y \to S$.
\end{definition}

\begin{lemma}\label{faith_fact}
Assume that $T\to S$ is a strict base change faithful for a pair $(X,Y)$.
If $T \to S' \to S$ is a factorization of $T\to S$ such that $T\to S'$
is strict and finite and $S'\to S$ is strict and smooth, then $S'\to S$
is faithful for the pair $(X,Y)$ and $T \to S'$ is faithful for the pair
$(X\times_S S',Y\times_S S')$.
\end{lemma}
\begin{proof}
Follows from corollary~\ref{ec_fact}.
\end{proof}

\begin{lemma}
If $\phi:T\to S$ is a strict base change faithful for a pair $(X,Y)$
then both squares in the commutative diagram
$$
\xymatrix{
X_T \ar[d]^\phi &
X_T\times_T Y_T \ar[d]^\phi \ar[r]^-{q_T} \ar[l]_-{p_T}
&
Y_T \ar[d]^\phi \\
X &
X\times_S Y \ar[r]^-q \ar[l]_-p &
Y
}
$$
are exact cartesian.
\end{lemma}
\begin{proof}
Apply the ``only if'' part of lemma~\ref{squ3}
to the following cartesian diagrams
$$
\vcenter{\xymatrix{
X_T\times_T Y_T \ar[d]^\phi \ar[r]^-{p_T} &
X_T \ar[d]^\phi \ar[r]^-{f_T} &
T \ar[d]^\phi \\
X\times_S Y \ar[r]^-p &
X \ar[r]^-f &
S
}}
\qquad\text{and}\qquad
\vcenter{\xymatrix{
X_T\times_T Y_T \ar[d]^\phi \ar[r]^-{q_T} &
Y_T \ar[d]^\phi \ar[r]^-{g_T} &
T \ar[d]^\phi \\
X\times_S Y \ar[r]^-q &
Y \ar[r]^-g &
S
}}
$$
\end{proof}

\begin{lemma}\label{phit}
If $\phi:T \to S$ is a strict base change faithful for a pair $(X,Y)$,
and $f$ is projective then we have
$$
\Phi_{i_{T*}K_T}\phi^* = \phi^*\Phi_{i_*K},\quad
\Phi_{i_*K}\phi_* = \phi_*\Phi_{i_{T*}K_T},\quad
\Phi_{i_{T*}K_T}^!\phi^* = \phi^*\Phi_{i_*K}^!,\quad\text{and}\quad
\Phi_{i_*K}^!\phi_* = \phi_*\Phi_{i_{T*}K_T}^!.
$$
\end{lemma}
\begin{proof}
Put for short $\Phi = \Phi_{i_*K}$, $\Phi^! = \Phi_{i_*K}^!$,
$\Phi_T = \Phi_{i_{T*}K_T}$, $\Phi_T^! = \Phi_{i_{T*}K_T}^!$.
Note that $\phi$ is strict, hence $\phi^* = \phi_\circ^*$.
For any $F\in \D^-(X)$, $G\in\D^-(X_T)$ we have
\begin{multline*}
\Phi_T\phi^*(F) =
q_{T*}({p_{T\circ}}^*\phi^*F\otimes_{\CA_X} K_T) \cong
q_{T*}(\phi^*p_\circ^*F\otimes_{\CA_X} \phi^*K) \cong \\ \cong
q_{T*}\phi^*(p_\circ^*F\otimes_{\CA_X} K) \cong
\phi^*q_*(p_\circ^*F\otimes_{\CA_X} K) =
\phi^*\Phi(F),
\end{multline*}
\begin{multline*}
\Phi\phi_*(G) =
q_*({p_\circ}^*\phi_*G\otimes_{\CA_X} K) \cong
q_*(\phi_*{p_{T\circ}}^*G\otimes_{\CA_X} K) \cong \\ \cong
q_*\phi_*({p_{T\circ}}^*G\otimes_{\CA_X} \phi^*K) \cong
\phi_*q_{T*}({p_{T\circ}}^*G\otimes_{\CA_X} K_T) =
\phi_*\Phi_T(G).
\end{multline*}
Note that $\phi$ has finite $\Tor$-dimension since $S$ is smooth,
and $q$ is projective because $f$ is. Hence for any
$F\in \D^+(X)$, $G\in\D^+(X_T)$ using lemma~\ref{sbpb} we deduce
\begin{multline*}
\Phi^!_T\phi^*(F) =
p_{T*}\RCHom_{\CA_Y}(K_T,{q_{T\circ}}^!\phi^*F) \cong
p_{T*}\RCHom_{\CA_Y}(\phi^*K,\phi^*q_\circ^!F) \cong \\ \cong
p_{T*}\phi^*\RCHom_{\CA_Y}(K,q_\circ^!F) \cong
\phi^*p_*\RCHom_{\CA_Y}(K,q_\circ^!F) =
\phi^*\Phi^!(F),
\end{multline*}
and using proposition~\ref{sbpf} we deduce
\begin{multline*}
\Phi^!\phi_*(G) =
p_*\RCHom_{\CA_Y}(K,{q_\circ}^!\phi_*F) \cong
p_*\RCHom_{\CA_Y}(K,\phi_*{q_{T\circ}}^!F) \cong \\ \cong
p_*\phi_*\RCHom_{\CA_Y}(\phi^*K,{q_{T\circ}}^!F) \cong
\phi_*p_{T*}\RCHom_{\CA_Y}(K_T,q_{T\circ}^!F) =
\phi_*\Phi^!_T(F).
\end{multline*}
\end{proof}

\begin{lemma}
Assume that $\phi$ is faithful for the pair $(X,Y)$ and $K\in\D^b(X\times_S Y)$.

\noindent$(a)$
If $K$ has finite $\Tor$-amplitude over $Y$ then
$K_T$ has finite $\Tor$-amplitude over $Y_T$.

\noindent$(b)$
If $K$ has finite $\Ext$-amplitude over $X$ then
$K_T$ has finite $\Ext$-amplitude over $X_T$.

\noindent$(c)$
If $\supp K$ is projective over $X$ then
If $\supp K_T$ is projective over $X_T$.

\noindent$(d)$
If $\supp K$ is projective over $Y$ then
If $\supp K_T$ is projective over $Y_T$.
\end{lemma}
\begin{proof}
$(a)$ and $(b)$ follow from corollary~\ref{fted_bc};
$(c)$ and $(d)$ are evident.
\end{proof}

Let $X_1$, \dots, $X_n$, and $Y$ be smooth Azumaya varieties
projective over the same smooth algebraic variety~$S$. Consider kernels
$K_1\in\D^-(X_1\times_S Y,\CA_{X_1}^\opp\otimes\CA_Y)$, \dots,
$K_n\in\D^-(X_n\times_S Y,\CA_{X_n}^\opp\otimes\CA_Y)$ and
denote by $\Phi_1:\D^b(X_1)\to\D^b(Y)$, \dots, $\Phi_n:\D^b(X_n)\to\D^b(Y)$
the corresponding kernel functors.
Let $\phi:T\to S$ be a strict base change. We use the notations
$X_{kT}$, $Y_T$, $f_{kT}$, $g_T$, $\Phi_{kT}$, $\Phi_{kT}^!$
in an evident sense.
Finally, assume that $K_i$ have finite $\Tor$-amplitude over $Y$,
finite $\Ext$-amplitude over $X_i$, and that morphisms $g$ and $f_i$
are projective for all $i$. Then the above lemma
together with lemma~\ref{phi_bounded} imply that
all functors preserve bounded derived categories and that
$\Phi^!_i$, $\Phi^!_{iT}$ are right adjoint to $\Phi_i$, $\Phi_{iT}$.

\begin{proposition}\label{fbc_sod}
Assume that $\phi$ is faithful for the pairs $(X_1,Y)$, $(X_2,Y)$.

\noindent$(i)$
If the subcategory $\Phi_2(\D^b(X_2))$ is right orthogonal
to $\Phi_1(\D^b(X_1))$ then
the subcategory $\Phi_{2T}(\D^b(X_{2T}))$ is right orthogonal
to $\Phi_{1T}(\D^b(X_{1T}))$.

\noindent$(ii)$
If $\Phi_1$ is fully faithful then $\Phi_{1T}$ is fully faithful.

\noindent$(iii)$
If $\Phi_1^!$ is fully faithful then $\Phi_{1T}^!$ is fully faithful.

\noindent$(iv)$
If $\Phi_1$ is an equivalence then $\Phi_{1T}$ is an equivalence.
\end{proposition}
\begin{proof}
First we note that $\Phi_2(\D^b(X_2))$ is right orthogonal
to $\Phi_1(\D^b(X_1))$ iff $\Phi_1^!\Phi_2 = 0$ and
$\Phi_{2T}(\D^b(X_{2T}))$ is right orthogonal to
$\Phi_{1T}(\D(X_{1T}))$ iff $\Phi_{1T}^!\Phi_{2T} = 0$.
So, in $(i)$ we must check that
$\Phi_1^!\Phi_2 = 0$ implies $\Phi_{1T}^!\Phi_{2T} = 0$.

Similarly, $\Phi_1$ is fully faithful iff the adjunction morphism
$\id_{\D^b(X_1)} \to \Phi_1^!\Phi_1$ is an isomorphism and
$\Phi_{1T}$ is fully faithful iff the adjunction morphism
$\id_{\D^b(X_{1T})} \to \Phi_{1T}^!\Phi_{1T}$ is an isomorphism.
So, in $(ii)$ we must check that $\id_{\D^b(X_1)} \cong \Phi_1^!\Phi_1$
implies $\id_{\D^b(X_{1T})} \cong \Phi_{1T}^!\Phi_{1T}$.

By lemma~\ref{faith_fact} it suffices to consider the case when
$\phi$ is smooth and the case when $\phi$ is finite.

If $\phi$ is smooth then $X_{1T}$, $X_{2T}$, and $Y_T$ are smooth varieties,
hence $\Phi_1^!$ and $\Phi_{1T}^!$ are usual kernel functors by
lemma~\ref{keronsm}.
Moreover, if $\CF_1$ is a perfect spanning class for $X_1$,
$\CF_2$ is a perfect spanning class for $X_2$, and $\CG$
is a perfect spanning class for $T$ then
$f_{2T\circ}^*\CG\otimes_{\CO_{X_{2T}}} \CF_2$ is a perfect spanning class
for $X_{2T}$ by lemma~\ref{psc_fp}. Further,
for all $F\in \CF_2$, $G\in\CG$
by lemma~\ref{phit} and lemma~\ref{fp_slin} we have
\begin{multline*}
\Phi_{1T}^!\circ\Phi_{2T}(f_{2T\circ}^*G\otimes_{\CO_{X_{2T}}} F) \cong
\Phi_{1T}^!(g_{T\circ}^*G\otimes_{\CO_{Y_{T}}} \Phi_{2T}(\phi^*F)) \cong
\\ \cong
f_{1T\circ}^*G\otimes_{\CO_{X_{1T}}} \Phi_{1T}^!(\phi^*\Phi_2(F)) \cong
f_{1T\circ}^*G\otimes_{\CO_{X_{1T}}} \phi^*\Phi_1^!(\Phi_2(F)).
\end{multline*}
We deduce that $\Phi_1^!\Phi_2 = 0$ implies that $\Phi_{1T}^!\Phi_{2T}$
vanishes on $f_{2T\circ}^*\CG\otimes_{\CO_{X_{2T}}} \CF_2$, hence
$\Phi_{1T}^!\Phi_{2T} = 0$ by proposition~\ref{kfz}.

Similarly, $f_{1T\circ}^*\CG\otimes_{\CO_{X_{1T}}} \CF_1$ is a perfect
spanning class for $X_{1T}$ and for all $F\in \CF_1$, $G\in\CG$
by lemma~\ref{phit} and lemma~\ref{fp_slin} we have
\begin{multline*}
\Phi_{1T}^!\circ\Phi_{1T}(f_{1T\circ}^*G\otimes_{\CO_{X_{1T}}} \phi^*F) \cong
\Phi_{1T}^!(g_{T\circ}^*G\otimes_{\CO_{Y_{T}}} \Phi_{1T}(\phi^*F)) \cong \\ \cong
f_{1T\circ}^*G\otimes_{\CO_{X_{1T}}} \Phi_{1T}^!(\phi^*\Phi_1(F)) \cong
f_{1T\circ}^*G\otimes_{\CO_{X_{1T}}} \phi^*\Phi_1^!(\Phi_1(F)).
\end{multline*}
We deduce that $\id_{\D^b(X_1)} \cong \Phi_1^!\Phi_1$ implies
that $\Phi_{1T}^!\Phi_{1T}$ is isomorphic to $\id_{\D^b(X_{1T})}$ on
$f_{1T\circ}^*\CG\otimes_{\CO_{X_{1T}}} \CF_1$, hence
$\id_{\D^b(X_{1T})} \cong \Phi_{1T}^!\Phi_{1T}$
by corollary~\ref{mkfi}.

If $\phi$ is finite then lemma~\ref{phit} implies that
$$
\phi_*\Phi_{1T}^!\Phi_{2T} \cong
\Phi_1^!\phi_*\Phi_{2T} \cong
\Phi_1^!\Phi_2\phi_*,
\quad\text{and}\quad
\phi_*\Phi_{1T}^!\Phi_{1T} \cong
\Phi_1^!\phi_*\Phi_{1T} \cong
\Phi_1^!\Phi_1\phi_*.
$$
We deduce that $\Phi_1^!\Phi_2 = 0$ implies $\phi_*\Phi_{1T}^!\Phi_{2T} = 0$,
hence $\Phi_{1T}^!\Phi_{2T} = 0$ by lemma~\ref{pff}~a), and that
$\id_{\D^b(X_1)} \cong \Phi_1^!\Phi_1$ implies
$\phi_* \cong \phi_*\Phi_{1T}^!\Phi_{1T}$, hence
$\id_{\D^b(X_{1T})} \cong \Phi_{1T}^!\Phi_{1T}$ by lemma~\ref{pff}~b).

The claim $(iii)$ can be proved by the similar arguments as $(ii)$
using the fact that $\Phi_1^!$ is fully faithful iff
$\Phi_1\Phi^!_1 \cong \id_{\D^b(Y)}$ and $\Phi_{1T}^!$
is fully faithful iff $\Phi_{1T}\Phi^!_{1T} \cong \id_{\D^b(Y_T)}$.
Finally, $(iv)$ follows from $(ii)$ combined with $(iii)$.
\end{proof}

\begin{corollary}\label{phitperp}
In the conditions of the above proposition
if $g_*\RCHom_{\CA_Y}(E,\Phi_2(\D^b(X_2))) = 0$ for some $E\in\D^b(Y)$
then $g_{T*}\RCHom_{\CA_{Y_T}}(\phi^*E,\Phi_{2T}(\D^b(X_{2T}))) = 0$.
\end{corollary}
\begin{proof}
We take $X_1 = S$, consider $E$ as a kernel on
$X_1\times_S Y = S\times_S Y \cong Y$, put $\Phi_1 = \Phi_E$,
and apply proposition~\ref{fbc_sod}~(i) and lemma~\ref{sod_f}.
\end{proof}

\begin{theorem}\label{phitsod}
In the conditions of the proposition if
$$
\D^b(Y) = \langle\Phi_1(\D^b(X_1)),\dots,\Phi_n(\D^b(X_n))\rangle
$$
is a semiorthogonal decomposition for $Y$ then
$$
\D^b(Y_T) =
\langle\Phi_{1T}(\D^b(X_{1T})),\dots,\Phi_{nT}(\D^b(X_{nT}))\rangle
$$
is a semiorthogonal decomposition for $Y_T$.
\end{theorem}
\begin{proof}
By proposition~\ref{fbc_sod} it suffices to check that the orthogonal
to the subcategory $\D'\subset\D^b(Y_T)$ generated by
$\Phi_{1T}(\D^b(X_{1T}))$, \dots, $\Phi_{nT}(\D^b(X_{nT}))$ is zero.
Choose a perfect spanning class $\CF$ for $Y$ and a perfect spanning class
$\CG$ for $T$. Then any object $F\in\CF$ can be decomposed with respect
to the initial semiorthogonal decomposition.
Applying lemma~\ref{phit} and lemma~\ref{fp_slin} we deduce that
$g_{T\circ}^*G \otimes_{\CO_{Y_T}} \phi^*(F)$ is contained in $\D'$
for all $F\in\CF$, $G\in\CG$. Hence $\D'$ contains a perfect spanning class
for $Y_T$, hence ${\D'}^\perp = 0$ by lemma~\ref{psc_ort}.
\end{proof}

\subsection{Relative Serre functor and Bridgeland's trick}\label{ssbt}

Assume that $X$ and $Y$ are Azumaya varieties, $X$ is smooth, and
$f:X\to S$, $g:Y\to S$ are projective morphisms into a smooth
algebraic variety $S$. If $\Phi:\D^b(X) \to \D^b(Y)$ is an $S$-linear
functor then for all $F,F'\in\D^b(X)$, $G\in\D^b(S)$ we have
a diagram
$$
\xymatrix@C=0pt{
\Hom_{\CO_S}(G,f_*\RCHom_{\CA_X}(F,F')) \ar@{=}[r] \ar@{-->}[d] &
\Hom_{\CO_X}(f_\circ^*G,\RCHom_{\CA_X}(F,F')) \ar@{=}[r] &
\Hom_{\CA_X}(f_\circ^*G\otimes_{\CO_X} F,F')
\ar[d]^{\Phi_{f_\circ^*G\otimes_{\CO_X} F,F'}}
\\
\Hom_{\CO_S}(G,g_*\RCHom_{\CA_Y}(\Phi F,\Phi F')) \ar@{=}[r] &
\Hom_{\CA_Y}(g_\circ^*G\otimes_{\CO_Y} \Phi F,\Phi F')  \ar@{=}[r] &
\Hom_{\CA_Y}(\Phi(f_\circ^*G\otimes_{\CO_X} F),\Phi F')
}
$$
where the dashed arrow is the unique morphism making the diagram commutative.
Note that we have $f_*\RCHom_{\CA_X}(F,F') \in \D^b(S)$, since $X$ is smooth
and $f$ is projective and the dashed arrow is evidently functorial in $G$.
Hence by the Yoneda lemma it is induced by a unique morphism
$$
\Phi^S_{F,F'} : f_*\RCHom_{\CA_X}(F,F') \to g_*\RCHom_{\CA_Y}(\Phi F,\Phi F').
$$

\begin{lemma}
Morphisms $\Phi^S_{F,F'}$ are bifunctorial in $F,F'$.
\end{lemma}
\begin{proof}
Use uniqueness of $\Phi^S_{F,F'}$
and functorial properties of the above commutative diagram.
\end{proof}

\begin{definition}
(cf.\ \cite{BK,BO4})
A covariant additive functor $\FS_{X/S}:\D^b(X)\to\D^b(X)$ is called
a {\sf relative Serre functor}, if it is an $S$-linear category equivalence
and for all objects $F,F'\in\D^b(X)$ there are given bifunctorial isomorphisms
$$
\varphi_{F,F'}:
f_*\RCHom_{\CA_X}(F,F') \to
(f_*\RCHom_{\CA_X}(F',\FS_{X/S}F))^*.
$$
\end{definition}

Let $\omega_{X/S}$ denote the relative canonical class
of the underlying algebraic variety of $X$ over $S$.

\begin{lemma}\label{rsf}
{\bf(}cf.\ \cite{BK}{\bf)}
The functor $\FS_{X/S}(F) = \omega_{X/S}[\dim X - \dim S]\otimes_{\CO_X} F$
is a relative Serre functor, and any relative Serre functor is
canonically isomorphic to it.
\end{lemma}
\begin{proof}
Using the duality theorem we get
\begin{multline*}
(f_*\RCHom_{\CA_X}(F,F'))^* \cong
\RCHom_{\CO_S}(f_*\RCHom_{\CA_X}(F,F'),\CO_S) \cong \\ \cong
f_*\RCHom_{\CO_X}(\RCHom_{\CA_X}(F,F'),f_\circ^!\CO_S) \cong
f_*\RCHom_{\CA_X}(F',f_\circ^!\CO_S\otimes_{\CO_X} F).
\end{multline*}
But $f_\circ^!\CO_S \cong \omega_{X/S}[\dim X - \dim S]$, and it is
clear that $\FS_{X/S}(F) = \omega_{X/S}[\dim X - \dim S]\otimes_{\CO_X} F$
is an $S$-linear category equivalence. Dualizing the above isomorphism we
deduce that it is a relative Serre functor. The second part is proved
in the same way as the uniqueness of a usual Serre functor.
\end{proof}

\begin{corollary}\label{relserre}
If $X$ is smooth, $\Phi:\D^b(Y) \to \D^b(X)$ is
an $S$-linear equivalence and $\omega_{X/S} \cong f_\circ^*L$
for a line bundle $L$ on $S$, then $Y$ is smooth, $\dim Y = \dim X$,
and $\omega_{Y/S} \cong g_\circ^*L$.
\end{corollary}
\begin{proof}
If $X$ is smooth then $\D^b(X)$ is $\Ext$-bounded, hence $\D^b(Y)$
is  $\Ext$-bounded, hence $Y$ is smooth (see lemma~\ref{snav_sm}).
Further, it is clear that $\FS_{Y/S} = \Phi^!\circ\FS_{X/S}\circ\Phi$
is a relative Serre functor on $Y$. But denoting $d = \dim X - \dim S$ we get
$$
(\Phi^!\circ\FS_{X/S}\circ\Phi)(F) =
\Phi^!(f_\circ^*L[d] \otimes_{\CO_X} \Phi(F)) \cong
g_\circ^*L[d] \otimes_{\CO_Y} \Phi^!(\Phi(F)) \cong
g_\circ^*L[d] \otimes_{\CO_Y} F,
$$
and it remains to apply the uniqueness of the relative Serre functor.
\end{proof}

\begin{proposition}\label{ffisequ}
If $X$ is smooth and connected, $\omega_X \cong f_\circ^*L$
for a line bundle $L$ on $S$, and $Y$ is not empty, then
any fully faithful $S$-linear functor $\Phi:\D^b(Y) \to \D^b(X)$
is an equivalence. Furthermore, in this case $Y$ is smooth, connected,
$\dim Y = \dim X$ and $\omega_Y \cong g_\circ^*L$.
\end{proposition}
\begin{proof}
The functor $\Phi$ is an equivalence of $\D^b(Y)$ with a strictly full
right admissible $S$-linear triangulated subcategory $\CC\subset\D^b(X)$.
Consider a semiorthogonal decomposition $(\CC^\perp,\CC)$ for $\D^b(X)$.
Note that the Serre functor
$\FS_{X/S}(F) = f_\circ^*L\otimes_{\CO_X} F[\dim X - \dim S]$
preserves any $S$-linear subcategory of $\D^b(X)$.
Therefore for all $F\in\CC$, $F'\in\CC^\perp$ we have
$f_*\RCHom_{\CA_X}(F',F) \cong
(f_*\RCHom_{\CA_X}(F,\FS_{X/S}F'))^* = 0$,
which implies
$\Hom_{\CA_X}(F',F) \cong \Hom_{\CO_S}(\CO_S,f_*\RCHom_{\CA_X}(F',F)) = 0$.
Thus the subcategories $\CC$ and $\CC^\perp$ are totally orthogonal.
But the category $\D^b(X)$ is indecomposable by \cite{Br1}, hence
$\CC^\perp = 0$ and $\Phi$ is an equivalence. The rest of the claim
follows from corollary~\ref{relserre}.
\end{proof}

\section{The universal hyperplane section}\label{uhps}

Though $Y$ is in general an Azumaya variety, for unburdening of the notation
we drop all the Azumaya stuff in this section. However, a patient reader
can easily fill in the details and check that section~\ref{sect1}
and appendix~D contain a verification of all necessary facts
concerning Azumaya varieties.

The goal of this section is to make the first step in the proof
of the theorem~\ref{themain} as indicated in the introduction.
So assume that we have data (D.1)--(D.5) satisfying
the conditions (C.1)--(C.6).
In fact, the condition (C.8) is never used in this section,
while the condition (C.7) will be explicitly stated in theorem~\ref{ff}.

\subsection{Notation}

Let $\BP = \PP(V^*)\setminus \BZ$ and denote by $\CL$ the restriction to $\BP$
of $\CO_{\PP(V^*)}(-1)$. Recall that~(C.6) implies that we have the following
exact sequence on $X\times Y$
\begin{equation}\label{equ_e}
0 \to E_1\boxtimes F_1 \exto{e} E_2\boxtimes F_2 \to i_*\CE \to 0,
\end{equation}
where $\CE$ is a coherent sheaf on $Q(X,Y)$ and
$i:Q(X,Y) \to X\times Y$ denotes the embedding.

Consider $\CX = X\times\BP$, and let $\CX_1 \subset X\times\BP$
be the {\sf universal hyperplane section of $X$},
i.e.\ the zero locus of the canonical section of the line bundle
$\CO_X(1)\boxtimes\CL^*$.
Let $\alpha:\CX_1\to X\times\BP$ denote the embedding.
Then we have a resolution
\begin{equation}\label{ex1}
0 \to \CO_X(-1)\boxtimes\CL \to \CO_{X\times\BP} \to \alpha_*\CO_{\CX_1} \to 0.
\end{equation}

Denote by $\beta$ the embedding of $Y$ to $\BP\times Y$
given by the graph of $g$. We will show below that the composition
$$
\xymatrix{
Q(X,Y) \ar[r]^i &
X\times Y \ar[rr]^{\id_X\times\beta} &&
X\times \BP \times Y
}
$$
factors through $\CX_1\times Y$. Let $j:Q(X,Y) \to \CX_1\times Y$
denote the corresponding closed embedding, so that
$(\id_X\times\beta)\circ i = (\alpha\times\id_Y)\circ j$
and we have a commutative square
\begin{equation}\label{abij}
\vcenter{\xymatrix{
Q(X,Y) \ar[r]^i \ar[d]^j &
X\times Y \ar[d]^{\beta} \\
\CX_1\times Y \ar[r]^-{\alpha} &
X\times \BP \times Y
}}
\end{equation}
where we write $\alpha$ instead of $\alpha\times\id_Y$ and $\beta$
instead of $\id_X\times\beta$ for brevity.
Further, consider an object
$\CE^* = \RCHom_{Q(X,Y)}(\CE,\CO_{Q(X,Y)}) \in \D(Q(X,Y))$
and denote
\begin{equation}\label{ce1d}
\CE_1 = j_*\CE,
\qquad
\CE^{\#t} = \CE^*(1+t-i,1),
\qquad
\CE_1^{\#t} = j_*\CE^{\#t}\otimes{\CL^*}^{\otimes t}[\dim X - 1],
\end{equation}
where $\CO(a,b)$ stands for $\CO_X(a)\otimes\CO_Y(b)$.
Let
\begin{equation}\label{b}
B := \dim X + \dim\BP - \codim_{Y\times Y}(Y\times_\BP Y).
\end{equation}
Note that the condition (C.7) is equivalent to $B = 2i - 1$.

The main result of this section is the following

\begin{theorem}\label{ff}
We have $B \ge 2i - 1$.
If $B = 2i - 1$ then the functor $\Phi_{\CE_1}:\D^b(Y) \to \D^b(\CX_1)$
is fully faithful.
Moreover, in this case $Y$ is smooth.
\end{theorem}

The proof goes as follows. First of all we note that the functor
$\Phi_{\CE_1^{\#0}}$ is left adjoint to $\Phi_{\CE_1}$, so the claim
of the theorem reduces to the computation of the convolution
$\CE_1\circ\CE_1^{\#0}$. However, a straightforward computation
is not as easy. Roughly speaking, we have to compute $\RCHom$ between
two sheaves on a singular divisor, when we only know resolutions of these
sheaves on the ambient variety. So we use another argument to circumvent this.
First of all we consider the convolutions $\CE_1\circ\CE_1^{\#t}$
for all $t\in\ZZ$. Second, we consider objects
$\TE = i^*i_*\CE \in \D(Q(X,Y))$, $\TE_1 = j_*\TE \in \D(\CX_1\times Y)$
and the convolutions $\TE_1\circ\CE_1^{\#t}$. Since $i$ is a divisorial
embedding ($Q(X,Y)$ is a divisor of bidegree $(1,1)$) we have an exact
triangle $i^*i_*\CE \to \CE \to \CE(-1,-1)[2]$ which gives an exact
triangle of convolutions
$\TE_1\circ\CE_1^{\#t} \to
\CE_1\circ\CE_1^{\#t} \to
\CE_1\circ\CE_1^{\#(t-1)}[2]$.
It turns out that the first term can be easily computed for $0\le t\le i$.
Comparing this result with uniform (in $t$) cohomological boundedness of
$\CE_1\circ\CE_1^{\#t}$ we deduce the theorem.

Now we give a detailed proof.

\subsection{Preparations}

\begin{lemma}\label{x1}
The universal hyperplane section $\CX_1$ is smooth,
projective over $\BP$ and its relative dimension over $\BP$
equals $\dim X - 1$.
\end{lemma}
\begin{proof}
It is easy to see that the projection $\CX_1 \to X$ is smooth
(in fact it is a projectivization of a vector bundle), hence $\CX_1$ is smooth.
On the other hand, the fibers of $\CX_1$ over $\BP$ are hyperplane sections
of $X$, hence all of them have dimension $\dim X - 1$ by (C.1) and projective
since $X$ is.
\end{proof}

\begin{lemma}
We have $\CX_1\times_\BP Y \cong Q(X,Y)$.
If $j$ denotes the induced embedding $Q(X,Y) \to \CX_1\times Y$ then
$(\id_X\times\beta)\circ j = (\alpha\times\id_Y)\circ i$
and the square $(\ref{abij})$ is commutative.
\end{lemma}
\begin{proof}
It is clear that $\CX_1\times_\BP Y$ coincides with the zero locus
of the canonical section of the line bundle $\CO_X(1)\otimes\CL^*$
on $\CX\times_\BP Y = (X\times\BP)\times_\BP Y = X\times Y$.
On the other hand, on $\CX\times_\BP Y$ we have
$\CL^* \cong \CO_Y(1)$ hence $\CX_1\times_\BP Y$ is the zero locus
of the canonical section of the line bundle $\CO_X(1)\otimes\CO_Y(1)$,
hence coincides with $Q(X,Y)$. It remains to note that the canonical
embedding $(X\times\BP)\times_\BP Y \to X\times\BP\times Y$ coincides with
$\id_X\times\beta$, while the embedding $\CX_1\times Y \to X\times\BP\times Y$
coincides with $\alpha\times\id_Y$.
\end{proof}

\begin{lemma}\label{ces_sh}
An object $\CE^* = \RCHom_{Q(X,Y)}(\CE,\CO_{Q(X,Y)}) \in \D(Q(X,Y))$ is
a coherent sheaf and we have the following exact sequence on $X\times Y$
\begin{equation}\label{equ_estar}
0 \to
E_2^*(-1)\boxtimes F_2^*(-1) \exto{e^*}
E_1^*(-1)\boxtimes F_1^*(-1) \to
i_*\CE^* \to
0.
\end{equation}
\end{lemma}
\begin{proof}
Applying the functor
$\RCHom(-,\CO_{X\times Y}(-1,-1))$ to~(\ref{equ_e}) we obtain a triangle
$$
E_2^*(-1)\boxtimes F_2^*(-1) \exto{e^*}
E_1^*(-1)\boxtimes F_1^*(-1) \to
\RCHom(i_*\CE,\CO_{X\times Y}(-1,-1))[1]
$$
on $X\times Y$. Since $e$ is an isomorphism at the generic point
of $X\times Y$, we deduce that $e^*$ is an embedding, hence the third term
of the triangle is a sheaf. On the other hand, $i$ is a divisorial embedding,
hence $\CO_{Q(X,Y)} = i^*\CO_{X\times Y} = i^!\CO_{X\times Y}(-1,-1)[1]$,
and by the duality theorem we have
$$
i_*\RCHom(\CE,\CO_{Q(X,Y)}) =
i_*\RCHom(\CE,i^!\CO_{X\times Y}(-1,-1))[1] =
\RCHom(i_*\CE,\CO_{X\times Y}(-1,-1))[1].
$$
Since $i_*$ is exact this implies that $\CE^* = \RCHom(\CE,\CO_{Q(X,Y)})$
is a coherent sheaf on $Q(X,Y)$ and (\ref{equ_estar}) is exact.
\end{proof}

\begin{lemma}\label{e1perf}
The sheaf $\CE_1$ is a perfect complex on $\CX_1\times Y$.
\end{lemma}
\begin{proof}
Note that the sheaf $\CE_1 = j_*\CE$ on $\CX_1\times Y$ has
finite $\Tor$-amplitude over $Y$. Indeed, consider the commutative
square~(\ref{abij}). For any sheaf $\CF$ on $Y$, denoting by the same symbol
its pull-backs to $\CX_1\times Y$, $X\times\BP\times Y$ and $X\times Y$,
and using the projection formula, we obtain
$$
\alpha_*(j_*\CE\otimes\CF) =
\alpha_*j_*\CE\otimes\CF =
\beta_*i_*\CE\otimes\CF =
\beta_*(i_*\CE\otimes\CF).
$$
Since $\alpha_*$ and $\beta_*$ are exact it follows that
the $\Tor$-amplitude of $j_*\CE$ over $Y$ is equal to
the $\Tor$-amplitude of $i_*\CE$ over $Y$, and the latter is finite
since $i_*\CE$ admits a finite locally free resolution (\ref{equ_e})
on $X\times Y$. Since $\CX_1\times Y$ is smooth over $Y$ by lemma~\ref{x1},
we conclude by lemma~\ref{snav_sm1}.
\end{proof}

\begin{lemma}\label{ed_la}
The functor $\Phi_{\CE_1^{\#0}}$ is left adjoint to $\Phi_{\CE_1}$.
\end{lemma}
\begin{proof}
Note that $\CE_1$ is a perfect complex by lemma~\ref{e1perf},
and its support $Q(X,Y) \cong \CX_1\times_\BP Y$ is projective
over $\CX_1$ by~(C.4) and over $Y$ by lemma~\ref{x1}. Therefore,
according to lemma~\ref{ladj} it suffices to verify that
$\CE_1^{\#0} \cong \RCHom(\CE_1,\omega_{\CX_1}[\dim \CX_1])$.
Let $q_1$ and $q$ be the projections $\CX_1\times Y\to Y$ and
$X\times Y \to  Y$. It is clear that $q\circ i = q_1\circ j$.
Using the duality theorem and the functoriality of the twisted pullback
we deduce
\begin{multline*}
\RCHom(\CE_1,\omega_{\CX_1}[\dim \CX_1]) \cong
\RCHom(j_*\CE,q_1^!\CO_Y) \cong
j_*\RCHom(\CE,j^!q_1^!\CO_Y) \cong
j_*\RCHom(\CE,i^!q^!\CO_Y) \cong \\ \cong
j_*\RCHom(\CE,i^!\omega_X[\dim X]) \cong
j_*\RCHom(\CE,\omega_X(1,1)[\dim X-1]) \cong
j_*\CE^*(1-i,1)[\dim X-1],
\end{multline*}
since $\omega_X \cong \CO_X(-i)$ by (C.2).
\end{proof}

Twisting (\ref{equ_estar}) by $\CO(1+t-i,1)$
and taking into account~(\ref{ce1d}),
we obtain an exact triangle
\begin{equation}\label{equ_ediez}
E_2^*(t-i)\boxtimes F_2^* \exto{e^*}
E_1^*(t-i)\boxtimes F_1^* \to
i_*\CE^{\#t}
\end{equation}

Consider the following diagram
\begin{equation}\label{diag_1}
\vcenter{\xymatrix{
Y \times \CX_1 \ar[d]^\alpha &
Y \times \CX_1 \times Y \ar[d]^\alpha \ar[l]_{p_{12}} \ar[r]^{p_{23}} &
\CX_1 \times Y \ar[d]^\alpha \\
Y \times (X \times \BP)  &
Y \times (X \times \BP) \times Y \ar[d]^q \ar[l]_{p_{12}} \ar[r]^{p_{23}} &
(X \times \BP) \times Y \\
& Y \times \BP \times Y \ar[d]^\pi \\
& Y \times Y
}}
\end{equation}
where $\pi$ is the projection along $\BP$.
Consider the following objects in $\D(Y\times \BP \times Y)$:
\begin{equation}\label{cct}
\CC_t = q_*\alpha_*(p_{12}^*\CE_1^{\#t}\otimes p_{23}^*\CE_1),
\qquad\qquad
\TCC_t = q_*(\alpha_*p_{12}^*\CE_1^{\#t}\otimes \alpha_*p_{23}^*\CE_1)
\end{equation}

\begin{lemma}\label{k}
We have $\CE_1\circ\CE_1^{\#t} \cong \pi_*\CC_t \in \D(Y\times Y)$.
\end{lemma}
\begin{proof}
Use the definition of the convolution and note that
$p_{13} = \pi\circ q\circ \alpha$.
\end{proof}

\begin{lemma}\label{lbcc}
There exists $N\in\ZZ$ such that
$\CC_t \in \D^{[-N,0]}(Y\times\BP\times Y)$ for all $t\in\ZZ$.
\end{lemma}
\begin{proof}
Since $\CE_1$ is perfect by lemma~\ref{e1perf}, it follows
that $p_{23}^*\CE_1$ is perfect as well.
Hence $p_{12}^*\CE_1^{\#t}\otimes p_{23}^*\CE_1$ is cohomologically bounded.
Moreover, since
$\CE_1^{\#t} \cong \CE_1^{\#0}\otimes ({\CL^*}^{\otimes t}\otimes\CO_X(t))$,
the left bound doesn't depend on $t$. Thus, there exists $N\in\ZZ$, such that
$p_{12}^*\CE_1^{\#t} \otimes p_{23}^*\CE_1 \in
\D^{\ge -N}(Y\times\CX_1\times Y)$.
On the other hand,
$p_{12}^*\CE_1^{\#t}\in \D^{1-\dim X}(Y\times\CX_1\times Y)$
by~(\ref{ce1d}) and lemma~\ref{ces_sh},
and
$p_{23}^*\CE_1\in \D^{0}(Y\times\CX_1\times Y)$ by~(\ref{ce1d}).
Hence $p_{12}^*\CE_1^{\#t} \otimes p_{23}^*\CE_1 \in
\D^{[-N,1-\dim X]}(Y\times\CX_1\times Y)$.
Since $q\circ\alpha$ is projective and has relative dimension $\dim X - 1$
by lemma~\ref{x1}, we are done.
\end{proof}

\begin{lemma}\label{cctt}
We have an exact triangle
\begin{equation}\label{tcccct}
\TCC_t \to
\CC_t \to
\CC_{t-1}[2]
\end{equation}
in $\D^b(Y\times\BP\times Y)$.
\end{lemma}
\begin{proof}
Since $\alpha:\CX_1\to X\times\BP$ is a divisorial embedding,
and $\CX_1$ is a zero locus of a section of the bundle
$\CO_X(1)\boxtimes\CL^*$, we have an exact triangle
$$
\alpha^*\alpha_*\CF \to \CF \to \CF\otimes(\CO_X(-1)\boxtimes\CL)[2]
$$
for any object $\CF$ on $Y\times\CX_1\times Y$.
Taking $\CF = p_{23}^*\CE_1$, tensoring with $p_{12}^*\CE_1^{\#t}$,
applying $q_*\alpha_*$ and taking into account the projection formula
$\alpha_*(p_{12}^*\CE_1^{\#t}\otimes \alpha^*\alpha_*p_{23}^*\CE_1) \cong
\alpha_*p_{12}^*\CE_1^{\#t}\otimes \alpha_*p_{23}^*\CE_1$,
and the definition of $\CC_t$ and $\TCC_t$ (\ref{cct}),
we obtain (\ref{tcccct}).
\end{proof}

Denote by $\beta':Y \to Y\times\BP$ the composition of $\beta$
with the transposition $\BP\times Y \to Y\times\BP$. Consider the maps
$\id_Y\times\beta,\beta'\times\id_Y:Y\times Y \to Y\times\BP\times Y$
and denote them by $\beta$ and $\beta'$ for brevity. Let
\begin{equation}\label{fd}
\FD := \beta'_*\CO_{Y\times Y} \otimes \beta_*\CO_{Y\times Y}
\in \D(Y\times\BP\times Y).
\end{equation}

\begin{lemma}\label{h0fd}
We have $\FD \in \D^{[\codim Y\times_\BP Y - \dim \BP,0]}(Y\times\BP\times Y)$.
Moreover, $\FD$ is supported scheme-theoretically on a infinitesimal
neighborhood of $Y\times_\BP Y \subset Y\times\BP\times Y$ and
$\CH^0(\FD) \cong \CO_{Y\times_\BP Y}$.
\end{lemma}
\begin{proof}
Consider a diagram
$$
\xymatrix{
Y\times_\BP Y \ar[r] \ar[d] &
Y\times Y \ar[r]^{g_2} \ar[d]_{\beta} &
Y\times\BP \ar[d]^{\id_Y\times\Delta_\BP} \\
Y\times Y \ar[r]^{\beta'} &
Y\times\BP\times Y \ar[r]^{g_3} &
Y\times\BP\times\BP
}
$$
where $g_2 = \id_Y\times g$, $g_3 = \id_Y\times\id_\BP\times g$ and
$\Delta_\BP$ is the diagonal embedding $\BP \to \BP\times\BP$.
Both squares are cartesian. Moreover, since $\BP$ is smooth the right
vertical arrow is a locally complete intersection embedding. Since
$\codim_{Y\times\BP\times Y}(Y\times Y) =
\codim_{Y\times\BP\times\BP}(Y\times\BP) = \dim \BP$,
the right square is exact cartesian by corollary~\ref{lci}.
Further, locally over $Y\times\BP\times\BP$ we can represent
$Y\times\BP$ as the zero locus of a regular section $s$ of a vector
bundle on $Y\times\BP\times\BP$ of rank equal to $\dim\BP$,
hence we have a local representation
$(\id_Y\times\Delta_\BP)_*\CO_{Y\times\BP} \cong
\Kosz_{Y\times\BP\times\BP}(s)$.
It follows that
$$
\beta_*\CO_{Y\times Y} \cong
\beta_*g_2^*\CO_{Y\times\BP} \cong
g_3^*(\id_Y\times\Delta_\BP)_*\CO_{Y\times\BP} \cong
g_3^*\Kosz_{Y\times\BP\times\BP}(s)
$$
and
$$
\FD =
\beta'_*\CO_{Y\times Y} \otimes \beta_*\CO_{Y\times Y} \cong
\beta'_*{\beta'}^*\beta_*\CO_{Y\times Y} \cong
\beta'_*{\beta'}^*g_3^*\Kosz_{Y\times\BP\times\BP}(s) \cong
\beta'_*\Kosz_{Y\times Y}({\beta'}^*g_3^*s).
$$
It remains to note that the zero locus of the section
${\beta'}^*g_3^*s$ on $Y\times Y$ coincides with $Y\times_\BP Y$
and the proof concludes with lemma~\ref{koszul}.
\end{proof}

Consider on $Y\times Y$ the objects $\CT,\CT^*\in\D(Y\times Y)$
defined by the following exact triangles
\begin{equation}\label{tts}
\vcenter{\xymatrix@R=5pt{
F_2^*\boxtimes F_2 \oplus F_1^*\boxtimes F_1 \ar[rr]^-{\ad(\phi)} &&
W^*\otimes F_1^*\boxtimes F_2 \ar[r] & \CT,
\\
W\otimes F_2^*\boxtimes F_1 \ar[rr]^-{\ad(\phi)} &&
F_2^*\boxtimes F_2 \oplus F_1^*\boxtimes F_1 \ar[r] & \CT^*.
}}
\end{equation}

\begin{lemma}\label{tcc}
We have $\TCC_t \in \D^{[-1 - B,1]}(Y\times\BP\times Y)$
for all $t\in\ZZ$. Moreover,
$$
\TCC_t = \begin{cases}
\pi^*\CT^*\otimes \FD,          & \text{for $t=0$}\\
0,                              & \text{for $1\le t \le i-1$}\\
\pi^*\CT\otimes {\CL^*}^{\otimes i}\otimes\FD[\dim X - 1],
                                & \text{for $t=i$}
\end{cases}
$$
In particular, $\TCC_i \in \D^{[-B,1-\dim X]}(Y\times\BP\times Y)$
and $\CH^0(\TCC_0) \ne 0$.
\end{lemma}
\begin{proof}
Consider the diagram
\begin{equation}\label{diag_2}
\vcenter{\xymatrix{
&
Y \times \CX_1 \ar[d]^\alpha \ar@{}[dr]|\ecart &
Y \times \CX_1 \times Y \ar[d]^\alpha \ar[l]_{p_{12}} \ar[r]^{p_{23}} &
\CX_1 \times Y \ar[d]^\alpha \ar@{}[dl]|\ecart \\
Q(X,Y) \ar[d]^i \ar[ur]^j &
Y \times (X \times \BP)  &
Y \times (X \times \BP) \times Y \ar[l]_{p_{12}} \ar[r]^{p_{23}} \ar[d]^q &
(X \times \BP) \times Y &
Q(X,Y) \ar[d]^i \ar[ul]_j \\
Y \times X \ar[ur]^{\beta'} &
Y \times X \times Y \ar[ur]^{\beta'} \ar[l]_{p_{12}} \ar@{}[u]|\ecart &
Y\times\BP\times Y \ar[d]^\pi &
Y \times X \times Y \ar[ul]_\beta    \ar[r]^{p_{23}} \ar@{}[u]|\ecart &
X \times Y \ar[ul]_\beta \\
&& Y\times Y
}}
\end{equation}
and note that the squares marked with $\ecart$ are exact cartesian
by corollary~\ref{lci}.
Hence we have,
$$
\alpha_*p_{23}^*\CE_1 \cong
p_{23}^*\alpha_*\CE_1 =
p_{23}^*\alpha_*j_*\CE \cong
p_{23}^*\beta_*i_*\CE \cong
\beta_*p_{23}^*i_*\CE
$$
and similarly
$\alpha_*p_{12}^*\CE_1^{\#t} \cong
\beta'_*p_{12}^*i_*\CE^{\#t}\otimes{\CL^*}^{\otimes t}[\dim X - 1]$.
Therefore,
$$
\TCC_t \cong
q_*(\beta'_*p_{12}^*i_*\CE^{\#t}\otimes{\CL^*}^{\otimes t}[\dim X - 1] \otimes
\beta_*p_{23}^*i_*\CE).
$$
Applying $\beta_*p_{23}^*$ to~(\ref{equ_e}) and
$\beta'_*p_{12}^*$ to~(\ref{equ_ediez}) we obtain exact triangles
\begin{equation}\label{bpie}
\arraycolsep=1pt
\begin{array}{lllll}
\beta_*(\CO\boxtimes E_1\boxtimes F_1) & \exto{e} &
\beta_*(\CO\boxtimes E_2\boxtimes F_2) & \to &
\beta_*p_{23}^*i_*\CE, \\
\beta'_*(F_2^*\boxtimes E_2^*(t-i)\boxtimes \CO) & \exto{e^*} &
\beta'_*(F_1^*\boxtimes E_1^*(t-i)\boxtimes \CO) & \to &
\beta'_*p_{12}^*i_*\CE^{\#t}.
\end{array}
\end{equation}
So, roughly speaking we have a resolution
\begin{equation}\label{tccres}
\left\lbrace \TCC_t^{2,1} \to
\TCC_t^{1,1} \oplus \TCC_t^{2,2} \to
\TCC_t^{1,2} \right\rbrace \cong \TCC_t,
\end{equation}
where
\begin{multline*}
\TCC_t^{k,l} :=
q_*(\beta'_*(F_k^*\boxtimes E_k^*(t-i)\boxtimes \CO)
\otimes{\CL^*}^{\otimes t}[\dim X - 1] \otimes
\beta_*(\CO\boxtimes E_l\boxtimes F_l)) \cong  \\ \cong
q_*((F_k^* \boxtimes (E_k^*(t-i)\otimes E_l)
\boxtimes {\CL^*}^{\otimes t} \boxtimes F_l) \otimes
\beta'_*\CO_{Y\times X\times Y} \otimes
\beta_*\CO_{Y\times X\times Y})[\dim X - 1]
\end{multline*}
for $k,l = 1,2$ (more rigorously, we should say that there are three
exact triangles, obtained by tensoring the first triangle of~(\ref{bpie})
by terms of the second triangle). To compute $\TCC_t^{k,l}$
we note that the squares in the diagram
$$
\xymatrix{
Y \times X \times Y \ar[r]^{\beta'} \ar[d]_q &
Y \times (X \times \BP) \times Y \ar[d]_q &
Y \times X \times Y \ar[l]_{\beta} \ar[d]_q \\
Y \times Y \ar[r]^{\beta'} &
Y \times \BP \times Y &
Y \times Y \ar[l]_{\beta}
}
$$
are exact cartesian by corollary~\ref{flat_ec}, hence
$$
\beta'_*\CO_{Y\times X\times Y} \otimes \beta_*\CO_{Y\times X\times Y} \cong
\beta'_*q^*\CO_{Y\times Y} \otimes \beta_*q^*\CO_{Y\times Y} \cong
q^*\beta'_*\CO_{Y\times Y} \otimes q^*\beta_*\CO_{Y\times Y} \cong
q^*\FD.
$$
Substituting this into the formula for $\TCC_t^{k,l}$ we deduce
$$
\TCC^{k,l}_t \cong
(F_k^* \boxtimes
(H^\bullet(X,E_k^*(t-i)\otimes E_l) \otimes {\CL^*}^{\otimes t})
\boxtimes F_l) \otimes \FD)[\dim X - 1]
$$
Since
$\FD \in
\D^{[\codim_{Y\times Y} Y\times_\BP Y - \dim \BP,0]}(Y\times\BP\times Y)$
by lemma~\ref{h0fd},
$H^\bullet(X,E_k^*(t-i)\otimes E_l) \in \D^{[0,\dim X]}(\Ab)$, and
$(\codim_{Y\times Y} Y\times_\BP Y - \dim \BP) - (\dim X - 1) = 1 - B$,
by~(\ref{b}), we deduce
$\TCC_t^{k,l} \in \D^{[1-B,1]}(Y\times\BP\times Y)$,
for all $t\in\ZZ$, $k,l = 1,2$.
Hence, looking at~(\ref{tccres}) we see that
$\TCC_t \in \D^{[-1-B,1]}(Y\times\BP\times Y)$
for all $t\in\ZZ$. Moreover, for $1\le t\le i-1$ we have
$H^\bullet(X,E_k^*(t-i)\otimes E_l) \cong \Hom^\bullet(E_k,E_l(t-i)) = 0$,
hence $\TCC_t^{k,l} = 0$ for $k,l=1,2$,
therefore $\TCC_t = 0$ for $1\le t\le i-1$.

Finally, for $t = 0$ we have
$$
H^\bullet(X,E_k^*(t-i)\otimes E_l) \cong
\Hom^\bullet(E_k,E_l(-i)) \cong
\begin{cases}
\kk[-\dim X], & \text{for $k=l$} \\
W[-\dim X], & \text{for $l=1$, $k=2$} \\
0, & \text{for $l=2$, $k=1$}
\end{cases}
$$
hence
$\TCC_0^{2,1} = \pi^*(W\otimes F_2^*\boxtimes F_1)\otimes\FD[-1]$,
$\TCC_0^{1,1} = \pi^*(F_1^*\boxtimes F_1)\otimes\FD[-1]$,
$\TCC_0^{2,2} = \pi^*(F_2^*\boxtimes F_2)\otimes\FD[-1]$,
and (\ref{tccres}) gives a triangle
$$
\pi^*(W\otimes F_2^*\boxtimes F_1) \otimes \FD \to
\pi^*(F_2^*\boxtimes F_2 \oplus F_1^*\boxtimes F_1) \otimes \FD \to \TCC_0.
$$
whereof we deduce $\TCC_0 \cong \pi^*\CT^* \otimes \FD$.
Note also that $\CH^0(\CT^*) = \Coker \ad(\phi) \cong {\Delta_Y}_*\CO_Y$
by condition (C.5),
$\CH^0(\FD) \cong \CO_{Y\times_\BP Y}$
by lemma~\ref{h0fd},
so that $\CH^0(\TCC_0) \ne 0$.

Similarly, for $t=i$ we have
$$
H^\bullet(X,E_k^*(t-i)\otimes E_l) \cong
\Hom^\bullet(E_k,E_l) \cong
\begin{cases}
\kk, & \text{for $k=l$} \\
W^*, & \text{for $l=2$, $k=1$} \\
0, & \text{for $l=1$, $k=2$}
\end{cases}
$$
hence
$\TCC_i^{1,1} =
\pi^*(F_1^*\boxtimes F_1)\otimes{\CL^*}^{\otimes i}\otimes\FD[\dim X - 1]$,
$\TCC_i^{2,2} =
\pi^*(F_2^*\boxtimes F_2)\otimes{\CL^*}^{\otimes i}\otimes\FD[\dim X - 1]$,
and
$\TCC_i^{1,2} =
\pi^*(W\otimes F_1^*\boxtimes F_2)\otimes{\CL^*}^{\otimes i}\otimes\FD[\dim X - 1]$.
So, (\ref{tccres}) gives a triangle
$$
\pi^*(F_2^*\boxtimes F_2 \oplus F_1^*\boxtimes F_1) \otimes
{\CL^*}^{\otimes i}\otimes \FD [\dim X - 1] \exto{(e^*,e)}
\pi^*(W^*\otimes F_1^*\boxtimes F_2) \otimes
{\CL^*}^{\otimes i}\otimes \FD [\dim X - 1] \to \TCC_i.
$$
whereof we deduce
$\TCC_i \cong \pi^*\CT \otimes {\CL^*}^{\otimes i}\otimes \FD[\dim X - 1]$.
In particular, $\TCC_i \in \D^{[-B,1-\dim X]}(Y\times\BP\times Y)$ since
$\TCC_i^{2,1} = 0$.
\end{proof}

\subsection{Proof of theorem~\ref{ff}}

\begin{proposition}\label{cc0}
We have $B \ge 2i - 1$. Moreover, if $B = 2i - 1$, then
$$
\CC_0 = \CC_1[-2] = \dots = \CC_{i-1}[2-2i] \cong
\CH^0(\TCC_0)
$$
\end{proposition}
\begin{proof}
First of all, we note that for all $t\in \ZZ$ we have
\begin{equation}\label{lbcc1}
\CC_t \in \D^{[-B,0]}(Y\times\BP\times Y).
\end{equation}
Indeed, it follows from triangle~(\ref{tcccct}) and lemma~\ref{tcc} that
for any $l \le - 1 - B$ we have an exact sequence
$$
0 = \CH^{l-2}(\TCC_t) \to
\CH^{l-2}(\CC_t) \to
\CH^{l-2}(\CC_{t-1}[2]) \to
\CH^{l-1}(\TCC_t) = 0.
$$
Therefore,
$\CH^l(\CC_{t-1}) = \CH^{l-2}(\CC_{t-1}[2]) \cong \CH^{l-2}(\CC_t)$.
So, if $\CH^l(\CC_t) = 0$ for $l \le -1-B$ then
$\CH^{l-2s}(\CC_{t+s}) \ne 0$ for all $s$ which
contradicts the claim of lemma~\ref{lbcc}.

Further, using triangle~(\ref{tcccct}) and
lemma~\ref{tcc} for $t=i-1,\dots,1$ we see that
$$
\CC_{i-1} = \CC_{i-2}[2] = \dots = \CC_0[2i-2].
$$
Therefore, $\CC_{i-1}[2] = \CC_0[2i]$, hence
\begin{equation}\label{lbcci}
\CC_{i-1}[2] \in \D^{\le -2i}
\end{equation}
by lemma~\ref{lbcc}.
Using triangle~(\ref{tcccct}) for $t=i$:
$\{\TCC_i \to \CC_i \to \CC_{i-1}[2]\}$
and recalling that $\TCC_i \in \D^{\ge -B}$ by lemma~\ref{tcc}
and $\CC_i \in \D^{\ge -B}$ by~(\ref{lbcc1}), we deduce
\begin{equation}\label{rbcci}
\CC_{i-1}[2] \in \D^{\ge -1-B}.
\end{equation}
Now, if $B < 2i-1$, then $-2i < -1-B$,
and comparing (\ref{lbcci}) with (\ref{rbcci}) we deduce $\CC_{i-1} = 0$.
Hence $\CC_0 = \CC_{i-1}[2-2i] = 0$ as well.
But then the triangle~(\ref{tcccct}) for $t=0$:
$\{\TCC_0 \to \CC_0 \to \CC_{-1}[2]\}$ implies
$\CC_{-1} = \TCC_0[-1]$, hence
$\CH^1(\CC_{-1}) = \CH^0(\TCC_0) \ne 0$ by lemma~\ref{tcc},
which contradicts lemma~\ref{lbcc}. Therefore, $B \ge 2i-1$.

If $B = 2i - 1$, then $-2i = -1-B$, and comparing (\ref{lbcci})
with (\ref{rbcci}) we deduce that $\CC_{i-1}[2] \in \D^{-2i}$.
Hence $\CC_0 = \CC_{i-1}[2-2i] \in D^0$. Since by lemma~\ref{lbcc} we have
$\CC_{-1}[2] \in \D^{\le -2}$ the triangle~(\ref{tcccct}) for $t=0$
implies that $\CC_0 = \CH^0(\TCC_0)$.
\end{proof}

\begin{corollary}\label{picc}
If $B = 2i - 1$ then we have
$$
\pi_*\CC_0 = \pi_*\CC_1[-2] = \dots = \pi_*\CC_{i-1}[2-2i] \cong
{\Delta_Y}_*\CO_Y.
$$
\end{corollary}
\begin{proof}
Note that
$$
\pi_*\TCC_0 \cong
\pi_*(\pi^*\CT^*\otimes\FD) \cong
\CT^*\otimes\pi_*\FD.
$$
On the other hand, by lemma~\ref{h0fd} the object $\FD$ is supported
scheme-theoretically on a infinitesimal neighborhood of
$Y\times_\BP Y \subset Y\times\BP\times Y$. But the restriction of
$\pi:Y\times\BP\times Y \to Y\times Y$ to the infinitesimal neighborhood
of $Y\times_\BP Y \subset Y\times\BP\times Y$ is finite, hence
$$
\pi_*\CC_0 \cong
\pi_*\CH^0(\TCC_0) \cong
\CH^0(\pi_*\TCC_0) \cong
\CH^0(\CT^*\otimes \pi_*\FD).
$$
Finally, we have
$\CH^0(\CT^*) =
\Coker(\xymatrix@1{W\otimes F_2^*\boxtimes F_1 \ar[rr]^{\ad(\phi)} &&
F_2^*\boxtimes F_2 \oplus F_1^*\boxtimes F_1}) \cong
{\Delta_Y}_*\CO_Y$
by~(\ref{tts}), and
$\CH^0(\pi_*\FD) \cong
\CH^0(\pi_*\CO_{Y\times_\BP Y}) \cong
\CO_{Y\times_\BP Y}$
by lemma~\ref{h0fd},
hence
$$
\pi_*\CC_0 \cong
\CH^0({\Delta_Y}_*\CO_Y \otimes \CO_{Y\times_\BP Y}) \cong
{\Delta_Y}_*\CO_Y.
$$
since ${\Delta_Y}(Y) \subset Y\times_\BP Y$.
\end{proof}

\begin{corollary}\label{convt}
If $B = 2i - 1$ then the functor
$\Phi_{\CE_1^{\#t}}\circ\Phi_{\CE_1}$ is isomorphic to
the shift by $[2t]$ for any $t = 0, 1, \dots, i-1$.
\end{corollary}
\begin{proof}
Apply lemma~\ref{k} and corollary~\ref{picc}.
\end{proof}

\noindent{\bf Proof of theorem~\ref{ff}.}

\noindent
We have $B\ge 2i-1$ by proposition~\ref{cc0}. Assume that $B = 2i - 1$.
Since $\Phi_{\CE_1^{\#0}}$ is left adjoint to $\Phi_{\CE_1}$
by lemma~\ref{ed_la}, it follows from corollary~\ref{convt}
that $\Phi_{\CE_1}$ is fully faithful and it remains to show
that $Y$ is smooth.
To this end we note that $\CX_1$ is smooth by lemma~\ref{x1},
hence $\D^b(\CX_1)$ is $\Ext$-bounded by lemma~\ref{snav_sm}.
On the other hand, as we have shown above $\D^b(Y)$ embeds
fully and faithfully into $\D^b(\CX_1)$, hence $\D^b(Y)$ is also
$\Ext$-bounded. Applying again lemma~\ref{snav_sm} we deduce
that $Y$ is smooth.
\qed

\medskip

We will need also the following

\begin{proposition}\label{homephi}
If $B = 2i - 1$ then
$$
\Phi_1(\D^b(Y)) \subset
\langle E_1(1)\otimes\D^b(\BP),\dots,E_2(i-1)\otimes\D^b(\BP)\rangle^\perp.
$$
\end{proposition}
\begin{proof}
We must check that $\Hom(E_s(k)\otimes F,\Phi_{\CE_1}(G)) = 0$
for all $F\in\D^b(\BP)$, $G\in\D^b(Y)$. Note that
$E_s(k)\otimes F \cong \alpha^*(E_s(k)\boxtimes F)$ and
$\Hom(\alpha^*(E_s(k)\boxtimes F),\Phi_{\CE_1}(G)) \cong
\Hom(E_s(k)\boxtimes F,\alpha_*\Phi_{\CE_1}(G))$.
Further, $\alpha_*\Phi_{\CE_1}$ is a kernel functor
with kernel $\alpha_*j_*\CE \cong \beta_*i_*\CE$ and
\begin{multline*}
\Hom_{X\times\BP}(E_s(k)\boxtimes F,\Phi_{\beta_*i_*\CE}(G)) =
\Hom_{X\times\BP}(E_s(k)\boxtimes F,p_*(\beta_*i_*\CE\otimes q^*G)) \cong
\\ \cong
\Hom_{X\times\BP\times Y}(E_s(k)\boxtimes F \boxtimes G^*,\beta_*i_*\CE) \cong
\Hom_{X\times Y}(E_s(k)\boxtimes (g^*F\otimes G^*),i_*\CE)
\end{multline*}
and the last space is zero because
$i_*\CE \in \langle E_1\boxtimes\D^b(Y),E_2\boxtimes\D^b(Y)\rangle$
and by condition (C.3) we have
$E_1,E_2 \in \langle E_1(1),\dots,E_2(i-1)\rangle^\perp$.
\end{proof}

\section{Semiorthogonal decompositions}

In this section we conclude the proof of theorem~\ref{themain}.
In fact, instead of individual linear sections of $X$ and $Y$
we consider the {\em universal families}\/ of linear sections
and prove a similar statements for them. After that theorem~\ref{themain}
follows by a faithful base change argument.

\subsection{Universal families of linear sections}

All $r$-dimensional subspaces $L \subset V^*$ are parameterized
by the Grassmannian $\Gr(r,V^*)$. Those of them, for which
the linear sections $X_L = X\times_{\PP(V)}\PP(L^\perp)$ and
$Y_L = Y\times_{\PP(V^*)}\PP(L)$ are compact are parameterized
by an open subset of the Grassmannian $\BP_r \subset \Gr(r,V^*)$
consisting of all $r$-dimensional subspaces $L\subset V^*$ such that
$L\cap\BZ = \emptyset$.
Let $\CL_r$ be the pullback of the tautological rank $r$ subbundle
of the Grassmannian $\Gr(r,V^*)$ to $\BP_r$ and let
$\CL_r^\perp:=(V^*\otimes\CO_{\BP_r}/\CL_r)^* \subset V\otimes\CO_{\BP_r}$,
be the orthogonal subbundle.
Then the subvarieties
$$
\arraycolsep=2pt
\begin{array}{lll}
\CX_r & =
(X\times\BP_r)\times_{\PP(V)\times\BP_r}\PP_{\BP_r}(\CL_r^\perp) &
\subset X\times\BP_r, \\
\CY_r & =
(Y\times\BP_r)\times_{\PP(V^*)\times\BP_r}\PP_{\BP_r}(\CL_r) &
\subset Y\times\BP_r,
\end{array}
$$
are the universal families of compact linear sections of $X$ and $Y$.

Consider the fiber product $\CX_r\times_{\BP_r}\CY_r$ and the projection
$\pi_r:\CX_r\times_{\BP_r}\CY_r \to X\times Y$. Since for any vector subspace
$L\subset V^*$ the product $\PP(L^\perp)\times\PP(L)$ is contained
in the incidence quadric $Q \subset \PP(V)\times\PP(V^*)$ it follows
that $\pi_r$ factors via a map
$\zeta_r:\CX_r\times_{\BP_r}\CY_r \to Q(X,Y) \subset X\times Y$.

Consider the object $\CE_r = \zeta_r^*\CE \in \D(\CX_r\times_{\BP_r}\CY_r)$
as a kernel on $\CX_r\times\CY_r$. It gives the following kernel functors
$\Phi_r = \Phi_{\CE_r}:\D^b(\CY_r) \to \D^b(\CX_r)$ and
$\Phi_r^! = \Phi_{\CE_r}^!:\D^b(\CX_r) \to \D^b(\CY_r)$.
The goal of this section is to show that
the functor $\Phi_r$ is fully faithful for $r\le i$,
the functor $\Phi_r^!$ is fully faithful for $r\ge i$
and that they give the following semiorthogonal decompositions
$$
\arraycolsep = 2pt
\begin{array}{ll}
\D^b(\CX_{r}) & = \langle
\Phi_{r}(\D^b(\CY_{r})),
E_1(1)\otimes\D^b(\BP_{r}),
\dots,
E_2(i-r)\otimes\D^b(\BP_{r})
\rangle,\smallskip\\
\D^b(\CY_{r}) & = \langle
\Phi^!_{r}(\D^b(\CX_{r})),
F^*_2(i-r)\otimes\D^b(\BP_{r}),
\dots,
F^*_1(-1)\otimes\D^b(\BP_{r})
\rangle.
\end{array}
$$
After that we deduce from this the main theorem~\ref{themain}
using the faithful base change theorem~\ref{phitsod}.

To prove the fullness and faithfulness of the functors
we use induction in $r$ (the base is given by theorem~\ref{ff}).
To compare the universal families $\CX_{r-1}$, $\CY_{r-1}$ and
$\CX_r$, $\CY_r$ we take
for a base scheme
$$
\BS_r = \Fl(r-1,r;V^*)\cap (\BP_{r-1}\times\BP_r),
$$
where $\Fl(r-1,r;V^*) \subset \Gr(r-1,V^*)\times\Gr(r,V^*)$
is the partial flag variety. The scheme $\BS_r$ parameterizes
flags $L_{r-1}\subset L_r\subset V^*$ such that $\dim L_{r-1} = r-1$,
$\dim L_r = r$ and $\PP(L_r) \cap \BZ = \emptyset$.

Let $\phi:\BS_r\to\BP_{r-1}$ and
$\psi:\BS_r\to\BP_r$ denote the natural projections.
Let $\TCL_{r-1} = \phi^*\CL_{r-1}$, $\TCL_{r} = \psi^*\CL_{r}$,
$\TCL^\perp_{r-1} = \phi^*\CL_{r-1}^\perp$,
$\TCL^\perp_{r} = \psi^*\CL^\perp_{r}$.
Then we have
\begin{equation}\label{tclrss}
\TCL_{r-1} \subset \TCL_r \subset V^*\otimes\CO_{\BS_r}, \qquad
\TCL_r^\perp \subset \TCL_{r-1}^\perp \subset V\otimes\CO_{\BS_r},
\end{equation}
Denote
$$
\arraycolsep=2pt
\begin{array}{llllll}
\TCX_{r-1} & =
\CX_{r-1}\times_{\BP_{r-1}}\BS_r  &
\subset X\times\BS_r, \qquad\qquad &
\TCX_r & =
\CX_r\times_{\BP_r}\BS_r &
\subset X\times\BS_r, \smallskip\\
\TCY_{r-1} & =
\CY_{r-1}\times_{\BP_{r-1}}\BS_r  &
\subset Y\times\BS_r, &
\TCY_r & =
\CY_r\times_{\BP_r}\BS_r  &
\subset Y\times\BS_r.
\end{array}
$$
Note that
$$
\arraycolsep=2pt
\begin{array}{llll}
\TCX_{r-1} &
= (X\times\BS_r) \times_{\PP(V)\times\BS_r}\PP_{\BS_r}(\TCL_{r-1}^\perp),
\qquad\qquad &
\TCX_r &
= (X\times\BS_r) \times_{\PP(V)\times\BS_r}\PP_{\BS_r}(\TCL_{r}^\perp),
\smallskip\\
\TCY_{r-1} &
= (Y\times\BS_r) \times_{\PP(V^*)\times\BS_r}\PP_{\BS_r}(\TCL_{r-1}),
&
\TCY_r &
= (Y\times\BS_r) \times_{\PP(V^*)\times\BS_r}\PP_{\BS_r}(\TCL_{r}).
\end{array}
$$
Therefore the embeddings~(\ref{tclrss}) induce embeddings
$\xi:\TCX_r \to \TCX_{r-1}$ and $\eta:\TCY_{r-1} \to \TCY_r$.
Consider the following commutative diagrams
(the squares marked with $\ecart$ are exact cartesian,
because the maps $\phi$ and $\psi$ are flat)
\begin{equation}\label{yayb}
\vcenter{\xymatrix{
\BS_r \ar[d]_\psi \ar@{}[dr]|{\ecart} &
\TCY_r \ar[l]_{g_r} \ar[d]_\psi &
\TCY_{r-1} \ar[d]^\phi \ar[l]_\eta \ar[dl]_{\hat{\psi}} \ar[r]^{g_{r-1}} &
\BS_r \ar[d]_\phi \ar@{}[dl]|{\ecart} \\
\BP_r &
\CY_r \ar[l]_{g_r} &
\CY_{r-1} \ar[r]^{g_{r-1}}  &
\BP_{r-1}
}}
\qquad\qquad
\vcenter{\xymatrix{
\BS_r \ar[d]_\psi \ar@{}[dr]|{\ecart} &
\TCX_r \ar[l]_{f_r} \ar[d]_\psi \ar[r]^\xi \ar[dr]^{\hat{\phi}} &
\TCX_{r-1} \ar[d]^\phi \ar[r]^{f_{r-1}} &
\BS_r \ar[d]_\phi \ar@{}[dl]|{\ecart}  \\
\BP_r &
\CX_r \ar[l]_{f_r} &
\CX_{r-1} \ar[r]^{f_{r-1}}  &
\BP_{r-1}
}}
\end{equation}
where $f_{r-1}$, $f_r$, $g_{r-1}$ and $g_r$ are the natural projections
and $\hat{\phi} = \phi\circ\xi$, $\hat{\psi} = \psi\circ\eta$.

Let $\TE_{r-1}\in\D(\TCX_{r-1}\times_{\BS_r}\TCY_{r-1})$ and
$\TE_r \in \D(\TCX_r\times_{\BS_r}\TCY_r)$ denote the pullbacks
of the objects $\CE_{r-1}$ and $\CE_r$ via the projections
$\TCX_{r-1}\times_{\BS_r}\TCY_{r-1} \to \CX_{r-1}\times_{\BP_{r-1}}\CY_{r-1}$,
$\TCX_r\times_{\BS_r}\TCY_r \to \CX_r\times_{\BP_r}\CY_r$.
Then we have the corresponding kernel functors
$\TPhi_{r-1}$, $\TPhi_r$ e.t.c between the derived categories
of $\TCX_{r-1}$, $\TCY_{r-1}$, $\TCX_r$ and $\TCY_r$.

The induction step is based on relation of the functors
$\Phi_{r-1}$, $\Phi_r$, $\TPhi_{r-1}$ and $\TPhi_r$
to the base change functors $\psi^*$, $\psi_*$
and to the functors of the pushforward and pullback
via $\xi$ and $\eta$. The relation to $\psi^*$ and $\psi_*$ is given
by lemma~\ref{phit}. The relation to $\xi$ and $\eta$
in a sense is the key point of the proof.
We prove that $\xi^*\TPhi_{r-1} \cong \TPhi_r\eta_*$
and that the ``difference'' between $\xi_*\TPhi_r$ and $\TPhi_{r-1}\eta^!$
is given by a very simple functor, the kernel of which has a resolution
of the form $E_1\boxtimes F_1 \to E_2\boxtimes F_2$ (up to a twist and a shift)
on $\TCX_{r-1}\times_{\BS_r}\TCY_r$.

Other results in this section (e.g.\ the above semiorthogonal
decompositions) are proved by similar arguments using (either
ascending or descending) induction in $r$.

The section is organized as follows. We start with some preparations
concluding with a description of the relation of the functors
$\TPhi_{r-1}$ and $\TPhi_r$ to the pushforward and pullback
via $\xi$ and $\eta$. Then we use induction in $r$ to prove that
$\Phi_r$ is fully faithful for $r\le i$. Then we use a relative
Bridgeland's trick~\ref{ffisequ} to establish a semiorthogonal
decomposition for $r=i$. Then we use descending induction in $r$
to establish semiorthogonal decompositions for $\CX_r$ when $r\le i$.
Then we show that the collection $(F_1(1),F_2(1),\dots,F_1(N-i),F_2(N-i))$
on $Y$ is an exceptional collection in a certain sense (if $\BZ=\emptyset$
then it is exceptional in the usual sense). Then we use induction in $r$
to establish semiorthogonal decompositions for $\CY_r$ when $r\ge i$.
Finally we deduce theorem~\ref{themain} and show that the set of
critical values of the morphism $g:Y\to\PP(V^*)\setminus\BZ$ coincides
with $X^\vee\setminus\BZ$.

\subsection{Preparations}\label{indsetup}

Recall that $\BP_r\subset\Gr(r,V^*)$ is an open subset
parameterizing the compact linear sections of $X$ and $Y$,
and $\CX_r$, $\CY_r$ are the universal families over $\BP_r$
of linear sections.

\begin{lemma}\label{xyksm}
Assume that $\BP_r \ne \emptyset$.
Then $\CX_r$ and $\CY_r$ are smooth,
$$
\begin{array}{l}
\dim\CX_r = \dim X + \dim\BP_r - r,\qquad
\dim\CY_r = \dim Y + \dim\BP_r + r - N,\\
\dim\CX_r\times_{\BP_r}\CY_r = \dim X + \dim Y + \dim\BP_r - N,
\end{array}
$$
and the maps $f_r:\CX_r \to \BP_r$ and $g_r:\CY_r \to \BP_r$ are projective.
\end{lemma}
\begin{proof}
Note that we have open embeddings into the relative Grassmannians
$$
\CX_r \subset \Gr_X(r,\CV_X),\qquad
\CY_r \subset \Gr_Y(r-1,\CV_Y),\qquad
\CX_r\times_{\BP_r}\CY_r \subset \Gr_{Q(X,Y)}(r-1,\CV_Q),
$$
where the bundles $\CV_X$, $\CV_Y$ are defined from exact sequences
$$
0 \to \CV_X \to V^*\otimes\CO_X \to \CO_X(1) \to 0,\qquad
0 \to \CO_Y(-1) \to V^*\otimes\CO_Y \to \CV_Y \to 0,
$$
and $\CV_Q$ is the middle cohomology bundle of the complex
$$
\CO_{Q(X,Y)}(0,-1) \to V^*\otimes\CO_{Q(X,Y)} \to \CO_{Q(X,Y)}(1,0).
$$
From this we easily deduce the smoothness and compute the dimensions.
It is also clear that the fibers of the projections $\CX_r \to \BP_r$, and
$\CY_r \to \BP_r$ are linear sections of $X$ and $Y$ corresponding
to subspaces $L \in \BP_r$, so they are projective.
\end{proof}

Recall that we have defined the objects $\CE_r$ on
$\CX_r\times_{\BP_r}\CY_r$ as the pullbacks of $\CE\in\D^b(Q(X,Y))$
via the map $\zeta_r:\CX_r\times_{\BP_r}\CY_r \to Q(X,Y)$,
and the objects $\TE_{r-1}$ and $\TE_r$ as the pullbacks
of $\CE_{r-1}$ and $\CE_r$ via the maps
$\phi:\TCX_{r-1}\times_{\BS_r}\TCY_{r-1} \to
\CX_{r-1}\times_{\BP_{r-1}}\CY_{r-1}$ and
$\psi:\TCX_r\times_{\BS_r}\TCY_r \to \CX_r\times_{\BP_r}\CY_r$.
The functors $\Phi_{r-1}$, $\Phi_r$, $\Phi_{r-1}^!$, $\Phi_r^!$,
$\TPhi_{r-1}$, $\TPhi_r$, $\TPhi_{r-1}^!$, and $\TPhi_r^!$ are the
kernel functors of the first and second type corresponding
to the kernels $\CE_{r-1}$, $\CE_r$, $\TE_{r-1}$ and $\TE_r$ respectively.

\begin{lemma}\label{b_adj}
The functors $\Phi_{r-1}$, $\Phi_r$, $\Phi_{r-1}^!$, $\Phi_r^!$,
$\TPhi_{r-1}$, $\TPhi_r$, $\TPhi_{r-1}^!$, and $\TPhi_r^!$,
take the bounded derived categories to the bounded derived categories.
Moreover, the functors with the shriek are right adjoint
to the corresponding functors without the shriek.
\end{lemma}
\begin{proof}
Since $\CX_r$ and $\CY_r$ are smooth it follows that the pushforward
of $\CE_r$ to $\CX_r\times\CY_r$ is a perfect complex. Therefore,
$\CE_r$ has finite $\Tor$ and $\Ext$-amplitude over $\CX_r$ and $\CY_r$
by lemma~\ref{isfted}. Similarly, $\CE_{r-1}$ has finite $\Tor$ and
$\Ext$-amplitude over $\CX_{r-1}$ and $\CY_{r-1}$. Since $\phi$ and $\psi$ are
flat it follows from corollary~\ref{fted_bc} that $\TE_{r-1}$ and $\TE_r$ also
have finite $\Tor$ and $\Ext$-amplitude. On the other hand, the projections
of $\CX_{r-1}\times_{\BP_{r-1}}\CY_{r-1}$, $\CX_r\times_{\BP_r}\CY_r$,
$\TCX_{r-1}\times_{\BS_r}\TCY_{r-1}$ and $\TCX_r\times_{\BS_r}\TCY_r$
to the factors are projective because the projections of $\CX_{r-1}$,
$\CY_{r-1}$ to $\BP_{r-1}$ and of $\CX_r$, $\CY_r$ to $\BP_r$ are.
It remains to apply lemma~\ref{phi_bounded}.
\end{proof}

\begin{lemma}\label{zl1}
Let $k=r-1$ or $k=r$. We have

\noindent
$(i)$ $\TCX_k$ is the zero locus of a section
of vector bundle $\CO_X(1)\boxtimes\TCL_k^*$ on $X\times\BS_r$;

\noindent
$(ii)$
$\TCY_k$ is the zero locus of a section
of vector bundle $\CO_Y(1)\boxtimes\TCL_l^{\perp*}$ on $Y\times\BS_r$;

\noindent
$(iii)$
$\TCX_r$ is the zero locus of a section of line bundle
$\CO_X(1)\otimes(\TCL_r/\TCL_{r-1})^*$ on $\TCX_{r-1}$;

\noindent
$(iv)$
$\TCY_{r-1}$ is the zero locus of a section of line bundle
$\CO_Y(1)\otimes(\TCL_{r-1}^\perp/\TCL_r^\perp)^* \cong
\CO_Y(1)\otimes(\TCL_r/\TCL_{r-1})$ on $\TCY_r$.\\
All these sections are regular.
\end{lemma}
\begin{proof}
The parts $(i)$ and $(ii)$ evidently follow from the definition
of $\TCX_k \subset X\times\BS_r$. The parts $(iii)$ and $(iv)$
follow from the exact sequences
$$
0 \to (\TCL_r/\TCL_{r-1})^* \to \TCL_r^* \to \TCL_{r-1}^* \to 0,
\qquad
0 \to (\TCL_{r-1}^\perp/\TCL_r^\perp)^* \to
(\TCL_{r-1}^\perp)^* \to (\TCL_r^\perp)^* \to 0.
$$
Finally, it follows from lemma~\ref{xyksm} that
$\dim\TCX_{r-1} = \dim X + \dim\BS_r - (r-1)$,
$\dim\TCX_r = \dim X + \dim\BS_r - r$. Therefore the sections
in the parts~$(i)$ and $(iii)$ are regular. The sections
in the parts~$(ii)$ and $(iv)$ are regular by similar reasons.
\end{proof}

Now we describe the maps $\psi$ and $\hat{\psi}$.

\begin{lemma}\label{im1}
The maps $\psi$ and $\hat{\psi}$ are projectivizations of vector bundles.
Explicitly, $\psi$ is the projectivization of $\CL_r^*$ and
$\hat\psi$ is the projectivization of $(\CL_r/\CO_{\PP(V^*)}(-1))^*$,
where the embedding $\CO_{\PP(V^*)}(-1) \to \CL_r$ is induced
by the projection $\CY_r \to \PP_{\BP_r}(\CL_r) \to \PP(V^*)$.
\end{lemma}
\begin{proof}
By definition of $\BS_r$ the fiber of $\psi$ is the set of all
hyperplanes in $L_r$. Similarly, the fiber of $\hat\psi$ is the set
of all hyperplanes in $L_r$ passing through a point $y\in\PP(L_r)$.
\end{proof}

Applying results of~\cite{O2} we deduce the following.

\begin{corollary}\label{cim1}
The functors $\psi^*$ and ${\hat\psi}^*$ are fully faithful and we have
$\psi_*\psi^* \cong \id$, ${\hat\psi}_*{\hat\psi}^* \cong \id$.
\end{corollary}

The situation with morphisms $\phi$ and $\hat\phi$ is
a little bit more involved because of the condition
$L\cap\BZ = \emptyset$ entering in the definition of $\BP_r$.
In fact, it is easy to see that both maps factor as a composition
of an open embedding followed by the projectivization of a vector bundle.
More precisely, we have

\begin{lemma}\label{im2}
The map $\phi$ is a composition of an open embedding
with the projectivization of the vector bundle $V^*/\CL_{r-1}$.
Moreover, any divisorial component of the complement
$\PP_{\BP_{r-1}}(V^*/\CL_{r-1}) \setminus \BS_r$
is ample over $\BP_{r-1}$.
Similarly, the map $\hat\phi$ is a composition of an open embedding
with the projectivization of the vector bundle $\CV_X/\CL_{r-1}$,
where $\CV_X$ was defined in lemma~$\ref{xyksm}$.
Moreover, the map $\hat\phi$ is surjective if $r\le i$.
\end{lemma}
\begin{proof}
By definition of $\BS_r$ the fiber of $\phi$ is the set of all
$L_r$ containing the given $L_{r-1}$ and nonintersecting with $\BZ$.
Thus it is an open subset of $\PP(V^*/L_{r-1})$. Moreover, the complement
coincides with the preimage of the set of all $L_r \in \Gr(r,V^*)$ such that
$\PP(L_r)\cap\BZ\ne\emptyset$ and it remains to note that any effective
divisor in $\Gr(r,V^*)$ is ample.

Similarly, the fiber of $\hat{\phi}$ is the set of all
$L_r$ containing the given $L_{r-1}$, contained in the fiber $\CV_x$
of the bundle $\CV_X$ at the given point $x$ and nonintersecting with $\BZ$.
Thus it is an open subset of $\PP(\CV_X/L_{r-1})$.
It remains to check the surjectivity.
Assume that a point $(x,L_{r-1})\in\CX_{r-1}$ doesn't lie
in the image of~$\TCX_r$. This means that for any subspace
$L_r\subset V^*$ such that $L_{r-1}\subset L_r\subset \CV_x$ we have
$\PP(L_r)\cap\BZ \ne \emptyset$. Therefore the linear projection
of $\BZ\cap\PP(\CV_x) \subset \PP(\CV_x)\setminus\PP(L_{r-1})$ to
$\PP(\CV_x/L_{r-1})$ is surjective. But then we have
$\dim(\BZ\cap\PP(\CV_x)) \ge
\dim(\PP(\CV_x/L_{r-1})) =
(N-1) - (r-1) - 1 =
N - r - 1 \ge
N - i -1$
contradicting~(C.8) since $\CV_x$ is a hyperplane in $V^*$.
\end{proof}

The functors $\phi^*$ and ${\hat\phi}^*$ are not fully faithful.
However, they enjoy the following properties

\begin{corollary}\label{cim2}
If ${\phi}^*\CL \cong \CO_{\TCY_{r-1}}$
for some line bundle $\CL$ on $\CY_{r-1}$ and $r\le i$ then
$\CL \cong \CO_{\CY_{r-1}}$.\\
If $\hat{\phi}^*F = 0$ for some $F\in\D(\CX_{r-1})$
and $r\le i$ then $F = 0$.
\end{corollary}
\begin{proof}
For the first claim note that the kernel of the map
$\Pic(\PP_{\CY_{r-1}}(V^*/\CL_{r-1})) \to \Pic(\TCY_{r-1})$ is
generated by divisorial components of the complement
$\PP_{\CY_{r-1}}(V^*/\CL_{r-1}) \setminus \TCY_{r-1}$,
which are ample over $\TCY_{r-1}$. Therefore, this kernel
intersects the image of $\Pic(\CY_{r-1})$ in
$\Pic(\PP_{\CY_{r-1}}(V^*/\CL_{r-1}))$ only at zero.
The second claim is evident.
\end{proof}

Now we go to the relation of the functors $\TPhi_{r-1}$ and $\TPhi_r$
to the pushforward and pullback via $\xi$ and $\eta$.

Consider the following diagram
$$
\xymatrix{
\TCX_r\times_{\BS_r}\TCY_{r-1} \ar[r]^{\eta} \ar[d]_{\xi} &
\TCX_r\times_{\BS_r}\TCY_r \ar[d]_{\xi} \\
\TCX_{r-1}\times_{\BS_r}\TCY_{r-1} \ar[r]^{\eta} &
\TCX_{r-1}\times_{\BS_r}\TCY_r
}
$$
and the projection
$\tpi:\TCX_{r-1}\times_{\BS_r}\TCY_r \subset
(X\times\BS_r)\times(Y\times\BS_r) \to
X\times Y$.

\begin{lemma}
The maps $\xi$ and $\eta$ in the above diagram are divisorial embeddings,
we have the following scheme-theoretical equalities
\begin{equation}\label{capcup}
\begin{array}{l}
(\TCX_{r-1}\times_{\BS_r}\TCY_{r-1}) \cap (\TCX_r\times_{\BS_r}\TCY_r) =
\TCX_r\times_{\BS_r}\TCY_{r-1},\\
(\TCX_{r-1}\times_{\BS_r}\TCY_{r-1}) \cup (\TCX_r\times_{\BS_r}\TCY_r) =
\tpi^{-1}(Q(X,Y)).
\end{array}
\end{equation}
and the following square is exact cartesian
\begin{equation}\label{xaybs}
\vcenter{\xymatrix{
(\TCX_{r-1}\times_{\BS_r}\TCY_{r-1}) \cup (\TCX_r\times_{\BS_r}\TCY_r)
\ar[rr]^-i \ar[d]^{\tzeta} &&
\TCX_{r-1}\times_{\BS_r}\TCY_r \ar[d]^{\tpi} \\
Q(X,Y) \ar[rr]^i &&
X\times Y \ar@{}[ull]|{\ecart}
}}
\end{equation}
\end{lemma}
\begin{proof}
Consider the projections of
$\TCX_r\times_{\BS_r}\TCY_{r-1}$,
$\TCX_{r-1}\times_{\BS_r}\TCY_{r-1}$,
$\TCX_r\times_{\BS_r}\TCY_r$ and
$\TCX_{r-1}\times_{\BS_r}\TCY_r$
to $X\times Y$. It is easy to check that
their fibers over a point $(x,y) \in (X,Y)$ are open subsets
in the subsets of the flag variety $\Fl(r-1,r;V^*)$
consisting of all flags $L_{r-1} \subset L_r$ satisfying
the following incidence conditions
$$
\arraycolsep=2pt
\begin{array}{ccccc}
y & \subset & L_{r-1} \\
&& \cap \\
&& L_r & \subset & \CV_x
\end{array}
\qquad,\qquad
\begin{array}{ccccc}
y & \subset & L_{r-1} & \subset & \CV_x \\
&& \cap \\
&& L_r
\end{array}
\qquad,\qquad
\begin{array}{ccccc}
&& L_{r-1} \\
&& \cap \\
y & \subset & L_r & \subset & \CV_x
\end{array}
\qquad\text{and}\qquad
\begin{array}{ccccc}
&& L_{r-1} & \subset & \CV_x \\
&& \cap \\
y & \subset & L_r
\end{array}
$$
respectively. In particular, the first three fibers are empty
if $(x,y)\not\in Q(X,Y)$. On the other hand, over $Q(X,Y)$
the first three fibers are irreducible, have dimension
$(r-1)(N-r)-1$, $(r-1)(N-r)$ and $(r-1)(N-r)$ respectively,
and the first of them is the intersection of the other two.
On the contrary, the fourth fiber is irreducible and $(r-1)(N-r)$-dimensional
if $(x,y)\not\in Q(X,Y)$ and for $(x,y)\in Q(X,Y)$ it coincides with the union
of the second and the third fibers (if $y\subset \CV_x$ and
$y\not\subset L_{r-1}$ then $L_r = \langle y, L_{r-1}\rangle \subset \CV_x$).
It follows that images of $\xi$ and $\eta$ have pure codimension~$1$.
Since they are also zero loci of line bundles by lemma~\ref{zl1}~$(iii)$
and $(iv)$, we conclude that $\xi$ and $\eta$ are divisorial embeddings.

The above arguments also prove the first equality of~(\ref{capcup})
on the scheme-theoretical level and the second equality on the set-theoretical
level. Taking into account that the LHS of the second equality is the zero
locus of the line bundle $\CO_X(1)\otimes\CO_Y(1)$ by definition of $Q(X,Y)$,
and that the RHS of the equality is the zero locus of the line bundle
$(\CO_X(1)\otimes(\TCL_r/\TCL_{r-1})^*)\otimes
(\CO_Y(1)\otimes(\TCL_r/\TCL_{r-1}))$
by lemma~\ref{zl1}~$(iii)$ and $(iv)$, and noting that the bundles
are isomorphic, we deduce that the second equality is also true
on the scheme-theoretical level.
Finally, we note that the square~(\ref{xaybs}) is exact cartesian
by lemma~\ref{lci} since $X\times Y$ and $Q(X,Y)$ are Cohen-Macaulay.
\end{proof}

Consider the pullback $\HE = \tzeta^*\CE$ of $\CE$ from $Q(X,Y)$ to
$(\TCX_{r-1}\times_{\BS_r}\TCY_{r-1}) \cup (\TCX_r\times_{\BS_r}\TCY_r)$.
Let us also denote $\CO(k,l) := \CO_X(k)\otimes\CO_Y(l)$ for brevity.
The following lemma gives a relation of $\TE_{r-1}$ and $\TE_r$.

\begin{lemma}
We have the following exact sequences on
$\TCX_{r-1}\times_{\BS_r}\TCY_r$:
\begin{equation}\label{ieab}
0 \to E_1\otimes F_1 \to E_2\otimes F_2 \to i_*\HE \to 0.
\end{equation}
\begin{equation}\label{ieab1}
0 \to
\eta_*\TE_{r-1}(-1,0)\otimes(\TCL_r/\TCL_{r-1}) \to
i_*\HE \to
\xi_*\TE_r \to 0,
\end{equation}
\begin{equation}\label{ieab2}
0 \to
\xi_*\TE_r(0,-1)\otimes(\TCL_r/\TCL_{r-1})^* \to
i_*\HE \to
\eta_*\TE_{r-1} \to 0.
\end{equation}
Moreover, we have an isomorphism on $\TCX_r\times_{\BS_r}\TCY_{r-1}$:
\begin{equation}\label{eaeb}
\eta^*\TE_r \cong \xi^*\TE_{r-1}.
\end{equation}
\end{lemma}
\begin{proof}
Since the square~(\ref{xaybs}) is exact cartesian we have
$i_*\HE = i_*\tzeta^*\CE = \tpi^*i_*\CE$
and applying the functor $\tpi^*$ to (\ref{equ_e}) we deduce~(\ref{ieab}).
Sequences (\ref{ieab1}) and (\ref{ieab2}) can be obtained by tensoring
resolutions
$$
0 \to \CO_{\TCX_{r-1}\times_{\BS_r}\TCY_{r-1}}(-1,0)\otimes(\TCL_r/\TCL_{r-1})
\to
\CO_{(\TCX_{r-1}\times_{\BS_r}\TCY_{r-1}) \cup (\TCX_r\times_{\BS_r}\TCY_r)}
\to
\CO_{\TCX_r\times_{\BS_r}\TCY_r} \to 0,
$$
$$
0 \to \CO_{\TCX_r\times_{\BS_r}\TCY_r}(0,-1)\otimes(\TCL_r/\TCL_{r-1})^*
\to
\CO_{(\TCX_{r-1}\times_{\BS_r}\TCY_{r-1}) \cup (\TCX_r\times_{\BS_r}\TCY_r)}
\to
\CO_{\TCX_{r-1}\times_{\BS_r}\TCY_{r-1}} \to 0
$$
with $\HE$ and applying $i_*$, since the pullback of $\HE$ to
$\TCX_{r-1}\times_{\BS_r}\TCY_{r-1}$ and $\TCX_r\times_{\BS_r}\TCY_r$
coincides with $\TE_{r-1}$ and $\TE_r$ respectively.
Finally, (\ref{eaeb}) is evident, because both sides are isomorphic
to the pullback of $\CE$.
\end{proof}

\begin{corollary}\label{phixieta1}
We have the following exact triangles of functors between
$\D^b(\TCX_{r-1})$ and $\D^b(\TCY_r)$:
\begin{equation}\label{ft1}
\TPhi_{r-1}\eta^! \to \xi_*\TPhi_r \to
\Phi_{i_*\HE(0,1)\otimes(\TCL_r/\TCL_{r-1})},
\end{equation}
\begin{equation}\label{ft2}
\TPhi_r^!\xi^* \to \eta_*\TPhi_{r-1}^! \to
\Phi^!_{i_*\HE(1,0)\otimes(\TCL_r/\TCL_{r-1})^*},
\end{equation}
and the following canonical isomorphisms of functors between
$\D^b(\TCX_r)$ and $\D^b(\TCY_{r-1})$:
$$
\xi^*\TPhi_{r-1} \cong \TPhi_r\eta_*,\qquad
\eta^!\TPhi_r^! \cong \TPhi_{r-1}^!\xi_*.
$$
\end{corollary}
\begin{proof}
Apply lemma~\ref{etf} to the exact triangles~(\ref{ieab1})
and~(\ref{ieab2}) and lemma~\ref{mkf} to the isomorphisms~(\ref{eaeb}).
\end{proof}

\begin{lemma}\label{phiise}
We have

\noindent$(i)$
$\Phi_{i_*\HE(0,1)\otimes(\TCL_r/\TCL_{r-1})}
(\langle F^*_1(-1)\otimes\D^b(\BS_r),F^*_2(-1)\otimes\D^b(\BS_r)\rangle^\perp) = 0$;

\noindent$(ii)$
$\Phi_{i_*\HE(0,1)\otimes(\TCL_r/\TCL_{r-1})}(\D^b(\TCY_r)) \subset
\langle E_1\otimes\D^b(\BS_r),E_2\otimes\D^b(\BS_r)\rangle$;

\noindent$(ii)$
$\Phi^!_{i_*\HE(1,0)\otimes(\TCL_r/\TCL_{r-1})^*}
(\langle E_1(1)\otimes\D^b(\BS_r),E_2(1)\otimes\D^b(\BS_r)\rangle^\perp) = 0$.
\end{lemma}
\begin{proof}
Note that for $k=1$ or $k=2$ we have
$$
\Phi_{E_k\boxtimes F_k(1)\otimes(\TCL_r/\TCL_{r-1})}(G)
\cong
E_k\otimes H^\bullet(\TCY_r,F_k(1)\otimes G)\otimes(\TCL_r/\TCL_{r-1})
\cong
E_k \otimes \RHom(F^*_k(-1),G) \otimes (\TCL_r/\TCL_{r-1}).
$$
Tensoring exact triangle~(\ref{ieab}) by $\CO(0,1)\otimes(\TCL_r/\TCL_{r-1})$
we deduce from this claims $(i)$ and $(ii)$. Similarly, we have
$$
\Phi^!_{E_k(1)\boxtimes F_k\otimes(\TCL_r/\TCL_{r-1})^*}(G)
\cong
\RHom(E_k(1),G)\otimes F^*_k\otimes(\TCL_r/\TCL_{r-1})
\otimes\omega_{\TCY_r}[\dim\TCY_r]
$$
whereof we deduce $(iii)$.
\end{proof}

\subsection{Semiorthogonal collections for linear sections of $X$}\label{sodlsx}

We prove in this subsection that the derived category $\D^b(\CX_r)$
admits a semiorthogonal collection of the form
$$
(\Phi_r(\D^b(\CY_r)),
E_1(1)\otimes\D^b(\BP_r),  E_2(1)\otimes\D^b(\BP_r),
\dots,
E_1(i-r)\otimes\D^b(\BP_r),  E_2(i-r)\otimes\D^b(\BP_r) ).
$$
We prove it by induction in $r$. We start with the following

\begin{lemma}\label{eec}
If $r < i$ then the functors $\D^b(\BP_r) \to \D^b(\CX_r)$,
$G\mapsto E_s(k)\otimes f_r^*G$ are fully faithful for all~$k$.
Moreover, the collections
$$
\begin{array}{l}
\left\langle
E_1(1)\otimes\D^b(\BP_r),  E_2(1)\otimes\D^b(\BP_r),
\dots,
E_1(i-r)\otimes\D^b(\BP_r),  E_2(i-r)\otimes\D^b(\BP_r)
\right\rangle
\qquad\text{and}
\smallskip\\
\left\langle
E_1\otimes\D^b(\BP_r),  E_2\otimes\D^b(\BP_r),
\dots,
E_1(i-r-1)\otimes\D^b(\BP_r),  E_2(i-r-1)\otimes\D^b(\BP_r)
\right\rangle
\end{array}
$$
in $\D^b(\CX_r)$ are semiorthogonal.
\end{lemma}
\begin{proof}
Take any pair $E_s(k)$, $E_t(l)$ such that
either $0 < k - l < i - r$, or $k = l$ and $s \ge t$. We have
\begin{multline*}
\Hom(E_s(k)\otimes f_r^*G,E_t(l)\otimes f_r^*G') \cong
H(\CX_r,E_s^*(-k)\otimes E_t(l)\otimes f_r^*G^*\otimes f_r^*G') \cong
\\ \cong
H(X\times\BP_r,(E_s^*\otimes E_t(l-k))\boxtimes
(G^*\otimes G')\otimes i_*\CO_{\CX_r}).
\end{multline*}
But $i_*\CO_{\CX_r}$ by lemma~\ref{zl1}~$(i)$ admits a Koszul resolution
$i_*\CO_{\CX_r} \cong \Lambda^\bullet(\CO_X(-1)\boxtimes\CL_r)$.
Moreover, we have
\begin{multline*}
H(X\times\BP_r,(E_s^*\otimes E_t(l-k))\boxtimes
(G^*\otimes G')\otimes \Lambda^p(\CO_X(-1)\boxtimes\CL_r)) \cong
\\ \cong
H(X\times\BP_r,(E_s^*\otimes E_t(l-k-p))\boxtimes
(G^*\otimes G'\otimes \Lambda^p\CL_r)) \cong
\\ \cong
H(X,E_s^*\otimes E_t(l-k-p))\otimes
H(\BP_r,G^*\otimes G'\otimes \Lambda^p\CL_r)) \cong
\\ \cong
\Hom_X(E_s,E_t(l-k-p))\otimes
\Hom_{\BP_r}(G,G'\otimes \Lambda^p\CL_r).
\end{multline*}
It follows from condition (C.3) that
$\Hom_X(E_s,E_t(l-k-p)) = 0$ for $0\le p\le r$ in all cases
with the exception of $k=l$, $s=t$ and $p=0$, when we have
$\Hom_X(E_s,E_t(l-k-p)) = \kk$.
Therefore the collections are semiorthogonal and we have
$\Hom(E_s(k)\otimes f_r^*G,E_s(k)\otimes f_r^*G') \cong \Hom_{\BP_r}(G,G')$,
that is the functor $G\mapsto E_s(k)\otimes f_r^*G$ is fully faithful.
\end{proof}

\begin{proposition}\label{phirff}
For all $1\le r \le i$ we have

\noindent
$(i)$
$\Phi_r(\D^b(\CY_r)) \subset
\langle E_{1r}(1)\otimes\D^b(\BP_r),\dots,
E_{2r}(i-r)\otimes\D^b(\BP_r)\rangle^\perp$;

\noindent
$(ii)$
the functor $\Phi_r:\D^b(\CY_r) \to \D^b(\CX_r)$ is fully faithful.
\end{proposition}
\begin{proof}
We use induction in $r$. The base of induction, $r=1$ is given
by theorem~\ref{ff} and proposition~\ref{homephi}.

Assume that $(i)$ and $(ii)$ are true for $r-1$ and consider
the first diagram of~(\ref{yayb}).
First of all, note that
\begin{equation}\label{philrm1}
\TPhi_{r-1}(\D^b(\TCY_{r-1})) \subset
\langle E_1(1)\otimes\D^b(\BS_r),\dots,
E_2(i-r+1)\otimes\D^b(\BS_r)\rangle^\perp
\end{equation}
by the induction hypothesis and proposition~\ref{phitperp}.
On the other hand,
$$
\Phi_{i_*\HE(0,1)\otimes(\TCL_r/\TCL_{r-1})}(\D^b(\TCY_{r-1})) \subset
\langle E_1\boxtimes\D^b(\BS_r),E_2\boxtimes\D^b(\BS_r)\rangle \subset
\langle E_1(1)\otimes\D^b(\BS_r),\dots,E_2(i-r)\otimes\D^b(\BS_r)\rangle^\perp
$$
by lemma~\ref{phiise} and lemma~\ref{eec}.
Taking into account exact triangle~(\ref{ft1}) we deduce that
$$
\xi_*\TPhi_r(\D^b(\TCY_{r})) \subset
\langle E_1(1)\otimes\D^b(\BS_r),\dots,E_2(i-r)\otimes\D^b(\BS_r)\rangle^\perp,
$$
hence by adjunction we have
\begin{multline*}
\TPhi_r(\D^b(\TCY_r)) \subset
\xi^*\langle E_1(1)\otimes\D^b(\BS_r),\dots,
E_2(i-r)\otimes\D^b(\BS_r)\rangle^\perp =
\\ =
\langle E_1(1)\otimes\D^b(\BS_r),\dots,
E_2(i-r)\otimes\D^b(\BS_r)\rangle^\perp \subset
\psi^*\langle E_1(1)\otimes\D^b(\BP_r),\dots,
E_2(i-r)\otimes\D^b(\BP_r)\rangle^\perp.
\end{multline*}
Again by adjunction we have
$$
\psi_*\TPhi_r(\D^b(\TCY_r)) \subset
\langle E_1(1)\otimes\D^b(\BP_r),\dots,E_2(i-r)\otimes\D^b(\BP_r)\rangle^\perp.
$$
By corollary~\ref{cim1} we have $\psi_*\psi^*\cong \id_{\CY_r}$,
hence applying lemma~\ref{phit} we deduce
$$
\Phi_r(\D^b(\CY_r)) =
\psi_*\psi^*\Phi_r(\D^b(\CY_r)) =
\psi_*\TPhi_r\psi^*(\D^b(\CY_r)) \subset
\psi_*\TPhi_r(\D^b(\TCY_r)),
$$
which combined with the previous inclusion proves $(i)$ for $r$.

To prove $(ii)$ we note that by~(\ref{philrm1}) and \ref{phiise}~$(iii)$
we have
$$
\Phi^!_{i_*\HE(1,0)\otimes(\TCL_r/\TCL_{r-1})^*}\circ\TPhi_{r-1} = 0.
$$
Therefore, (\ref{ft2}) implies that
$\TPhi^!_r\xi^*\TPhi_{r-1} \cong \eta_*\TPhi^!_{r-1}\TPhi_{r-1}$.
Combining this with corollary~\ref{phixieta1} and with the induction hypothesis
we obtain
$$
\TPhi^!_r\TPhi_r\eta_* \cong
\TPhi^!_r\xi^*\TPhi_{r-1} \cong
\eta_*\TPhi^!_{r-1}\TPhi_{r-1} \cong \eta_*.
$$
Finally, applying corollary~\ref{cim1} and lemma~\ref{phit}
we deduce
$$
\Phi_r^!\Phi_r \cong
\Phi_r^!\Phi_r\hat\psi_*\hat\psi^* \cong
\Phi_r^!\Phi_r\psi_*\eta_*\hat\psi^* \cong
\psi_*\TPhi_r^!\TPhi_r\eta_*\hat\psi^* \cong
\psi_*\eta_*\hat\psi^* \cong
\hat\psi_*\hat\psi^* \cong \id_{\CY_r}.
$$
Thus $\Phi_r$ is fully faithful and $(ii)$ is proved for $r$.
\end{proof}

From the above results we obtain the following

\begin{corollary}\label{phiaff}
If $r \le i$ then the collection
\begin{equation}\label{xasod}
\hspace{-1cm}
\langle
\Phi_r(\D^b(\CY_r)),
E_1(1)\otimes\D^b(\BP_r),E_2(1)\otimes\D^b(\BP_r),
\dots,
E_1(i-r)\otimes\D^b(\BP_r),E_2(i-r)\otimes\D^b(\BP_r)
\rangle
\end{equation}
is semiorthogonal in $\D^b(\CX_r)$.
\end{corollary}

\subsection{Fullness}\label{fullness}

In this subsection we establish fullness of the semiorthogonal
collection~(\ref{xasod}). First of all we use a relative
Bridgeland's trick~\ref{ffisequ} to check that
$\Phi_i:\D^b(\CY_i) \to \D^b(\CX_i)$ is an equivalence of categories.
Then we use descending induction in $r$ to check fullness for $r\le i$.
We also compute here the canonical class of $Y$.

First of all we need the following

\begin{lemma}\label{bpine0}
We have $\BP_i \ne \emptyset$.
\end{lemma}
\begin{proof}
Follows from lemma~\ref{im1} and lemma~\ref{im2} by induction in $r$.
\end{proof}

\begin{proposition}\label{phiiequ}
The functor $\Phi_i:\D^b(\CY_i) \to \D^b(\CX_i)$ is
an equivalence of categories.
Moreover, we have $\dim\CY_i - \dim\BP_i = \dim X - i$ and
$\omega_{\CY_i} \cong \omega_{\BP_i}\otimes\det\CL_i^*$.
\end{proposition}
\begin{proof}
Note that $\CX_i$ is smooth and connected by lemma~\ref{xyksm},
$\Phi_i$ is fully faithful by proposition~\ref{phirff}~(ii), and
by adjunction formula and lemma~\ref{zl1}~(i) we have
$$
\omega_{\CX_i} \cong
(\omega_{X\times\BP_i} \otimes \det(\CO_X(1)\boxtimes\CL_i^*))_{|\CX_i} \cong
((\CO_X(-i)\boxtimes\omega_{\BP_i}) \otimes
(\CO_X(i)\boxtimes\det\CL_i^*))_{|\CX_i} \cong
\omega_{\BP_i}\otimes\det\CL_i^*,
$$
and $\dim\CX_i - \dim\BP_i = \dim X - i$.
Note that the line bundle $\omega_{\BP_i}\otimes\det\CL_i^*$
is a pullback from $\BP_i$, hence, by proposition~\ref{ffisequ}
the functor $\Phi_i$ is an equivalence,
$\dim\CY_i - \dim\BP_i = \dim X - i$ and
$\omega_{\CY_i} \cong \omega_{\BP_i}\otimes\det\CL_i^*$.
\end{proof}

Now we compute the relative canonical classes of $\CY_r$ over $\BP_r$
by descending induction in $r$. This allows to compute the canonical
class of $Y$.

\begin{lemma}
For all $r \le i$ we have
$\omega_{\CY_r/\BP_r} \cong \det\CL_r^*(i - r)$.
\end{lemma}
\begin{proof}
We use descending induction in $r$. The base of induction, $r = i$,
follows from proposition~\ref{phiiequ}. Indeed,
$\omega_{\CY_i/\BP_i} \cong
\omega_{\CY_i}\otimes g^*\omega_{\BP_i}^{-1} \cong
\det\CL_i^*$.
Now assume that the claim is true for $r$ and consider the
diagrams~(\ref{yayb}). Using the adjunction formula and lemma~\ref{zl1}~$(iv)$
we deduce
\begin{multline*}
\omega_{\TCY_{r-1}/\BS_r} \cong
\omega_{\TCY_r/\BS_r} \otimes (\TCL_r/\TCL_{r-1})(1) \cong
\det\TCL_r^*(i - r) \otimes (\TCL_r/\TCL_{r-1})(1) \cong
\\ \cong
\det\TCL_{r-1}^*(i - (r - 1)) \cong
\phi^*\det\CL_{r-1}^*(i - (r - 1)).
\end{multline*}
On the other hand,
$\omega_{\TCY_{r-1}/\BS} \cong \phi^*\omega_{\CY_{r-1}/\BP_{r-1}}$
since $\phi$ is flat. Therefore
$\omega_{\CY_{r-1}/\BP_{r-1}} \cong \det\CL_{r-1}^*(i - (r - 1))$
by corollary~\ref{cim2}.
\end{proof}

\begin{corollary}\label{omegay}
We have $\omega_Y \cong \CO_Y(i-N)$ and $\dim Y = \dim X + N - 2i$.
\end{corollary}
\begin{proof}
Note that $\CY_1 = Y$, $\BP_1 = \BP$ and
$\det\CL_1 \cong \CL_1 \cong \CO_Y(-1)$.
Hence
$$
\omega_Y \cong
\omega_{Y/\BP} \otimes \omega_{\BP} \cong
\omega_{\CY_1/\BP_1} \otimes \omega_{\BP} \cong
\left(\CO_Y(1)\otimes\CO_Y(i-1)\right) \otimes \CO_Y(-N) \cong
\CO_Y(i-N).
$$
On the other hand, by lemma~\ref{xyksm} we have
$\dim\CY_i - \dim\BP_i = \dim Y + i - N$, hence
by proposition~\ref{phiiequ} we have $\dim Y + i - N = \dim X - i$,
whereof $\dim Y = \dim X + N - 2i$.
\end{proof}

Now we establish semiorthogonal decompositions for the derived categories
$\D^b(\CX_r)$ with $r\le i$.

\begin{theorem}\label{dxa}
For any $r \le i$ we have a semiorthogonal decomposition
\begin{equation}\label{xrm1}
\hspace{-1cm}
\D^b(\CX_{r}) = \langle
\Phi_{r}(\D^b(\CY_{r})),
E_1(1)\otimes\D^b(\BP_{r}),
E_2(1)\otimes\D^b(\BP_{r}),
\dots,
E_1(i-r)\otimes\D^b(\BP_{r}),
E_2(i-r)\otimes\D^b(\BP_{r})
\rangle
\end{equation}
\end{theorem}
\begin{proof}
We use descending induction in $r$.
The base of induction, $r = i$, is given by proposition~\ref{phiiequ}.
Assume that~(\ref{xrm1}) is true for $r$. Then by the faithful base change
theorem~\ref{phitsod} we have
\begin{equation}\label{xtclr}
\D^b(\TCX_r) = \langle
\TPhi_r(\D^b(\TCY_r)),
E_1(1)\otimes\D^b(\BS_r),
\dots,
E_2(i-r)\otimes\D^b(\BS_r)
\rangle
\end{equation}
Take arbitrary
$$
G \in \langle
\Phi_{r-1}(\D^b(\CY_{r-1})),
E_1(1)\otimes\D^b(\BP_{r-1}),\dots,E_2(i-r+1)\otimes\D^b(\BP_{r-1})
\rangle^\perp.
$$
Consider the diagram~(\ref{yayb}).
By lemma~\ref{sod_f} we have
${f_{r-1}}_*\RCHom(E_s(k),G) = 0$ for $s=1,2$ and $1\le k\le i-r+1$.
Therefore
$$
{f_{r-1}}_*\RCHom(\phi^*E_s(k),\phi^*G) =
{f_{r-1}}_*\phi^*\RCHom(E_s(k),G) =
\phi^*{f_{r-1}}_*\RCHom(E_s(k),G) = 0
$$
since the right square in~(\ref{yayb}) is exact cartesian.
Therefore
\begin{equation}\label{phisg}
\phi^*G \in
\langle
E_1(1)\otimes\D^b(\BS_r),\dots,E_2(i-r+1)\otimes\D^b(\BS_r)
\rangle^\perp
\end{equation}
On the other hand, it follows from lemma~\ref{zl1}~$(iii)$ that
we have a resolution
$$
0 \to E_s(k)\otimes(\TCL_{r}/\TCL_{r-1}) \to E_s(k+1) \to \xi_*E_s(k+1) \to 0,
$$
hence
$$
\phi^*G \in
\langle
\xi_*E_1(2)\otimes\D^b(\BS_r),\dots,\xi_*E_2(i-r+1)\otimes\D^b(\BS_r)
\rangle^\perp,
$$
and by adjunction
$$
\xi^!\phi^*G \in
\langle
E_1(2)\otimes\D^b(\BS_r),\dots,E_2(i-r+1)\otimes\D^b(\BS_r)
\rangle^\perp.
$$
Further,
$\xi^!\phi^*G \cong
\xi^*\phi^*G\otimes(\TCL_{r}/\TCL_{r-1})^*(1)[-1] \cong
\hat{\phi}^*G\otimes(\TCL_{r}/\TCL_{r-1})^*(1)[-1]$,
therefore
\begin{equation}\label{hphisg}
\hat{\phi}^*G \in
\langle E_1(1)\otimes\D^b(\BS_r),\dots,E_2(i-r)\otimes\D^b(\BS_r)\rangle^\perp.
\end{equation}
On the other hand, (\ref{phisg}) and lemma~\ref{phiise}~$(iii)$ imply that
$\Phi^!_{i_*\HE(1,0)\otimes(\TCL_r/\TCL_{r-1})^*}(\phi^*G) = 0$.
Therefore, applying exact triangle~(\ref{ft2}) to $\phi^*G$ we deduce that
$\TPhi_r^!\xi^*\phi^*G \cong \eta_*\TPhi_{r-1}^!\phi^*G$.
But
$\TPhi_{r-1}^!\phi^*(G) \cong \phi^*\Phi_{r-1}^!(G) = 0$,
since $G \in \Phi_{r-1}(\D^b(\CY_{r-1}))^\perp$ by assumption.
Thus we deduce that $\TPhi_r^!\hat{\phi}^*(G) = 0$,
hence
$$
\hat{\phi}^*G \in \TPhi_r(\D^b(\TCY_r))^\perp.
$$
Combining this with~(\ref{hphisg}) and (\ref{xtclr}) we deduce
that $\hat{\phi}^*G = 0$. Therefore $G = 0$ by corollary~\ref{cim2}
and we obtain decomposition~(\ref{xrm1}) for $r-1$.
\end{proof}

\subsection{Exceptional collection for $Y$}\label{ecfory}

In this subsection we show that  in a certain sense the collection
$(F_1(1),F_2(1),\dots,F_1(N-i),F_2(N-i))$ on $Y$ is exceptional.
The idea behind the proof is very simple.
It is easy to see that if a collection $(E(1),\dots,E(i))$ on a Fano
variety $X$ of index $i$ is exceptional then the restriction of the bundle $E$
to a Calabi--Yau linear section of $X$ is also exceptional. It turns out
that the inverse is also true: if the restrictions of a vector bundle $E$
to all Calabi--Yau linear sections of $X$ are exceptional then
the collection $(E(1),\dots,E(i))$ is exceptional on $X$.
We use this idea as follows. First of all we consider
an exceptional pair $(E_0(1),E_1(1))$ on $X$ obtained from $(E_1,E_2)$
by a mutation and a twist, and show that the functor $\Phi_i^!$ takes
this pair to the pair $(F_2^*,F_1^*)$ on $\CY_i$ up to a twist and a shift.
Since $\Phi_i$ is an equivalence, its adjoint functor $\Phi_i^!$ is
an equivalence as well, hence we deduce that the pair $(F_2^*,F_1^*)$
is exceptional on $\CY_i$. But this pair is formed by the restrictions
of a pair $(F_2^*,F_1^*)$ on $Y$ to the Calabi--Yau sections, hence
the above argument shows that the collection
$(F_2^*(1),F_1^*(1),\dots,F_2^*(N-i),F_1^*(N-i))$ on $Y$ is exceptional.

Now we give a detailed proof.
We define an object $E_0 \in \D(X)$ as the mutation of $E_2$ through $E_1$
(see~\cite{B}). Thus we have the following exact triangle
\begin{equation}\label{e0}
E_0 \to \RHom(E_1,E_2)\otimes E_1 \to E_2.
\end{equation}

\begin{lemma}[cf.~\cite{B}]\label{homtoe0}
We have
$$
\RHom(E_1,E_0(k)) = 0,\quad\text{if $1-i \le k \le 0$},\qquad
\RHom(E_2,E_0(k)) = \begin{cases}
0, & \text{if $1-i \le k \le -1$,}\\
\kk[-1], & \text{if $k=0$.}
\end{cases}
$$
\end{lemma}
\begin{proof}
Apply $\RHom(E_1,-)$ and $\RHom(E_2,-)$ to~(\ref{e0}).
\end{proof}

\begin{proposition}\label{phises}
For all $1\le r\le i$ and $s = 0,1$ we have
$$
\Phi_r^!(E_s(k)) =
\begin{cases}
0, & \text{for $r+1-i \le k \le 0$;}\\
F_{2-s}^*(i-r-1)\otimes\det\CL_r^*[\dim Y - N + r], &
\text{for $k = 1$.}
\end{cases}
$$
\end{proposition}
\begin{proof}
We use induction in $r$. First, let us check the case $r=1$.
Consider the diagram
$$
\xymatrix{
& \CX_1\times Y \ar[d]_p \ar[drr]_(.7){q} \ar[dl]_{p_1} &
Q(X,Y) \ar[l]_j \ar[r]^i &
X\times Y \ar[dll]^(.7)p \ar[d]^{q} \\
\CX_1 \ar[r]^\pi & X && Y
}
$$
Note that the bundle $E_s(k)$ on $\CX_1$ is a pullback via $\pi$. But
$\Phi_1^!(\pi^*G) \cong
\Phi_1^!(\pi^!G\otimes\omega_{\CX_1/X}^{-1}[\dim X - \dim\CX_1]))$
and since by adjunction formula and lemma~\ref{zl1}~$(i)$
$\omega_{\CX_1/X} \cong \omega_\BP\otimes(\CO_X(1)\boxtimes\CL_1^*)$,
$\dim\CX_1 = \dim X + N - 2$,
and since the functor $\Phi_1^!$ is $\BP$-linear we have
\begin{equation}\label{phishpis}
\Phi_1^!(\pi^*G) \cong
\Phi_1^!(\pi^!(G(-1))\otimes\omega_{\BP}^{-1}\otimes\CL_1[2 - N].
\end{equation}
On the other hand, it is clear that the functor $\Phi_1^!\circ\pi^!$
is a kernel functor of the second type with kernel
$(\pi\times\id_Y)_*\CE_1 \cong (\pi\times\id_Y)_*j_*\CE \cong i_*\CE$.
Therefore
$$
\Phi_1^!(\pi^!(G(-1)) \cong
\Phi^!_{i_*\CE}(G(-1)) =
q_*\RCHom(i_*\CE,p^!(G(-1))) \cong
q_*\RCHom(i_*\CE,G(-1)\boxtimes\omega_Y[\dim Y]).
$$
Using resolution~(\ref{equ_e}) combined with lemma~\ref{homtoe0}
we obtain $\Phi_1^!(\pi^!(E_s(k-1)) = 0$ for $s = 0, 1$ and
$2-i\le k\le 0$ and for $k=1$ we obtain the following exact triangles
$$
\begin{array}{l}
\Phi_1^!(\pi^!(E_0(-1)) \to F_2^*\otimes\omega_Y[\dim Y - 1] \to 0,\\
\Phi_1^!(\pi^!(E_1(-1)) \to 0 \to F_1^*\otimes\omega_Y[\dim Y].
\end{array}
$$
Finally, we get
$$
\Phi_1^!(\pi^!(E_s(k-1)) \cong
\begin{cases}
0, & \text{for $2-i \le k \le 0$, $s = 0, 1$;}\\
F_{2-s}^*\otimes\omega_Y[\dim Y - 1], & \text{for $k = 1$, $s = 0, 1$.}
\end{cases}
$$
Substituting this into~(\ref{phishpis}) and taking into account
the isomorphisms $\omega_Y \cong \CO_Y(i-N)$ (see corollary~\ref{omegay}),
$\omega_\BP^{-1} \cong \CO_Y(N)$ and
$\det\CL_1 \cong \CL_1 \cong \CO_Y(-1)$ we deduce the claim for $r=1$.

Now, assume that the claim is true for $r-1$.
Consider the diagram~(\ref{yayb}).
We have
\begin{multline}\label{phitclrm1}
\TPhi_{r-1}^!(E_s(k)) \cong
\TPhi_{r-1}^!\phi^*(E_s(k)) \cong
\phi^*\Phi_{r-1}^!(E_s(k)) \cong
\\ \cong
\begin{cases}
0, & \text{for $r-i \le k \le 0$, $s = 0,1$;}\\
\phi^*F_{2-s}^*(i-r)\otimes\det\CL_{r-1}^*[\dim Y - N + r - 1], &
\text{for $k = 1$, $s = 0,1$.}
\end{cases}
\end{multline}
by lemma~\ref{phit} and the induction hypothesis.
On the other hand,
\begin{multline*}
\hat{\psi}^*\Phi_r^!(E_s(k)) \cong
\eta^*\psi^*\Phi_r^!(E_s(k)) \cong
\eta^*\TPhi_r^!(\psi^*E_s(k)) \cong
\\ \cong
\eta^!\TPhi_r^!(E_s(k)) \otimes \CO_Y(-1)\otimes (\TCL_r/\TCL_{r-1})^*[1] \cong
\TPhi_{r-1}^!(\xi_*E_s(k)) \otimes \CO_Y(-1)\otimes (\TCL_r/\TCL_{r-1})^*[1]
\end{multline*}
by lemma~\ref{phit}, lemma~\ref{zl1}~$(iv)$, and lemma~\ref{phixieta1}.
Taking into account resolution
$$
0 \to E_s(k-1)\otimes(\TCL_r/\TCL_{r-1}) \to E_s(k) \to \xi_*E_s(k) \to 0,
$$
and~(\ref{phitclrm1}) we deduce
$$
\hat{\psi}^*\Phi_r^!(E_s(k)) \cong
\begin{cases}
0, & \text{for $r+1-i \le k \le 0$, $s = 0,1$;}\\
F_{2-s}^*(i-r-1)\otimes\det\TCL_r^*[\dim Y - N + r], &
\text{for $k = 1$, $s = 0,1$.}
\end{cases}
$$
Finally, applying corollary~\ref{cim1} we deduce the claim.
\end{proof}

\begin{corollary}\label{fsifti}
Let $s,t=1,2$. Then
$$
{g_i}_*\RCHom_{\CY_i}(F_s,F_t) \cong
\begin{cases}
\CO_{\BP_i} \oplus \det\CL_i[i - \dim X], & \text{for $s = t$}\\
W\otimes\CO_{\BP_i}, & \text{for $s=1$, $t=2$}\\
W^*\otimes\det\CL_i[i - \dim X], & \text{for $s=2$, $t=1$}
\end{cases}
$$
\end{corollary}
\begin{proof}
By proposition~\ref{phises}, proposition~\ref{slinadj}
and proposition~\ref{phiiequ} we have
\begin{multline*}
{g_i}_*\RCHom_{\CY_i}(F_s,F_t) \cong \\ \cong
{g_i}_*\RCHom_{\CY_i}(F_t^*(-1)\otimes\det\CL_i^*[\dim Y - N + i],
F_s^*(-1)\otimes\det\CL_i^*[\dim Y - N + i]) \cong \\ \cong
{g_i}_*\RCHom_{\CY_i}(\Phi_i^!(E_{2-t}(1)),\Phi_i^!(E_{2-s}(1))) \cong
{f_i}_*\RCHom_{\CX_i}(\Phi_i\Phi_i^!(E_{2-t}(1)),E_{2-s}(1)) \cong \\ \cong
{f_i}_*\RCHom_{\CX_i}(E_{2-t}(1),E_{2-s}(1)) \cong
{f_i}_*\RCHom_{\CX_i}(E_{2-t},E_{2-s}).
\end{multline*}
On the other hand,
${f_i}_*\RCHom_{\CX_i}(E_{2-t},E_{2-s}) \cong
{f_i}_*(E_{2-t}^*\otimes E_{2-s})$.
Tensoring the Koszul resolution of $\CX_i$ in $X\times\BP_i$
(see lemma~\ref{zl1}~$(i)$) by $E_{2-t}^*\otimes E_{2-s}$
we see that the pushforward of $E_{2-t}^*\otimes E_{2-s}$
from $\CX_i$ to $X\times\BP_i$ is quasiisomorphic to the complex
$$
(E_{2-t}^*\otimes E_{2-s}(-i))\boxtimes\det\CL_i \to
(E_{2-t}^*\otimes E_{2-s}(1-i))\boxtimes\Lambda^{i-1}\CL_i \to
\dots \to
(E_{2-t}^*\otimes E_{2-s}(-1))\boxtimes\CL_i \to
E_{2-t}^*\otimes E_{2-s}
$$
Applying the functor ${f_i}_*$ and using the K\"unneth formula on
$X\times\BP_i$ and (C.3) we deduce the claim.
\end{proof}

\begin{proposition}\label{fsrftr}
If $r\ge i$ and $\BP_r\ne\emptyset$ then
$$
\CH^p({g_r}_*\RCHom_{\CY_r}(F_s(k),F_t(l))) \cong
\begin{cases}
\CO_{\BP_r}, &
        \text{if $s=t$, $k=l$ and $p=0$,}\\
\det\CL_r, &
        \text{if $s=t$, $k-l = r-i$}\\
&\quad  \text{and $p=\dim X + r - 2i$,}\\
W\otimes\CO_{\BP_r}, &
        \text{if $s=1$, $t=2$, $k=l$ and $p=0$,}\\
W^*\otimes\det\CL_r, &
        \text{if $s=2$, $t=1$, $k-l = r-i$}\\
&\quad  \text{and $p=\dim X + r - 2i$,}\\
0, &
        \text{otherwise, if $0 \le k-l \le r-i$.}
\end{cases}
$$
\end{proposition}
\begin{proof}
We use induction in $r$. The base of induction, $r=i$,
is given by corollary~\ref{fsifti}. Now assume that
the claim is true for $r-1$ and consider the diagram~(\ref{yayb}).
Let us denote
$$
\CF_{s,k,t,l}^r = {g_r}_*\RCHom_{\CY_r}(F_s(k),F_t(l)),\qquad
\TCF_{s,k,t,l}^r = {g_r}_*\RCHom_{\TCY_r}(F_s(k),F_t(l)).
$$
First of all we note that
\begin{multline*}
\TCF_{s,k,t,l}^{r-1} =
{g_{r-1}}_*\RCHom_{\TCY_{r-1}}(F_s(k),F_t(l)) \cong
\\ \cong
{g_{r-1}}_*\phi^*\RCHom_{\CY_{r-1}}(F_s(k),F_t(l)) \cong
\phi^*{g_{r-1}}_*\RCHom_{\CY_{r-1}}(F_s(k),F_t(l)) =
\phi^*\CF_{s,k,t,l}^{r-1}.
\end{multline*}
and similarly
$$
\TCF_{s,k,t,l}^{r} \cong \psi^*\CF_{s,k,t,l}^{r}.
$$
On the other hand, by~\ref{zl1}~$(iv)$ we have a resolution
$$
0 \to F_t(l-1)\otimes(\TCL_r/\TCL_{r-1})^* \to F_t(l) \to \eta_*F_t(l) \to 0.
$$
Applying to it the functor ${g_r}_*\RCHom_{\TCY_r}(F_s(k),-)$
and taking into account that
$$
{g_r}_*\RCHom_{\TCY_r}(F_s(k),\eta_*F_t(l)) \cong
{g_r}_*\eta_*\RCHom_{\TCY_{r-1}}(\eta^*F_s(k),F_t(l)) \cong
{g_{r-1}}_*\RCHom_{\TCY_{r-1}}(F_s(k),F_t(l)) =
\TCF_{s,k,t,l}^{r-1},
$$
we deduce the following exact triangle on $\BS$
$$
\TCF_{s,k,t,l-1}^{r}\otimes(\TCL_r/\TCL_{r-1})^* \to
\TCF_{s,k,t,l}^{r} \to
\TCF_{s,k,t,l}^{r-1}.
$$
Since
$\TCF_{s,k,t,l}^{r-1} \cong \phi^*\CF_{s,k,t,l}^{r-1}$,
$\TCF_{s,k,t,l}^{r} \cong \psi^*\CF_{s,k,t,l}^{r}$
we can rewrite this as
\begin{equation}\label{cftr}
\psi^*(\CF_{s,k,t,l-1}^{r}) \otimes (\TCL_r/\TCL_{r-1})^* \to
\psi^*(\CF_{s,k,t,l}^{r}) \to
\phi^*\CF_{s,k,t,l}^{r-1}.
\end{equation}
Finally, note that
\begin{equation}\label{psis}
\psi_*(\TCL_r/\TCL_{r-1})^{\otimes -k} =
\begin{cases}
\CO_{\BP_r}, & \text{if $k=0$,}\\
0, & \text{if $0 \le k\le r-1$,}\\
\det\CL_r^*[1-r], & \text{if $k=r$}
\end{cases}
\end{equation}
by lemma~\ref{im1} (because $\TCL_r/\TCL_{r-1}$ is
the Grothendieck line bundle on $\BS_r \subset \PP_{\BP_r}(\CL_r^*)$).
Therefore, applying the functor $\psi_*$ to the triangle~(\ref{cftr}) we get
$$
\CF_{s,k,t,l}^{r} \cong \psi_*\phi^*\CF_{s,k,t,l}^{r-1},
$$
by corollary~\ref{cim1},
and applying the functor $\psi_*$ to the triangle~(\ref{cftr})
tensored by $(\TCL_r/\TCL_{r-1})^{\otimes 1-r}$ we get
$$
\CF_{s,k,t,l-1}^{r}\otimes\det\CL_r^*[2-r] \cong
\psi_*(\phi^*\CF_{s,k,t,l}^{r-1}\otimes(\TCL_r/\TCL_{r-1})^{\otimes 1-r}).
$$
Now we use the first formula and the induction assumption to compute
$\CF_{s,k,t,l}^{r}$ for $0\le k-l \le r-1-i$ and
the second formula and the induction assumption to compute
$\CF_{s,k,t,l}^{r}$ for $k-l = r-i$.
\end{proof}

\begin{corollary}
If\/ $\BZ = \emptyset$ then $(F_1(1),F_2(1),\dots,F_1(N-i),F_2(N-i))$
is an exceptional collection on~$Y$ and $\RHom(F_1,F_2) = W$.
\end{corollary}
\begin{proof}
If $\BZ = \emptyset$ then $\BP_N = \Spec\kk$, $\CY_N = Y$,
${g_N}_*\RCHom_{\CY_N}(F_s(k),F_t(l)) = \RHom_Y(F_s(k),F_t(l))$.
\end{proof}

\subsection{Semiorthogonal decompositions for linear section of $Y$}\label{sodlsy}

In this subsection we establish semiorthogonal decompositions
for $\CY_r$ when $r\ge i$.

First of all we have to show that $\Phi_r^!(\D^b(\CX_r))$ is orthogonal
to $F_s^*(k)$ for $i-r\le k\le -1$. This follows from the

\begin{proposition}\label{phirfsr}
If $r\ge i$ and $\BP_r\ne\emptyset$ then
$$
\Phi_r(F_s^*(k)) \cong
\begin{cases}
0, & \text{for $i-r\le k\le -1$,}\\
E_{2-s}(1)\otimes\det\CL_r[N-r-\dim Y], & \text{for $k=i-r-1$.}
\end{cases}
$$
\end{proposition}
\begin{proof}
We use the same arguments as in the proof of proposition~\ref{fsrftr}.
The base of the induction, $r=i$, follows from proposition~\ref{phises},
since $\Phi_i^!$ is an equivalence by proposition~\ref{phiiequ}.
Now assume that the claim is true for $r-1$ and consider
the diagram~(\ref{yayb}). Let us denote
$$
\CG_{s,k}^r = \Phi_r(F_s^*(k)),
\qquad
\TCG_{s,k}^r = \TPhi_r(F_s^*(k)).
$$
First of all we note that
$$
\TCG_{s,k}^{r-1} =
\TPhi_{r-1}(F_s^*(k)) \cong
\TPhi_{r-1}(\phi^*F_s^*(k)) \cong
\phi^*\Phi_{r-1}(F_s^*(k)) =
\phi^*\CG_{s,k}^{r-1},
$$
and similarly
$$
\TCG_{s,k}^{r} \cong \psi^*\CG_{s,k}^{r}.
$$
On the other hand, by~\ref{zl1}~$(iv)$ we have a resolution
$$
0 \to F^*_s(k-1)\otimes(\TCL_r/\TCL_{r-1})^* \to
F^*_s(k) \to \eta_*F^*_s(k) \to 0.
$$
Applying to it the functor $\TPhi_r$ and taking into account
lemma~\ref{phixieta1} we deduce the following exact triangle on $\BS_r$
$$
\TCG_{s,k-1}^{r}\otimes(\TCL_r/\TCL_{r-1})^* \to
\TCG_{s,k}^{r} \to
\xi^*\TCG_{s,k}^{r-1}.
$$
Since
$\TCG_{s,k}^{r-1} \cong \phi^*\CG_{s,k}^{r-1}$,
$\TCG_{s,k}^{r} \cong \psi^*\CG_{s,k}^{r}$
we can rewrite this as
\begin{equation}\label{cgtr}
\psi^*(\CG_{s,k-1}^{r}) \otimes (\TCL_r/\TCL_{r-1})^* \to
\psi^*(\CG_{s,k}^{r}) \to
\phi^*\CG_{s,k}^{r-1}.
\end{equation}
Finally, recall~(\ref{psis}).
Applying the functor $\psi_*$
to the triangle~(\ref{cftr}) we get
$$
\CG_{s,k}^{r} \cong \psi_*\phi^*\CG_{s,k}^{r-1},
$$
by corollary~\ref{cim1},
and applying the functor $\psi_*$ to the triangle~(\ref{cgtr})
tensored by $(\TCL_r/\TCL_{r-1})^{\otimes 1-r}$ we get
$$
\CG_{s,k-1}^{r}\otimes\det\CL_r^*[2-r] \cong
\psi_*(\phi^*\CG_{s,k}^{r-1}\otimes(\TCL_r/\TCL_{r-1})^{\otimes 1-r}).
$$
Now we use the first formula and the induction assumption to compute
$\CG_{s,k}^{r}$ for $i-r\le k \le -1$ and
the second formula and the induction assumption to compute
$\CG_{s,k}^{r}$ for $k = i-r-1$.
\end{proof}

Now we can check semiorthogonality.

\begin{proposition}\label{dya1}
For all $r\ge i$ such that $\BP_r\ne\emptyset$ the following collection
$$
\langle
\Phi^!_{r}(\D^b(\CX_{r})),
F^*_2(i-r)\otimes\D^b(\BP_{r}),
\dots,
F^*_1(-1)\otimes\D^b(\BP_{r})
\rangle
$$
is semiorthogonal in $\D^b(\CY_r)$.
\end{proposition}
\begin{proof}
First of all, for all $G,G'\in\D^b(\BP_r)$ we have
\begin{multline*}
\Hom_{\CY_r}({g_r}^*G\otimes F^*_t(l),{g_r}^*G'\otimes F^*_s(k)) \cong
\Hom_{\CY_r}({g_r}^*G,{g_r}^*G'\otimes F^*_s(k)\otimes F_t(l)) \cong
\\ \cong
\Hom_{\BP_r}(G,{g_r}_*({g_r}^*G'\otimes F^*_s(k)\otimes F_t(l))) \cong
\Hom_{\BP_r}(G,G'\otimes {g_r}_*\RCHom_{\CY_r}(F_s(k),F_t(l))))
\end{multline*}
and it follows immediately from proposition~\ref{fsrftr} that
the functors $G\mapsto {g_r}^*G\otimes F^*_t(l)$ are fully faithful
and that the collection
$\langle
F^*_2(i-r)\otimes\D^b(\BP_{r}),\dots,F^*_1(-1)\otimes\D^b(\BP_{r})
\rangle$
in $\D^b(\CY_{r})$ is semiorthogonal.
Further, for any $G\in\D^b(\BP_r)$, $H\in\D^b(\CX_r)$,
since $\Phi_r^!$ is right adjoint to $\Phi_r$ and
$\Phi_r$ is $\BS$-linear we have
$$
\Hom_{\CY_r}(g_r^*G\otimes F_s^*(k),\Phi_r^!(H)) \cong
\Hom_{\CX_r}(\Phi_r(g_r^*G\otimes F_s^*(k)),H) \cong
\Hom_{\CX_r}(f_r^*G\otimes \Phi_r(F_s^*(k)),H)
$$
which is zero for $i-r \le k \le -1$ by proposition~\ref{phirfsr}.
Hence the whole collection is semiorthogonal.
\end{proof}

Now we establish semiorthogonal decompositions for the derived categories
$\D^b(\CY_r)$ with $r\ge i$.

\begin{theorem}\label{dya}
If $r \ge i$ and $\BP_r\ne\emptyset$
then we have a semiorthogonal decomposition
\begin{equation}\label{yrm1}
\D^b(\CY_{r}) =
\langle
\Phi^!_r(\D^b(\CX_{r})),
F^*_2(i-r)\otimes\D^b(\BP_{r}),\dots,F^*_1(-1)\otimes\D^b(\BP_{r})
\rangle
\end{equation}
\end{theorem}
\begin{proof}
By proposition~\ref{dya1} it suffices to show that $\Phi_r^!$ is
fully faithful and that the RHS in~(\ref{yrm1}) generates $\D^b(\CY_r)$.
For this we use induction in $r$. The base of induction, $r = i$,
is given by proposition~\ref{phiiequ}. Assume that the induction
hypothesis is true for $r-1$ and let us check that it
is also true for $r$. Consider the diagram~(\ref{yayb}).
Then by the faithful base change theorem~\ref{phitsod}
we have a semiorthogonal decomposition
\begin{equation}\label{dym1}
\hspace{-1cm}
\D^b(\TCY_{r-1}) =
\langle
\TPhi^!_{r-1}(\D^b(\TCX_{r-1})),
F^*_2(i-r+1)\otimes\D^b(\BS_r),\dots,F^*_1(-1)\otimes\D^b(\BS_r)
\rangle
\end{equation}
By proposition~\ref{dya1} we have
$\TPhi^!_r(\D^b(\TCX_r)) \subset
\langle F^*_2(-1)\otimes\D^b(\BS_r),F^*_1(-1)\otimes\D^b(\BS_r)\rangle^\perp$
Using lemma~\ref{phiise}~$(i)$ we deduce that
$\Phi_{i_*\HE(0,1)\otimes(\TCL_r/\TCL_{r-1})} \circ \TPhi^!_r = 0$.
Therefore, composing exact triangle~(\ref{ft1}) with $\TPhi^!_r$
we obtain an isomorphism
$\TPhi_{r-1}\eta^!\TPhi^!_r \cong \xi_*\TPhi_r\TPhi^!_r$.
Applying also lemma~\ref{phixieta1} and the induction hypothesis we deduce
$$
\xi_*\TPhi_r\TPhi^!_r \cong
\TPhi_{r-1}\eta^!\TPhi^!_r \cong
\TPhi_{r-1}\TPhi^!_{r-1}\xi_* \cong
\xi_*.
$$
Then by lemma~\ref{pff}~$b)$ we have
$\TPhi_r\TPhi^!_r \cong \id_{\TCX_r}$, and
by lemma~\ref{phit} and corollary~\ref{cim1} we have
$$
\Phi_{r}\Phi^!_{r} \cong
\psi_*\psi^*\Phi_{r}\Phi^!_{r} \cong
\psi_*\TPhi_r\psi^*\Phi^!_{r} \cong
\psi_*\TPhi_r\TPhi^!_r\psi^* \cong
\psi_*\psi^* \cong
\id_{\CX_r},
$$
hence $\Phi^!_{r}$ is fully faithful.

Finally, assume that $G\in\D^b(\CY_r)$ is in the left orthogonal
to the RHS of~(\ref{yrm1}), so that
\begin{equation}\label{exp_g}
\Hom_{\CY_r}(G,F^*_s(k)) = 0
\qquad\text{for $i-r\le k\le -1$}\qquad\text{and}\qquad
\Phi_r(G) = 0.
\end{equation}
Since $\hat{\psi}^*$ is fully faithful
by corollary~\ref{cim1} we have for all $i-r+1\le k\le -1$
$$
\Hom_{\TCY_{r-1}}(\hat{\psi}^*G(1),F^*_s(k)) \cong
\Hom_{\TCY_{r-1}}(\hat{\psi}^*G,\hat{\psi}^*F^*_s(k-1)) \cong
\Hom_{\CY_r}(G,F^*_s(k-1)) = 0.
$$
On the other hand, since
$$
\hat{\psi}^*G(1) \cong
\eta^*\psi^*G(1) \cong
\eta^!\psi^*G\otimes(\TCL_r/\TCL_{r-1})^*[1]
$$
by lemma~\ref{zl1}~(iv) and
$\psi^*G \in \langle F^*_2(-1)\otimes\D^b(\BS_r),
F^*_1(-1)\otimes\D^b(\BS_r) \rangle^\perp$
by fullness and faithfulness of $\psi^*$ (corollary~\ref{cim1}),
we have
$\Phi_{i_*\HE(0,1)\otimes(\TCL_r/\TCL_{r-1})}(\psi^*G) = 0$
by lemma~\ref{phiise}. Taking into account~(\ref{ft1}) we deduce
\begin{multline*}
\TPhi_{r-1}(\hat{\psi}^*G(1)) \cong
\TPhi_{r-1}(\eta^!\psi^*G\otimes(\TCL_r/\TCL_{r-1})^*[1]) \cong
(\TPhi_{r-1}\eta^!\psi^*G)\otimes(\TCL_r/\TCL_{r-1})^*[1] \cong
\\ \cong
(\xi_*\TPhi_r\psi^*G)\otimes(\TCL_r/\TCL_{r-1})^*[1] \cong
(\xi_*\psi^*\Phi_r(G)\otimes(\TCL_r/\TCL_{r-1})^*[1] = 0
\end{multline*}
by~(\ref{exp_g}). Thus $\hat{\psi}^*G(1)$ is in the left orthogonal
to the RHS of~(\ref{dym1}), hence $\hat{\psi}^*G(1) = 0$.
But then $G(1) = \hat{\psi}_*\hat{\psi}^*G(1) = 0$ and $G = 0$.
Thus the induction claim is proved.
\end{proof}

\subsection{Proof of the main theorem}

Let $L\subset V^*$ be an admissible subspace, $\dim L = r$.

\begin{lemma}
The map $\varphi:\Spec\kk \to \BP_r$ induced by $L$ is
a faithful base change for the pair $(\CX_r,\CY_r)$.
\end{lemma}
\begin{proof}
Use lemma~\ref{lci}, \ref{zl1}, \ref{xyksm} and
the definition of admissible subspace.
\end{proof}

Now we apply the faithful base change theorem~\ref{phitsod}.
Then theorem~\ref{themain} follows from theorem~\ref{dxa} and theorem~\ref{dya}.

\begin{corollary}
Assume that $L\subset V^*$ is an admissible subspace.
Then $Y_L$ is smooth iff $X_L$ is smooth.
\end{corollary}
\begin{proof}
Note that lemma~\ref{snav_sm} implies that
$X_L$ is smooth iff the category $\D^b(X_L)$ is $\Ext$-bounded,
and $Y_L$ is smooth iff the category $\D^b(Y_L,\CA_Y)$ is $\Ext$-bounded.
On the other hand, one of these categories is embedded fully and faithfully
into another, and the orthogonal is generated by an exceptional collection.
Hence $\D^b(X_L)$ is $\Ext$-bounded iff $\D^b(Y_L,\CA_Y)$ is $\Ext$-bounded.
\end{proof}

Let $\sing(g) \subset \PP(V^*)\setminus\BZ$ denote
the set of critical values of $g:Y \to \PP(V^*)\setminus\BZ$.
Let $X^\vee \subset \PP(V^*)$ be the projectively dual variety of $X$.

\begin{corollary}\label{singg}
If $g$ is flat then $\sing(g) = X^\vee \setminus \BZ$.
\end{corollary}
\begin{proof}
If $g$ is flat then every 1-dimensional subspace $H\subset V^*$ is admissible
hence
$$
\sing(g) =
\{H\in\PP(V^*)\setminus\BZ\ |\ Y_H\text{ is singular}\} =
\{H\in\PP(V^*)\setminus\BZ\ |\ X_H\text{ is singular}\} =
X^\vee \setminus \BZ.
$$
\end{proof}

%
%
%


\section{General construction of $Y$}\label{s5}

In this section we will give a general approach to construction
of the Azumaya variety $Y$ for a given $X$. The approach uses
moduli spaces of quiver representations. Our general reference
for this subject is~\cite{Ki}. We will freely use conventions
and notations thereof.

So, assume that we are given the data (D.1) and (D.2) satisfying
the conditions (C.1)--(C.3).
Assume additionally that
\begin{equation}\label{addc}
\Ext^{>0}(E_2(-1),E_1) = \Ext^{>0}(E_2(-1),E_2) = 0.
\end{equation}

Consider a quiver $\SQ = \xymatrix@1{ \bullet \ar[r]^{W} & \bullet }$
with two vertices and $\dim W$ arrows from the first vertex to the second,
and a slope function $\theta:K_0(\SQ) \to \ZZ$, given by the formula
\begin{equation}\label{chi}
\theta(d_1,d_2) = d_2\rank(E_2) - d_1\rank(E_1).
\end{equation}
For every dimension vector $(d_1,d_2)\in K_0(\SQ)$,
take a pair of vector spaces $R_1$, $R_2$ of dimensions $d_1$ and $d_2$
respectively, and consider an affine space
$\BR = \BR_\SQ(d_1,d_2) = \Hom(R_1\otimes W,R_2)$
parameterizing all $(d_1,d_2)$-dimensional representations of $\SQ$.
We denote the trivial vector bundles on $\BR$ with fibers $R_1$ and $R_2$
by $\CR_1$, $\CR_2$.
Then on $\BR$ we have a universal representation of the quiver
$\urho:W\otimes\CR_1 \to \CR_2$.
Consider on $X\times\BR$ a map $e:E_1\boxtimes\CR_1\to E_2\boxtimes\CR_2$
defined as the composition
$$
e:\xymatrix@1{E_1\boxtimes\CR_1 \ar[r]^-{\urho^*} &
E_1\boxtimes (W^*\otimes\CR_2) \cong
(W^*\otimes E_1)\boxtimes\CR_2 \ar[r]^-{\ev} &
E_2\boxtimes\CR_2},
$$
where $\ev:W^*\otimes E_1 = \Hom(E_1,E_2)\otimes E_1 \to E_2$
is the evaluation homomorphism.
For each $\rho\in\BR$ we denote by $e_\rho$ the restriction of $e$
to $X \cong X\times\{\rho\} \subset X\times\BR$.
Let
$$
C = \Coker e,\qquad
C_\rho = \Coker e_\rho \cong C_{|X\times\{\rho\}}.
$$
An embedding
$V^* \subset H^0(X,\CO_X(1)) \subset
\Hom(E_2(-1)\boxtimes \CR_2,E_2\boxtimes \CR_2)$
combined with the projection $E_2\boxtimes \CR_2 \to C$
induces a morphism
$V^*\otimes\CO_{X\times\BR} \to (E_2^*(1)\boxtimes \CR_2^*)\otimes C$.
Taking a pushforward with respect to the projection
$q:X\times\BR \to \BR$ we obtain a homomorphism
\begin{equation}\label{vsto}
\nu:V^*\otimes\CO_{\BR} \to R^0q_*((E_2^*(1)\boxtimes \CR_2^*)\otimes C).
\end{equation}

From now on assume that
$$
\theta(d_1,d_2) = 0.
$$
This condition determines the dimension vector $(d_1,d_2)$ up to a constant.
Let $\gcd(d_1,d_2) = m$ and denote the corresponding representation space
$\BR_\SQ(d_1,d_2)$ by $\BR_m$. Consider the following subschemes of $\BR_m$:
\begin{itemize}
\item the $\theta$-stable locus $\BR^{s}_m$;
\item the $\theta$-semistable locus $\BR^{ss}_m$;
\item an open subscheme
$\BRO{m} = \{\rho\in\BR_m\ |\ \Ker e_\rho = 0\} \subset \BR_m$;
\item a subset
$\BRP{m} = \{\rho\in\BR_m\ |\ \Ker\nu_\rho \ne 0\} \subset \BR_m$;
\item a subset
$\BRTP{m} = \{(\rho,H)\in\BR_m\times\PP(V^*)\ |\ \nu_\rho(H) = 0\} \subset
\BR_m\times\PP(V^*)$.
\end{itemize}

\begin{proposition}\label{mainqu}
We have $\BRTP{m} \subset \BRP{m}\times\PP(V^*)$ and
$\BRP{m} \subset \BRO{m} \subset \BR^{ss}_m$.
Moreover, the sheaf $R^0q_*((E_2^*(1)\boxtimes \CR_2^*)\otimes C)$
is locally free on $\BRO{m}$.
\end{proposition}
\begin{proof}
The first embedding is evident.
Assume that $\rho\in\BRP{m}$ and $0\ne H\in\Ker\nu_\rho$.
Note that the composition of $\nu_\rho$ with the canonical map
$$
(R^0q_*((E_2^*(1)\boxtimes \CR_2^*)\otimes C))_\rho \to
H^0(X,(E_2^*(1)\boxtimes \CR_2^*)\otimes C_\rho) \cong
\Hom(E_2(-1)\otimes R_2,C_\rho)
$$
coincides with the composition
$V^* \to \Hom(E_2(-1)\otimes R_2,E_2\otimes R_2) \to
\Hom(E_2(-1)\otimes R_2,C_\rho)$.
Therefore, the assumption implies that the dashed arrow in the diagram
$$
\xymatrix{
E_2(-1)\otimes R_2 \ar[r] \ar[d]^H \ar@{-->}[dr] & C_\rho(-1) \ar[d]^H \\
E_2\otimes R_2 \ar[r] & C_\rho
}
$$
vanishes. But since the upper arrow is surjective it follows
that the map $C_\rho(-1) \exto{H} C_\rho$ is zero, hence
$C$ is supported scheme-theoretically on the hyperplane
section $X_H$ of $X$.
\begin{equation}\label{brd}
\rho\in\BRP{m} \qquad\text{if and only if}\qquad
\supp(C_\rho)\subset H \in |\CH|,
\end{equation}
where $|\CH|$ is the linear system of hyperplane sections
of $X \subset \PP(V)$.
In particular, we have $\rank(C_\rho) = 0$, hence
$\rank\Ker e_\rho =
\rank (E_1\otimes R_1) - \rank(E_2\otimes R_2) =
-\theta(d_1,d_2) = 0$,
hence $e_\rho$ is a monomorphism, so $\BRP{m}\subset\BRO{m}$.
Moreover, we have the following exact sequence
\begin{equation}\label{equ_crho}
0 \to E_1\otimes R_1 \exto{e_\rho} E_2\otimes R_2 \to C_\rho \to 0
\end{equation}
for any $\rho\in\BRP{m}$.
Note also that $e_{|X\times\BRO{m}}$ is also a monomorphism
and we have on $X\times\BRO{m}$ the following exact sequence
\begin{equation}\label{equ_c}
0 \to E_1\boxtimes \CR_1 \exto{e} E_2\boxtimes \CR_2 \to C \to 0.
\end{equation}

Now let us check that $\BRO{m}\subset\BR^{ss}_m$.
Assume that $\rho\in\BRO{m}$ is unstable.
Then there exist subspaces $R'_1\subset R_1$, $R'_2\subset R_2$,
such that $\theta(\dim R'_1,\dim R'_2) < 0$. Then
$$
\rank\Ker e_{\rho'} \ge
\rank(E_1\boxtimes R'_1) - \rank(E_2\boxtimes R'_2) =
-\theta(\dim R'_1,\dim R'_2) > 0,
$$
hence $\Ker(e_{\rho'}:E_1\boxtimes R'_1 \to E_2\boxtimes R'_2) \ne 0$.
On the other hand, it is clear that $\Ker e_{\rho'} \subset \Ker e_\rho = 0$,
which gives a contradiction.

Comparing (\ref{equ_crho}) and (\ref{equ_c}) we see that
$C\otimes\CO_{X\times\rho}\cong C_\rho$. It follows that the sheaf
$C$ (and therefore $(E_2^*(1)\boxtimes \CR_2^*)\otimes C$ as well)
is flat over $\BRO{m}$.
Furthermore, using (\ref{equ_crho}) and
taking into account~(\ref{addc}) we deduce that
$H^p(X,(E_2^*(1)\otimes R_2^*)\otimes C_\rho) = 0$ for $p>0$
and the dimension of the zero cohomology doesn't depend on $\rho$.
Finally, for any closed point $\rho\in\BRO{m}$ we have
$$
q_*((E_2^*(1)\boxtimes \CR_2^*)\otimes C)\otimes\CO_\rho \cong
H^\bullet(X,(E_2^*(1)\otimes R_2^*)\otimes C\otimes\CO_{X\times \rho}).
$$
Therefore,
$q_*((E_2^*(1)\boxtimes \CR_2^*)\otimes C) =
R^0q_*((E_2^*(1)\boxtimes \CR_2^*)\otimes C)$ is a vector bundle.
\end{proof}

\begin{lemma}
If
\begin{equation}\label{codimbro}
\codim_{\BR_m}(\BR_m\setminus\BRO{m}) \ge 2
\end{equation}
then $\BRP{m}$ is closed in $\BR_m$ and
$\BRTP{m}$ is closed in $\BR_m\times\PP(V^*)$.
\end{lemma}
\begin{proof}
It is clear that $\BRTP{m}$ is closed in $\BRP{m}\times\PP(V^*)$ and
$\BRP{m}$ is the degeneration locus of a the morphism
$\nu:V^*\otimes\CO_{\BR_m} \to R^0q_*((E_2^*(1)\boxtimes \CR_2^*)\otimes C)$,
so it suffices to check that the sheaf
$R^0q_*((E_2^*(1)\boxtimes \CR_2^*)\otimes C)$
is torsion free.
It follows from~(\ref{addc}) and (\ref{equ_c}) that we have
the following exact sequence on $\BR_m$:
$$
0 \to \Hom(E_2(-1),E_1)\otimes\CR_2^*\otimes\CR_1 \to
\Hom(E_2(-1),E_2)\otimes\CR_2^*\otimes\CR_2 \to
R^0q_*((E_2^*(1)\boxtimes \CR_2^*)\otimes C) \to 0,
$$
so the sheaf admits a length 2 locally free resolution.
Therefore, it is torsion free if and only if it is locally free
in codimension $2$ (see \cite{OSS}). But we have proved
in proposition~\ref{mainqu} that it is locally free on~$\BRO{m}$.
\end{proof}

The following lemma is useful for verification of~(\ref{codimbro}).

\begin{lemma}\label{brombro1}
If $\codim_{\BR_1}(\BR_1\setminus\BRO{1}) \ge 2$ then
$\codim_{\BR_m}(\BR_m\setminus\BRO{m}) \ge 2$.
\end{lemma}
\begin{proof}
Consider the embedding $(\BR_1)^m \to \BR_m$,
$(\rho_1,\dots,\rho_m) \mapsto \rho_1\oplus\dots\oplus\rho_m$.
Since $\BR_m$ is a vector space and $\BRO{m}$ is dilations-invariant,
it suffices to check that
$\codim_{(\BR_1)^m}((\BR_1)^m\setminus(\BRO{m}\cap(\BR_1)^m)) \ge 2$.
But it is clear that $(\BRO{m}\cap(\BR_1)^m) = (\BRO1)^m$
and the claim follows.
\end{proof}

For any divisor class $D$ such that both linear systems
$|D|$ and $|\CH-D|$ are nonempty let $Z^{|D|} \subset \PP(V^*)$
denote the image of $|D|\times|\CH-D|$ in $|\CH| = \PP(V^*)$.
Let $\tg_m$ denote the projection $\BRTP{m} \to \PP(V^*)$ and
let $p$ denote the projection $\BRTP{m} \to \BRP{m}$.
Consider the open subset
$$U :=
\BRTP{m} \setminus \mathop{\cup}\limits_{0 < D < \CH} \tg_m^{-1}(Z^{|D|})
\subset \BRTP{m}.
$$

\begin{proposition}\label{brpmd}
The projection $p:U \to p(U) \subset \BRP{m}$ is an isomorphism.
\end{proposition}
\begin{proof}
Let us check that $\Ker\nu_{|p(U)}$ is a line subbundle
in $V^*\otimes\CO_{p(U)}$. Since $\BRP{m}$ is the degeneration locus
of a morphism of vector bundles
$\nu:V^*\otimes\CO_{\BRO{m}} \to R^0q_*((E_2^*(1)\boxtimes \CR_2^*)\otimes C)$
on $\BRO{m}$ it suffices to check that $\dim\Ker\nu_\rho = 1$
for all $\rho\in p(U)$. Indeed, assume that $\dim\Ker\nu_\rho \ge 2$
and let $\langle H,H'\rangle$ be a $2$-dimensional subspace in $\Ker\nu_\rho$.
Then (\ref{brd}) implies that $\supp(C_\rho) \subset H \cap H'$.
On the other hand, it follows from~(\ref{equ_crho}) that
either $C_\rho = 0$, which would imply that $\Hom(E_2,E_1)\ne 0$
thus contradicting (C.3), or that $D = \supp(C_\rho)$ is a divisor in $X$.
Since $D \subset H\cap H'$, we have $H \in Z^{|D|}$
and $(\rho,H)\in \tg_m^{-1}(Z^{|D|})$.

Consider a map $p(U) \to \PP(V^*)$ given by the line subbundle
$\Ker\nu_{|p(U)} \subset V^*\otimes\CO_{|p(U)}$.
Combining it with the embedding $p(U)\subset\BRP{m}$ we obtain a map to
$p(U) \to U \subset \BRP{m}\times\PP(V^*)$, which is clearly
the inverse to the projection $p$.
\end{proof}

Consider the group $G = \left(\GL(R_1)\times\GL(R_2)\right)/\kk^*$.
Let $\chi:G \to \kk^*$ denote the character
$$
\chi(g_1,g_2) = \det(g_1)^{-\rank(E_1)}\det(g_2)^{\rank(E_2)},
$$
corresponding to the slope function~(\ref{chi}).
Consider the action of $G$ on $\BR_m$ by isomorphisms of representations
linearized by the character $\chi$, and
the trivial action of $G$ on $\PP(V^*)$.
Note that $\BRP{m}\subset\BR^{ss}_m$
and $\BRTP{m} \subset \BR^{ss}_m\times\PP(V^*)$
are $G$-invariant closed subschemes (if (\ref{codimbro}) is satisfied).
Consider the following GIT quotients:
\begin{itemize}
\item
$\CM_m^{ss} = \BR^{ss}_m\git G$,
the moduli space of semistable representations;
\item
$\CM_m^{s} = \BR^{s}_m\git G \subset \CM_m^{ss}$,
the moduli space of stable representations;
\item
$\BY_m = \BRTP{m}\git G \subset \CM_m^{ss}\times\PP(V^*)$;
\end{itemize}
Let $\bg_m:\BY_m\to\PP(V^*)$ denote the projection to $\PP(V^*)$.
Let
$$
Z_m = \bg_m(\BY_m) \subset \PP(V^*),\qquad
\OZ_m = (\bigcup_{k=1}^{[m/2]} Z_k) \cup (\bigcup_{0 < D < \CH} Z^{|D|}).
$$
Then $\OZ_m$ is closed in $\PP(V^*)$. Consider
$$
Y_m = \BY_m \setminus \bg_m^{-1}(\OZ_m) = \bg_m^{-1}(\PP(V^*)\setminus\OZ_m)
$$
and let $g_m:Y_m \to \PP(V^*)$ denote the restriction
of $\bg_m$ to $Y_m$.

\begin{lemma}\label{gmpro}
The map $g_m:Y_m \to \PP(V^*) \setminus \OZ_m$ is projective and
the projection $Y_m \to \CM^s_m$ is an embedding.
\end{lemma}
\begin{proof}
Since $\BRTP{m}$ is closed in $\BR^{ss}_m\times\PP(V^*)$ it follows that
$\BY_m$ is closed in $\CM_m^{ss}\times\PP(V^*)$, which is a projective variety
by results of~\cite{Ki}. Therefore, $\BY_m$ is projective
and the map $\bg_m$ is projective as well.
On the other hand, the map $g_m$ is obtained from $\bg_m$
by a change of base $\PP(V^*)\setminus\OZ_m \to \PP(V^*)$,
hence it is also projective.
Further note that any $(d_1,d_2)$-dimensional
strictly semistable representation $\rho$ of $\SQ$ with $\gcd(d_1,d_2)=m$
contains a $(d'_1,d'_2)$-dimensional grading factors with
$\gcd(d'_1,d'_2)\le [m/2]$ in its Jordan--G\"older filtration.
Hence $\bg_m(\rho)\in \OZ_m$ and $\rho\not\in Y_m$.
So, the projection of $Y_m$ in $\CM^{ss}_m$ is contained in $\CM^s_m$.
Moreover, by proposition~\ref{brpmd} this projection is an embedding
(compare the definition of $Y_m\subset\BY_m$ and $U\subset\BRTP{m}$).
\end{proof}

\begin{lemma}\label{resdia}
There exists a sheaf of Azumaya algebras $\CA_m$ on $\CM_m^s$,
such that the universal family $(\CR_1,\CR_2)$ of representations
of the quiver $\SQ$ on $\BR^s_m$ descends to a universal family
$(F_1,F_2)$ in the category of coherent $\CA_m$-modules on $\CM_m^s$.
Moreover, we have a right exact sequence
$$
W\otimes F_2^*\boxtimes F_1 \to
F_2^*\boxtimes F_2 \oplus F_1^*\boxtimes F_1 \to
\Delta_*\CA_\CM \to 0
$$
on $\CM_m^s\times \CM_m^s$.
\end{lemma}
\begin{proof}
Since the action of $G$ on $\BR^s_m$ is free
every $\GL(R_1)\times\GL(R_2)$-equivariant vector bundle on $\BR^{s}_m$
descends to $\CM_m^s$ iff the diagonal subgroup
$\kk^*\subset\GL(R_1)\times\GL(R_2)$ acts on it trivially.
However, the diagonal subgroup acts on both $\CR_1$ and $\CR_2$
with weight~$1$. To get rid of this action
we choose one more representation $R_0$ of the group
$\GL(R_1)\times\GL(R_2)$, on which the diagonal subgroup
acts with weight $1$ and consider the trivial vector bundle $\CR_0$
with fiber $R_0$ on $\BR$. Then the diagonal subgroup acts trivially on
$(\CR_0^*\otimes\CR_1,\CR_0^*\otimes\CR_2)$, hence this family descends
to $\CM_m^s = \BR^s_m//_\chi G$ and gives a family of $(F_1,F_2)$
of representations of the quiver $\SQ$ on $\CM_m^s$, which is quasiuniversal.
Moreover, $\CEnd(\CR_0)$ is a sheaf of $G$-equivariant Azumaya algebras,
which descends to a sheaf of Azumaya algebras $\CA_m$ over $\CM_m^s$.
It is clear that $(F_1,F_2)$ is a universal family of representations
of the quiver $\SQ$ in the category of coherent $\CA_m$-modules.

It remains to check that the above sequence is right exact.
The action of the group $G$ on $\BR^s_m$ gives a map
$\pi:\BR^s_m\times G \to \BR^s_m\times \BR^s_m$,
$(\rho,g)\mapsto (\rho,g\rho)$.
Let us check that the complex
$$
\CC =
\{
\xymatrix@1{
W\otimes \CR_2^*\boxtimes \CR_1 \ar[r]^-{\ad(\urho)} &
\CR_2^*\boxtimes \CR_2 \oplus \CR_1^*\boxtimes\CR_1 \ar[r] &
\pi_*\CO
}
\},
$$
on $\BR^s_m\times\BR^s_m$ is acyclic in the middle and right term.
For any point $(\rho_1,\rho_2) \in \BR_m\times\BR_m$ we have
an exact sequence
$$
0 \to \Hom_\SQ(\rho_2,\rho_1) \to
\CR_{2\rho_1}\otimes\CR_{2\rho_2}^* \oplus \CR_{1\rho_1}\otimes\CR_{1\rho_2}^* \to
W^*\otimes\CR_{2\rho_1}\otimes\CR_{1\rho_2}^* \to
\Ext^1_\SQ(\rho_2,\rho_1) \to 0.
$$
Since both $\rho_1$ and $\rho_2$ are stable,
we have $\Hom_\SQ(\rho_2,\rho_1) = 0$
for $(\rho_1,\rho_2)\not\in\pi(\BR^s_m\times G)$,
and $\Hom_\SQ(\rho_2,\rho_1) = \kk$
for $(\rho_1,\rho_2)\in\pi(\BR^s_m\times G)$.
Therefore, $\CH^0(\CC)$ and $\CH^{-1}(\CC)$ are supported set-theoretically
on $\pi(\BR^s_m\times G)$.

Now compare two spectral sequences computing $\CH^n(\pi^*\CC)$.
The first term of the first spectral sequence
$E_1^{p,q} = \CH^q(\pi^*(\CC^p)) \Longrightarrow \CH^{p+q}(\pi^*\CC)$ has form
$$
\xymatrix@R=2pt{
W\otimes \CR_2^*\otimes \CR_1 \ar[r]^-{\ad(\urho)} \ar@{-->}[drr] &
\CR_2^*\otimes \CR_2 \oplus \CR_1^*\otimes\CR_1 \ar[r] &
\CO \\
0 & 0 & \CH^{-1}(\pi^*\pi_*\CO) \\
0 & 0 & \CH^{-2}(\pi^*\pi_*\CO)
}
$$
Note that the top line is acyclic in the middle and right term, its cohomology
in the left term at a point $(\rho,g)\in\BR^s_m\times G$ is equal to
$\Ext_\SQ^1(\rho,\rho)^* \cong
T^*_\rho\CM_m^s \cong
(\CH^{-1}(\pi^*\pi_*\CO))_{(\rho,g)}$,
and it is clear that the dashed arrow is an isomorphism.
It follows that $\pi^*\CC \in \D^{\le -2}(\BR^s_m\times G)$.
Then the second spectral sequence
$E_2^{p,q} = \CH^{q}(\pi^*(\CH^p(\CC))) \Longrightarrow \CH^{p+q}(\pi^*\CC)$
implies that $\CH^{0}(\pi^*\CH^0(\CC)) = \CH^{0}(\pi^*\CH^{-1}(\CC)) = 0$,
hence $\CH^0(\CC) = \CH^{-1}(\CC) = 0$ since both are supported
set-theoretically on $\pi(\BR^s_m\times G)$. Tensoring $\CC$ with
$\CR_0\boxtimes\CR_0^*$ we obtain a complex
$$
W\otimes (\CR_0^*\otimes\CR_2)^*\boxtimes (\CR_0^*\otimes\CR_1)
\exto{\ad(\urho)}
(\CR_0^*\otimes\CR_2)^*\boxtimes (\CR_0^*\otimes\CR_2) \oplus
(\CR_0^*\otimes\CR_1)^*\boxtimes (\CR_0^*\otimes\CR_1) \to
\pi_*\CEnd(\CR_0)
$$
which is acyclic in the middle and right term. Therefore,
after descent to $\CM_m^s\times\CM_m^s$ it gives a complex
$W\otimes F_2^*\boxtimes F_1 \exto{\ad(\urho)}
F_2^*\boxtimes F_2 \oplus F_1^*\boxtimes F_1 \to \Delta_*\CA_\CM$
acyclic in the middle and in the right term.
\end{proof}

We restrict the sheaf of Azumaya algebras $\CA_m$
and $\CA_m$-modules $(F_1,F_2)$ from $\CM_m^s$ to $Y_m$
and denote the restrictions by the same symbols.

\begin{theorem}\label{qu_th}
Assume that we have a data {\rm(D.1)}, {\rm(D.2)} satisfying
the conditions {\rm(C.1)--(C.3)}. Assume that additional
conditions~$(\ref{addc})$ and $(\ref{codimbro})$ are satisfied.
Then the Azumaya variety $(Y_m,\CA_m)$ with the map
$g_m:Y_m \to \PP(V^*)\setminus\OZ_m$, the pair of $\CA_m$-modules $(F_1,F_2)$,
and the canonical morphism $W \to \Hom(F_1,F_2)$ satisfy
the conditions {\rm(C.4)--(C.6)}.
\end{theorem}
\begin{proof}
Condition~(C.4) is proved in lemma~\ref{gmpro}, and
condition~(C.5) is proved in lemma~\ref{resdia}.
It remains to check condition~(C.6).
To this end denote
$\TU = \BRTP{m} \setminus \tg_m^{-1}(\OZ_m) =
\BRP{m} \setminus q(\tg_m^{-1}(\OZ_m))$,
so that $Y_m = \TU\git G$, and recall that the composition
$$
\tg_m^*\CO_{\PP(V^*)}(-1) \to
V^*\otimes\CO_{U} \exto{\nu}
q_*((E_2^*(1)\boxtimes \CR_2^*)\otimes C)
$$
on $\TU$ is zero. By adjunction we deduce that the composition
$$
q^*\tg_m^*\CO_{\PP(V^*)}(-1) \to
V^*\otimes\CO_{X\times U} \exto{\nu}
(E_2^*(1)\boxtimes \CR_2^*)\otimes C
$$
on $X\times\TU$ is zero. Tensoring by $q^*\tg_m^*\CO_{\PP(V^*)}(1)$
we deduce that $\nu$ takes the image $\lambda$ of the canonical element
$\id_V \in V^*\otimes V =
H^0(X,\CO_X(1))\otimes H^0(\PP(V^*),\CO_{\PP(V^*)}(1))$ in
$H^0(X\times U,\CO_X(1)\boxtimes\tg_m^*\CO_{\PP(V^*)}(1))$
to
$0\in H^0(X\times U,
(E_2^*(1)\boxtimes \CR_2^*\otimes\tg_m^*\CO_{\PP(V^*)}(1))\otimes C)$.
This means that  the composition of the left and bottom arrows
in the following commutative diagram
$$
\xymatrix{
E_2(-1)\boxtimes (\CR_2\otimes\tg_m^*\CO_{\PP(V^*)}(-1))
\ar[r] \ar[d]^{\lambda} &
C\otimes(\CO_X(-1)\boxtimes\tg_m^*\CO_{\PP(V^*)}(-1))
\ar[d]^{\lambda} \\
E_2\boxtimes \CR_2 \ar[r] & C
}
$$
vanishes. But the upper arrow is an epimorphism. Hence the right
arrow also vanishes. This means, that $C$ is supported
scheme-theoretically on the zero locus of $\lambda$.
But it is clear that after going down to the quotient
$X\times Y_m = (X\times \TU)\git G$ the bundle
$\CO_X(1)\boxtimes\tg_m^*\CO_{\PP(V^*)}(1)$ descends to
$\CO_X(1)\boxtimes\CO_{Y_m}(1)$ and its section $\lambda$
descends to the canonical section of $\CO_X(1)\boxtimes\CO_{Y_m}(1)$.
\end{proof}

The problem thus reduces to the choice of $m$ such that
conditions~(C.7) and (C.8) are satisfied. I don't know
how to make this choice in general. See however the following examples
where the choices are made explicitly.

\section{Examples}

Examples in this section are ordered by the complexity of the variety $Y$.
In examples~\ref{g25} and~\ref{og510} the variety $Y$ can be described
explicitly and $\CA_Y = \CO_Y$. In examples~\ref{lg36}, \ref{g2g27}
and~\ref{iq6} we use the approach of section~\ref{s5}, consider
the corresponding varieties $Y_2$, $Y_3$ and $Y_4$ respectively,
and only after that give their description.

\subsection{Grassmannian}\label{g25}

Let $W = \kk^5$ and
$$
X = \Gr(2,W) = \Gr(2,5),\qquad
(E_1,E_2) = (\CO_X,\CU_X^*),
$$
where $\CU_X \subset W\otimes\CO_X$ is the tautological rank 2 subbundle.
Let $V = \Lambda^2 W = \kk^{10}$, and
let $f:X\to\PP(V)$ be the Pl\"ucker embedding.
Take
$$
Y = \Gr(2,W^*) = \Gr(2,5),\qquad
\CA_Y = \CO_Y,\qquad
(F_1,F_2) = (\CU_Y,\CO_Y),
$$
where $\CU_Y \subset W^*\otimes\CO_Y$ is the tautological rank 2 subbundle.
Let $g:Y\to\PP(V^*)$ be the Pl\"ucker embedding (i.e.\ $\BZ = \emptyset$).
Then $\Hom(E_1,E_2) = W^*$, $\Hom(F_1,F_2) = W$,
and we take $\phi = \id_W$.

It is easy to check that all conditions (C.1)--(C.8) are satisfied,
so theorem~\ref{themain} applies, and we have the following
semiorthogonal decompositions
(note that $\codim_{\PP(V^*)}(Y) = 3$, therefore for any admissible
subspace $L\subset V^*$ of dimension $\dim L\le 3$ we have $Y_L = \emptyset$)
$$
\begin{array}{lll}
\D^b(X) & = &
\langle E_1(1),E_2(1),E_1(2),E_2(2),E_1(3),E_2(3),E_1(4),E_2(4),E_1(5),E_2(5) \rangle \\
\D^b(X_1) & = &
\langle E_1(1),E_2(1),E_1(2),E_2(2),E_1(3),E_2(3),E_1(4),E_2(4) \rangle \\
\D^b(X_2) & = &
\langle E_1(1),E_2(1),E_1(2),E_2(2),E_1(3),E_2(3) \rangle \\
\D^b(X_3) & = &
\langle E_1(1),E_2(1),E_1(2),E_2(2) \rangle \\
\D^b(X_4) & = &
\langle \D^b(Y_4),E_1(1),E_2(1) \rangle \\
\D^b(X_5) & = &
\D^b(Y_5)
\end{array}
$$
where $Y_r = Y\cap\PP(L)$ is a linear section of $Y$ by an admissible
subspace $L\subset V^*$ of dimension~$r$, and
$X_r = X\cap\PP(L^\perp)$ is the corresponding
orthogonal linear section of $X$ of codimension~$r$.
Note that $X_1$, $X_2$ and $X_3$ are always smooth,
$X_3$ is the Fano threefold of index~$2$ and degree~$5$ ($V_5$~threefold),
$X_4$ is a del Pezzo surface of degree~$5$ (possibly singular),
$Y_4$ is a zero-dimensional scheme of length~$\deg Y = 5$,
and $X_5$, $Y_5$ are elliptic curves of degree~$5$ (possibly singular).

In particular, the results of Orlov~\cite{O1} about the derived category
of the Fano threefold $V_5$ follow from theorem~\ref{themain}.

\subsection{Orthogonal Grassmannian}\label{og510}

Let $W = \kk^{10}$ and
$$
X = \OGr_+(5,W),\qquad
(E_1,E_2) = (\CO_X,\CU_X^*),
$$
where $\OGr_+(5,W)$ is a connected component of the Grassmannian
of $5$-dimensional subspaces in~$W$, isotropic with respect to
a chosen nondegenerate quadratic form $\bq\in S^2W^*$, and
$\CU_X \subset W\otimes\CO_X$ is the tautological rank 5 subbundle.
Let $V \cong \kk^{16}$ be the corresponding half-spinor representation
of the spin-group $\Spin(W)$ (see \cite{Ch}), and let $f:X\to\PP(V)$ denote
the canonical embedding.
Take
$$
Y = \OGr_-(5,W),\qquad
\CA_Y = \CO_Y,\qquad
(F_1,F_2) = (\CU_Y,\CO_Y),
$$
where $\CU_Y \subset W\otimes\CO_Y$ is the tautological rank 5 subbundle.
We denote by $g:Y\to\PP(V^*)$ the composition of the embedding $Y\to\PP(V)$
and of the isomorphism $\PP(V)\to\PP(V^*)$ given by $\bq$.
Then
$\Hom(E_1,E_2) \cong W^*$,
$\Hom(F_1,F_2) \cong W^*$,
and we take $\phi = \bq^{-1}$.

It is easy to check that all conditions (C.1)--(C.8) are satisfied,
so theorem~\ref{themain} applies, and we have the following
semiorthogonal decompositions
(note that $\codim_{\PP(V^*)}(Y) = 5$, therefore for any admissible
subspace $L\subset V^*$ of dimension $\dim L\le 5$ we have $Y_L = \emptyset$)
$$
\footnotesize
\arraycolsep = 3pt
\begin{array}{lll}
\D^b(X) & = &
\langle E_1(1),E_2(1),E_1(2),E_2(2),E_1(3),E_2(3),E_1(4),E_2(4),
E_1(5),E_2(5),E_1(6),E_2(6),E_1(7),E_2(7),E_1(8),E_2(8) \rangle \\
\D^b(X_1) & = &
\langle E_1(1),E_2(1),E_1(2),E_2(2),E_1(3),E_2(3),E_1(4),E_2(4),
E_1(5),E_2(5),E_1(6),E_2(6),E_1(7),E_2(7) \rangle \\
\D^b(X_2) & = &
\langle E_1(1),E_2(1),E_1(2),E_2(2),E_1(3),E_2(3),E_1(4),E_2(4),
E_1(5),E_2(5),E_1(6),E_2(6) \rangle \\
\D^b(X_3) & = &
\langle E_1(1),E_2(1),E_1(2),E_2(2),E_1(3),E_2(3),E_1(4),E_2(4),
E_1(5),E_2(5) \rangle \\
\D^b(X_4) & = &
\langle E_1(1),E_2(1),E_1(2),E_2(2),E_1(3),E_2(3),E_1(4),E_2(4) \rangle \\
\D^b(X_5) & = &
\langle E_1(1),E_2(1),E_1(2),E_2(2),E_1(3),E_2(3) \rangle \\
\D^b(X_6) & = &
\langle \D^b(Y_6),E_1(1),E_2(1),E_1(2),E_2(2) \rangle \\
\D^b(X_7) & = &
\langle \D^b(Y_7),E_1(1),E_2(1) \rangle \\
\D^b(X_8) & = & \D^b(Y_8)
\end{array}
$$
where $Y_r = Y\cap\PP(L)$ is a linear section of $Y$ by an admissible
subspace $L\subset V^*$ of dimension~$r$, and
$X_r = X\cap\PP(L^\perp)$ is the corresponding
orthogonal linear section of $X$ of codimension~$r$.
Note that $X_1$, $X_2$, $X_3$, $X_4$ and $X_5$ are always smooth,
$X_6$ is a Fano fourfold of index $2$ (possibly singular),
$Y_6$~is a zero-dimensional scheme of length $\deg Y=12$,
$X_7$ is a Fano threefold $V_{12}$ (possibly singular),
$Y_7$ is a curve of genus $7$ (possibly singular),
and $X_8$, $Y_8$ are $K3$-surfaces of degree $12$ (possibly singular).

In particular, the results of~\cite{K3} about the derived category
of the Fano threefold $V_{12}$ follow from theorem~\ref{themain}.
Indeed when $X_7$ is a smooth Fano threefold $V_{16}$ then $Y_7$
is a smooth curve of genus~$7$, and we obtain a decomposition
$$
\D^b(V_{12}) = \langle \D^b(C_7),\CO,\CU^*\rangle,
$$
where $C_7$ is a curve of genus~7.


\subsection{Lagrangian Grassmannian}\label{lg36}

Let $W = \kk^{6}$ and
$$
X = \LGr(3,W),\qquad
(E_1,E_2) = (\CO_X,\CU_X^*),
$$
where $\LGr(3,W)$ is the Grassmannian of $3$-dimensional subspaces in $W$,
Lagrangian with respect to a chosen symplectic form $\sigma\in \Lambda^2W^*$,
and $\CU_X \subset W\otimes\CO_X$ is the tautological rank 3 subbundle.
The natural representation of the group $\SP(W)$ in the space
$\Lambda^3 W$ decomposes into the direct sum of representations
$$
\Lambda^3 W = W \oplus V.
$$
The Pl\"ucker embedding $\Gr(3,W) \subset \PP(\Lambda^3 W)$ restricts
to an embedding $f:X\to\PP(V)$.
Now we apply the construction of section~5 and take
$Y = Y_2$, $\CA_Y = \CA_2$, $\BZ = \OZ_2$, $g = g_2$,
the corresponding universal family $(F_1,F_2)$, and
the map $\phi$ induced by the $\SQ$-representation structure on $(F_1,F_2)$.

It is shown in Appendix A that all conditions (C.1)--(C.8) are satisfied,
and in fact
$$
\BZ = \mu(\Gr(3,W^*)) = \mu(\Gr(3,6)),
$$
where $\mu:\xymatrix@1{\PP(\Lambda^3 W^*) \ar@{-->}[r] & \PP(V^*)}$ is
the linear projection from $\PP(\sigma\wedge W^*)$, and
%
%
$$
Y = {X^\vee} \setminus \BZ,
$$
where ${X^\vee}$ is the projectively dual variety of $X$,
which is a quartic hypersurface, singular along $\BZ \subset \PP(V^*)$
(see~\cite{Ho}). Thus we have the following semiorthogonal decompositions
(note that $\codim_{\PP(V^*)}(Y) = 1$, therefore for any admissible
subspace $L\subset V^*$ of dimension $\dim L\le 1$ we have $Y_L = \emptyset$)
$$
\begin{array}{lll}
\D^b(X) & = &
\langle E_1(1),E_2(1),E_1(2),E_2(2),E_1(3),E_2(3),E_1(4),E_2(4) \rangle \\
\D^b(X_1) & = &
\langle E_1(1),E_2(1),E_1(2),E_2(2),E_1(3),E_2(3) \rangle \\
\D^b(X_2) & = &
\langle \D^b(Y_2),E_1(1),E_2(1),E_1(2),E_2(2) \rangle \\
\D^b(X_3) & = &
\langle \D^b(Y_3,\CA_Y),E_1(1),E_2(1) \rangle \\
\D^b(X_4) & = & \D^b(Y_4,\CA_Y)
\end{array}
$$
where $Y_r = Y\cap\PP(L)$ is a linear section of $Y$ by an admissible
subspace $L\subset V^*$ of dimension~$r$, and
$X_r = X\cap\PP(L^\perp)$ is the corresponding
orthogonal linear section of $X$ of codimension~$r$.
Note that $X_1$ is always smooth,
$X_2$ is a Fano fourfold of index $2$ (possibly singular),
$Y_2$ is a zero-dimensional scheme of length $\deg Y = 4$,
$X_3$ is a Fano threefold $V_{16}$ (possibly singular),
$Y_3$ is a plane quartic (possibly singular),
and $X_4$, $Y_4$ are $K3$-surfaces of degree $16$ and $4$ respectively
(possibly singular).

In particular, the results of~\cite{S1,S2} about the derived category
of the Lagrangian Grassmannian~$X$ and of its hyperplane section $X_1$
follow from theorem~\ref{themain}. Moreover, when $X_3$ is a smooth
Fano threefold $V_{16}$ then $Y_3$ is a smooth curve of genus~$3$,
hence $\D^b(Y_3,\CA_Y) \cong \D^b(Y_3)$ (the Brauer group
of a smooth curve is trivial), and we obtain a decomposition
$$
\D^b(V_{16}) = \langle \D^b(C_3),\CO,\CU^*\rangle,
$$
where $C_3$ is a curve of genus~3.


\subsection{$G_2$ Grassmannian}\label{g2g27}

Let $W = \kk^7$ and
$$
X = \GTGr(2,W),\qquad
(E_1,E_2) = (\CO_X,\CU_X^*),
$$
where $\GTGr(3,W)$ is the Grassmannian of the simple Lie group $G_2$,
realized as the zero locus of the section
$s_\lambda \in H^0(\Gr(2,W),\CU^\perp(1)) \cong \Lambda^3W^*$,
corresponding to a unique $G_2$-invariant $3$-form $\lambda$
on the tautological $G_2$-representation $W$,
$\CU \subset W\otimes\CO_{\Gr(2,W)}$ is the tautological rank 2 subbundle,
and $\CU_X$ is the restriction of $\CU$ to $X$.
The natural representation of the group $G_2$ in the space
$\Lambda^2 W$ decomposes into the direct sum of representations
$$
\Lambda^2 W = W^* \oplus V
$$
(the projection $\Lambda^2W \to W^*$ is given by $\lambda$).
The Pl\"ucker embedding $\Gr(2,W) \subset \PP(\Lambda^2 W)$ restricts
to an embedding $f:X\to\PP(V)$.
Now we apply the construction of section~5 and take
$Y = Y_3$, $\CA_Y = \CA_3$, $\BZ = \OZ_3$, $g = g_3$,
the corresponding universal family $(F_1,F_2)$, and
the map $\phi$ induced by the $\SQ$-representation structure on $(F_1,F_2)$.

It is shown in Appendix B that all conditions (C.1)--(C.8) are satisfied,
and in fact
$$
\BZ = \mu(\Gr(2,W^*)) = \mu(\Gr(2,7)),
$$
where $\mu:\xymatrix@1{\PP(\Lambda^2 W^*) \ar@{-->}[r] & \PP(V^*)}$ is
the linear projection from $\PP(\lambda(W))$, and
$$
\xymatrix@1{g:Y \ar[r]^-{2:1} & \PP(V^*) \setminus \BZ}
\quad\text{is the double covering ramified at}\quad
{X^\vee} \setminus \BZ \subset \PP(V^*) \setminus \BZ,
$$
where ${X^\vee}$ is the projectively dual variety of $X$ which is
a sextic hypersurface, singular along $\BZ \subset \PP(V^*)$
(see~\cite{Ho}). Thus we have the following semiorthogonal
decompositions
$$
\begin{array}{lll}
\D^b(X) & = &
\langle E_1(1),E_2(1),E_1(2),E_2(2),E_1(3),E_2(3) \rangle \\
\D^b(X_1) & = &
\langle \D^b(Y_1),E_1(1),E_2(1),E_1(2),E_2(2) \rangle \\
\D^b(X_2) & = &
\langle \D^b(Y_2,\CA_Y),E_1(1),E_2(1) \rangle \\
\D^b(X_3) & = & \D^b(Y_3,\CA_Y)
\end{array}
$$
where $Y_r = Y\cap\PP(L)$ is a linear section of $Y$ by an admissible
subspace $L\subset V^*$ of dimension~$r$, and
$X_r = X\cap\PP(L^\perp)$ is the corresponding
orthogonal linear section of $X$ of codimension~$r$.
Note that
$X_1$ is a Fano fourfold of index $2$ (possibly singular),
$Y_1$ is a zero-dimensional scheme of length~$\deg Y = 2$,
$X_2$ is a Fano threefold $V_{18}$ (possibly singular),
$Y_2$ is a curve of genus~$2$ (possibly singular),
and $X_3$,~$Y_3$ are $K3$-surfaces of degree $18$ and $2$ respectively
(possibly singular).

In particular, the results of~\cite{Ra} about the derived category
of the $G_2$-Grassmannian~$X$ follow from theorem~\ref{themain}.
Moreover, when $X_2$ is a smooth Fano threefold $V_{18}$
then $Y_2$ is a smooth curve of genus~$2$ (by projective duality),
hence $\D^b(Y_2,\CA_Y) \cong \D^b(Y_2)$ (the Brauer group of a smooth
curve is trivial), and we obtain a decomposition
$$
\D^b(V_{18}) = \langle \D^b(C_2),\CO,\CU^*\rangle
$$
where $C_2$ is a curve of genus $2$.


\subsection{Intersection of quadrics}\label{iq6}

Let $W = \kk^{6}$ and
$$
X = \PP(W) = \PP^5,\qquad
(E_1,E_2) = (\CO_X,\CO_X(1)).
$$
Let $V = S^2W = \kk^{21}$ and let $f:X \to \PP(V)$ be
the double Veronese embedding.

Now we apply the construction of section~5 and take
$Y = Y_4$, $\CA_Y = \CA_4$, $\BZ = \OZ_4$, $g = g_4$,
the corresponding universal family $(F_1,F_2)$, and
the map $\phi$ induced by the $\SQ$-representation structure on $(F_1,F_2)$.

It is shown in Appendix C that all conditions (C.1)--(C.8) are satisfied,
and in fact $\BZ\subset\PP(S^2W^*)$ is the locus of rank~4 (corank~2)
quadratic forms and
$$
\xymatrix@1{g:Y \ar[r]^-{2:1} & \PP(V^*) \setminus \BZ}
\quad\text{is the double covering ramified at}\quad
{X^\vee} \setminus \BZ \subset \PP(V^*) \setminus \BZ,
$$
where ${X^\vee}$ is the projectively dual variety of $X$ which is
the determinantal sextic hypersurface (the locus of degenerate
quadratic forms), singular along $\BZ \subset \PP(V^*)$.
Thus we have the following semiorthogonal decompositions
$$
\begin{array}{lll}
\D^b(X) & = &
\langle \CO_X(1),\CO_X(2),\CO_X(3),\CO_X(4),\CO_X(5),\CO_X(6) \rangle \\
\D^b(X_1) & = &
\langle \D^b(Y_1),\CO_X(1),\CO_X(2),\CO_X(3),\CO_X(4) \rangle \\
\D^b(X_2) & = &
\langle \D^b(Y_2,\CA_Y),\CO_X(1),\CO_X(2) \rangle \\
\D^b(X_3) & = & \D^b(Y_3,\CA_Y)
\end{array}
$$
where $Y_r = Y\cap\PP(L)$ is a linear section of $Y$ by an admissible
subspace $L\subset V^*$ of dimension~$r$, and
$X_r = X\cap\PP(L^\perp)$ is the corresponding
orthogonal linear section of $X$ of codimension~$r$.
Note that
$X_1$ is a four-dimensional quadric (possibly singular),
$Y_1$ is a zero-dimensional scheme of length $2$,
$X_2$ is a Fano threefold of index $2$ and degree $4$
(an intersection of two four-dimensional quadrics) (possibly singular),
$Y_2$ is a curve of genus~$2$ (possibly singular),
and $X_3$,~$Y_3$ are $K3$-surfaces of degree $8$ and $2$ respectively
(possibly singular).

In particular, the results of~\cite{BO1} about the derived category
of an intersection of two four-dimensional quadrics follow from
theorem~\ref{themain}.



It is also worth to note that the sheaf of Azumaya algebras $\CA_Y$
on $Y$ can be constructed in quite another way. More precisely,
one can check that the sheaf of algebras $g_*\CA_Y$ is isomorphic
to the even part of the universal Clifford algebra (see the introduction).


\section*{Appendix A. Lagrangian Grassmannian}

\setcounter{section}{0}
\refstepcounter{section}
\renewcommand{\thesection}{\Alph{section}}

In this appendix we show that a quartic hypersurface in $\PP^{13}$
is homologically projectively dual to the Lagrangian Grassmannian $\LGr(3,6)$.

Recall the notation. Let $W = \kk^{6}$ and
$$
X = \LGr(3,W),\qquad
(E_1,E_2) = (\CO_X,\CU_X^*),
$$
where $\LGr(3,W)$ is the Grassmannian of $3$-dimensional subspaces in $W$,
Lagrangian with respect to a chosen symplectic form $\sigma\in \Lambda^2W^*$,
and $\CU_X \subset W\otimes\CO_X$ is the tautological rank 3 subbundle.
The natural representation of the group $\SP(W)$ in the space
$\Lambda^3 W$ decomposes into the direct sum of representations
$$
\Lambda^3 W = W \oplus V.
$$
The Pl\"ucker embedding $\Gr(3,W) \subset \PP(\Lambda^3 W)$ restricts
to an embedding $f:X\to\PP(V)$.

\begin{lemma}\label{conds_a}
We have $\dim X = 6$, $\omega_X \cong \CO_X(-4)$.
Moreover, conditions {\rm(C.1)--(C.3)}\/ as well as the additional
conditions {\rm(\ref{addc})} and {\rm(\ref{codimbro})} are satisfied for $X$.
\end{lemma}
\begin{proof}
Note that $X\subset\Gr(3,W)$ is a zero locus of a regular section
$s_\sigma\in H^0(\Gr(3,W),\Lambda^2\CU^*)$. Hence
$$
\begin{array}{l}
\dim X = \dim\Gr(3,W) - \rank(\Lambda^2\CU^*) = 9 - 3 = 6,\\
\omega_X \cong \omega_{\Gr(3,W)|X} \otimes \det(\Lambda^2\CU^*) \cong
\CO_X(-6) \otimes \CO_X(2) \cong \CO_X(-4).
\end{array}
$$
In particular, (C.2) is true and $i_X = 4$.
Further, using the Koszul resolution
$\Lambda^\bullet(\Lambda^2\CU) \cong \CO_X$ of $X$ in $\Gr(3,W)$
it is easy to check that $H^0(X,\CO_X(1)) \cong V^*$, hence
$f(X)$ is not contained in a hyperplane, i.e.\ (C.1) is true.
Using again the Koszul resolution we compute $\Hom(E_s(k),E_t(l))$
via Borel--Bott--Weil theorem and check that (C.3) and
(\ref{addc}) is true. Finally, we note that
$\BR_1 = \Hom(\kk^3,W^*)$ and $\BRO1\subset\BR_1$
corresponds to the subset of embeddings $\kk^3\subset W^*$.
Therefore, $\codim_{\BR_1}(\BR_1\setminus\BRO1) = 4$
and $(\ref{codimbro})$ is satisfied by lemma~\ref{brombro1}.
\end{proof}

Now we are going to apply the construction of section~5 and
describe the space $Y = Y_2$ explicitly. The question turns out
to be rather cumbersome. However, it seems easier to describe instead
$\PP_{Y_2,\CA_2}(F_2^*)$, the moduli space of quotient $\CA_2$-modules
of rank~$1$. Indeed, the fiber of
$\PP_{Y_2,\CA_2}(F_2^*)$ at the point of $Y_2$ corresponding to
a $(6,2)$-dimensional representation $\rho = (R_1,R_2)$ of $\SQ$
coincides with $\PP(R_2)$. We show that any line $R'_2\subset R_2$
can be extended in a unique way to a $(1,1)$-dimensional subrepresentation
$\rho' = (R'_1,R'_2)$. The moduli space of $(1,1)$-dimensional
representations of $\SQ$ is $\PP(W)$. We claim that defined in this way map
$\PP_{Y_2,\CA_2}(F_2^*) \to \PP(W)$ identifies $\PP_{Y_2,\CA_2}(F_2^*)$
with an open subset in a $\PP^8$-bundle over $\PP(W)$, explicitly
$\TY = \PP_{\PP(W)}(\Omega_{\PP(W)}^3(3)/\CO_{\PP(W)}(-1))$.
To prove this we go the other way round. We consider
the space $\TY$, and identify its open subset with
$\PP_{Y_2,\CA_2}(F_2^*)$.

We start with several useful lemmas.

\begin{lemma}\label{yz1_a}
We have $\BY_1 = \OZ_1 = \mu(\Gr(3,W^*))$, where
$\mu:\xymatrix@1{\PP(\Lambda^3W^*) \ar@{-->}[r] & \PP(V^*)}$
is the linear projection from
$\PP(\sigma\wedge W^*) \subset \PP(\Lambda^3W^*)$.
\end{lemma}
\begin{proof}
If $m=1$ then $(d_1,d_2) = (3,1)$. A $(3,1)$-dimensional representation
$\rho:W\otimes R_1 \to R_2$ is $\chi$-semistable iff the corresponding map
$R_1\otimes R_2^* \to W^*$ is an embedding. Thus $\CM_1^{ss}$ gets
identified with $\Gr(3,W^*)$. Moreover, it is easy to see that
in this case every hyperplane section of $X$ is irreducible,
and we have $\BRTP1 = \BRP1 = \BRO1 = \BR_1^{ss}$ and the lemma follows.
\end{proof}

Consider a product $X\times\PP(W)$ and a subscheme
$\PP_X(\CU_X) \subset X\times\PP(W)$.

\begin{lemma}\label{ppxu}
The following sequences on $X\times\PP(W)$ are exact
{\small
$$
0 \to
\CO_X(-1)\boxtimes\CO_{\PP(W)}(-3) \to
\Lambda^2\CU_X\boxtimes\CO_{\PP(W)}(-2) \to
\CU_X\boxtimes\CO_{\PP(W)}(-1) \to
\CO_X\boxtimes\CO_{\PP(W)} \to
\CO_{\PP_X(\CU_X)} \to 0,
$$
$$
0 \to
\CO_X(-1)\boxtimes\CO_{\PP(W)}(-3) \to
\CU^*_X(-1)\boxtimes\CO_{\PP(W)}(-2) \to
\CO_X\boxtimes\Omega_{\PP(W)} \to
\CU^*_X\boxtimes\CO_{\PP(W)}(-1) \to
\CO_{\PP_X(\CU_X)} \to 0.
$$
}
\end{lemma}
\begin{proof}
The first sequence is just the Koszul resolution of
$\PP_X(\CU_X) \subset X\times\PP(W)$ (note that $\CU_X^\perp \cong \CU_X$,
the isomorphism is given by $\sigma$).
Further, it follows from exact sequences
$$
\xymatrix@R=12pt{
0 \ar[r] &
\CU_X\boxtimes\CO_{\PP(W)}(-1) \ar[r] &
W^*\otimes\CO_X\boxtimes\CO_{\PP(W)}(-1) \ar[r] \ar@{=}[d]&
\CU_X^*\boxtimes\CO_{\PP(W)}(-1) \ar[r] &
0 \\
0 \ar[r] &
\CO_X\boxtimes\Omega_{\PP(W)} \ar[r] &
W^*\otimes\CO_X\boxtimes\CO_{\PP(W)}(-1) \ar[r] &
\CO_X\boxtimes\CO_{\PP(W)} \ar[r] &
0
}
$$
that the kernel and the cokernel of the composition maps
$\CU_X\boxtimes\CO_{\PP(W)}(-1) \to \CO_X\boxtimes\CO_{\PP(W)}$
and
$\CO_X\boxtimes\Omega_{\PP(W)} \to \CU^*_X\boxtimes\CO_{\PP(W)}(-1)$
coincide. Taking also into account that
$\Lambda^2\CU_X \cong \CU_X^*(-1)$ we obtain
the second sequence from the first one.
\end{proof}

Consider the truncated first sequence of lemma~\ref{ppxu}
as a resolution of the sheaf of ideals $J_{\PP_X(\CU_X)}$.
Tensoring it with $\CO_X(1)$ and pushing forward to $\PP(W)$
we obtain exact sequence
$$
0 \to \CO_{\PP(W)}(-3) \to
W^*\otimes\CO_{\PP(W)}(-2) \to
(\Lambda^2W^*/\kk\sigma)\otimes\CO_{\PP(W)}(-1) \to
{p_2}_*(J_{\PP_X(\CU_X)}\otimes\CO_X(1)) \to 0,
$$
where $p_2:X\times\PP(W) \to \PP(W)$ is the projection.
Taking into account exact sequence
$$
0 \to \CO_{\PP(W)}(-3) \to
W^*\otimes\CO_{\PP(W)}(-2) \to
\Lambda^2W^*\otimes\CO_{\PP(W)}(-1) \to
\Omega^3_{\PP(W)}(3) \to 0
$$
we can rewrite it as
$$
0 \to \CO_{\PP(W)}(-1) \exto{\sigma} \Omega^3_{\PP(W)}(3) \to
{p_2}_*(J_{\PP_X(\CU_X)}\otimes\CO_X(1)) \to 0.
$$

\begin{lemma}
Morphism $\CO_{\PP(W)}(-1) \exto{\sigma} \Omega^3_{\PP(W)}(3)$ is an embedding
of vector bundles on $\PP(W)$.
\end{lemma}
\begin{proof}
This morphism is $\SP(W)$-equivariant
and $\SP(W)$ acts transitively on $\PP(W)$.
\end{proof}


Consider the projectivization of the quotient bundle
$$
\TY = \PP_{\PP(W)}(\Omega_{\PP(W)}^3(3)/\CO_{\PP(W)}(-1)) =
\PP_{\PP(W)}({p_2}_*(J_{\PP_X(\CU_X)}\otimes\CO_X(1))).
$$
The natural embedding
${p_2}_*(J_{\PP_X(\CU_X)}\otimes\CO_X(1)) \to
H^0(X,\CO_X(1))\otimes\CO_{\PP(W)} =
V^*\otimes\CO_{\PP(W)}$
induces a morphism $\TY \exto{\tg} \PP(V^*)$.
Let $\pi:\TY\to\PP(W)$ denote the canonical projection.
Denote also
$$
\CO_\TY(H) := \tg^*\CO_{\PP(V)}(1),\quad
\tg:\TY \to \PP(V^*),\qquad\text{and}\qquad
\CO_{\TY}(h) := \pi^*\CO_{\PP(W)}(1),\quad
\pi:\TY\to\PP(W).
$$
%
%
Consider the product $X\times\TY$, the projection
$\id_X\times\pi:X\times\TY \to X\times\PP(W)$,
the subschemes
$$
\Sigma = (\id_X\times\pi)^{-1}(\PP_X(\CU_X)),\qquad
\TQ = Q(X,\TY) = (X\times\TY)\times_{\PP(V)\times\PP(V^*)}Q,
$$
and the equation
$s_\TQ\in
H^0(X\times\TY,\CO_X(1)\boxtimes\CO_\TY(H)) =
H^0(X\times\TY,(f\times\tg)^*\CO_{\PP(V)\times\PP(V^*)}(1,1)) \cong
V^*\otimes V$
of $\TQ$ in $X\times\TY$.

\begin{lemma}\label{sqsigma}
The section $s_\TQ$ vanishes on the scheme $\Sigma$.
\end{lemma}
\begin{proof}
We have on $\TY$ a canonical embedding
$$
\CO_\TY(-H) \to \pi^*{p_2}_*(J_{\PP_X(\CU_X)}\otimes\CO_X(1)) \cong
{p_2}_*(J_{\Sigma}\otimes\CO_X(1)).
$$
Tensoring it by $\CO_\TY(H)$ we obtain a canonical morphism
$\CO_{X\times\TY} \to J_{\Sigma}\otimes(\CO_X(1)\boxtimes\CO_\TY(H))$.
It is clear that the composition of this morphism with embedding
$J_{\Sigma}\otimes(\CO_X(1)\boxtimes\CO_\TY(H)) \to
\CO_X(1)\boxtimes\CO_\TY(H)$
coincides with $s_\TQ$.
\end{proof}

Consider the pullback of the second exact sequence of lemma~\ref{ppxu}
to $X\times\TY$.
Consider the cone $G'$ of its first morphism twisted by
$\CO_X(1)\boxtimes\CO_\TY(H)$,
and the cone $G''$ of its third morphism as objects of $\D^b(X\times\TY)$,
so that we have exact triangles
$$
\xymatrix@R=10pt{
E_1\boxtimes\CO_{\TY}(H-3h) \ar[r]^-{\te'} &
E_2\boxtimes\CO_{\TY}(H-2h) \ar[r] &
G' \\
E_1\boxtimes\pi^*\Omega_{\PP(W)} \ar[r]^-{\te''} &
E_2\boxtimes\CO_{\TY}(-h) \ar[r] &
G''
}
$$
Then the pullback of the second exact sequence of lemma~\ref{ppxu}
can be rewritten as exact triangle
\begin{equation}\label{gpgpp}
G'\otimes(\CO_X(-1)\boxtimes\CO_\TY(-H))[1] \to G'' \to \CO_\Sigma.
\end{equation}

It is clear that morphisms $\te'$ and $\te''$ determine a structure
of $\SQ$-representations on families $(\CO_{\TY}(H-3h),\CO_{\TY}(H-2h))$
and $(\pi^*\Omega_{\PP(W)},\CO_{\TY}(-h))$ on $\TY$.

\begin{lemma}\label{fam16}
There exist a unique extension
\begin{equation}\label{tf}
0 \to (\CO_{\TY}(H-3h),\CO_{\TY}(H-2h)) \to
(\TF_1,\TF_2) \to (\pi^*\Omega_{\PP(W)},\CO_{\TY}(-h)) \to 0
\end{equation}
in the category of families of $\SQ$-representations on $\TY$
such that the cone of the induced morphism
$\te:E_1\boxtimes\TF_1 \to E_2\boxtimes\TF_2$ is a coherent sheaf
supported scheme-theoretically at $\TQ$.
\end{lemma}
\begin{proof}
Since $(E_1,E_2)$ is an exceptional pair, $\Hom(E_1,E_2)=W^*$,
and the sheaf of homomorphisms of $\SQ$-representations
from $(\pi^*\Omega_{\PP(W)},\CO_{\TY}(-h))$ to
$(\CO_{\TY}(H-3h),\CO_{\TY}(H-2h))$ equals zero,
we have a bijection between the set of all extensions~(\ref{tf})
in the category of families of $\SQ$-representations on $\TY$
and $\Ext^1(G'',G') = \Hom(G''[-1],G')$. Moreover, it is clear
that the cone of $\te$ is isomorphic to the cone of the corresponding
morphism $\eps:G''[-1] \to G'$. Thus we have to construct a morphism
$\eps:G''[-1] \to G'$, such that the cone of $\eps$ is a sheaf
supported scheme-theoretically at $\TQ$. This can be done as follows.

Let $\ti:\TQ\to X\times\TY$ denote the embedding.
Consider the diagram
$$
\xymatrix{
& \ti_*\ti^*G' \ar[d] \\
\CO_\Sigma[-1] \ar[r] \ar@{-->}[dr] \ar@{..>}[ur]^{\eps'} &
G'\otimes(\CO_X(-1)\boxtimes\CO_\TY(-H))[1] \ar[r] \ar[d]_{s_\TQ} &
G'' \ar[r] \ar@{..>}[dl]^{\eps[1]} &
\CO_\Sigma \\
& G'[1]
}
$$
(the horizontal line is the exact triangle~(\ref{gpgpp})
and the vertical line is the exact triangle obtained by
tensoring the resolution of $\ti_*\CO_\TQ$ with $G'$).
It follows from lemma~\ref{sqsigma} that the dashed arrow is zero.
Hence there exist dotted arrows making the diagram commutative.
Moreover, since $\codim_{X\times\TY}\Sigma = 3$ it follows that
the functor $\Ext^{\le 2}(\CO_\Sigma,-)$ is zero
on any locally free sheaf on $X\times\TY$, hence
$\Hom(\CO_\Sigma,G'[1]) = \Ext^{1}(\CO_\Sigma,G') = 0$,
since $G'$ admits a locally free resolution of length~$2$.
Therefore the dotted arrows are unique.

Further, the octahedron axiom implies that the cone of $\eps$
is isomorphic to the cone of $\eps'$. Thus it remains to check
two things:
that $\eps' = \ti_*\eps''$ for some $\eps''\in\Hom(\CO_\Sigma[-1],\ti^*G')$,
and that $\ti^*G'$ is a coherent sheaf.
But by definition $\eps'$ is obtained from an element of the space
\begin{multline*}
\Hom(\ti_*\CO_\Sigma[-1],G'\otimes(\CO_X(-1)\boxtimes\CO_\TY(-H))[1]) = \\ =
\Hom(\CO_\Sigma[-1],\ti^!(G'\otimes(\CO_X(-1)\boxtimes\CO_\TY(-H))[1])) \cong
\Hom(\CO_\Sigma[-1],\ti^*G'),
\end{multline*}
so the first is true.
For the second we must check that $\Ker\ti^*\te' = 0$.
But the second exact sequence of lemma~\ref{ppxu} shows that
$\Ker\ti^*\te' \cong
\CH^{-3}(\ti^*(\CO_\Sigma\otimes(\CO_X(1)\boxtimes\CO_\TY(H))))$
which is zero because the $\Tor$-dimension of $\ti$ equals $1$.
\end{proof}

Lemma~\ref{fam16} says that the family $(\TF_1,\TF_2)$ on $\TY$
is a family of $\SQ$-representations and that  on $X\times\TY$
we have the following exact sequence:
$$
0 \to E_1\boxtimes \TF_1 \exto{\te} E_2\boxtimes \TF_2 \to \ti_*\TE \to 0,
$$
where $\TE$ is a coherent sheaf on $\TQ$.

\begin{lemma}\label{c_phi}
The family $(\TF_1,\TF_2)$ induces a morphism $\phi:\TY\to \BY_2$
such that $\tg = \bg_2\circ\phi$.
\end{lemma}
\begin{proof}
For any point $y\in\TY$ we choose a trivializations of $F_1$ and $F_2$
in a small neighborhood of~$y$. This induces a morphism
$\tphi$ of this neighborhood into a representation space
of the quiver~$\SQ$. Since $\rank(\TF_1) = 6$, $\rank(\TF_2) = 2$,
and $\gcd(6,2) = 2$ this is the space $\BR_2$.
Let $y'$ be an arbitrary point of the neighborhood.
Let us check that $\tphi(y')\in\BRO2$.
Indeed, this is equivalent to the injectivity of the map $\eta^*\te$,
where $\eta:X\times y' \to X\times\TY$ is the embedding, or equivalently
to $\eta^*\ti_*\TE \in \D^{\ge0}$.
But it was shown in the proof of lemma~\ref{fam16} that
we have an exact triangle
$$
\xymatrix@1{\CO_\Sigma[-1] \ar[r]^-{\eps'} & \ti_*\ti^*G' \ar[r] & \ti_*\TE}.
$$
Thus it suffices to check that
$\eta^*\ti_*\ti^* G' \in \D^{\ge0}$ and $\eta^*\CO_\Sigma \in \D^{\ge 0}$.
Note that the following diagrams are exact cartesian by~\ref{lci}
$$
\xymatrix{
X_{y'} \ar[d]_{\eta} \ar[r]^{i'} & X\times y' \ar[d]^{\eta} \\
Q(X,\TY) \ar[r]^{\ti} & X\times\TY
}
\qquad
\qquad
\xymatrix{
\Sigma\cap X \ar[d] \ar[r] & X\times y' \ar[d]^{\eta} \\
\Sigma \ar[r] & X\times\TY
}
$$
and the top arrows are closed embeddings. The second inclusion follows
immediately and for the first we note that by~(\ref{gpgpp}) we have
$\CH^{-1}(\eta^*\ti_*\ti^* G'\otimes\CO_X(-1)\otimes\CO_\TY(-H)) =
\CH^{-3}(\eta^*\ti_*\ti^* \CO_\Sigma)$
and on the other hand
$\eta^*\ti_*\ti^* \CO_\Sigma =
i'_*\eta^*\ti^* \CO_\Sigma =
i'_*{i'}^*\eta^* \CO_\Sigma =
i'_*{i'}^* \CO_{\Sigma \cap X} \in \D^{\ge -1}$.
%
%
%
%

Further, $\supp\eta^*\ti_*\TE = \supp i'_*\eta^*\TE = X_{y'}$,
is a hyperplane section of $X$, hence $\tphi(y')\in\BRP2$.

Composing $\tphi$ with the factorization map
$\BRP2\to\BY_2 = \BRP2\git G$ we obtain a map $\phi_y$ from
the neighborhood of $y$ in $\TY$ to $\BY_2$. It is clear that
$\phi_y$ doesn't depend on a choice of trivializations,
hence constructed morphisms $\phi_y$ glue into a morphism $\phi:\TY\to\BY_2$.
Finally, since $\ti_*\TE$ is supported on $\TQ$ it follows that
$\tg = \bg_2\circ\phi$.
\end{proof}

Let $\TYO = \TY \setminus\tg^{-1}(Z_1)$ and denote the restriction
of $\tg$ and $\phi$ to $\TYO$ by the same symbols.

\begin{lemma}
The map $\tg:\TYO\to\PP(V^*)\setminus Z_1$ factors through
${X^\vee}\setminus Z_1$, where ${X^\vee} \subset \PP(V^*)$
is the projectively dual quartic hypersurface. Moreover, the general
fiber of $\tg:\TYO \to {X^\vee}\setminus Z_1$ is a conic.
\end{lemma}
\begin{proof}
Let $H_{\tg(y')}\subset\PP(V)$ be the hyperplane corresponding to a point
$\tg(y')\in\PP(V^*)$. It is clear from the proof of lemma~\ref{c_phi}
that for any point $y'\in\TY$ we have
$$
X\cap H_{\tg(y')} =
X_{y'} =
\eta^{-1}(\TQ) \supset \eta^{-1}(\Sigma) = \LGr(2,4),
$$
But a smooth hyperplane section of $X = \LGr(3,6)$ cannot contain
a $\LGr(2,4)$ by Lefschetz theorem. Hence $X\cap H_{\tg(y')}$
is singular for any $y'\in\TY$, hence
$\tg(\TY)\subset {X^\vee}\subset \PP(V^*)$.
Moreover, it is easy to compute the intersection index $H^{12}\cdot h = 8$
on $\TY$, which implies that the map $\tg:\TY\to {X^\vee}$ is surjective
and its generic fiber is a conic in $\PP(W)$, because ${X^\vee}$
is a quartic hypersurface by~\cite{Ho}.
\end{proof}

Recall that we have a sheaf of Azumaya algebras $\CA_2$ on $Y_2$
and that $F_2$ is a locally projective $\CA_2$-module on $Y_2$ of rank~$2$.
Let $\PP_{Y_2,\CA_2}(F_2^*)$ denote the moduli space of locally projective
quotient $\CA_2$-modules of rank~$1$.
Then $\PP_{Y_2,\CA_2}(F_2^*)$ is a conic bundle over $Y_2$.

\begin{proposition}\label{tyoispf2}
We have $\TYO = \PP_{Y_2,\CA_2}(F_2^*)$.
\end{proposition}
\begin{proof}
Let $y \in \TYO$. Then $\phi(y)\in\BY_2$ and
$\bg_2(\phi(y)) = \tg(y) \not\in Z_1$, hence $\phi(y)\in Y_2$.
Therefore $\phi$ takes $\TYO$ to $Y_2$.
Moreover, it follows from the construction of the map $\phi$
in lemma~\ref{c_phi} and from the construction of the family
$(F_1,F_2)$ on $Y_2$ in lemma~\ref{resdia} that
$\phi^*F_2 \cong \TF_0^*\otimes\TF_2$ and
$\phi^*\CA_2 \cong \CEnd(\TF_0)$
for a certain vector bundle $\TF_0$ on $\TYO$.
The epimorphism $\TF_2 \to \CO_\TY(-h)$ from~(\ref{tf})
gives an epimorphism of $\CEnd(\TF_0)$-modules
$\TF_0^*\otimes\TF_2 \to \TF_0^*(-h)$, hence induces a morphism
$\TYO \to \PP_{Y_2,\CA_2}(F_2^*)$.


Vice versa, choose a point $\ty\in\PP_{Y_2,\CA_2}(F_2^*)$.
Let $y \in Y_2$ be its projection and let $\rho:W\otimes R_1 \to R_2$
be the corresponding $(6,2)$-dimensional representation of the quiver $\SQ$,
and $R_2\to\kk$ --- the corresponding one-dimensional quotient.
Let $R_2'\subset R_2$ be its kernel and consider the composition
$E_1\otimes R_1 \exto{e_\rho} E_2\otimes R_2 \to E_2$.
Applying $\Hom(E_1,-)$ we obtain a map $R_1 \to W^*$.
Let $R'_1$ denote the kernel of this map. Then we have a diagram
\begin{equation}\label{ty}
\vcenter{\xymatrix@R=10pt{
&& 0 \ar[d] && 0 \ar[d]
\\
&
0 \ar[r] &
E_1\otimes R'_1 \ar[rr]^{e'_\rho} \ar[d] &&
E_2\otimes R'_2 \ar[r] \ar[d] &
C' \ar[r] & 0
\\
& 0 \ar[r] &
E_1\otimes R_1 \ar[rr]^{e_\rho} \ar[d] &&
E_2\otimes R_2 \ar[d] \ar[r] &
C \ar[r] &
0
\\
0 \ar[r] &
K'' \ar[r] &
E_1\otimes R_1/R'_1 \ar[rr]^{e''_\rho} \ar[d] &&
E_2 \ar[r] \ar[d] &
C'' \ar[r] &
0
\\
&& 0 && 0
}}
\end{equation}
which induces an exact sequence
\begin{equation}\label{cpc}
0 \to K'' \to C' \to C \to C'' \to 0.
\end{equation}
Note that the sheaf $C$ is supported on the hyperplane section
$X_y = H_{g_2(y)}\cap X\subset X$ and
its rank at the generic point of $X_y$ equals
$c_1(E_2\otimes R_2) - c_1(E_1\otimes R_1) = 2$.

We claim that $\dim R'=1$ and we are going to prove it
by a case by case analysis.

Since $(R'_1,R'_2)$ is a subrepresentation of $(R_1,R_2)$
and $(R_1,R_2)$ is stable by lemma~\ref{gmpro} we have
$$
\dim R'_1 < 3,
$$
so we need to reject two cases: $\dim R'_1=2$ and $\dim R'_1 = 0$.

Assume that $\dim R'_1=2$, $\dim R_1/R'_1=4$.
Then the bottom line of the diagram~(\ref{ty}) takes form
$$
0 \to K'' \to (R_1/R'_1)\otimes\CO_X \exto{e''} \CU_X^* \to C'' \to 0.
$$
Since $C''$ is a quotient of $C$ by~(\ref{cpc}) we have $\rank(C'')=0$,
hence $\rank(K'')=4-3=1$. Thus $K''$ is a reflexive sheaf of rank $1$,
hence $K''$ is a line bundle by~\cite{OSS}. Since $\Pic X$ is generated
by $\CO_X(1)$ we have $K''\cong\CO_X(-k)$, $k\ge 0$. Moreover,
$k\ne 0$ since the map $H^0(e'')$ is an embedding (by definition of $R'_1$).
On the other hand,
$-k = c_1(K'') = c_1(C'') - c_1(\CU_X^*) \ge -c_1(\CU_X)^* = -1$,
hence $k=1$ and $c_1(C'') = c_1(K'')+c_1(\CU_X^*) = 0$.
Therefore, the rank of $C''$ at the generic point of $X_y$ is zero.

Further, the top line of the diagram~(\ref{ty}) takes form
$$
0 \to R'_1\otimes\CO_X \exto{e'} \CU_X^* \to C' \to 0.
$$
We have $\rank(C') = \rank(\CU_X^*) - \dim R_1' = 1$,
and $c_1(C') = c_1(\CU_X^*) = 1$. If $T(C')$ denotes the torsion part of $C'$
then $C'/T(C')$ is a torsion free rank 1 sheaf, hence
$C'/T(C') \cong J_S(k)$ (see \cite{OSS}), where $S\subset X$ is a subscheme
of codimension $\ge 2$ in $X$. Therefore $c_1(C') = c_1(T(C')) + k \ge k$,
so $k\le 1$. On the other hand, $k\ge 1$ since we have a nontrivial
morphism $\CU_X^* \to C' \to J_S(k) \to \CO_X(k)$.
Therefore, $k=1$, $c_1(T(C')) = c_1(C')-k = 0$,
hence the torsion part of $C'$ is supported in codimension $\ge 2$
and the rank of $C'$ at the generic point of $X_y$ equals $\rank(C') = 1$.

Now we have shown that the rank of $C''$ at the generic point of $X_y$
equals $0$ and the rank of $C'$ at the generic point of $X_y$ equals $1$.
Looking at~(\ref{cpc}) we see that this contradicts to the fact that
the rank of $C$ at the generic point of $X_y$ equals $2$.
So the case $\dim R'_1=2$ is impossible.

Assume that $R'_1 = 0$. Then the map $R_1\to W^*$ is injective,
hence an isomorphism since $\dim R_1 = 6 = \dim W^*$.
Thus $R_1/R'_1 \cong W^*$ and the map $e''_\rho$ gets identified
with the canonical morphism $W^*\otimes\CO_X \to \CU_X^*$.
Therefore, $C''=0$ and $K'' = \CU_X^\perp \cong \CU_X$.
On the other hand, it is clear that $C' = R_2'\otimes\CU_X^* \cong \CU_X^*$.
Thus the sequence~(\ref{cpc}) takes form
$$
0 \to \CU_X \to \CU_X^* \to C \to 0.
$$
Note that $\Hom(\CU_X,\CU_X^*) \cong S^2W^* \oplus \Lambda^2W^*/\kk$.
If a homomorphism has a nontrivial component in $S^2W^*$ then
its cokernel is supported on a hyperquadric section of $X$,
and if a homomorphism lies in $\Lambda^2W^*/\kk$ then it is
skew-symmetric, hence degenerate (since $\CU_X$ has odd rank).
Thus $C$ cannot be a rank $2$ sheaf on a hyperplane section of $X$.
So the case $\dim R'_1=2$ is impossible as well.

Therefore, $\dim R'_1 = 1$, $\dim R_1/R'_1 = 5$.
The top and the bottom lines of the diagram~(\ref{ty}) take form
\begin{equation}\label{wp}
0 \to \CO_X \exto{e'} \CU_X^* \to C' \to 0,
\end{equation}
$$
0 \to K'' \to \CO_X^{\oplus 5} \exto{e''} \CU_X^* \to C'' \to 0.
$$
The map $e'$ is given by an element $w'\in\PP(\Hom(\CO_X,\CU_X^*)) = \PP(W^*)$
and the map $e''$ is given by a hyperplane in $\PP(W^*)$, that is by a point
$w''\in\PP(W)$. Moreover, it follows from the canonical exact sequence
$0 \to \CU_X \to W^*\otimes\CO_X \to \CU_X^* \to 0$ that we have
the following exact sequence
\begin{equation}\label{wpp}
0 \to K'' \to \CU_X \exto{w''} \CO_X \to C'' \to 0.
\end{equation}
In particular, $K''$ is reflexive and $\supp C'' = \LGr(2,4)$
and has codimension $3$.
Moreover, it follows from~(\ref{wp}) that $C'$ is reflexive and $\rank(C')=2$.
Therefore $C'_{|X_y}$ is torsion free of rank $2$ on $X_y$ and we have
the following exact sequence
$$
0 \to C'(-1) \to C' \to C'_{|X_y} \to 0.
$$
Now, note that since $C$ is supported on $X_y$ the map $C' \to C$
factors through $C'_{|X_y}$. Since the ranks of both $C'_{|X_y}$
and $C$ at generic points of $X_y$ equal to $2$ and the rank of $C''$
at generic point of $X_y$ is zero, it follows that the map
$C'_{|X_y} \to C$ is an embedding. Therefore, the composition
$K'' \to C' \to C'_{|X_y}$ is zero, hence the map $K'' \to C'$
factors through the map $K'' \to C'(-1)$. But this map, being an embedding
of reflexive sheaves of equal rank and $c_1$ is an isomorphism.
An isomorphism $K''\cong C'(-1)$ allows to glue the sequence~(\ref{wp})
twisted by $\CO_X(-1)$ with the sequence~(\ref{wpp}) into
the following exact sequence
$$
0 \to \CO_X(-1) \exto{w'} \CU_X^*(-1) \to \CU_X \exto{w''} \CO_X \to C'' \to 0,
$$
Then it is easy to see that this sequence is the Koszul complex of $w''$,
hence $w' = w''$ and $C''\cong \CO_{Z(w'')}$, where $Z(w'')$ is the zero locus
of  $w''\in H^0(X,\CU_X^*)$.
It is easy to check that the space of hyperplane section of $X$
passing through $Z(w'')$ equals to the fiber of the bundle
$\Omega_{\PP(W)}^3(3)/\CO_{\PP(W)}(-1) \subset V^*\otimes\CO_{\PP(W)}$
at $w''$. On the other hand, the support hyperplane $X_y$ of $C$
contains $\supp(C'') = Z(w'')$, therefore the map
$$
\PP_{Y_2,\CA_2}(F_2^*) \to \PP(W)\times\PP(V^*),\qquad
\ty \mapsto (w'',g_2(y))
$$
factors through
$\TY = \PP_{\PP(W)}(\Omega_{\PP(W)}^3(3)/\CO_{\PP(W)}(-1)) \subset
\PP(W)\times\PP(V^*)$. Moreover, it is clear that its image lies in $\TYO$,
so we have constructed a map $\PP_{Y_2,\CA_2}(F_2^*) \to \TYO$.

It remains to note that the constructed maps
are mutually inverse.
\end{proof}

\begin{theorem}\label{th_a}
We have $Y_2 = {X^\vee} \setminus Z_1$.
\end{theorem}
\begin{proof}
Since $\TYO$ is smooth it follows from proposition~\ref{tyoispf2} that
$Y_2$ is smooth and moreover $\dim Y_2 = \dim\TYO - 1 = 12$.
On the other hand, since $\tg = \bg_2\circ\phi$, $\tg$ and $\phi$
are proper and dominant, and $\tg$ takes $\TYO$ to ${X^\vee} \setminus Z_1$,
it follows that $g_2$ is a proper dominant map
$Y_2 \to {X^\vee} \setminus Z_1$. On the other hand,
by~\cite{Ho} we have $Z_1 = \sing({X^\vee})$, hence
${X^\vee}\setminus Z_1$ is also smooth and
$\dim({X^\vee}\setminus Z_1) = \dim \PP(V^*) - 1 = 12$.
So, $g_2:Y_2\to{X^\vee}\setminus Z_1$ is a proper dominant map
of smooth varieties of equal dimension.
Therefore, either $Y_2 = {X^\vee} \setminus Z_1$, or
$\rank\Pic Y_2 = \rank \Pic {X^\vee} + 1 = 2$.
Assume that the latter is true, and note that since $\phi$ is
a conic bundle, we have $\rank\Pic\TYO = \rank\Pic Y_2 + 1 = 3$.
But on the other hand, it is easy to see that
$\rank\Pic\TYO \le \rank\Pic\TY = 2$.
Thus we get a contradiction, hence $Y_2 = {X^\vee} \setminus Z_1$.
\end{proof}

\begin{corollary}
If the data {\rm(D.3)--(D.5)} is given by
$Y = Y_2$, $\CA_Y = \CA_2$, $\BZ = \OZ_2 = Z_1$, $g = g_2$,
with $(F_1,F_2)$ being the universal family, and
with $\phi$ being induced by the $\SQ$-representation structure on $(F_1,F_2)$,
then all conditions {\rm(C.1)--(C.8)} are satisfied.
\end{corollary}
\begin{proof}
Since conditions~(\ref{addc}) and (\ref{codimbro}) are
satisfied by lemma~\ref{conds_a}, it follows from theorem~\ref{qu_th}
that conditions (C.4)--(C.6) are true. On the other hand,
(C.8) is true because $\BZ = \mu(\Gr(3,W^*))$ by lemma~\ref{yz1_a}
and $\dim\BZ\cap H = \dim\Gr(3,W^*) - 1 = 8$, while
$N - i - 2 = 14 - 4 - 2 = 8$.
Finally, by theorem~\ref{th_a} the map $Y = Y_2 \to \PP(V^*)$ is an embedding
and $\codim_{Y\times Y}(Y\times_{\PP(V^*)}Y) = \dim Y = \dim{X^\vee} =
\dim\PP(V^*) - 1 = 12$, hence
$\dim X + N - \codim_{Y\times Y}(Y\times_{\PP(V^*)}Y) = 6 + 14 - 12 = 8 = 2i$,
and (C.7) is true.
\end{proof}

\section*{Appendix B. $G_2$ Grassmannian}

\refstepcounter{section}

In this appendix we show that a double covering of $\PP^{13}$
ramified in a sextic hypersurface is homologically projectively dual
to the Grassmannian $\GTGr(2,7)$ of the simple Lie group $G_2$.

Recall the notation. Let $W = \kk^{7}$ be an irreducible representation
of the simple Lie group $G_2$ and
$$
X = \GTGr(2,W),\qquad
(E_1,E_2) = (\CO_X,\CU_X^*),
$$
where $\GTGr(2,W)$ is the Grassmannian of the Lie group $G_2$,
realized as the zero locus of the section
$s_\lambda \in H^0(\Gr(2,W),\CU^\perp(1)) \cong \Lambda^3W^*$,
corresponding to the unique $G_2$-invariant $3$-form $\lambda$ on~$W$,
$\CU \subset W\otimes\CO_{\Gr(2,W)}$ is the tautological rank 2 subbundle,
and $\CU_X$ is the restriction of $\CU$ to~$X$.
The natural representation of the group $G_2$ in the space
$\Lambda^2 W$ decomposes into the direct sum of representations
$$
\Lambda^2 W = W^* \oplus V
$$
(the projection $\Lambda^2W \to W^*$ is given by the 3-form $\lambda$).
The Pl\"ucker embedding $\Gr(2,W) \subset \PP(\Lambda^2 W)$ restricts
to an embedding $f:X\to\PP(V)$.

\begin{lemma}\label{conds_b}
We have $\dim X = 5$, $\omega_X \cong \CO_X(-3)$.
Moreover, conditions {\rm(C.1)--(C.3)}\/ as well as the additional
conditions {\rm(\ref{addc})} and {\rm(\ref{codimbro})} are satisfied for $X$.
\end{lemma}
\begin{proof}
Note that $X\subset\Gr(2,W)$ is the zero locus of a regular section
$s_\lambda\in H^0(\Gr(2,W),\CU^\perp(1))$. Hence
$$
\begin{array}{l}
\dim X = \dim\Gr(2,W) - \rank(\CU^\perp(1)) = 10 - 5 = 5,\\
\omega_X \cong \omega_{\Gr(2,W)|X} \otimes \det(\CU^\perp(1))^* \cong
\CO_X(-7) \otimes \CO_X(4) \cong \CO_X(-3).
\end{array}
$$
In particular, (C.2) is true and $i_X = 3$.
Further, using the Koszul resolution
$\Lambda^\bullet(\CU^\perp(1)) \cong \CO_X$ of $X$ in $\Gr(2,W)$
it is easy to check that $H^0(X,\CO_X(1)) \cong V^*$, hence
$f(X)$ is not contained in a hyperplane, i.e.\ (C.1) is true.
Using again the Koszul resolution we compute $\Hom(E_s(k),E_t(l))$
via Borel--Bott--Weil theorem and check that (C.3) and
(\ref{addc}) are true. Finally, we note that
$\BR_1 = \Hom(\kk^2,W^*)$ and $\BRO1\subset\BR_1$
corresponds to the subset of embeddings $\kk^2\subset W^*$.
Therefore, $\codim_{\BR_1}(\BR_1\setminus\BRO1) = 6$
and $(\ref{codimbro})$ is satisfied by lemma~\ref{brombro1}.
\end{proof}

Now we are going to apply the construction of section~5 and
to describe the space $Y = Y_3$ explicitly. We use the same
approach as in Appendix~A. The arguments in most cases are
the same. So we will not repeat them, but describe the differences.

Instead of $Y_3$ we will describe $\PP_{Y_3,\CA_3}(F_2^*)$,
the moduli space of quotient $\CA_3$-modules of rank $1$.
The fiber of $\PP_{Y_3,\CA_3}(F_2^*)$ over a point of $Y_3$
corresponding to a $(6,3)$-dimensional representation
$\rho = (R_1,R_2)$ of $\SQ$ coincides with $\PP(R^*_2)$.
We show that any plane $R'_2\subset R_2$ can be extended in a unique way
to a $(1,2)$-dimensional subrepresentation $\rho' = (R'_1,R'_2)$.
The moduli space of $(1,2)$-dimensional representations of $\SQ$
is $\Gr(2,W^*)$. We claim that defined in this way map
$\PP_{Y_3,\CA_3}(F_2^*) \to \Gr(2,W^*)$ identifies $\PP_{Y_3,\CA_3}(F_2^*)$
with an open subset in a $\PP^8$-bundle over $\OGr(2,W^*) \subset \Gr(2,W^*)$,
the isotropic Grassmannian.

\begin{lemma}\label{yz1_b}
We have $\BY_1 = \OZ_1 = \mu(\Gr(2,W^*))$, where
$\mu:\xymatrix@1{\PP(\Lambda^2W^*) \ar@{-->}[r] & \PP(V^*)}$
is the linear projection from $\PP(\lambda(W)) \subset \PP(\Lambda^2W^*)$.
\end{lemma}

Recall that the group $G_2$ is a subgroup of the group $\SO(W)$.
Consider the corresponding $G_2$-invariant nondegenerate
quadratic form on $W$ and let $\CQ\subset\PP(W)$,
$\OGr(2,W) \subset \Gr(2,W)$ and $\OF(1,2;W) \subset \Fl(1,2;W)$
denote the corresponding quadric, the isotropic Grassmannian and
the isotropic partial flag variety. Let $\tp:\OF(1,2;W) \to \OGr(2,W)$ and
$\tq:\OF(1,2;W) \to \CQ$ denote the projections. Let $S_\CQ$ denote
the spinor bundle on $\CQ$, let $\CU_\OGr$ denote the tautological
rank $2$ subbundle in $W\otimes\CO_{\OGr(2,W)}$ and put
$M := \tp_*\tq^* S_\CQ \in \D(\OGr(2,W))$.

\begin{lemma}\label{g2gr}
We have the following commutative diagram
$$
\xymatrix{
&& \PP_X(\CU_X) \ar[dl]_p \ar[d] \ar[dr]^q \\
&
X \ar[dl] \ar[d] &
\OF(1,2;W) \ar[dl]_\tp \ar[dr]^\tq &
\CQ \ar@{=}[d] \ar[dr] \\
\Gr(2,W) &
\OGr(2,W) \ar[l] &&
\CQ \ar[r] &
\PP(W)
}
$$
Moreover, $M$ is a rank $2$ vector bundle on $\OGr(2,W)$ and we have
an isomorphism $M_{|X} \cong \CU_X$ and an exact sequence
$$
0 \to \CO_{\OGr(2,W)}(-1) \to M \to \CO_{\OGr(2,W)} \to \CO_X \to 0.
$$
\end{lemma}
\begin{proof}
Note that $\PP_X(\CU_X)$ is the flag variety of the group $G_2$.
The embedding $G_2 \subset \SO(W)$ induces an embedding
of the flag varieties $\PP_X(\CU_X) \to \OF(1,2;W)$,
and of the Grassmannians $X \to \OGr(2,W)$ such that
the diagram is commutative (note that the quadric $\CQ$
is also a Grassmannian of the group $G_2$).

The second claim follows from the fact that the restriction
of the spinor bundle $S_\CQ$ to any line $\PP^1\subset\CQ$ splits as
$\CO_{\PP^1}\oplus\CO_{\PP^1}\oplus\CO_{\PP^1}(-1)\oplus\CO_{\PP^1}(-1)$.

To check the third claim we note that $\PP_X(\CU_X) \cong \PP_\CQ(K)$,
where $K$ is the kernel of the restriction to $\CQ$ of the morphism
$T_{\PP(W)}(-2) \to \Omega_{\PP(W)}(1)$ given by the form $\lambda$.
On the other hand, the image of this morphism is just the spinor bundle,
so we have an exact sequence
\begin{equation}\label{kts}
0 \to K \to T_{\PP(W)}(-2)_{|\CQ} \to S_\CQ \to 0.
\end{equation}
Taking into account the sequences
$$
\begin{array}{l}
0 \to p^*\CO_X(-1) \to q^*K \to p^*\CO_X(1)\otimes q^*\CO_\CQ(-3) \to 0,\\
0 \to q^*\CO_\CQ(-1) \to p^*\CU_X \to q^*\CO_\CQ(1)\otimes p^*\CO_X(-1) \to 0,
\end{array}
$$
we see that $\CH^0(p_*q^*K) \cong \CO_X(-1)$, $\CH^1(p_*q^*K) \cong \CU_X$
and $p_*q^*T_{\PP(W)}(-2)_{|\CQ} \cong \CO_X(-1)$, whereof we deduce
$p_*q^*S_\CQ \cong \CU_X$. Thus
$M_{|X} \cong (\tp_*\tq^*S_\CQ)_{|X} \cong
p_*((\tq^*S_\CQ)_{|\PP_X(\CU_X)}) \cong p_*q^*S_\CQ \cong \CU_X$.

Finally, we note that (\ref{kts}) implies that
$\PP_X(\CU_X) = \PP_\CQ(K)$ is the zero locus of a section
of the bundle $\tq^*S_\CQ\otimes\tp^*\CO_{\OGr(2,W)}(1)$ on
$\PP_\CQ(T_{\PP(W)}(-2)_{|\CQ})\subset\PP_{\PP(W)}(T_{\PP(W)}(-2))=\Fl(1,2;W)$.
Restricting to $\OF(1,2;W) = \PP_{\OGr(2,W)}(\CU_{\OGr})$ we see that
$\PP_X(\CU_X)$ is the zero locus of a section
of the bundle $\tq^*S_\CQ\otimes\tp^*\CO_{\OGr(2,W)}(1)$,
hence $X$ is the zero locus of a section of the vector bundle
$\tp_*(\tq^*S_\CQ\otimes\tp^*\CO_{\OGr(2,W)}(1)) =
M\otimes\CO_{\OGr(2,W)}(1) \cong M^*$.
Therefore, its structure sheaf admits a Koszul resolution
which takes form of the exact sequence of the lemma.
\end{proof}

Consider $X\times\OGr(2,W)$ and the subscheme
$\Sigma'\subset X\times \OGr(2,W) \subset \Gr(2,W) \times \Gr(2,W)$
of pairs of intersecting two-dimensional subspaces.

\begin{lemma}\label{sigmab}
The following sequences on $X\times\OGr(2,W)$ are exact
{\small
$$
0 \to
\CU_X(-1)\boxtimes\CO_{\OGr(2,W)}(-1) \to
\CO_X(-1)\boxtimes(W/\CU_{\OGr})(-1) \to
\CU_X\boxtimes M \to
\CO_X\boxtimes\CO_{\OGr(2,W)} \to
\CO_{\Sigma'} \to 0.
$$
$$
0 \to
\CO_X(-1)\boxtimes\CO_{\OGr(2,W)}(-1) \to
\CU^*_X(-1)\boxtimes M \to
\CO_X\boxtimes \CU_{\OGr}^\perp \to
\CU_X^*\boxtimes\CO_{\OGr(2,W)} \to
\CF_{\Sigma'} \to 0.
$$
}
where
$\CF_{\Sigma'} = \CExt^3(\CO_{\Sigma'},\CO_X(-1)\boxtimes\CO_{\OGr(2,W)}(-1))$.
\end{lemma}
\begin{proof}
Consider the decomposition with respect to the exceptional collection
$$
\langle \CO_\CQ(-3),\CO_\CQ(-2),\CS_\CQ(-1),
\CO_\CQ(-1),\CO_\CQ,\CO_\CQ(1)\rangle = \D^b(\CQ),
$$
of the structure sheaf of the variety
$\OF(1,2;W) \subset \OGr(2,W)\times\CQ$.
It is easy to check that the decomposition takes form
\begin{multline}\label{ofres}
0 \to
\CO_{\OGr(2,W)}(-1)\boxtimes \CO_\CQ(-3) \to
(W/\CU_{\OGr})(-1)\boxtimes \CO_\CQ(-2) \to
M\boxtimes \CS_\CQ(-1) \to
\\ \to
\CU_{\OGr}^\perp\boxtimes \CO_\CQ(-1) \to
\CO_{\OGr(2,W)} \boxtimes \CO_\CQ \to
\CO_{\OF(1,2;W)} \to 0.
\end{multline}
Now consider the product $\OGr(2,W)\times\CQ\times\OGr(2,W)$
and the subvariety
$$
\Sigma'' = \OF(1,2;W)\times_\CQ\OF(1,2;W) \subset
\OGr(2,W)\times\CQ\times\OGr(2,W).
$$
It is clear that
$\Sigma''' = p_{13}(\Sigma'') \subset \OGr(2,W)\times\OGr(2,W)$
is just the subscheme of pairs of intersecting two-dimensional subspaces.
Tensoring the pullbacks of the resolutions~(\ref{ofres})
via the projections
$\OGr(2,W)\times\CQ\times\OGr(2,W) \to \OGr(2,W)\times\CQ$ and
$\OGr(2,W)\times\CQ\times\OGr(2,W) \to \CQ\times\OGr(2,W)$,
and applying ${p_{13}}_*$ we obtain
\begin{multline*}
0 \to
\CU_{\OGr}(-1)\boxtimes\CO_{\OGr(2,W)}(-1) \to
\CO_{\OGr(2,W)}(-1)\boxtimes(W/\CU_{\OGr})(-1) \to \\ \to
M\boxtimes M \to
\CO_{\OGr(2,W)}\boxtimes\CO_{\OGr(2,W)} \to
\CO_{\Sigma'''} \to 0.
\end{multline*}
Restricting to $X\times\OGr(2,W)$ and taking into account
that $\CU_{\OGr|X}\cong M_{|X}\cong\CU_X$, we obtain
the first sequence. Applying the functor
$\RCHom(-,\CO_X(-1)\boxtimes\CO_{\OGr(2,W)}(-1))$
we obtain the second sequence.
\end{proof}

Consider the truncated first sequence of lemma~\ref{sigmab}
as a resolution of the sheaf of ideals $J_{\Sigma'}$.
Tensoring it with $\CO_X(1)$ and pushing forward
via the projection $p_2:X\times\OGr(2,W) \to \OGr(2,W)$
we obtain exact sequence
\begin{equation}\label{p2jb}
0 \to (W/\CU_\OGr)(-1) \to W^*\otimes M \to
{p_2}_*(J_{\Sigma'}\otimes\CO_X(1)) \to 0.
\end{equation}

\begin{lemma}\label{fpbundle}
${p_2}_*(J_{\Sigma'}\otimes\CO_X(1))$ is a vector bundle
of rank $9$ on $\OGr(2,W)$. Moreover, its restriction to $X$
is isomorphic to the bundle $\fp$ on $X = G_2/P$
{\rm(}here $P$ is a parabolic subgroup of $G_2$
and $\fp$ is the corresponding parabolic Lie subalgebra{\rm)}.
\end{lemma}
\begin{proof}
For the first claim it suffices to check that
${p_2}_*(\CO_{\Sigma'}\otimes\CO_X(1))$ is a vector bundle.
But it is easy to see that the fiber of $\Sigma'$ over a point $u\in\OGr(2,W)$
is either a Hirzebruch surface $F_1$ (if $u\not\in X$),
or a cone over a rational twisted cubic curve (if $u\in X$).
In both cases $H^0(\Sigma'_u,\CO_X(1))$ is 5-dimensional,
so ${p_2}_*(\CO_{\Sigma'}\otimes\CO_X(1))$ is a vector bundle of rank $5$
and the first claim follows. Moreover, restricting this bundle to $X$
we obtain a $G_2$-equivariant rank $5$ quotient bundle of the bundle
$V^*\otimes\CO_X \cong \fg\otimes\CO_X$. But the only such bundle
is the tangent bundle $T_X \cong \fg/\fp$, hence
${p_2}_*(J_{\Sigma'}\otimes\CO_X(1))_{|X} \cong \fp$.
\end{proof}

Consider the projectivization
$$
\TY = \PP_{\OGr(2,W)}({p_2}_*(J_{\Sigma'}\otimes\CO_X(1))).
$$
The natural embedding
${p_2}_*(J_{\Sigma'}\otimes\CO_X(1)) \to
H^0(X,\CO_X(1))\otimes\CO_{\OGr(2,W)} =
V^*\otimes\CO_{\OGr(2,W)}$
induces a morphism $\TY \exto{\tg} \PP(V^*)$.
Let $\pi:\TY\to\OGr(2,W)$ denote the canonical projection.
So we have
$$
\CO_\TY(H) := \tg^*\CO_{\PP(V^*)}(1),\quad\!\!
\tg:\TY \to \PP(V^*),\quad\text{and}\quad
\CO_{\TY}(h) := \pi^*\CO_{\OGr(2,W)}(1),\quad\!\!
\pi:\TY\to\OGr(2,W).
$$
%
%
Consider the product $X\times\TY$, the projection
$\id_X\times\pi:X\times\TY \to X\times\OGr(2,W)$,
the subschemes
$$
\Sigma = (\id_X\times\pi)^{-1}(\Sigma'),\qquad
\TQ = Q(X,\TY) = (X\times\TY)\times_{\PP(V)\times\PP(V^*)}Q,
$$
and the equation
$s_\TQ\in
H^0(X\times\TY,\CO_X(1)\boxtimes\CO_\TY(H)) =
H^0(X\times\TY,(f\times\tg)^*\CO_{\PP(V)\times\PP(V^*)}(1,1)) \cong
V^*\otimes V$
of $\TQ$ in $X\times\TY$.

\begin{lemma}\label{sqsigmb}
The section $s_\TQ$ vanishes on the scheme $\Sigma$.
\end{lemma}

Consider the pullback of the second exact sequence of lemma~\ref{sigmab}
to $X\times\TY$.
Consider the cone $G'$ of its first morphism twisted by
$\CO_X(1)\boxtimes\CO_\TY(H)$,
and the cone $G''$ of its third morphism as objects of $\D^b(X\times\TY)$,
so that we have exact triangles
$$
\xymatrix@R=10pt{
E_1\boxtimes\CO_{\TY}(H-h) \ar[r]^-{\te'} &
E_2\boxtimes M(H) \ar[r] &
G' \\
E_1\boxtimes \CU_\OGr^\perp \ar[r]^-{\te''} &
E_2\boxtimes\CO_{\TY} \ar[r] &
G''
}
$$
Then the pullback of the second exact sequence of lemma~\ref{sigmab}
can be rewritten as exact triangle
\begin{equation}\label{gpgppb}
G'\otimes(\CO_X(-1)\boxtimes\CO_\TY(-H))[1] \to G'' \to \CF_\Sigma.
\end{equation}

It is clear that morphisms $\te'$ and $\te''$ determine a structure
of $\SQ$-representations on families $(\CO_{\TY}(H-h),M(H))$
and $(\CU_\OGr^\perp,\CO_{\TY})$ on $\TY$.

\begin{lemma}\label{fam16b}
There exist a unique extension
\begin{equation}\label{tfb}
0 \to (\CO_{\TY}(H-h),M(H)) \to
(\TF_1,\TF_2) \to (\CU_\OGr^\perp,\CO_{\TY}) \to 0
\end{equation}
in the category of families of $\SQ$-representations on $\TY$
such that the cone of the induced morphism
$\te:E_1\boxtimes\TF_1 \to E_2\boxtimes\TF_2$ is a sheaf
supported at $\TQ$.
\end{lemma}

Lemma~\ref{fam16b} says that the family $(\TF_1,\TF_2)$ on $\TY$
is a family of $\SQ$-representations and that  on $X\times\TY$
we have the following exact sequence:
$$
0 \to E_1\boxtimes \TF_1 \exto{\te} E_2\boxtimes \TF_2 \to \ti_*\TE \to 0,
$$
where $\TE$ is a coherent sheaf on $\TQ$.

\begin{lemma}\label{c_phib}
The family $(\TF_1,\TF_2)$ induces a morphism $\phi:\TY\to \BY_3$
such that $\tg = \bg_3\circ\phi$.
\end{lemma}

Let $\TYO = \TY \setminus\tg^{-1}(Z_1)$ and denote the restriction
of $\tg$ and $\phi$ to $\TYO$ by the same symbols.

Recall that we have a sheaf of Azumaya algebras $\CA_3$ on $Y_3$
and that $F_2$ is a locally projective $\CA_3$-module on $Y_3$ of rank~$3$.
Let $\PP_{Y_3,\CA_3}(F_2^*)$ denote the moduli space of locally projective
quotient $\CA_3$-modules of rank~$1$. Then $\PP_{Y_2,\CA_2}(F_2^*)$ is
a twisted $\PP^2$-bundle over $Y_3$.

\begin{proposition}\label{tyoispf2b}
We have $\TYO = \PP_{Y_3,\CA_3}(F_2^*)$.
\end{proposition}
\begin{proof}
The map $\TYO \to \PP_{Y_3,\CA_3}(F_2^*)$ can be constructed
like in proposition~\ref{tyoispf2}. Thus it remains to
construct the inverse map.
%
Choose a point $\ty\in\PP_{Y_3,\CA_3}(F_2^*)$.
Let $\rho:W\otimes R_1 \to R_2$ be the corresponding $(6,3)$-dimensional
representation of the quiver $\SQ$, and let $R_2\to\kk$ be the corresponding
one-dimensional quotient. Let $R_2'\subset R_2$ be its kernel
and consider the composition
$E_1\otimes R_1 \exto{e_\rho} E_2\otimes R_2 \to E_2$.
Taking $\Hom(E_1,-)$ we obtain a map $R_1 \to W^*$.
Let $R'_1$ denote the kernel of this map. Then we have
a diagram~(\ref{ty}) and an exact sequence~(\ref{cpc})
Note that the sheaf $C$ is supported on the hyperplane section
$X_y = H_{g_3(y)}\cap X\subset X$ and
its rank at the generic point of $X_y$ equals
$c_1(E_2\otimes R_2) - c_1(E_1\otimes R_1) = 3$.

We claim that $\dim R'=1$ and we are going to prove it
by a case by case analysis.

Since $(R'_1,R'_2)$ is a subrepresentation of $(R_1,R_2)$
and $(R_1,R_2)$ is stable by lemma~\ref{gmpro} we have
$$
\dim R'_1 < 4,
$$
so we need to reject three cases: $\dim R'_1=3$, $\dim R'_1=2$
and $\dim R'_1 = 0$.

The cases $\dim R'_1 = 2$ and $\dim R'_1=3$ are treated along the same lines
as the case $\dim R'_1=2$ in the proof of proposition~\ref{tyoispf2}.
The only difference is that in the case $\dim R'_1 = 2$ the sheaf
$T(C')$ can have rank 1 at the generic point of $X_y$,
but it is easy to see that this is possible only in the case
$C' \cong \CU_X^* \oplus \CU_X^*/(\CO_X\oplus\CO_X)$,
which means that $g_3(y)\in Z_1$.

Assume that $R'_1 = 0$. Then the map $R_1\to W^*$ is injective,
so we can consider $R_1$ as a hyperplane in $\PP(W^*)$, or equivalently
as a point $w_y\in\PP(W)$.
If the cokernel of the morphism $K'' \to C' = \CU_X^*\oplus\CU_X^*$
is supported scheme-theoretically on a hyperplane section $X_y\subset X$,
then the morphism
$\CU_X \oplus \CU_X \cong
(\CU_X^*\oplus\CU_X^*)\otimes\CO_X(-1) \to
\CU_X^*\oplus\CU_X^*$
factors through $K''$ and the cokernel of the
corresponding morphism $\CU_X \oplus \CU_X \to K''$ is supported
on the same hyperplane section $X_y\subset X$. On the other hand,
it is easy to compute that $\dim\Hom(\CU_X,K'') = 1$ for $w_y\not\in\CQ$,
and $\dim\Hom(\CU_X,K'') = 4$ for $w_y\in\CQ$. Therefore we must have
$w_y\in\CQ$. Moreover, it is clear that the dimension of the set of such
hyperplane sections $X_y\subset X$ has dimension less than
$\dim\CQ + \dim\Gr(2,4) = 5 + 4 = 9$. On the other hand,
the set of such hyperplane sections is evidently $G_2$-invariant
in $\PP(V^*)$, and it follows from the description of all $G_2$-invariant
closed subsets of $\PP(V^*)$ that this set is contained in $Z_1$.

Thus we have shown that if $g_3(y)\not\in Z_1$ then
$\dim R'_1 = 1$, $\dim R_1/R'_1 = 5$.
The top and the bottom lines of the diagram~(\ref{ty}) take form
\begin{equation}\label{wpb}
0 \to \CO_X \exto{e'} \CU_X^*\oplus\CU_X^* \to C' \to 0,
\end{equation}
\begin{equation}\label{wppb}
0 \to K'' \to \CO_X^{\oplus 5} \exto{e''} \CU_X^* \to C'' \to 0.
\end{equation}
The map $e''$ is given by a 5-dimensional subspace of $W^*$, that is by a
point $u \in \Gr(2,W)$. The same arguments as in the proof
of proposition~\ref{tyoispf2} show that we can glue the sequence~(\ref{wpb})
twisted by $\CO_X(-1)$ with the sequence~(\ref{wppb}) into
the following exact sequence
$$
0 \to \CO_X(-1) \exto{e'} \CU_X^*(-1)\oplus\CU_X^*(-1) \to
\CO_X^{\oplus 5} \exto{u} \CU_X^* \to C'' \to 0.
$$
Then it is easy to see that this sequence must coincide
with the second sequence of lemma~\ref{sigmab} restricted
to a fiber of $X\times\OGr(2,W)$ over a point of $\OGr(2,W)$.
In particular, $u\in\OGr(2,W)$. Associating to the point $\ty$
the point $u\in\OGr(2,W)$ we obtain a map
$\PP_{Y_3,\CA_3}(F_2^*) \to \OGr(2,W)$.
On the other hand, the support hyperplane $X_y$ of $C$
contains $\supp(C'') = \Sigma'_u$, therefore the map
$$
\PP_{Y_3,\CA_3}(F_2^*) \to \OGr(2,W)\times\PP(V^*),\qquad
\ty \mapsto (u,g_3(y))
$$
factors through
$\TY = \PP_{\OGr(2,W)}({p_2}_*(J_{\Sigma'}\otimes\CO_X(1))) \subset
\OGr(2,W)\times\PP(V^*)$. Moreover, it is clear that its image lies in $\TYO$,
so we have constructed a map $\PP_{Y_3,\CA_3}(F_2^*) \to \TYO$.

It remains to note that the constructed maps
are mutually inverse.
\end{proof}

\begin{theorem}\label{th_b}
The map $g_3:Y_3\to\PP(V^*)\setminus Z_1$ is a double covering.
\end{theorem}
\begin{proof}
Since $\TYO$ is smooth it follows from proposition~\ref{tyoispf2b} that
$Y_3$ is smooth and moreover $\dim Y_3 = \dim\TYO - 2 = 13$.
On the other hand, since $\tg = \bg_3\circ\phi$, $\tg$ and $\phi$
are proper and dominant, it follows that $g_3$ is a proper dominant map
$Y_3 \to \PP(V^*) \setminus Z_1$.
Since $Y_3$ and $\PP(V^*)\setminus Z_1$ are smooth varieties
of equal dimension and
$\rank\Pic Y_3 = \rank\Pic\TYO - 1 = 1 = \rank\Pic(\PP(V^*)\setminus Z_1)$
it follows that $g_3$ is a finite covering. It remains to show that
$\deg g_3 = 2$. This can be done by a direct calculation using
the intersection theory of $\OGr(2,W)$.
Alternatively, we can argue as follows.

Consider the fiber product $\TYO\times_{\OGr(2,W)}X$.
By lemma~\ref{fpbundle} we have $\TYO\times_{\OGr(2,W)}X$
is an open subset of $\PP_X(\fp)$, so and it is clear that the map
$\tg_{|\PP_X(\fp)}:\PP_X(\fp) \to \PP(V^*) = \PP(\fg)$
is the projectivization of the partial Springer--Grothendieck map
$\Tot(\fp) \to \fg$. Therefore, by~\cite{Hu} its degree equals $12/2 = 6$.
On the other hand, since $X\subset\OGr(2,W)$ is the zero locus of
a section of the vector bundle $M^*$ and the restriction of $M$
to the fiber $\PP^2$ of $\TYO$ over $Y_3$ is isomorphic to $\Omega_{\PP^2}$
by~(\ref{tfb}) we see that the degree of the map
$\phi_{|\PP_X(\fp)}:\PP_X(\fp) \to Y_3$ equals $c_2(T_{\PP^2}) = 3$.
Hence $\deg g_3 = 6/3 = 2$.
\end{proof}

\begin{corollary}
If the data {\rm(D.3)--(D.5)} is given by
$Y = Y_3$, $\CA_Y = \CA_3$, $\BZ = \OZ_3 = Z_1$, $g = g_3$,
with $(F_1,F_2)$ being the universal family, and
with $\phi$ being induced by the $\SQ$-representation structure on $(F_1,F_2)$,
then all conditions {\rm(C.1)--(C.8)} are satisfied.
Moreover, the branching locus of the map $g:Y \to \PP(V^*)\setminus\BZ$
is the hypersurface ${X^\vee}\setminus\BZ$, where
${X^\vee} \subset \PP(V^*)$ is the projectively dual sextic hypersurface.
\end{corollary}
\begin{proof}
Since the conditions~(C.1)--(C.3), (\ref{addc}) and (\ref{codimbro}) are
satisfied by lemma~\ref{conds_b}, it follows from theorem~\ref{qu_th}
that the conditions (C.4)--(C.6) are true. On the other hand,
(C.8) is true because $\BZ = \mu(\Gr(2,W))$ by lemma~\ref{yz1_b}
and $\dim\BZ\cap H = \dim\Gr(2,W) - 1 = 9$, while
$N - i - 2 = 14 - 3 - 2 = 9$.
Finally, by theorem~\ref{th_b} the map $Y = Y_3 \to \PP(V^*)$ is a double
covering and
$\codim_{Y\times Y}(Y\times_{\PP(V^*)}Y) = \dim Y = \dim\PP(V^*) = 13$,
hence
$\dim X + N - \codim_{Y\times Y}(Y\times_{\PP(V^*)}Y) = 5 + 14 - 13 = 6 = 2i$,
and (C.7) is true.

Finally, we note that the branching locus of $g$
is nothing but the singular locus of $g$, hence
it coincides with ${X^\vee}\setminus\BZ$ by corollary~\ref{singg}.
\end{proof}

\section*{Appendix C. Intersection of quadrics in $\PP^5$}

\refstepcounter{section}

In this appendix we show that a double covering of $\PP^{20}$
ramified in a sextic hypersurface is homologically projectively dual
to the double Veronese subvariety $\PP^5 \subset \PP^{20}$.

Recall the notation. Let $W = \kk^{6}$ and
$$
X = \PP(W) = \PP^5,\qquad
(E_1,E_2) = (\CO_X,\CO_X(1)).
$$
Let $V = S^2W = \kk^{21}$ and let $f:X \to \PP(V)$ be
the double Veronese embedding.

\begin{lemma}\label{conds_c}
We have $\dim X = 5$, $\omega_X \cong \CO_X(-6) \cong f^*\CO_{\PP(V)}(-3)$.
Moreover, conditions {\rm(C.1)--(C.3)}\/ as well as the additional
conditions {\rm(\ref{addc})} and {\rm(\ref{codimbro})} are satisfied for $X$.
\end{lemma}

Now we are going to apply the construction of section~5 and
to describe the space $Y = Y_4$ explicitly. We use the same
approach as in Appendix~A. The arguments in most cases are
the same. So we will not repeat them, but describe the differences.

Instead of $Y_4$ we will describe $\PP_{Y_4,\CA_4}(F_2^*)$,
the moduli space of quotient $\CA_4$-modules of rank $1$.
The fiber of $\PP_{Y_4,\CA_4}(F_2^*)$ over a point of $Y_4$
corresponding to a $(4,4)$-dimensional representation $\rho = (R_1,R_2)$
of the quiver $\SQ$ coincides with $\PP(R^*_2)$.
We show that any hyperplane $R'_2\subset R_2$ can be extended in a unique way
to a $(1,3)$-dimensional subrepresentation $\rho' = (R'_1,R'_2)$.
The moduli space of $(1,3)$-dimensional representations of $\SQ$
is $\Gr(3,W)$. We claim that defined in this way map
$\PP_{Y_4,\CA_4}(F_2^*) \to \Gr(3,W)$ identifies $\PP_{Y_4,\CA_4}(F_2^*)$
with an open subset in a $\PP^{14}$-bundle over $\Gr(3,W^*)$.

\begin{lemma}\label{yz1_c}
We have $Z_2 \subset \PP(V^*)$ is the locus of rank~$4$ quadrics
and $\cup_{0 < D < \CH} Z^{|D|} \subset \PP(V^*)$
is the locus of rank~$2$ quadrics.
\end{lemma}

Consider the product $X\times\Gr(3,W)$ and the subscheme
$\Sigma' = \Fl(1,3;W)\subset X\times \Gr(3,W)$.

\begin{lemma}\label{sigmac}
The following sequence on $X\times\Gr(3,W)$ is exact
{\small
$$
0 \to
\CO_X(-3)\boxtimes\CO_{\Gr(3,W)}(-1) \to
\CO_X(-2)\boxtimes\Lambda^2\CU_{\Gr(3,W)}^\perp \to
\CO_X(-1)\boxtimes\CU_{\Gr(3,W)}^\perp \to
\CO_X\boxtimes\CO_{\Gr(3,W)} \to
\CO_{\Sigma'} \to 0.
$$
}
\end{lemma}
\begin{proof}
This is just the Koszul resolution of $\Sigma'$.
\end{proof}

Consider the truncated sequence of lemma~\ref{sigmac}
as a resolution of the sheaf of ideals $J_{\Sigma'}$.
Tensoring it with $\CO_X(2)=f^*\CO_{\PP(V)}(1)$ and
pushing forward via the projection $X\times\Gr(3,W) \to \Gr(3,W)$
we obtain exact sequence
\begin{equation}\label{p2jc}
0 \to
\Lambda^2\CU_{\Gr(3,W)}^\perp \to
W^*\otimes\CU_{\Gr(3,W)}^\perp \to
{p_2}_*(J_{\Sigma'}\otimes\CO_X(2)) \to 0.
\end{equation}

\begin{lemma}\label{fpbundlec}
${p_2}_*(J_{\Sigma'}\otimes\CO_X(1))$ is a vector bundle
of rank $15$ on $\Gr(3,W)$. Moreover, its restriction to
$\LGr(3,W) \subset \Gr(3,W)$ is isomorphic to the bundle
$\fp$ on $\LGr(3,W) = \SP(W)/P$
{\rm(}here $P$ is a parabolic subgroup of $\SP(W)$
and $\fp$ is the corresponding parabolic Lie subalgebra{\rm)}.
\end{lemma}
\begin{proof}
It is clear that the morphism
$\Lambda^2\CU_{\Gr(3,W)}^\perp \to W^*\otimes\CU_{\Gr(3,W)}^\perp$
is an embedding of vector bundles.
Moreover, restricting this bundle to $\LGr(3,W)$
we obtain a $\SP(W)$-equivariant rank $15$ subbundle of the bundle
$V^*\otimes\CO_X \cong \fg\otimes\CO_X$. But the only such bundle
is the bundle $\fp$.
\end{proof}

Consider the projectivization
$$
\TY = \PP_{\Gr(3,W)}({p_2}_*(J_{\Sigma'}\otimes\CO_X(2))).
$$
The natural embedding
${p_2}_*(J_{\Sigma'}\otimes\CO_X(2)) \to
H^0(X,\CO_X(2))\otimes\CO_{\Gr(3,W)} =
V^*\otimes\CO_{\Gr(3,W)}$
induces a morphism $\TY \exto{\tg} \PP(V^*)$.
Let $\pi:\TY\to\Gr(3,W)$ denote the canonical projection.
So we have
$$
\CO_\TY(H) = \tg^*\CO_{\PP(V^*)}(1),\quad\!\!
\tg:\TY \to \PP(V^*),\qquad\text{and}\qquad
\CO_{\TY}(h) = \pi^*\CO_{\Gr(3,W)}(1),\quad\!\!
\pi:\TY\to\Gr(3,W).
$$
%
%
Consider the product $X\times\TY$, the projection
$\id_X\times\pi:X\times\TY \to X\times\Gr(3,W)$,
the subschemes
$$
\Sigma = (\id_X\times\pi)^{-1}(\Sigma'),\qquad
\TQ = Q(X,\TY) = (X\times\TY)\times_{\PP(V)\times\PP(V^*)}Q,
$$
and the equation
$s_\TQ\in
H^0(X\times\TY,\CO_X(2)\boxtimes\CO_\TY(H)) =
H^0(X\times\TY,(f\times\tg)^*\CO_{\PP(V)\times\PP(V^*)}(1,1)) \cong
V^*\otimes V$
of $\TQ$ in $X\times\TY$.

\begin{lemma}\label{sqsigmc}
The section $s_\TQ$ vanishes on the scheme $\Sigma$.
\end{lemma}

Consider the pullback of the exact sequence of lemma~\ref{sigmac}
to $X\times\TY$.
Consider the cone $G'$ of its first morphism twisted by
$\CO_X(3)\boxtimes\CO_\TY(H)$,
and the cone $G''$ of its third morphism twisted by $\CO_X(1)$
as objects of $\D^b(X\times\TY)$,
so that we have exact triangles
$$
\xymatrix@R=10pt{
E_1\boxtimes\CO_{\TY}(H-h) \ar[r]^-{\te'} &
E_2\boxtimes \Lambda^2\CU_{\Gr(3,W)}^\perp(H) \ar[r] &
G' \\
E_1\boxtimes \CU_{\Gr(3,W)}^\perp \ar[r]^-{\te''} &
E_2\boxtimes\CO_{\TY} \ar[r] &
G''
}
$$
Then the pullback of the exact sequence of lemma~\ref{sigmac}
can be rewritten as exact triangle
\begin{equation}\label{gpgppc}
G'\otimes(\CO_X(-2)\boxtimes\CO_\TY(-H))[1] \to
G'' \to
\CO_\Sigma\otimes\CO_X(1).
\end{equation}

It is clear that morphisms $\te'$ and $\te''$ determine a structure
of $\SQ$-representations on families
$(\CO_{\TY}(H-h),\Lambda^2\CU_{\Gr(3,W)}^\perp(H))$ and
$(\CU_{\Gr(3,W)}^\perp,\CO_{\TY})$ on $\TY$.

\begin{lemma}\label{fam16c}
There exist a unique extension
\begin{equation}\label{tfc}
0 \to (\CO_{\TY}(H-h),\Lambda^2\CU_{\Gr(3,W)}^\perp(H)) \to
(\TF_1,\TF_2) \to (\CU_{\Gr(3,W)}^\perp,\CO_{\TY}) \to 0
\end{equation}
in the category of families of $\SQ$-representations on $\TY$
such that the cone of the induced morphism
$\te:E_1\boxtimes\TF_1 \to E_2\boxtimes\TF_2$ is a sheaf
supported at $\TQ$.
\end{lemma}

Lemma~\ref{fam16c} says that the family $(\TF_1,\TF_2)$ on $\TY$
is a family of $\SQ$-representations and that  on $X\times\TY$
we have the following exact sequence:
$$
0 \to E_1\boxtimes \TF_1 \exto{\te} E_2\boxtimes \TF_2 \to \ti_*\TE \to 0,
$$
where $\TE$ is a coherent sheaf on $\TQ$.

\begin{lemma}\label{c_phic}
The family $(\TF_1,\TF_2)$ induces a morphism $\phi:\TY\to \BY_4$
such that $\tg = \bg_4\circ\phi$.
\end{lemma}

Let $\TYO = \TY \setminus\tg^{-1}(Z_2)$ and denote the restriction
of $\tg$ and $\phi$ to $\TYO$ by the same symbols.

Recall that we have a sheaf of Azumaya algebras $\CA_4$ on $Y_4$
and that $F_2$ is a locally projective $\CA_4$-module on $Y_4$ of rank~$4$.
Let $\PP_{Y_4,\CA_4}(F_2^*)$ denote the moduli space of locally projective
quotient $\CA_4$-modules of rank~$1$. Then $\PP_{Y_4,\CA_4}(F_2^*)$ is
a twisted $\PP^3$-bundle over $Y_4$.

\begin{proposition}\label{tyoispf2c}
We have $\TYO = \PP_{Y_4,\CA_4}(F_2^*)$.
\end{proposition}
\begin{proof}
The map $\TYO \to \PP_{Y_4,\CA_4}(F_2^*)$ can be constructed
like in proposition~\ref{tyoispf2}. Thus it remains to
construct the inverse map.
%
Choose a point $\ty\in\PP_{Y_4,\CA_4}(F_2^*)$.
Let $\rho:W\otimes R_1 \to R_2$ be the corresponding $(4,4)$-dimensional
representation of the quiver $\SQ$, and let $R_2\to\kk$ be the corresponding
one-dimensional quotient. Let $R_2'\subset R_2$ be its kernel
and consider the composition
$E_1\otimes R_1 \exto{e_\rho} E_2\otimes R_2 \to E_2$.
Taking $\Hom(E_1,-)$ we obtain a map $R_1 \to W^*$.
Let $R'_1$ denote the kernel of this map. Then we have a diagram~(\ref{ty})
and an exact sequence~(\ref{cpc})
%
Note that the sheaf $C$ is supported on the quadric
$X_y = H_{g_3(y)}\cap X\subset X$ and
its rank at the generic point of $X_y$ equals
$(c_1(E_2\otimes R_2) - c_1(E_1\otimes R_1))/2 = 2$.

We claim that $\dim R'=1$ and we are going to prove it
by a case by case analysis.

Since $(R'_1,R'_2)$ is a subrepresentation of $(R_1,R_2)$
and $(R_1,R_2)$ is stable by lemma~\ref{gmpro} we have
$$
\dim R'_1 < 3,
$$
so we need to reject two cases: $\dim R'_1=2$ and $\dim R'_1 = 0$.

The case $\dim R'_1 = 2$ is  treated along the same lines
as the case $\dim R'_1=2$ in the proof of proposition~\ref{tyoispf2}.
The only difference is that in the case $\dim R'_1 = 2$ the sheaf
$T(C')$ can have rank 1 at the generic point of $X_y$,
but it is easy to see that this is possible only in the case
$C' \cong \CO_X(1) \oplus (\CO_X(1)\oplus\CO_X(1))/(\CO_X\oplus\CO_X)$,
which means that $g_4(y)\in Z_2$.

Assume that $R'_1 = 0$. Then the map $R_1\to W^*$ is injective,
so we can consider $R_1$ as a 4-dimensional subspace in $W^*$.
If the cokernel of the morphism $K'' \to C' = \CO_X(1)^{\oplus3}$
is supported on a quadric $X_y\subset X$, then the morphism
$\CO_X(-1)^{\oplus3} \cong \CO_X(1)^{\oplus3}\otimes\CO_X(-2) \to
\CO_X(1)^{\oplus3}$ factors through $K''$ and the cokernel of the
corresponding morphism $\CO_X(-1)^{\oplus3} \to K''$ is supported
on the same quadric $X_y\subset X$. On the other hand,
it is easy to compute that $\dim\Hom(\CO_X(-1),K'') = 6$.
Therefore the dimension of the set of such
quadrics $X_y\subset X$ has dimension less than
$\dim\Gr(4,W^*) + \dim\Gr(3,6) = 8 + 9 = 17$. On the other hand,
the set of such quadrics is evidently $\GL(W)$-invariant
in $\PP(V^*) = \PP(S^2W^*)$, and it follows from the description
of all $\GL(W)$-invariant closed subsets of $\PP(S^2W^*)$ that
this set is contained in $Z_2$.

Thus we have shown that if $g_4(y)\not\in Z_2$ then
$\dim R'_1 = 1$, $\dim R_1/R'_1 = 3$.
The top and the bottom lines of the diagram~(\ref{ty}) take form
\begin{equation}\label{wpc}
0 \to \CO_X \exto{e'} \CO_X(1)^{\oplus3} \to C' \to 0,
\end{equation}
\begin{equation}\label{wppc}
0 \to K'' \to \CO_X^{\oplus3} \exto{e''} \CO_X(1) \to C'' \to 0.
\end{equation}
The map $e''$ is given by a 3-dimensional subspace of $W$, that is by a
point $u \in \Gr(3,W)$. The same arguments as in the proof
of proposition~\ref{tyoispf2} show that we can glue the sequence~(\ref{wpc})
twisted by $\CO_X(-3)$ with the sequence~(\ref{wppc}) twisted by $\CO_X(-1)$
into the following exact sequence
$$
0 \to \CO_X(-3) \exto{e'} \CO_X(-2)^{\oplus3} \to
\CO_X(-1)^{\oplus3} \exto{e''} \CO_X \to C'' \to 0.
$$
Then it is easy to see that this sequence must coincide
with the sequence of lemma~\ref{sigmac} restricted
to a fiber of $X\times\Gr(3,W)$ over the point of $u\in\Gr(3,W)$.
Associating to the point $\ty$ the point $u\in\Gr(3,W)$ we obtain a map
$\PP_{Y_4,\CA_4}(F_2^*) \to \Gr(3,W)$.
On the other hand, the support quadric $X_y$ of $C$
contains $\supp(C'') = \Sigma'_u$, therefore the map
$$
\PP_{Y_4,\CA_4}(F_2^*) \to \Gr(3,W)\times\PP(V^*),\qquad
\ty \mapsto (u,g_4(y))
$$
factors through
$\TY = \PP_{\Gr(3,W)}({p_2}_*(J_{\Sigma'}\otimes\CO_X(2))) \subset
\Gr(3,W)\times\PP(V^*)$. Moreover, it is clear that its image lies in $\TYO$,
so we have constructed a map $\PP_{Y_4,\CA_4}(F_2^*) \to \TYO$.

It remains to note that the constructed maps
are mutually inverse.
\end{proof}

\begin{theorem}\label{th_c}
The map $g_4:Y_4\to\PP(V^*)\setminus Z_2$ is a double covering.
\end{theorem}
\begin{proof}
Since $\TYO$ is smooth it follows from proposition~\ref{tyoispf2c} that
$Y_4$ is smooth and moreover $\dim Y_4 = \dim\TYO - 3 = 20$.
On the other hand, since $\tg = \bg_4\circ\phi$, $\tg$ and $\phi$
are proper and dominant, it follows that $g_4$ is a proper dominant map
$Y_4 \to \PP(V^*) \setminus Z_2$.
Since $Y_4$ and $\PP(V^*)\setminus Z_2$ are smooth varieties
of equal dimension and
$\rank\Pic Y_4 = \rank\Pic\TYO - 1 = 1 = \rank\Pic(\PP(V^*)\setminus Z_2)$
it follows that $g_4$ is a finite covering. It remains to show that
$\deg g_4 = 2$. This can be done by a direct calculation using
the intersection theory of $\Gr(3,W)$.
Alternatively, we can argue as follows.

Take a nondegenerate 2-form on $W^*$ and consider the zero
locus of the corresponding section of the vector bundle $W/\CU_{\Gr(3,W)}$
on $\Gr(3,W)$. It is clear that this is the Lagrangian Grassmannian
$\LGr(3,W)$. Consider the fiber product $\TYO\times_{\Gr(3,W)}\LGr(3,W)$.
By lemma~\ref{fpbundlec} we have $\TYO\times_{\Gr(3,W)}\LGr(3,W)$
is an open subset of $\PP_{\LGr(3,W)}(\fp)$, so and it is clear that the map
$\tg_{|\PP_{\LGr(3,W)}(\fp)}:\PP_{\LGr(3,W)}(\fp) \to \PP(V^*) = \PP(\fg)$
is the projectivization of the partial Springer--Grothendieck map
$\Tot(\fp) \to \fg$. Therefore, by~\cite{Hu} its degree equals $16/2 = 8$.
On the other hand, since the restriction of $W/\CU_{\Gr(3,W)}$ to the fiber
$\PP^3$ of $\TYO$ over $Y_4$ is isomorphic to $\Omega_{\PP^3}$ by~(\ref{tfc})
we see that the degree of the map
$\phi_{|\PP_{\LGr(3,W)}(\fp)}:\PP_{\LGr(3,W)}(\fp) \to Y_4$ equals
$c_3(T_{\PP^3}) = 4$. Hence $\deg g_4 = 8/4 = 2$.
\end{proof}

\begin{corollary}
If the data {\rm(D.3)--(D.5)} is given by
$Y = Y_4$, $\CA_Y = \CA_4$, $\BZ = \OZ_4 = Z_2$, $g = g_4$,
with $(F_1,F_2)$ being the universal family, and
with $\phi$ being induced by the $\SQ$-representation structure on $(F_1,F_2)$,
then all conditions {\rm(C.1)--(C.8)} are satisfied.
Moreover, the branching locus of the map $g:Y \to \PP(V^*)\setminus\BZ$
is the hypersurface ${X^\vee}\setminus\BZ$, where
${X^\vee} \subset \PP(V^*)$ is the projectively dual sextic hypersurface.
\end{corollary}
\begin{proof}
Since the conditions~(C.1)--(C.3), (\ref{addc}) and (\ref{codimbro}) are
satisfied by lemma~\ref{conds_b}, it follows from theorem~\ref{qu_th}
that the conditions (C.4)--(C.6) are true. On the other hand,
(C.8) is true because $\BZ$ is the locus of rank 4 quadrics
by lemma~\ref{yz1_c} and $\dim\BZ\cap H = 17 - 1 = 16$, while
$N - i - 2 = 21 - 3 - 2 = 16$.
Finally, by theorem~\ref{th_c} the map $Y = Y_4 \to \PP(V^*)$ is a double
covering and
$\codim_{Y\times Y}(Y\times_{\PP(V^*)}Y) = \dim Y = \dim\PP(V^*) = 20$,
hence
$\dim X + N - \codim_{Y\times Y}(Y\times_{\PP(V^*)}Y) = 5 + 21 - 20 = 6 = 2i$,
and (C.7) is true.

Finally, we note that the branching locus of $g$
is nothing but the singular locus of $g$, hence
it coincides with ${X^\vee}\setminus\BZ$ by corollary~\ref{singg}.
\end{proof}

\section*{Appendix D. Azumaya algebraic varieties}

\refstepcounter{section}
\nc{\mysubsection}[1]{{\bf{#1}.\/}}

In this section we work out some basic facts about Azumaya varieties.
All algebraic varieties in this section are assumed to be embeddable.

\mysubsection{Spaces, categories and functors}

\begin{definition}
An {\sf Azumaya algebraic variety}\/ over a field $\kk$ is a pair
$(X,\CA_X)$, where $X$ is an algebraic variety of finite type over $\kk$
and $\CA_X$ is a sheaf of semisimple $\CO_X$-algebras which is locally free
of finite rank over $\CO_X$. A {\sf morphism of Azumaya varieties}\/
$f:(X,\CA_X) \to (Y,\CA_Y)$ is a pair $(f_\circ,f_\CA)$,
where $f_\circ:X\to Y$ is a morphism of algebraic varieties, and
$f_\CA:f_\circ^*\CA_Y \to \CA_X$ is a homomorphism of
$f_\circ^*\CO_Y\cong\CO_X$-algebras. The composition of morphisms
is defined naturally.
\end{definition}

It is clear that Azumaya algebraic varieties over $\kk$
form a category. Every algebraic variety $X$ can be considered
as an Azumaya variety by taking $\CA_X = \CO_X$. Thus the category
of algebraic varieties can be considered as a full subcategory of
the category of Azumaya varieties. For each Azumaya variety $(X,\CA_X)$
we have a canonical morphism $\pi_X:(X,\CA_X) \to X$,
the {\sf structure morphism} of $(X,\CA_X)$.

With each Azumaya algebraic variety $(X,\CA_X)$ we associate
abelian categories $\Qcoh(X,\CA_X)$ and $\Coh(X,\CA_X)$
of quasicoherent and coherent sheaves of right $\CA_X$-modules on $X$
and their bounded and unbounded derived categories $\D^?_{qc}(X,\CA_X)$
and $\D^?(X,\CA_X)$, where $? = b$, $-$, $+$ or nothing.

\begin{remark}
Recall that the Morita-equivalent algebras $\CA_X$ and $\CA'_X$ give rise
to equivalent categories $\Coh(X,\CA_X)$ and $\Coh(X,\CA'_X)$.
In particular, if $\CA_X \cong \CEnd(\CV)$, the algebra
of endomorphisms of a vector bundle $\CV$,
then $\Coh(X) \cong \Coh(X,\CA_X)$, the equivalence
is given by tensoring with $\CV$.
However, we prefer to distinguish between
Morita-equivalent Azumaya varieties.
\end{remark}

\begin{definition}
A morphism of Azumaya varieties $f:(X,\CA_X) \to (Y,\CA_Y)$
is called {\sf strict}, if $\CA_X \cong f_0^*\CA_Y$ and
$f_\CA$ is the identity.
Similarly, morphism $f$ is called {\sf an extension},
if $X \cong Y$ and $f_\circ$ is the identity.
\end{definition}

Every morphism $f:(X,\CA_X) \to (Y,\CA_Y)$ admits the following
{\sf canonical decomposition}:
$$
(X,\CA_X) \exto{f^e} (X,f_\circ^*\CA_Y) \exto{f^s} (Y,\CA_Y),
$$
where $f^e$ is an extension, and $f^s$ is strict.

If $f:(X,\CA_X) \to (Y,\CA_Y)$ is a morphism of Azumaya varieties
and $F\in\Coh(X,\CA_X)$, we define $R^0f_*F\in\Coh(Y,\CA_Y)$ as
the sheaf $R^0{f_\circ}_*F$ with an $\CA_Y$-module structure
induced by the homomorphism
$$
R^0{f_\circ}_*F \otimes \CA_Y \cong
R^0{f_\circ}_*(F \otimes f_\circ^*\CA_Y) \exto{f_\CA}
R^0{f_\circ}_*(F \otimes \CA_X) \to
R^0{f_\circ}_*F.
$$
It is clear that $R^0f_*$ is a functor and that exactness properties
of $R^0f_*$ are the same as those of $R^0{f_\circ}_*$.
If $f:(X,\CA_X) \to \Spec\kk$ is the projection to a point,
we denote the pushforward functor by $\Gamma(X,-)$.
Note that the pushforward with respect to the structure morphism
$\pi_X$ is just the forgetting of the $\CA_X$-module structure.

Similarly, for $G\in\Coh(Y,\CA_Y)$ we define $L_0f^*G\in\Coh(X,\CA_X)$ as
$$
L_0f^*G = L_0f_\circ^*G \otimes_{f_\circ^*\CA_Y}\CA_X.
$$
Then $L_0f^*$ is also a functor.
Note that if $f$ is a strict morphism then we have $L_0f^* = L_0f_\circ^*$.

If $F\in\Coh(X,\CA_X)$ and $F'\in\Coh(X,\CA_X^\opp)$, where
$\CA_X^\opp$ is the opposite algebra, we can consider their
tensor product over $\CA_X$, $F\otimes_{\CA_X} F' \in \Coh(X)$.
Note, that if $F$ (resp.~$F'$) admits some additional module structure,
commuting with the structure of the right (resp.~left) $\CA_X$-module, then
the tensor product $F\otimes_{\CA_X} F'$ preserves this structure.
For example, if $F'$ is a sheaf of $\CA_X$-bimodules, then
$F\otimes_{\CA_X} F' \in \Coh(X,\CA_X)$.

Similarly, for $F,F'\in\Coh(X,\CA_X)$ we can consider the sheaf
$\CHom_{\CA_X}(F,F') \in \Coh(X)$ of local homomorphisms of $\CA_X$-modules.
Note, that if $F$ (resp.~$F'$) admits some additional module structure,
commuting with the structure of the right $\CA_X$-module, then
the sheaf $\CHom_{\CA_X}(F,F')$ preserves this structure.
For example, if $F'$ is a sheaf of $\CA_X$-bimodules, then
$\CHom_{\CA_X}(F,F') \in \Coh(X,\CA_X)$.

\begin{lemma}\label{locproj}
A sheaf $F\in\Coh(X,\CA_X)$ is locally projective over $\CA_X$
in the Zariski topology iff ${\pi_X}_*F \in \Coh(X)$
is locally free, where $\pi_X:(X,\CA_X) \to X$ is
the structure morphism.
\end{lemma}
\begin{proof}
Note that $\CA_X^\opp\otimes\CA_X$ is a sheaf of matrix algebras,
therefore $\CA_X$ considered as a $\CA_X^\opp\otimes\CA_X$-module
is locally projective (in the Zariski topology), i.e.\ locally is
a direct summand of $\CA_X^\opp\otimes\CA_X$. It follows that
$F \cong F\otimes_{\CA_X}\CA_X$ is locally a direct summand of
$F\otimes_{\CA_X}(\CA_X^\opp\otimes\CA_X)\cong F\otimes_{\CO_X}\CA_X$,
a free $\CA_X$-module.
\end{proof}

It is clear that the category of quasicoherent $\CA_X$-modules on
an Azumaya variety $(X,\CA_X)$ admits enough injectives and enough
locally free objects. Therefore, the pushforward, the pullback,
the tensor product and the $\CHom$ functors can be extended to the derived
categories in a natural way. From now on we use the notation $f_*$, $f^*$,
$\otimes_{\CA_X}$ and $\RCHom_{\CA_X}$ to denote the derived functors.
The derived functor of $\Gamma(X,-)$ is denoted by $\RGamma(X,-)$.

\begin{remark}
Note that
$$
F\otimes_{\CA_X} G \cong
(F\otimes_{\CO_X} G)\otimes_{\CA_X\otimes\CA_X^\opp}\CA_X
\quad\text{and}\quad
\RCHom_{\CA_X}(F,G) \cong
\CHom_{\CA_X\otimes\CA_X^\opp}(\CA_X,\RCHom_{\CO_X}(F,G)).
$$
Therefore the bar-resolution
$$
\dots \to
\CA_X\otimes\CA_X^{\otimes 2}\otimes\CA_X^{\opp} \to
\CA_X\otimes\CA_X\otimes\CA_X^{\opp} \to
\CA_X\otimes\CA_X^{\opp} \to
\CA_X \to 0
$$
gives the bar-resolutions of the functors
$\otimes_{\CA_X}$ and $\RCHom_{\CA_X}$
$$
\begin{array}{l}
\dots \to
F\otimes_{\CO_X}\CA_X^{\otimes 2}\otimes_{\CO_X} G \to
F\otimes_{\CO_X}\CA_X\otimes_{\CO_X} G \to
F\otimes_{\CO_X}G \to
F\otimes_{\CA_X}G \to 0
\smallskip\\
0 \to \RCHom_{\CA_X}(F,G) \to
\RCHom_{\CO_X}(F,G) \to
\RCHom_{\CO_X}(F\otimes_{\CO_X}\CA_X,G) \to \dots
\end{array}
$$
\end{remark}

\begin{remark}
If the algebra $\CA_X$ is commutative then the category $\Coh(X,\CA_X)$
is equivalent to the category of coherent sheaves on the algebraic variety
$\Spec_X\CA_X$ (the relative spectrum) which is \'etale over $X$. Following
this analogy we consider a pair $(X,\CA_X)$ in a noncommutative case
as a {\em noncommutative \'etale covering}\/ of $X$.
\end{remark}

Finally, we define the twisted pull-back functor
$f^!:\D^+(Y,\CA_Y) \to \D^+(X,\CA_X)$ extending
the usual twisted pull-back functor defined
in~\cite{H} to the category of Azumaya varieties.

The definition for a strict morphism is the same as in {\em loc.\ cit.}.
If $f$ is strict and smooth we define
$$
f^!(G) = f^*G\otimes_{\CO_X}\omega_{X/Y}[\dim X - \dim Y],
$$
and if $f_\circ$ is strict and finite we define
$$
f^!G = \RCHom_{\CA_Y}({f_\circ}_*\CA_X,G),
$$
where we identify the category $\D^+(X,\CA_X)$ with $\D^+(Y,{f_\circ}_*\CA_X)$
via the functor ${f_\circ}_*$. Now assume that $f_\circ$ is embeddable
(this is always the case when $X$ is embeddable). Decomposing $f$
into a composition $f_\circ = f_\circ'\circ f_\circ''$ with
$f_\circ'$ smooth and $f''_\circ$ finite we define $f^! = {f''}^!\circ{f'}^!$.
A standard verification ({\em loc.\ cit.}) shows that this definition
doesn't depend on a choice of a decomposition
$f_\circ = f_\circ'\circ f_\circ''$.

Now, for an arbitrary morphism $f:(X,\CA_X) \to (Y,\CA_Y)$ we consider
its canonical decomposition $f = f^s\circ f^e$ with $f^s$ strict and
$f^e$ extension, note that for any $G\in\D^+(Y,\CA_Y)$ we have
${f^s}^!G \in \D^+(X,f_\circ^*\CA_Y)$ and define $f^!G$ by the formula
\begin{equation}\label{snavshriek}
f^!G = \RCHom_{f_\circ^*\CA_Y}(\CA_X,{f^s}^!G).
\end{equation}

\begin{lemma}\label{fshriek}
If $f$ is strict then $f^! = f_\circ^!$. Moreover
$$
\begin{array}{ll}
f_\circ^!F = f_\circ^*F\otimes_{\CO_X}\det\CN_{X/Y}[\dim X - \dim Y], &
\text{if $f$ is a locally complete intersection}\\
& \hspace{4cm} \text{closed embedding,}\\
f_\circ^!F = f_\circ^*F\otimes_{\CO_X}\omega_{X/Y}[\dim X - \dim Y], &
\text{if $f_\circ$ is smooth.}
\end{array}
$$
\end{lemma}
\begin{proof}
Evident.
\end{proof}

\mysubsection{Relations between functors}

\begin{lemma}[Functoriality of the pushforward]
If $(X,\CA_X) \exto{f} (Y,\CA_Y) \exto{g} (Z,\CA_Z)$
are morphisms of Azumaya varieties
then $g_*f_* \cong (gf)_*$ on $\D(X,\CA_X)$,
a functorial isomorphism.
In particular,\\
$\RGamma(Y,f_*(-)) \cong \RGamma(X,-)$.
\end{lemma}
\begin{proof}
Since the pushforward functor for Azumaya varieties coincides with
the usual pushforward, the objects $R^0g_*R^0f_*(F)$ and $R^0(gf)_*(F)$
are identified in $\Coh(Z)$ by a natural isomorphism for any
$F\in\Coh(X,\CA_X)$. Moreover, it is clear that this isomorphism
identifies also the $\CA_Z$-module structures on them.
Therefore $R^0g_*R^0f_* \cong R^0(gf)_*$ as functors
$\Coh(X,\CA_X) \to \Coh(Z,\CA_Z)$ whereof we deduce an isomorphism
of derived functors as well.
\end{proof}

\begin{lemma}[Functoriality of the pullback]
If $(X,\CA_X) \exto{f} (Y,\CA_Y) \exto{g} (Z,\CA_Z)$
are morphisms of Azumaya varieties
then $f^*g^* \cong (gf)^*$ on $\D^-(Z,\CA_Z)$.
Moreover, if both $f$ and $g$ have finite $\Tor$-dimension
then $f^*g^* \cong (gf)^*$ on $\D(Z,\CA_Z)$.
\end{lemma}
\begin{proof}
Since
$L_0f^*L_0g^*(F) :=
L_0f_\circ^*(L_0g_\circ^*F\otimes_{g_\circ^*\CA_Z}\CA_Y)
\otimes_{f_\circ^*\CA_Y}\CA_X \cong
L_0f_\circ^*L_0g_\circ^*F\otimes_{f_\circ^*g_\circ^*\CA_Z}f_\circ^*\CA_Y
\otimes_{f_\circ^*\CA_Y}\CA_X \cong
L_0f_\circ^*L_0g_\circ^*F\otimes_{f_\circ^*g_\circ^*\CA_Z}\CA_X =
L_0(gf)^*(F)$ we have an isomorphism of functors
$L^0f^*L_0g^* \cong L_0(gf)_*$ whereof we deduce an isomorphism
of derived functors as well.
\end{proof}

\begin{lemma}[Associativity of the tensor product]
If $F\in\D^-(X,\CA_1)$, $G\in\D^-(X,\CA_1^\opp\otimes\CA_2)$,
and $H\in\D^-(X,\CA_2^\opp)$ then
$(F\otimes_{\CA_1} G) \otimes_{\CA_2} H \cong
F \otimes_{\CA_1} (G\otimes_{\CA_2} H)$.
\end{lemma}
\begin{proof}
Similarly.
\end{proof}

\begin{lemma}[Commutativity of the tensor product]
If $F\in\D^-(X,\CA_X)$ and $G\in\D^-(Y,\CA_X^\opp)$ then
$F\otimes_{\CA_X} G \cong G\otimes_{\CA^\opp_X} F$.
\end{lemma}
\begin{proof}
Similarly.
\end{proof}

\begin{lemma}[The projection formula]
If $f:(X,\CA_X) \to (Y,\CA_Y)$ is a morphism of Azumaya varieties,
$F\in\D^-(X,\CA_X)$, $G\in\D^-(Y,\CA_Y^\opp)$, and $H\in\D^-(X,\CO_X)$
then
$$
f_*(F\otimes_{\CA_X} f^*G) \cong f_*F \otimes_{\CA_Y} G,
\qquad\text{and}\qquad
f_*(f_\circ^*G\otimes_{\CO_X} H) \cong G \otimes_{\CO_Y} f_*H.
$$
Moreover, if $f$ is strict, $F\in\D^-(X,\CO_X)$ and $G\in\D^-(Y,\CA_Y)$ then
$f_*(F\otimes_{\CO_X} f^*G) \cong f_*F \otimes_{\CO_Y} G$.
\end{lemma}
\begin{proof}
First of all assume that $f$ is strict.
If $F$ is $f_*$-acyclic and $G$ is locally projective over $\CA_Y$ then
$f_*(F\otimes_{\CO_X} f^*G) \cong f_*F \otimes_{\CO_Y} G$
by the usual projection formula. Using the bar-resolution we then deduce that
$f_*(F\otimes_{f_\circ^*\CA_Y} f*G) \cong f_*F \otimes_{\CA_Y} G$,
the first isomorphism of the lemma. The standard argument \cite{H}, II.5.6
then proves this for all $F$ and $G$. On the other hand, if $f$ is an extension
the first isomorphism of the lemma follows from
$F\otimes_{\CA_X}(\CA_X\otimes_{\CA_Y} G) \cong F\otimes_{\CA_Y} G$.
Finally, if $f$ is an arbitrary morphism we consider its canonical
decomposition into the product of a strict morphism and an extension and
deduce the first isomorphism of the lemma from the above remarks.

The second isomorphism is just the usual projection formula,
and it remains to note that it is compatible with the $\CA_Y$-module
structures.
\end{proof}

\begin{lemma}
If $F\in\D^-(X,\CA_X)$, $G\in\D^+(X,\CA_X)$ then we have
a functorial isomorphism\\
$\RHom_{\CA_X}(F,G) \cong \RGamma(X,\RCHom_{\CA_X}(F,G))$.
\end{lemma}
\begin{proof}
We have evidently $\Hom_{\CA_X}(F,G) \cong \Gamma(X,\CHom_{\CA_X}(F,G))$
for all sheaves $F,G\in\Coh(X,\CA_X)$. The standard argument shows that
we still have this identity for the derived functors.
\end{proof}

\begin{lemma}
If $f:(X,\CA_X) \to (Y,\CA_Y)$ is a strict morphism, such that
$f_\circ$ has finite Tor-dimension and
$F\in\D^-(Y,\CA_Y)$, $G\in\D^+(Y,\CA_Y)$,
then
$f_\circ^*\RCHom_{\CA_Y}(F,G) \cong \RCHom_{\CA_X}(f^*F,f^*G)$.
\end{lemma}
\begin{proof}
For the strict morphism $f$ we have $f^* = f_\circ^*$,
so combining the usual formula for the pullback of $\RCHom$
with the bar-resolution of $\RCHom_{\CA_Y}$ we deduce the lemma.
\end{proof}

\begin{lemma}
If $f:(X,\CA_X) \to (Y,\CA_Y)$ is a strict morphism,
$F\in\D^-(Y,\CA_Y)$, $G\in\D^-(Y,\CA_Y^\opp)$, then
$f_\circ^*(F\otimes_{\CA_Y} G) \cong f^*F \otimes_{\CA_X} f^*G$.
\end{lemma}
\begin{proof}
As above.
\end{proof}

\begin{lemma}
If $f:(X,\CA_X) \to (Y,\CA_Y)$ is a morphism,
$F\in\D^-(Y,\CA_Y)$ and $G\in\D^-(Y,\CO_Y)$, then
$f^*(F\otimes_{\CO_Y} G) \cong f^*F \otimes_{\CO_X} f_\circ^*G$.
\end{lemma}
\begin{proof}
We have
$f^*(F\otimes_{\CO_Y} G) \cong
f_\circ^*(G\otimes_{\CO_Y} F) \otimes_{f_\circ^*\CA_Y} \CA_X \cong
(f_\circ^*G\otimes_{\CO_X} f_\circ^*F) \otimes_{f_\circ^*\CA_Y} \CA_X \cong \\
f_\circ^*G\otimes_{\CO_X} (f_\circ^*F \otimes_{f_\circ^*\CA_Y} \CA_X) \cong
f_\circ^*G\otimes_{\CO_X} f^*F \cong
f^*F \otimes_{\CO_X} f_\circ^*G$.
\end{proof}

\begin{lemma}
If $f:(X,\CA_X) \to (Y,\CA_Y)$ is a morphism,
$F\in\D^-(Y,\CA_Y)$ and $G\in\D^+(X,\CA_X)$, then
$f_*\RCHom(f^*F,G) \cong \RCHom(F,f_*G)$.
In particular, the functor $f^*$ is left adjoint to $f_*$.
\end{lemma}
\begin{proof}
If $f$ is a strict morphism then $f_* = {f_\circ}_*$ and
$f_*\RCHom_{\CO_X}(f^*F,G) \cong \RCHom_{\CO_Y}(F,f_*G)$
by the usual adjointness property. Using the bar-resolution
we then deduce that
$f_*\RCHom_{f_\circ^*\CA_Y}(f^*F,G) \cong \RCHom_{\CA_Y}(F,f_*G)$.
On the other hand, if $f$ is an extension the isomorphism of the lemma
follows from
$\RCHom_{\CA_Y}(F,f_*G) \cong
\RCHom_{\CA_Y}(F,G) \cong
\RCHom_{\CA_X}(F\otimes_{\CA_Y}\CA_X,G) =
\RCHom_{\CA_X}(f^*F,G)$.
If $f$ is an arbitrary morphism we consider its canonical decomposition
into the product of a strict morphism and an extension and deduce the
first isomorphism of the lemma from the functoriality of the pullback
and of the pushforward.
\end{proof}

\mysubsection{Flatness and smoothness}

\begin{definition}
A morphism of Azumaya varieties $f:(X,\CA_X) \to (Y,\CA_Y)$
is called {\sf flat}, if the pullback functor
$f^*:\D(Y,\CA_Y) \to \D(X,\CA_X)$ is exact
with respect to the standard t-structures.
\end{definition}

\begin{lemma}
Morphism $f:(X,\CA_X) \to (Y,\CA_Y)$ is flat iff
$f_\circ:X \to Y$ is flat.
\end{lemma}
\begin{proof}
Let $\pi_X:(X,\CA_X) \to X$, and $\pi_Y:(Y,\CA_Y)\to Y$ be
the structure morphisms. It is clear that $\pi_X^*$ and $\pi_Y^*$
are exact and conservative ($\pi_X^*F = 0$ implies $F=0$).
On the other hand, we have $\pi_Y\circ f = f_\circ\circ\pi_X$,
hence $f^*\pi_Y^* = \pi_X^* f_\circ^*$.
If $f_\circ$ is not flat, then there exists $G \in \Coh(Y)$ such that
$\CH^i(f_\circ^*G) \ne 0$ for some $i<0$.
Then
$\CH^i(f^*\pi_Y^*G) =
\CH^i(\pi_X^*f_\circ^*G) =
\pi_X^*\CH^i(f_\circ^*G) \ne 0$,
while $\pi_Y^*G \in \Coh(Y,\CA_Y)$. Hence $f^*$ is not exact.
On the other hand, assume that $f_\circ$ is flat. Then the strict part
$f^s$ of $f$ is flat. On the other hand, the extension part $f^e$ of $f$
is flat, because it coincides with the functor
$F \mapsto F\otimes_{f_\circ^*\CA_Y}\CA_X$, and $\CA_X$ is
locally projective over $f_\circ^*\CA_Y$ by lemma~\ref{locproj}.
\end{proof}

\begin{definition}
An object $F \in \D(X,\CA_X)$ is called a {\sf perfect complex},
if it is locally quasiisomorphic to a bounded complex
of projective $\CA_X$-modules of finite rank.
We denote by $\D^\perf(X,\CA_X)$ the fully faithful subcategory
of $\D(X,\CA_X)$ formed by all perfect complexes.
\end{definition}

\begin{definition}
A triangulated category $\D$ is called {\sf$\Ext$-bounded}, if
for any objects $F,F'\in\D$ there exist $p,q\in\ZZ$ such that
$\Hom_\D(F,F'[t]) = 0$ for $t\not\in[p,q]\subset\ZZ$.
\end{definition}

\begin{lemma}\label{snav_sm}
The following conditions for an Azumaya variety $(X,\CA_X)$ are equivalent:

$(i)$ the underlying algebraic variety $X$ is smooth;

$(ii)$ $\D^b(X,\CA_X) = \D^\perf(X,\CA_X)$.

$(iii)$ the bounded derived category $\D^b(X,\CA_X)$ is $\Ext$-bounded.
\end{lemma}
\begin{proof}
$(i) \Longrightarrow (ii)$:
Take any $F\in\D^b(X,\CA_X)$. Since $\Coh(X,\CA_X)$ has enough
locally free objects, we can construct a locally free over $\CA_X$
resolution $F_{n} \to F_{n-1} \to \dots \to F_1 \to F_0 \to F \to 0$ for
arbitrary large $n$. On the other hand, if $n$ is sufficiently large, then
$F_{n+1} := \Ker (F_{n} \to F_{n-1})$ is locally free over $\CO_X$, hence
locally projective over $\CA_X$ by lemma~\ref{locproj}. Thus,
$F$ is quasiisomorphic to a finite complex of locally projective
$\CA_X$-modules, so $F$ is a perfect complex.

$(ii) \Longrightarrow (iii)$:
Since $F$ is a perfect complex, we have $\RCHom(F,F') \in \D^{[p,q]}(X)$
for some $p,q\in\ZZ$, hence
$\RHom_{\CA_X}(F,F') =
\RGamma(X,\RCHom_{\CA_X}(F,F')) \in
\D^{[p,q+\dim X]}(\Ab)$.

$(iii) \Longrightarrow (i)$:
Let $\pi:(X,\CA_X) \to X$ denote the structure morphism. Then
$$
\RHom_{\CA_X}(\pi^*F,\pi^*F') \cong
\RHom_{\CO_X}(F,\pi_*\pi^*F') \cong
\RHom_{\CO_X}(F,F'\otimes\CA_X)
$$
for all $F,F'\in\D^b(X,\CO_X)$. If $X$ is not smooth, take $F = F'$
to be the structure sheaf of a singular point. Then
$F'\otimes\CA_X \cong {F'}^{\oplus \rank\CA_X}$ and
$\Hom_{\CO_X}(F,F'\otimes_{\CO_X}\CA_X[t]) \ne 0$ for arbitrary large $t$.
Hence $\D^b(X,\CA_X)$ is not $\Ext$-bounded.
\end{proof}

\begin{definition}
An Azumaya variety $(X,\CA_X)$ is called {\sf smooth}, if any
of the equivalent conditions~\ref{snav_sm} is satisfied.
A morphism $f:(X,\CA_X) \to (Y,\CA_Y)$ of Azumaya varieties
is called {\sf smooth} if the underlying morphism $f_\circ:X\to Y$ is smooth.
\end{definition}

\begin{lemma}
If morphism $f:(X,\CA_X) \to (Y,\CA_Y)$ is smooth and
$(Y,\CA_Y)$ is a smooth Azumaya variety, then $(X,\CA_X)$ is smooth.
A composition of two smooth morphisms is smooth.
\end{lemma}
\begin{proof}
Evident.
\end{proof}

\mysubsection{Further relations}

\begin{lemma}\label{f22}
If $F\in\D^-(X,\CA_1)$, $G\in\D^\perf(X,\CA_2)$ and
$H\in\D^+(X,\CA_2^\opp\otimes\CA_1)$ then we have\\
$\RCHom_{\CA_1}(F,G\otimes_{\CA_2} H) \cong
G \otimes_{\CA_2}\RCHom_{\CA_1}(F,H)$.
\end{lemma}
\begin{proof}
It is easy to construct a homomorphism from the RHS to the LHS.
It is clear that it is an isomorphism for $G$ being locally projective,
hence it is an isomorphism as well for $G$ being any perfect complex.
\end{proof}

\begin{lemma}\label{fdg}
If either $F\in\D^-(X,\CA_X)$, $G\in\D^\perf(X,\CA_X)$ or
$F\in\D^\perf(X,\CA_X)$, $G\in\D^+(X,\CA_X)$ then we have
$\RCHom_{\CA_X}(F,G) \cong G\otimes_{\CA_X} \RCHom_{\CA_X}(F,\CA_X)$.
\end{lemma}
\begin{proof}
It is easy to construct a homomorphism from the RHS to the LHS.
It is clear that it is an isomorphism when either $F$ or $G$ is
locally projective. Hence it is an isomorphism when either $F$
or $G$ is perfect.
\end{proof}


\begin{lemma}\label{f11}
If $F\in\D^-(X,\CA_1)$, $G\in\D^-(X,\CA_1^\opp\otimes\CA_2)$
and $H\in\D^+(X,\CA_2)$ then we have\\
$\RCHom_{\CA_1}(F,\RCHom_{\CA_2}(G,H)) \cong
\RCHom_{\CA_2}(F\otimes_{\CA_1} G, H)$.
\end{lemma}
\begin{proof}
Evident.
\end{proof}

For any $L \in \D^-(X,\CA_X)$ we define the dual
$L_{\CA_X}^* := \RCHom_{\CA_X}(L,\CA_X) \in \D^+(X,\CA_X^\opp)$.
It is clear that the dual of a bimodule is a bimodule,
the dual of a perfect complex is a perfect complex,
and $(L_{\CA_X}^*)_{\CA_X}^* \cong L$ (functorial isomorphism) for perfect complexes.

\begin{lemma}
If $F\in\D^-(X,\CA_1)$, $G\in\D^+(X,\CA_2)$, and
$L\in\D^\perf(X,\CA_2^\opp\otimes\CA_1)$ then we have
$\RCHom_{\CA_1}(F,G \otimes_{\CA_2} L) \cong
\RCHom_{\CA_2}(F\otimes_{\CA_1} L_{\CA_2}^*,G)$.
Moreover, if $F\in\D^-(X,\CA_X)$, $G\in\D^+(X,\CA_X)$, and
$L\in\D^\perf(X,\CO_X)$ then we have
$\RCHom_{\CA_X}(F,G \otimes_{\CO_X} L) \cong
\RCHom_{\CA_X}(F\otimes_{\CO_X} L_{\CO_X}^*,G)$.
\end{lemma}
\begin{proof}
By~\ref{f11} we have
$\RCHom_{\CA_2}(F\otimes_{\CA_1} L_{\CA_2}^*,G) \cong
\RCHom_{\CA_1}(F,\RCHom_{\CA_2}(L_{\CA_2}^*,G))$
and by~\ref{fdg} we have
$\RCHom_{\CA_2}(L_{\CA_2}^*,G) \cong
G\otimes_{\CA_2}\RCHom_{\CA_2}(L_{\CA_2}^*,\CA_2) \cong
G\otimes_{\CA_2} L$. The second claim can be proved similarly.
\end{proof}

\begin{lemma}\label{fsshriek}
If $f:(X,\CA_X) \to (Y,\CA_Y)$ is a morphism,
$F \in \D^-(Y,\CA_Y)$ and $G\in\D^+(Y,\CA_Y)$ then we have
$f_\circ^!\RCHom_{\CA_Y}(F,G) \cong
\RCHom_{f_\circ^*\CA_Y}(f_\circ^*F,f_\circ^!G)$.
\end{lemma}
\begin{proof}
It is clear that
$f_\circ^!\RCHom_{\CO_Y}(F,G) \cong \RCHom_{\CO_X}(f_\circ^*F,f_\circ^!G)$
by the usual properties of the twisted pullback. Using the bar-resolution
we then deduce the claim.
\end{proof}

\begin{lemma}\label{fsshriek1}
If $f:(X,\CA_X) \to (Y,\CA_Y)$ is a morphism,
$F \in \D^\perf(Y,\CA_Y)$ and $G\in\D^+(Y,\CA_Y^\opp)$ then we have
$f_\circ^!(F\otimes_{\CA_Y}G) \cong
f_\circ^*F\otimes_{f_\circ^*\CA_Y}f_\circ^!G$.
\end{lemma}
\begin{proof}
It is clear that
$f_\circ^!(F\otimes_{\CO_Y}G) \cong f_\circ^*F\otimes_{\CO_X}f_\circ^!G$
by the usual properties of the twisted pullback. Using the bar-resolution
we then deduce the claim.
\end{proof}

\begin{lemma}[Functoriality of the twisted pullback]
If $(X,\CA_X) \exto{f} (Y,\CA_Y) \exto{g} (Z,\CA_Z)$ are morphisms
of Azumaya varieties then $f^!g^! \cong (gf)^!$ on $\D^+(Z,\CA_Z)$.
\end{lemma}
\begin{proof}
First, consider the case when both $f$ and $g$ are extensions, so $X = Y = Z$.
Then
$f^!g^!(F) \cong \RCHom_{\CA_Y}(\CA_X,\RCHom_{\CA_Z}(\CA_Y,F))$
which is isomorphic to
$\RCHom_{\CA_Z}(\CA_X,F) \cong (gf)^!(F)$
by~\ref{f11}.

Further, consider the case when both $f$ and $g$ are strict.
Then $f^! = f_\circ^!$, $g^! = g_\circ^!$ and the desired
isomorphism follows from the functoriality of the twisted
pullback functor for morphisms of algebraic varieties.

Since $f^!g^! \cong (gf)^!$ if $f$ is an extension and $g$ is strict,
it remains to show that $f^!g^! \cong (gf)^!$ when $f$ is strict and
$g$ is an extension. Indeed, in this case $Z = Y$, $\CA_X = f_\circ^*\CA_Y$,
the canonical decomposition of $gf:(X,f_\circ^*\CA_Y) \to (Y,\CA_Z)$ takes form
$\xymatrix@1{(X,f_{\circ}^*\CA_Y) \ar[r]^{f_\circ^* g_\CA} &
(X,f_\circ^*\CA_Z) \ar[r]^{f_\circ} &
(Y,\CA_Z)}$,
we have
$f^!g^!(G) = f_\circ^!\RCHom_{\CA_Z}(\CA_Y,G)$,
$(gf)^!(G) = \RCHom_{f_\circ^*\CA_Z}(f_{\circ}^*\CA_Y,f_\circ^! G)$
for all $G\in\D^+(Z,\CA_Z)$ and it remains to apply~\ref{fsshriek}
with $F = \CA_Y$.
\end{proof}

\begin{lemma}
If $f:(X,\CA_X) \to (Y,\CA_Y)$ is a morphism of Azumaya varieties
$F \in \D^+(Y,\CA_Y)$, and $G\in\D^\perf(Y,\CO_Y)$ then we have
$f^!(F\otimes_{\CO_Y} G) \cong f^!F \otimes_{\CO_X} f_\circ^*G$.
\end{lemma}
\begin{proof}
By functoriality of the pullback and of the twisted pullback
it suffices to prove this only for a strict morphism
and for an extension. For a strict morphism $f$ we have
$f^! = f_\circ^!$, hence we can use properties of the usual
twisted pullback. Now assume that $f$ is an extension. Then
by~\ref{f22} we have
$f^!(F\otimes_{\CO_Y} G) \cong
\RCHom_{\CA_Y}(\CA_X,F\otimes_{\CO_Y} G) \cong
\RCHom_{\CA_Y}(\CA_X,F)\otimes_{\CO_Y} G \cong
f^!F \otimes_{\CO_X} f_\circ^*G$, and we are done.
\end{proof}

\begin{lemma}
If $f:(X,\CA_X) \to (Y,\CA_Y)$ is a morphism of Azumaya varieties,
$F \in \D^+(Y,\CA_Y)$, and $G\in\D^\perf(Y,\CO_Y)$ then we have
$f^!\RCHom_{\CO_Y}(F,G) \cong \RCHom_{\CO_X}(f^*F,f_\circ^!G)$.
\end{lemma}
\begin{proof}
By functoriality of the pullback and of the twisted pullback
it suffices to prove this only for a strict morphism
and for an extension. For a strict morphism $f$ we have
$f^* = f_\circ^*$ and $f^! = f_\circ^!$, hence we can use
properties of the usual twisted pullback. Now assume that $f$
is an extension. Then we have
$f^!\RCHom_{\CO_Y}(F,G) \cong
\RCHom_{\CA_Y}(\CA_X,\RCHom_{\CO_Y}(F,G)) \cong
\RCHom_{\CO_Y}(\CA_X\otimes_{\CA_Y} F,G) \cong
\RCHom_{\CO_Y}(f^*F,G) \cong
\RCHom_{\CO_X}(f^*F,f_\circ^!G)$,
since $X=Y$ and $f_\circ = \id$.
\end{proof}

\begin{lemma}[Duality]\label{duality}
If $f:(X,\CA_X) \to (Y,\CA_Y)$ is a projective morphism of Azumaya varieties,
$F \in \D^-(X,\CA_X)$, $G\in\D^+(Y,\CA_Y)$ then we have
$\RCHom_{\CA_Y}(f_*F,G) \cong f_*\RCHom_{\CA_X}(F,f^!G)$.
\end{lemma}
\begin{proof}
By functoriality of the pushforward and of the twisted pullback
it suffices to prove the duality only for a strict morphism
and for an extension. For a strict morphism $f$ we have $f_* = {f_\circ}_*$,
$f^! = f_\circ^!$, and
$\RCHom_{\CO_Y}(f_*F,G) \cong f_*\RCHom_{\CO_X}(F,f^!G)$ by the usual
duality theorem. Using the bar-resolution we then deduce that
$\RCHom_{\CA_Y}(f_*F,G) \cong f_*\RCHom_{\CA_X}(F,f^!G)$.
Now assume that $f$ is an extension. Then we have
$f_*\RCHom_{\CA_X}(F,f^!G) \cong
\RCHom_{\CA_X}(F,\RCHom_{\CA_Y}(\CA_X,G)) \cong
\RCHom_{\CA_Y}(F\otimes_{\CA_X}\CA_X,G) \cong
\RCHom_{\CA_Y}(F,G) \cong
\RCHom_{\CA_Y}(f_*F,G)$
and we are done.
\end{proof}

\begin{corollary}
If $f:(X,\CA_X) \to (Y,\CA_Y)$ is a projective morphism of Azumaya varieties
then the functor $f^!$ is right adjoint to $f_*$
on the bounded derived category.
\end{corollary}

\mysubsection{Fiber products and base changes}

Let $(X,\CA_X)$ and $(Y,\CA_Y)$ be Azumaya varieties.
Then $(X\times Y,\CA_X\boxtimes\CA_Y)$ is also an Azumaya variety
and the projections of $X\times Y$ to $X$ and $Y$ with embeddings
of $\CA_X\boxtimes\CO_Y$ and $\CO_X\boxtimes\CA_Y$ to
$\CA_X\boxtimes\CA_Y$ define maps from $(X\times Y,\CA_X\boxtimes\CA_Y)$
to $(X,\CA_X)$ and $(Y,\CA_Y)$. Note that $(X\times Y,\CA_X\boxtimes\CA_Y)$
is not a categorical product of $(X,\CA_X)$ and $(Y,\CA_Y)$.
Nevertheless, we will denote it by $(X,\CA_X)\times(Y,\CA_Y)$.

\begin{lemma}
If $p:(X,\CA_X)\times(Y,\CA_Y) \to (X,\CA_X)$ is the projection then we have
$p^*F \cong F\boxtimes\CA_Y$ and
$p^!F \cong F\boxtimes\RCHom_{\CO_Y}(\CA_Y,\omega_Y)[\dim Y]$.
\end{lemma}
\begin{proof}
Evident.
\end{proof}

It is easy to see that in general the fiber products don't exist in
the category of Azumaya varieties. However, there are important cases
when they do exist.

Assume that $f:(X,\CA_X) \to (S,\CA_S)$ and $g:(Y,\CA_Y) \to (S,\CA_S)$
are morphisms of Azumaya varieties. Let $p_\circ:X\times_S Y \to X$
and $q_\circ:X\times_S Y \to Y$ denote the projections.

\begin{lemma}\label{fpxy}
If $f$ is strict then $(X\times_S Y, q_\circ^*\CA_Y)$ is a fiber product
of $(X,\CA_X)$ and $(Y,\CA_Y)$ over $(S,\CA_S)$.
If $g$ is strict then $(X\times_S Y, p_\circ^*\CA_X)$ is a fiber product
of $(X,\CA_X)$ and $(Y,\CA_Y)$ over $(S,\CA_S)$.
\end{lemma}
\begin{proof}
If $\phi:(Z,\CA_Z)\to(X,\CA_X)$ and $\psi:(Z,\CA_Z) \to (Y,\CA_Y)$ are
morphisms of Azumaya varieties such that $f\circ\phi = g\circ\psi$
then $f_\circ\circ\phi_\circ = g_\circ\circ\psi_\circ$, hence
there exists a morphism $\xi_\circ:Z\to X\times_SY$ such that
$p_\circ\circ\xi_\circ = \phi_\circ$ and $q_\circ\circ\xi_\circ = \psi_\circ$.
Moreover, if $f$ is strict then the map
$\xi_\circ^*q_\circ^*\CA_Y \cong \psi_\circ^*\CA_Y \exto{\psi_\CA} \CA_Z$
induces a morphism of Azumaya varieties
$(Z,\CA_Z) \to (X\times_S Y, q_\circ^*\CA_Y)$.
Similarly, if $g$ is strict then the map
$\xi_\circ^*p_\circ^*\CA_X \cong \phi_\circ^*\CA_X \exto{\phi_\CA} \CA_Z$
induces a morphism of Azumaya varieties
$(Z,\CA_Z) \to (X\times_S Y, p_\circ^*\CA_X)$.
\end{proof}

We often consider strict morphism $f:(X,\CA_X) \to (Y,\CA_Y)$
as a base change.

\begin{lemma}
A base change preserves flatness, smoothness, projectivity, e.t.c.\
of a morphism.
\end{lemma}
\begin{proof}
Evident.
\end{proof}

If the algebra $\CA_X$ is noncommutative, the multiplication
$\CA_X\otimes\CA_X \to \CA_X$ is not a homomorphism of algebras.
Therefore, there is no ``diagonal embedding'' in the category
of Azumaya varieties. However, $\CA_X$ is a $\CA_X\otimes\CA^\opp_X$-module,
thus ``the structure sheaf of the diagonal'',
$\Delta_*\CA_X$ on $X\times X$ can be considered as an object of
$\Coh(X\times X,\CA_X\boxtimes\CA^\opp_X)$.

\mysubsection{$\Tor$ and $\Ext$-amplitude}

Let $f:(X,\CA_X) \to (Y,\CA_Y)$ be a morphism of Azumaya varieties.

\begin{definition}
An object $F \in \D(X,\CA_X)$ has {\sf finite $\Tor$-amplitude
over $(Y,\CA_Y)$}, if there exist integers $p,q$ such that
for any object $G \in \D^{[s,t]}(Y,\CA_Y)$ we have
$F\otimes_{\CA_X} f^*G \in \D^{[p+s,q+t]}(X,\CO_X)$.
Morphism $f$ has {\sf finite $\Tor$-dimension}, if the sheaf $\CA_X$ has
finite $\Tor$-amplitude over $(Y,\CA_Y)$.
Similarly, an object $F \in \D(X,\CA_X)$ has {\sf finite $\Ext$-amplitude
over $(Y,\CA_Y)$}, if there exist integers $p,q$ such that
for any object $G \in \D^{[s,t]}(Y,\CA_Y)$ we have
$\RCHom_{\CA_X}(F,f^!G) \in \D^{[p+s,q+t]}(X,\CO_X)$.
Morphism $f$ has {\sf finite $\Ext$-dimension}, if the sheaf $\CA_X$ has
finite $\Ext$-amplitude over $(Y,\CA_Y)$.
\end{definition}

The full subcategory of $\D(X,\CA_X)$ consisting of objects
of finite $\Tor$-amplitude over $(Y,\CA_Y)$ is denoted by
$\D_{fTd/(Y,\CA_Y)}(X,\CA_X)$.
The full subcategory of $\D(X,\CA_X)$ consisting of objects
of finite $\Ext$-amplitude over $(Y,\CA_Y)$ is denoted by
$\D_{fEd/(Y,\CA_Y)}(X,\CA_X)$.
Both are triangulated subcategories of $\D^b(X,\CA_X)$.

\begin{lemma}\label{isfted}
If $i:(X,\CA_X) \to (X',\CA_{X'})$ is a finite morphism
of Azumaya varieties over $(Y,\CA_Y)$ then
$F\in\D_{fTd/(Y,\CA_Y)}(X,\CA_X)$ $\LRA$
$i_*F\in\D_{fTd/(Y,\CA_Y)}(X',\CA_{X'})$ and
$F\in\D_{fEd/(Y,\CA_Y)}(X,\CA_X)$ $\LRA$
$i_*F\in\D_{fEd/(Y,\CA_Y)}(X',\CA_{X'})$.
\end{lemma}
\begin{proof}
Use exactness of $i_*$, the isomorphisms
$i_*F \otimes_{\CA_{X'}} {f'}^*G \cong
i_*(F \otimes_{\CA_X} i^*{f'}^*G) \cong
i_*(F \otimes_{\CA_X} f^*G)$
and
$\RCHom_{\CA_{X'}}(i_*F,{f'}^!G) \cong
i_*\RCHom_{\CA_X}(F,i^!{f'}^!G) \cong
i_*\RCHom_{\CA_X}(F,f^!G)$.
\end{proof}

\begin{lemma}\label{ftedperf}
If morphism $f:(X,\CA_X) \to (Y,\CA_Y)$ has finite $\Tor$-dimension
then any perfect complex on $(X,\CA_X)$ has finite $\Tor$-amplitude
over $(Y,\CA_Y)$.
If morphism $f:(X,\CA_X) \to (Y,\CA_Y)$ has finite $\Ext$-dimension
then any perfect complex on $(X,\CA_X)$ has finite $\Ext$-amplitude
over $(Y,\CA_Y)$.
\end{lemma}
\begin{proof}
Evident.
\end{proof}

Now we are going to prove the inverse statements for a smooth morphism.
For this we will need the following fact about perfect complexes on
usual algebraic varieties. This fact certainly must be well known,
however I don't know a reference, so the proof is included here
for completeness.

\begin{lemma}\label{ftdisp}
Assume that $f:X \to S$ is a smooth morphism of algebraic varieties
and $F$ is a coherent sheaf on $X$. If $F$ is either flat or
locally projective over $S$ then $F$ is a perfect complex.
\end{lemma}
\begin{proof}
The claim is local in $X$, so we can assume that $X$ is affine
and the ideal of the diagonal $X\subset X\times_S X$ is generated
by a regular sequence $s = (s_1,\dots,s_n)$ of functions on $X$. Then
$\Delta_*\CO_X \cong
\Kosz_{X\times_S X}(s) \cong
\Lambda^\bullet(\CO_{X\times_S X}^{\oplus n})$.
Consider the diagram
$$
\xymatrix{
X \ar[r]^-{\Delta} &
X\times_S X \ar[r]^q \ar[d]_p &
X \ar[d]^f \\
& X \ar[r]^f &
S
}
$$
and note that $q\circ\Delta = p\circ\Delta = \id_X$, hence
for any closed point $x\in X$ we have
$$
\begin{array}{l}
\CO_x \cong
q_*\Delta_*\Delta^*p^*\CO_x \cong
q_*(\Delta_*\CO_X\otimes p^*\CO_x) \cong
q_*(\Kosz_{X\times_S X}(s)\otimes p^*\CO_x),\\
\CO_x \cong
q_*\Delta_*\Delta^!p^!\CO_x \cong
q_*\RCHom(\Delta_*\CO_X,p^!\CO_x) \cong
q_*\RCHom(\Kosz_{X\times_S X}(s),p^!\CO_x).
\end{array}
$$
But the flat base change implies that
$$
\begin{array}{l}
q_*(\Lambda^r(\CO_{X\times_S X}^{\oplus n})\otimes p^*\CO_x) \cong
q_*(\CO_{X\times_S X}^{\oplus\binom{n}{r}}\otimes p^*\CO_x) \cong
q_*p^*\CO_x^{\oplus\binom{n}{r}} \cong
f^*f_*\CO_x^{\oplus\binom{n}{r}} \cong
f^*\CO_{f(x)}^{\oplus\binom{n}{r}},\\
q_*\RCHom(\Lambda^r(\CO_{X\times_S X}^{\oplus n}),p^!\CO_x) \cong
q_*\RCHom(\CO_{X\times_S X}^{\oplus\binom{n}{r}},p^!\CO_x) \cong
q_*p^!\CO_x^{\oplus\binom{n}{r}} \cong
f^!f_*\CO_x^{\oplus\binom{n}{r}} \cong
f^!\CO_{f(x)}^{\oplus\binom{n}{r}}.
\end{array}
$$
Thus every sheaf $\CO_x$ on $X$ admits a length $n$ resolution
either by sheaves, which are pullbacks of coherent sheaves on $S$,
or by twisted pullbacks of coherent sheaves on $S$.
It follows that $\Tor_{>n}(F,\CO_x) = 0$
for all $x\in X$ if $F$ is flat over $S$,
and $\Ext^{>n}(F,\CO_x) = 0$ for all $x\in X$
if $F$ is locally projective over $S$.
On the other hand, it is easy to construct
a resolution of $F$ on $X$ of the form
$0 \to F' \to F_{n - 1} \to \dots \to F_1 \to F_0 \to F \to 0$
with locally free $F_i$. Then it follows that $\Tor_{>0}(F',\CO_x) = 0$
for all $x\in X$ if $F$ is flat over $S$, and $\Ext^{>0}(F',\CO_x) = 0$
for all $x\in X$ if $F$ is locally projective over $S$.
Hence in both cases $F'$ is locally free, so $F$ is a perfect complex.
\end{proof}

\begin{lemma}\label{snav_sm1}
If $f:(X,\CA_X) \to (Y,\CA_Y)$ is a smooth morphism of Azumaya varieties
then we have
$\D_{fTd/(Y,\CA_Y)}(X,\CA_X) =
\D^\perf(X,\CA_X) =
\D_{fEd/(Y,\CA_Y)}(X,\CA_X)$.
\end{lemma}
\begin{proof}
Any perfect complex has finite $\Tor$-amplitude by lemma~\ref{ftedperf}.
Take $F\in\D_{fTd/(Y,\CA_Y)}(X,\CA_X)$.
Since $\Coh(X,\CA_X)$ has enough locally free objects,
we can construct a locally free over $\CA_X$
resolution $F_{n-1} \to \dots \to F_1 \to F_0 \to F \to 0$ for
arbitrary large $n$. On the other hand, if $n$ is sufficiently large, then
$F_{n} := \Ker (F_{n-1} \to F_{n-2})$ is flat over $(Y,\CA_Y)$,
hence $F_{n}$ is flat over $Y$ (since the canonical projection
$(Y,\CA_Y) \to Y$ is flat), hence $F_{n}$ is locally free over $\CO_X$,
hence $F_{n+1}$ is locally projective over $\CA_X$ by lemma~\ref{locproj}.
Therefore, any $F$ is quasiisomorphic to a finite complex of locally
projective $\CA_X$-modules, hence $F$ is a perfect complex.
The same arguments show $\D_{fEd/(Y,\CA_Y)}(X,\CA_X) = \D^\perf(X,\CA_X)$.
\end{proof}

\begin{corollary}\label{ftd_perf}
$\D_{fTd/(X,\CA_X)}(X,\CA_X) =
\D^\perf(X,\CA_X) =
\D_{fEd/(X,\CA_X)}(X,\CA_X)$.
\end{corollary}

\begin{corollary}\label{fted_crit}
Let $f = f'\circ f''$ be a decomposition of a morphism
$f:(X,\CA_X) \to (Y,\CA_Y)$ with smooth $f':(X',\CA_{X'}) \to (Y,\CA_Y)$
and finite $f'':(X,\CA_X) \to (X',\CA_{X'})$.
Then
$F\in\D_{fTd/(Y,\CA_Y)}(X,\CA_X)$
$\LRA$
$f''_*F \in \D^\perf(X',\CA_{X'})$
$\LRA$
$F\in\D_{fEd/(Y,\CA_Y)}(X,\CA_X)$.
\end{corollary}
\begin{proof}
Use lemma~\ref{snav_sm1} and lemma~\ref{isfted}.
\end{proof}

\begin{corollary}\label{fted_bc}
If $\phi:(W,\CA_W) \to (Y,\CA_Y)$ is a {\rm(}strict and{\rm)} faithful
base change for a morphism $f:(X,\CA_X) \to (Y,\CA_Y)$ and
$(Z,\CA_Z) = (X,\CA_X)\times_{(Y,\CA_Y)}(W,\CA_W)$ is the fiber product
then we have
$\phi^*(\D_{fTd/(Y,\CA_Y)}(X,\CA_X)) \subset \D_{fTd/(W,\CA_W)}(Z,\CA_Z)$,
and
$\phi^*(\D_{fEd/(Y,\CA_Y)}(X,\CA_X)) \subset \D_{fEd/(W,\CA_W)}(Z,\CA_Z)$.
\end{corollary}
\begin{proof}
Decompose $f:(X,\CA_X) \to (Y,\CA_Y)$ as a product of
smooth $f':(X',\CA_{X'}) \to (Y,\CA_Y)$ and
finite $f'':(X,\CA_X) \to (X',\CA_{X'})$, and denote
$(Z',\CA_{Z'}) = (X',\CA_{X'})\times_{(Y,\CA_Y)}(W,\CA_W)$, the fiber product.
Then we have a commutative diagram
$$
\xymatrix{
(Z,\CA_Z) \ar[r]^{f''} \ar[d]^{\phi} &
(Z',\CA_{Z'}) \ar[r]^{f'} \ar[d]^{\phi} &
(W,\CA_W) \ar[d]^{\phi} \\
(X,\CA_X) \ar[r]^{f''} &
(X',\CA_{X'}) \ar[r]^{f'} &
(Y,\CA_Y)
}
$$
Its right square is exact cartesian because $f'$ is smooth,
and its left square is exact cartesian by lemma~\ref{squ3}.
Therefore $f''_*\phi^* F \cong \phi^*f''_* F$ and we deduce
the claim from the criterion~\ref{fted_crit}.
\end{proof}

\begin{lemma}
If $(Y,\CA_Y)$ is a smooth Azumaya variety then
$\D_{fTd/(Y,\CA_Y)}(X,\CA_X) = \D^b(X,\CA_X) = \D_{fEd/(Y,\CA_Y)}(X,\CA_X)$.
\end{lemma}
\begin{proof}
If $F\in\D_{fTd/(Y,\CA_Y)}(X,\CA_X)$ then
$F \cong F\otimes_{\CA_X}\CA_X \cong F\otimes_{\CA_X}f^*\CA_Y$ is bounded,
therefore we have $\D_{fTd/(Y,\CA_Y)}(X,\CA_X) \subset \D^b(X,\CA_X)$.
On the other hand, take any $F\in\D^{[p,q]}(X,\CA_X)$.
Since $Y$ is smooth any $G \in \D^{[s,t]}(Y,\CA_Y)$
is locally quasiisomorphic to a complex of locally projective
$\CA_Y$-modules concentrated in degrees $[s-\dim Y,t]$.
Hence $F\otimes_{\CA_X}f^*G \in \D^{[p+s-\dim Y,q+t]}(X,\CA_X)$
and $F\in\D_{fTd/(Y,\CA_Y)}(X,\CA_X)$.
\end{proof}

\end{document}